\documentclass[a4paper,11pt]{article}
\usepackage[utf8]{inputenc}
\usepackage{enumerate}
\usepackage{lmodern}
\usepackage[T1]{fontenc}
\usepackage{verbatim}
\usepackage{textcomp}

\usepackage[english]{babel}
\usepackage[a4paper,vmargin={3.5cm,3.5cm},hmargin={2.5cm,2.5cm}]{geometry}
\usepackage[font=sf, labelfont={sf,bf}, margin=1cm]{caption}
\usepackage[pdftex]{hyperref}

\usepackage[pdftex]{color,graphicx}

\usepackage{amsmath,amsfonts,amssymb,amsthm,mathrsfs}

\newtheorem{theorem}{Theorem}
\newtheorem{definition}{Definition}
\newtheorem{lemma}[theorem]{Lemma}

\newtheorem{corollary}[theorem]{Corollary}
\newtheorem{proposition}[theorem]{Proposition}

\def\llbracket{[\hspace{-.10em} [ }
\def\rrbracket{ ] \hspace{-.10em}]}

\def\w{\mathrm{w}}
\def\ll{{\mathcal L}}
\def\t{{\mathcal T}}
\def\cc{\mathcal{C}}
\def\ve{\varepsilon}
\def\wt{\widetilde}
\def\wh{\widehat}

\def\bp{{\bf p}}
\def\bv{\mathbf{v}}
\def\bV{\mathbf{V}}
\def\S{\mathcal{S}}
\def\W{\mathcal{W}}

\def\nn{\mathcal{N}}
\def\mm{\mathcal{M}}
\def\dd{\mathcal{D}}

\def\qq{\mathrm{q}}
\def\QQ{\mathcal{Q}}
\def\rr{\mathcal{R}}

\def\z{\mathcal{Z}}

\def\la{\longrightarrow}

\def\M{{\mathbb M}}

\def\E{{\mathbb E}}
\def\P{{\mathbb P}}
\def\N{{\mathbb N}}
\def\T{{\mathbb T}}
\def\Z{{\mathbb Z}}
\def\R{{\mathbb R}}
\def\Q{{\mathbb Q}}
\def\F{{\mathbb F}}
\def\D{{\mathbb D}}

\def\dg{\mathrm{d}_{\mathrm{gr}}}
\def\ov{\overline}
\def\un{\underline}
\def\tr{\mathrm{tr}}

\def\build#1_#2^#3{\mathrel{
\mathop{\kern 0pt#1}\limits_{#2}^{#3}}}
\def\rem{\noindent{\bf Remark. }}

\title{Brownian disks and the Brownian snake}
\author{Jean-Fran\c cois Le Gall\footnote{Supported by the ERC Advanced Grant 740943 {\sc GeoBrown}}}

\date{\small 
Revised version, October 20, 2017}

\begin{document}
\maketitle

\begin{abstract}We provide a new construction of the Brownian disks, which have been 
defined by Bettinelli and Miermont as scaling limits of quadrangulations 
with a boundary when the boundary size tends to infinity. 
Our method is
very similar to the construction of the Brownian map, but it makes use 
of the positive excursion measure of the Brownian snake which has been
introduced recently. This excursion measure
involves a random continuous tree whose vertices are assigned nonnegative
labels, which correspond to distances
from the boundary in our approach to the Brownian disk.
 We provide several applications of our construction. In
particular, we prove that the uniform measure on the boundary can be obtained as
the limit of the suitably normalized volume measure on a small tubular neighborhood
of the boundary. We also prove that connected components of
the complement of the Brownian net are Brownian disks, as it
was suggested in the recent work of Miller and Sheffield.
Finally, we show that connected components of the complement
of balls centered at the distinguished point of the Brownian map 
are independent Brownian disks, conditionally on their
volumes and perimeters. 
\end{abstract}

\section{Introduction}

In the last ten years, much work has been devoted to the construction
and study of continuous models of random geometry in two 
dimensions. A popular approach has been to start from 
the discrete models called planar maps, which are just graphs drawn 
on the two-dimensional sphere and viewed up to orientation-preserving 
homeomorphisms of the sphere. Choosing a planar map
uniformly at random in a suitable class, for instance the class
of all triangulations of the sphere with a fixed number of faces, 
yields a discrete model of random geometry. In order to get a
continuous model, one equips the vertex 
set of a random planar map with the (suitably rescaled) graph distance, 
and studies the convergence in distribution of the resulting
(discrete) metric space as the size of the map tends to infinity. 
This requires a notion of convergence for a sequence 
of compact metric spaces, which is provided by the
Gromov-Hausdorff distance (see for instance \cite{BBI}). It was proved independently in \cite{Mie2}
for quadrangulations and in \cite{Uniqueness} for more general 
cases including triangulations, that several important classes 
of random planar maps converge in this sense toward a
limiting random compact metric space called the Brownian map,
which had been introduced previously in \cite{MM1}. Several recent papers
(see in particular \cite{Abr,AA,BJM,Mar}) have shown that many other 
classes of random planar maps converge to the Brownian map, 
which thus provides a universal continuous model of random geometry
in two dimensions. On the other hand, in a series of
recent papers \cite{MS0,MS1,MS2,MS3}, Miller and Sheffield have shown 
that the Brownian map can be equipped with a conformal structure, which is linked
to Liouville quantum gravity, and 
that this conformal structure is in fact determined by the Brownian
map viewed as a random metric space. An important step
in this program was to derive an axiomatic characterization of the
Brownian map \cite{MS0}. 

The Brownian map is known to be homeomorphic to the sphere, but other
models homeomorphic to the disk have been introduced under the 
name of Brownian disks, and are expected to correspond to scaling limits
of planar maps with a boundary. Following the earlier work of Bettinelli \cite{Bet},
Bettinelli and 
Miermont \cite{BM} constructed Brownian disks as scaling limits of quadrangulations with a boundary,
when the number of faces grows like a constant times the 
square of the size of the boundary
(see also  \cite{BMR} for a general discussion of possible scaling limits for quadrangulations 
with a boundary). A Brownian disk $\D$  is homeomorphic to the closed unit disk of the plane, 
and its boundary $\partial\D$ may be defined as the set of all points that have no neighborhood homeomorphic
to an open disk. Furthermore, Brownian disks are equipped with a natural volume measure
and are thus viewed as random measure metric spaces. The distribution of a Brownian disk
depends on two parameters, on one hand the 
perimeter or size of the boundary, and on the other hand the volume or total mass of the volume measure. 
If the perimeter $r$ is fixed, there is a natural way of choosing the volume at random, and this leads to
the so-called free Brownian disk with perimeter $r$, which is shown in \cite{BM}
to be the scaling limit of Boltzmann quadrangulations with a boundary (see Section \ref{sec:freeBdisk} below).
Although this has not yet been proved, it is very likely that the Brownian disk also
appears as the scaling limit of more general classes of random planar maps 
with a boundary. The case of quadrangulations with a {\it simple} boundary
was treated very recently by Gwynne and Miller \cite{GM2}, in view of applications 
to the study of percolation interfaces on random quadrangulations with a boundary
\cite{GM3}. Interestingly, Brownian disks, as defined in \cite{BM}, also
play an important role in the recent papers \cite{GM0,GM1}.

Our goal in this work is to develop a new construction of Brownian disks in terms
of the Brownian snake positive excursion measure, which has been introduced and studied
in \cite{ALG}. We hope that this construction will be a useful tool for further
investigations of Brownian disks and for new
applications of these random objects. Our construction is
in a sense similar to the one in \cite{BM} as it relies on a 
a random continuous tree whose vertices are assigned real labels,
but, in contrast with \cite{BM}, labels in our approach correspond to distances from the boundary.
This has several advantages, and in particular it makes it
possible to define the uniform measure on the boundary as
the limit of the volume measure restricted to the $\ve$-neighborhood of the 
boundary and scaled by the factor $\ve^{-2}$. Our construction
also enables us to prove that Brownian disks can be embedded in the
Brownian map in various ways. In particular, we show that connected components
of the complement of balls in the Brownian map
are Brownian disks. We also establish that connected
components of the complement of the so-called Brownian net, which
is a closed subset of the Brownian map playing an important
role in the axiomatic characterization of \cite{MS0}, are Brownian disks.
As a matter of fact, these connected components were already 
called Brownian disks in \cite[Section 4.2]{MS0}, but the equivalence with the 
definition of Brownian disks as scaling limits of random planar maps
with a boundary had not been proved (see the discussion in \cite[Section 1.6]{BM}).
We note that a different method aiming to connect the definition of Brownian 
disks in \cite{MS0} with that in \cite{BM} is developed in 
the forthcoming paper \cite{JM}. 

Let us turn to a more precise description of our main results. We start by
recalling that the usual construction of the Brownian map relies on
considering the Brownian snake excursion measure $\N_0$, which is
a convenient way of representing Brownian motion indexed by the
Brownian tree.  Under $\N_0$, the values of the Brownian snake form 
a collection $(W_s)_{0\leq s\leq \sigma}$, 
such that, for every $s\in[0,\sigma]$, $W_s=(W_s(t))_{0\leq t\leq \zeta_s}$
is a finite Brownian path
with initial point $0$ and lifetime $\zeta_s\geq 0$.
The lifetime process
$(\zeta_s)_{0\leq s\leq \sigma}$ is distributed under $\N_0$
according to the It\^o measure of positive excursions of linear
Brownian motion, see e.g. \cite{RY} . Furthermore, the genealogical structure
of the paths $W_s$, $0\leq s\leq \sigma$ is described by
the tree $\t_\zeta$, which is the tree coded by $(\zeta_s)_{0\leq s\leq \sigma}$
(see Section \ref{sec:snake} for a precise definition of $\t_\zeta$).
 Recall that $\t_\zeta$
is obtained as a quotient space of the interval $[0,\sigma]$, and write
$p_\zeta$ for the corresponding canonical projection, so that $\t_\zeta$ is 
equipped with a volume measure defined as the push forward of
Lebesgue measure under $p_\zeta$.  
One can define ``Brownian labels'' on the tree $\t_\zeta$
by setting $Z_a:=\wh W_s$ if $a=p_\zeta(s)$, where 
$\wh W_s=W_s(\zeta_s)$ is the terminal point of the path $W_s$. 
One may also define ``lexicographical intervals'' on the tree $\t_\zeta$:
Informally, if $a,b\in\t_\zeta$, $[a,b]$ corresponds to the set of
points of $\t_\zeta$ that are visited when going ``clockwise around the
tree'' from $a$ to $b$ (see Section \ref{sec:snake} below for a 
more precise definition). We then set, for every $a,b\in\t_\zeta$,
\begin{equation}
\label{Dzero-i}
D^\circ(a,b)=Z_a + Z_b -2\max\Big(\min_{c\in[a,b]} Z_c, \min_{c\in[b,a]} Z_c\Big).
\end{equation}
and 
\begin{equation}
\label{formulaD-i}
D(a,b) = \inf\Big\{ \sum_{i=1}^k D^\circ(a_{i-1},a_i)\Big\},
\end{equation}
where the infimum is over all choices of the integer $k\geq 1$ and of the
elements $a_0,a_1,\ldots,a_k$ of $\t_\zeta$ such that $a_0=a$
and $a_k=b$. We observe that $D$ defines a pseudo-metric on $\t_\zeta$
and let $\approx$ stand for the associated equivalence relation
($a\approx b$ if and only if $D(a,b)=0$, which turns out to be equivalent
to $D^\circ(a,b)=0$). The Brownian map is then defined
as the quotient $\mm:=\t_\zeta\,/\!\approx$ equipped with the distance 
induced by $D$ and with the volume measure $\bv$ which is the push forward of
the volume measure on $\t_\zeta$ under the canonical projection. More precisely, this is what may be called the ``free Brownian 
map'' under the ($\sigma$-finite) measure $\N_0$, and the same construction under 
the probability measure $\N^{(1)}_0:=\N_0(\cdot\mid \sigma=1)$ (under which 
the lifetime process is
a normalized Brownian excursion) yields the standard Brownian map. 

\begin{figure}[!h]
 \begin{center}
 \includegraphics[width=10cm]{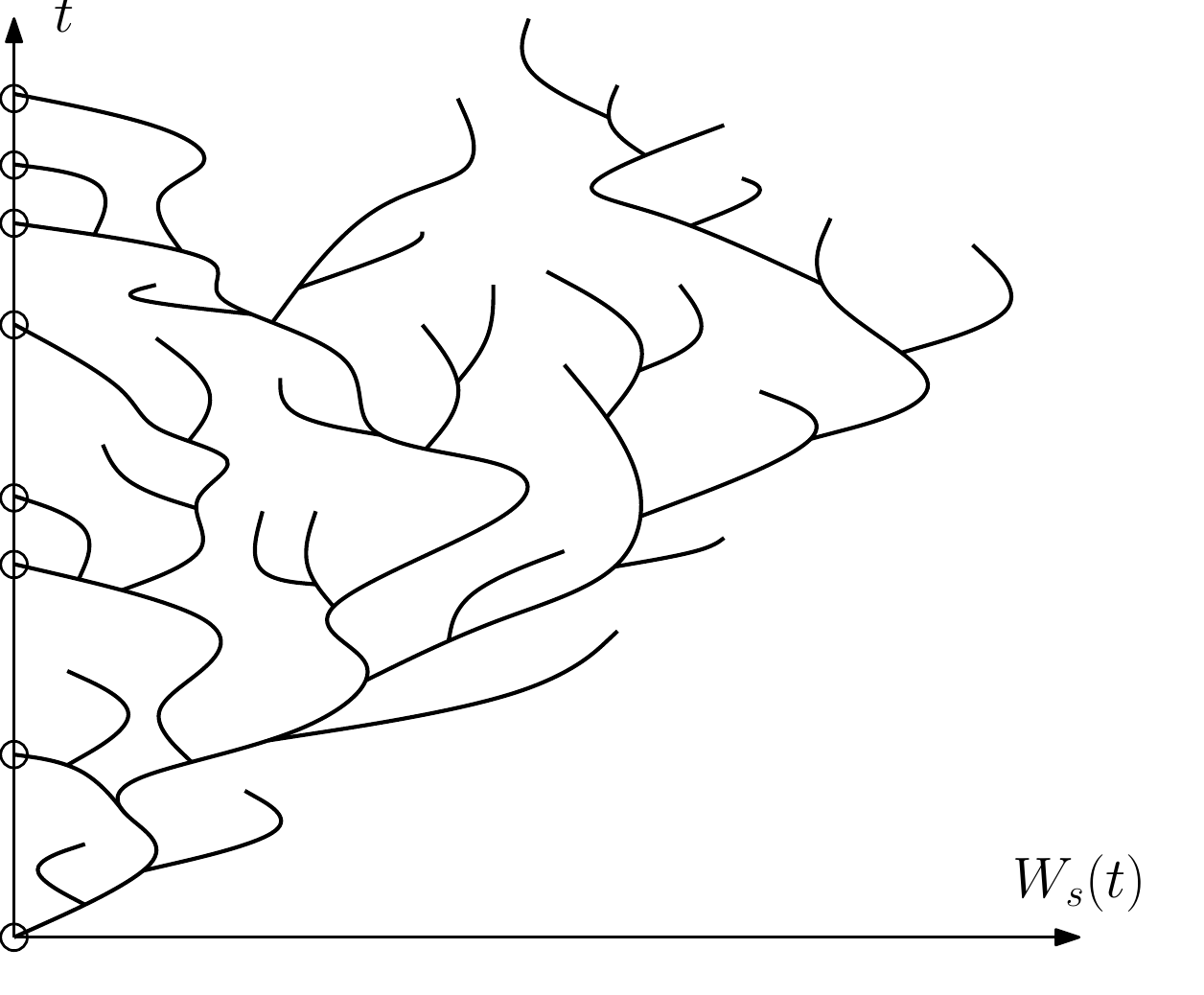}
 \caption{\label{tree}
 A schematic representation of the tree of paths $(W_s)_{0\leq s\leq \sigma}$ under $\N^*_0$.
 The time parameter $t$ for the paths $W_s$ is represented by the vertical coordinate.
 The exit measure $\z^*_0$ ``counts'' the number of circled points corresponding to those 
 paths $W_s$ that return to $0$. For the sake of clarity, the paths $W_s$ do not cross on the figure,
 although of course they should.}
 \end{center}
 \vspace{-5mm}
 \end{figure}

Quite remarkably, a close variant of the preceding construction will give the
Brownian disk. The main ingredient now is the positive excursion measure
$\N^*_0$ introduced in \cite{ALG}. Under $\N^*_0$, we still have a
collection of finite paths $W_s$, $0\leq s\leq \sigma$, with respective 
lifetimes $\zeta_s$, $0\leq s\leq \sigma$, whose genealogical structure is
determined by the tree $\t_\zeta$ coded by $(\zeta_s)_{0\leq s\leq \sigma}$. However the paths $W_s$
now behave quite differently from Brownian paths. Indeed, each path $W_s$, $0<s<\sigma$,
starts from $0$, then stays positive
on some nontrivial interval $(0,u)$ and if it returns to $0$ it is stopped at
that moment. To give a somewhat informal description of $\N^*_0$, let 
$\N_\ve$
be obtained from $\N_0$ by adding $\ve$ to the paths $W_s$.
Then, in a way similar to classical approximations
of the It\^o excursion measure of Brownian motion, the measure $\N^*_0$ 
appears as the limit when $\ve \downarrow 0$ of the distribution under the (rescaled) measure  $\N_\ve$ of the
collection of paths $(W_s(\cdot\wedge \tau_0(W_s)))_{0\leq s\leq \sigma}$, where $\tau_0(W_s)$ denotes the
first hitting time of $0$ by $W_s$
 (see \cite[Theorem 23]{ALG} for a more precise statement). 
Under $\N^*_0$, one can define a random quantity $\z^*_0$ called the
exit measure, which counts in some sense the number of paths $W_s$ that
return to $0$, and one can make sense of the conditional probability
measures $\N^{*,z}_0=\N^*_0(\cdot\mid \z^*_0=z)$. 
See Fig.\ref{tree} for a schematic illustration.

Arguing now under $\N^{*,z}_0$ for some fixed $z>0$, we consider
again the tree $\t_\zeta$, the labels $Z_a$, the functions
$D^\circ(a,b)$ and $D(a,b)$, and the 
equivalence relation $\approx$ defined exactly as above. Then, under $\N^{*,z}_0$,
the quotient space $\t_\zeta\,/\!\approx$ equipped with the distance induced by $D$ turns out to be the free Brownian 
disk with perimeter $z$ and boundary glued into a single point. Here,
if $(\D_z,\mathbf{d}_z)$ stands for a Brownian disk with perimeter
$z$, and if $\partial \D_z$ denotes its boundary, the Brownian disk 
with glued boundary is the quotient of $\D_z$ corresponding to
the pseudo-metric
$$\mathbf{d}^\dagger_z(x,y)=\min\{\mathbf{d}_z(x,y),\mathbf{d}_z(x,\partial\D_z)+\mathbf{d}_z(y,\partial\D_z)\},$$
which identifies all points of $\partial\D_z$. The reason why a naive imitation under $\N^*_0$
of the construction of the Brownian map yields an object with glued boundary
is easy to understand: If one defines the ``boundary'' of $\t_\zeta$ under $\N^{*,z}_0$
by $\partial\t_\zeta:=\{a\in\t_\zeta:Z_a=0\}$, it is immediate from 
\eqref{Dzero-i} that $D(a,b)=0$ for every $a,b\in \partial\t_\zeta$, $\N^{*,z}_0$ a.s.

The point is now to reconstruct the metric on the Brownian disk from
the one on the same object with glued boundary. Still arguing under
$\N^{*,z}_0$, we set $\Delta^\circ(a,b)=D^\circ(a,b)$ if 
$$\max\Big(\min_{c\in[a,b]} Z_c, \min_{c\in[b,a]} Z_c\Big)>0$$
and $\Delta^\circ(a,b)=\infty$ otherwise. Informally, the condition in the last 
display means that we can go from $a$ to $b$ around the tree without
visiting a vertex of $\partial\t_\zeta$. We then define $\Delta(a,b)$
for $a,b\in\t_\zeta\backslash\partial\t_\zeta$ by the analog
of formula \eqref{formulaD-i} where $D^\circ$ is replaced by $\Delta^\circ$.
It is not hard to verify that the mapping $(a,b)\mapsto \Delta(a,b)$ takes
finite values and is continuous on $(\t_\zeta\backslash\partial\t_\zeta)\times(\t_\zeta\backslash
\partial\t_\zeta)$.

\begin{theorem}
\label{main}
With probability one under $\N^{*,z}_0$, the function $(a,b)\mapsto \Delta(a,b)$
has a continuous extension to $\t_\zeta\times \t_\zeta$, which is a pseudo-metric on $\t_\zeta$.
We let $\Theta$ stand for the associated quotient space, and we equip $\Theta$ with the induced
metric, which is still denoted by $\Delta(a,b)$, and with the volume measure $\bV(\mathrm{d}x)$, which is
the image of the volume measure on $\t_\zeta$ under the canonical projection. Then, the random measure metric space 
$(\Theta,\Delta,\bV)$ is a free Brownian disk with perimeter $z$ under $\N^{*,z}_0$,
and its boundary $\partial\Theta$ is the image of $\partial\t_\zeta$
under the canonical projection from $\t_\zeta$ onto $\Theta$. Furthermore, 
if $x\in\Theta$ is the image of $a\in\t_\zeta$ under the canonical projection, we have
$$\Delta(x,\partial \Theta)=Z_a.$$
\end{theorem}

In order to recover a Brownian disk with perimeter $z$ and fixed volume $v$, it then suffices 
to condition $\N^{*,z}_0$ on the event $\{\bV(\Theta)=v\}$, noting that $\bV(\Theta)=\sigma$ is nothing 
but the duration of the excursion under $\N^{*,z}_0$.

\begin{proposition}
\label{measure-boundary}
Almost surely under $\N^{*,z}_0$, there exists a finite measure $\nu$
on $\partial\Theta$ with total mass $z$, such that, for every bounded
continuous function $\varphi$ on $\Theta$,
$$\langle \nu, \varphi\rangle = \lim_{\ve\to 0} \frac{1}{\ve^2}\int_\Theta \bV(\mathrm{d}x)\,\varphi(x)\,\mathbf{1}_{\{\Delta(x,\partial\Theta)<\ve\}}.$$
\end{proposition}

The measure $\nu$ is interpreted as the uniform measure on the boundary. Notice that the
construction in \cite{Bet,BM} also gives a natural way of defining a
measure on the boundary, which one should be able to identify with the measure $\nu$
of the preceding proposition. This identification does not seem to be easy because the construction in \cite{Bet,BM}
is rather different, and does not give access to distances from the boundary which are used in our approximations of $\nu$. We however note that there are other ways of defining
the same measure $\nu$ on the boundary: Using the so-called ``exit local time''
under $\N^{*,z}_0$, one can define a simple continuous loop whose range 
is $\partial\Theta$, and the measure $\nu$
coincides with the occupation time measure of this loop (Corollary \ref{bdry-disk}). Finally, we observe that \cite[Section 4.2]{MS0} 
constructs a local time measure on the boundary of a filled metric ball of the Brownian map (these filled metric balls are
called hulls in the present work), which can be identified with our measure $\nu$ once we know that the complement of a filled
metric ball is a Brownian disk --- see Theorem \ref{ccBm} below. In fact, the construction of \cite{MS0} relies on the exit local time 
of the Brownian snake, which we also use to construct $\nu$ in Section \ref{sec:unibdry}. 

As an application of Theorem \ref{main}, we show that connected components of the
complement of balls in the Brownian map are Brownian disks. To state this in a precise form,
we recall that (under $\N_0$ or under $\N^{(1)}_0$), the Brownian map
$\mm$ has a distinguished point $x_*$ corresponding to the point $a_*$ of $\t_\zeta$
with minimal label. The interest of considering $x_*$ comes from the 
property $D(a_*,a)=Z_a-Z_{a_*}$ for every $a\in\t_\zeta$, showing that 
labels correspond to distances from $x_*$, up to the shift by $-Z_{a_*}$
(on the other hand, the invariance of the Brownian map under uniform re-rooting
\cite[Theorem 8.1]{Geo} shows that $x_*$ plays no special role in the Brownian map).
We use the notation $B(x_*,r)$ for the closed ball of radius $r$
centered at $x_*$ in $\mm$. We also recall that $\bv$ stands for the volume measure on $\mm$
(which is a probability measure under $\N^{(1)}_0$). If $O$ is an open subset of $\mm$,
we can define an {\it intrinsic} metric $d^{O}_{\mathrm{intr}}$ on $O$ by declaring that
$d^{O}_{\mathrm{intr}}(x,y)$ is the minimal length of a continuous path connecting $x$ to $y$ in $O$
(see \cite[Chapter 2]{BBI}). 

For every $z>0$ and $v>0$, we let $\F_{z,v}$ be the distribution of the Brownian disk with perimeter
$z$ and volume $v$ (see \cite{BM}).

\begin{theorem}
\label{ccBm}
Let 
$r>0$. Then, $\N^{(1)}_0$ a.s. for every connected component $\mathbf{C}$ of $\mm\backslash B(x_*,r)$, the 
limit
\begin{equation}
\label{approx-bd}
|\partial\mathbf{C}|:=\lim_{\ve\to 0} \ve^{-2} \int_\mathbf{C} \bv(\mathrm{d}x)\,\mathbf{1}_{\{ D(x,\partial\mathbf{C})<\ve\}}
\end{equation}
exists and is called the boundary size of $\mathbf{C}$. On the event $\{\mm\backslash B(x_*,r)\not =\varnothing\}$, write $\mathbf{C}^{r,1},\mathbf{C}^{r,2},\ldots$
for the connected components of $\mm\backslash B(x_*,r)$ ranked in decreasing order of their boundary sizes.
Let $d^{r,j}_{\mathrm{intr}}$ be the intrinsic distance on $\mathbf{C}^{r,j}$, and let $\bv^{r,j}$ be the 
restriction of $\bv$ to $\mathbf{C}^{r,j}$.
Then, $\N^{(1)}_0$ a.s. on the event $\{\mm\backslash B(x_*,r)\not =\varnothing\}$, for every $j=1,2,\ldots$,
the metric $d^{r,j}_{\mathrm{intr}}$ has a continuous extension to the closure
$\ov{\mathbf{C}}^{r,j}$ of $\mathbf{C}^{r,j}$, and this extension is a metric on $\ov{\mathbf{C}}^{r,j}$. Furthermore, 
under $\N^{(1)}_0(\cdot\mid \mm\backslash B(x_*,r)\not =\varnothing)$ and conditionally on  the sequence 
$$(|\partial \mathbf{C}^{r,1}|,\bv(\mathbf{C}^{r,1})), (|\partial \mathbf{C}^{r,2}|,\bv(\mathbf{C}^{r,2})),\ldots$$
the measure metric spaces $(\ov{\mathbf{C}}^{r,j},d^{r,j}_{\mathrm{intr}},\bv^{r,j})$, $j=1,2,\ldots$, are independent
Brownian disks with respective distributions $\F_{|\partial \mathbf{C}^{r,j}|,\bv(\mathbf{C}^{r,j})}$, $j=1,2,\ldots$.
\end{theorem}

As a side remark, we note that the distribution of the collection of boundary lengths $(|\partial\mathbf{C}^{r,j}|)_{j\geq 1}$ as a function of $r$
should be related to the growth-fragmentation process obtained in \cite{BCK} as the limiting process for the lengths of cycles bounding 
the connected components of the complement of balls in large random triangulations.

As another application of Theorem \ref{main}, we study the connected components of
the complement of the metric net of the Brownian map. We now argue under
the measure $\N_0$, so that we deal with the free Brownian map as defined above. We note that, in
addition to $x_*$, 
$\mm$ has another distinguished point, namely  $x_0$ that corresponds to the root of $\t_\zeta$. The metric net
$\mathbf{N}$ is then defined as the closure of
$$\bigcup_{0<r<D(x_*,x_0)} \partial\mathfrak{H}_r,$$
where,
for every $r\in (0,D(x_*,x_0))$,  $\mathfrak{H}_r$ denotes 
the connected component containing $x_0$ of the complement of the closed ball $B(x_*,r)$ (the complement of $\mathfrak{H}_r$ is the so-called hull of radius $r$). The metric net plays an important role in the axiomatic characterization
of the Brownian map which is presented in \cite{MS0}. Write $\mathbf{D}^{(1)},\mathbf{D}^{(2)},\ldots$ for the sequence
of connected components of $\mm\backslash \mathbf{N}$, which
can be ordered by decreasing size of the boundaries (for every $j$, 
the size $|\partial \mathbf{D}^{(j)}|$ may be defined  by an approximation exactly similar to
\eqref{approx-bd}). For every
$j=1,2,\ldots$, let $\bv^{(j)}$ be the restriction  of the 
volume measure on $\mm$ to $\mathbf{D}^{(j)}$.
Then Theorem \ref{Bnet} below
shows that, for every $j$, the intrinsic metric $d^{(j)}_{\rm intr}$ on $\mathbf{D}^{(j)}$ has a continuous extension 
to $\ov{\mathbf{D}}^{(j)}$ and, conditionally on the sequence $(|\partial \mathbf{D}^{(j)}|,j=1,2,\ldots)$, the 
measure metric spaces $(\ov{\mathbf{D}}^{(j)},d^{(j)}_{\rm intr},\bv^{(j)})$, $j=1,2,\ldots$ are
independent free Brownian disks with respective perimeters $|\partial \mathbf{D}^{(j)}|$, $j=1,2,\ldots$. 
In particular, this makes it possible to identify the probability measure $\mu^L_{\rm DISK}$
introduced in \cite[Section 4.2]{MS0} with the law of the free Brownian disk with perimeter $L$. 

\medskip
Let us sketch the main steps of the proof of Theorem \ref{main}, which motivates much of the subsequent developments. We start from a (pointed) Boltzmann quadrangulation,
that is, a random rooted and pointed planar quadrangulation $Q$ defined under a probability measure $\P$, such that $\P(Q=q)=c\,12^{-k}$
for any fixed rooted and pointed planar quadrangulation $q$ with $k$ faces, with some normalizing constant $c>0$. Write $\dg$ for the graph distance on the vertex set of $Q$,
and note that this vertex set has two
distinguished elements, namely the root vertex $\rho$ and the distinguished vertex $\xi$.
Given $\delta>0$, it is not too hard 
to verify that the distribution of the metric space associated with $Q$ under $\P(\cdot\mid \dg(\rho,\xi)>\delta\sqrt{n})$, with distances rescaled by $\sqrt{3/2n}$,
converges as $n$ tends to infinity to the distribution of the free Brownian map under $\N_0(\cdot\mid W_*<-\delta\sqrt{3/2})$, where $W_*$ stands for the minimal label under $\N_0$
(see Corollary \ref{convBmapB} below).

The idea now is that we can ``embed'' a quadrangulation with a boundary in $Q$. To this end, still under
the conditional probability $\P(\cdot\mid \dg(\rho,\xi)>\delta\sqrt{n})$, we consider the hull of radius 
$\delta\sqrt{n}$ in $Q$, which roughly speaking contains those vertices of $Q$ that cannot be connected 
to $\xi$ by a path staying at distance greater than $\delta\sqrt{n}$ from $\rho$. It turns out that the quadrangulation 
$Q$ may be obtained by gluing to this hull a Boltzmann quadrangulation $Q_n$ with a boundary of size equal to the 
perimeter of the hull (the definition of a Boltzmann quadrangulation with a boundary
is similar to the above definition in the case without boundary). For the latter statement to hold (Proposition \ref{glue-Boltzmann}), we need to use the ``lazy hull'' whose
definition is recalled in Section \ref{sec:lazy}. Furthermore, in the limiting result connecting $Q$ under $\P(\cdot\mid \dg(\rho,\xi)>\delta\sqrt{n})$ to the free Brownian 
map, the boundary size of $Q_n$ (normalized by $n^{-1}$) converges to the exit measure at level $W_*+\delta\sqrt{3/2}$
(counting ``how many'' Brownian snake paths exit the interval $(W_*+\delta\sqrt{3/2},\infty)$). See Proposition \ref{conv-exit}
in Section \ref{subsec:asBol}.

We now take $\delta$ small but simultaneously condition on the event that the boundary size of $Q_n$
is close to $n$ (say between $(1-\ve)n$ and $(1+\ve)n$). We note that the graph distance (computed in $Q$) between any
two points of the boundary of $Q_n$ is smaller than $2\delta\sqrt{n}$ and thus will be small after rescaling by $\sqrt{3/2n}$. It follows that,
under the preceding conditioning, the (rescaled) metric space associated with $Q_n$ after
gluing its boundary in an appropriate sense
is close to the (rescaled) metric space associated with $Q$, which by the convergence
result for planar quadrangulations is close to the free Brownian map under the special conditioning 
that the exit measure at level $W_* + \delta\sqrt{3/2}$ is close to $1$ (see Section \ref{conv-glued-bdry} and 
in particular Proposition \ref{conv-truncated}).

We can analyse the behavior of the Brownian snake
under the latter conditioning and relate this behavior to the measure $\N^*_0$, thanks to the re-rooting representation theorem of \cite{ALG} (see
Theorem \ref{re-root-rep} and Proposition \ref{law-trunc} below). This indicates that the scaling limit of Boltzmann quadrangulations 
with a glued boundary can be expressed 
in terms of the Brownian snake under $\N^{*}_0$, which eventually leads to the representation of the Brownian disk
with glued boundary (Theorem \ref{identif-glued-unpointed}) that was already mentioned before Theorem \ref{main}. This representation in the glued case is simpler than the
one of Theorem \ref{main} as it involves only the (pseudo-)metric $D$ defined in \eqref{formulaD-i} and not the (pseudo-)metric $\Delta$ in Theorem \ref{main}.

On the other hand, to get the statement of Theorem \ref{main}, it is crucial to reconstruct the real Brownian disk metric from the one corresponding to the glued boundary
case (Section \ref{sec:cons-metric}). This is done by proving
that the pseudo-metric $\Delta$ of Theorem \ref{main} has a 
continuous extension to the boundary (Proposition \ref{extenD}, in combination with the
re-rooting representation of $\N^*_0$).  The idea of the proof
of Theorem \ref{main} in Section \ref{sec:ident} is then to observe that the collection of mutual 
distances in a sequence of independent uniformly distributed points has 
the same distribution in 
$\Theta$ and in the free Brownian disk. To this end, we rely on a result of
\cite{BM} stating that the geodesic path between two ``typical points'' of
a Brownian disk does not hit the boundary. It follows that we
can recover the distance between two such points from the information 
given by distances in the object with glued boundary, and this is basically what we need
to identify $\Delta$ with the Brownian disk metric.

\medskip

The present article is organized as follows. Section \ref{GHP} is devoted to some
preliminaries about the Gromov-Hausdorff-Prokhorov convergence for pointed measure metric spaces.
Section \ref{sec:Crq} recalls the convergence to the Brownian map for rescaled planar 
quadrangulations. For our purposes, it is convenient to define the Brownian snake
excursion measure within the framework of snake trajectories, which was first discussed in \cite{ALG}. In
particular, we explain how to associate a measure metric space with a snake trajectory in a deterministic setting
(Section \ref{sec:consmms}). After recalling the classical theorem of convergence to the
Brownian map for quadrangulations with a fixed size, in a form suitable for our applications, we 
derive the similar result for Boltzmann quadrangulations (with a random size), 
with convergence to the free Brownian map under a particular conditioning (Corollary \ref{convBmapB}).

In Section \ref{sec:quadbBd}, we discuss quadrangulations with a boundary 
and Brownian disks. We recall from \cite{BM} the main result of convergence of rescaled quadrangulations
with a boundary to the free Brownian disk, again
in a slightly more precise form in view of our applications. We also explain why this statement implies 
a similar convergence result for quadrangulations and Brownian disks with glued boundary.
Section \ref{sec:resu-snake} is devoted to some preliminary results about the Brownian snake truncated 
at a level of the form $W_* + \delta$. We
also recall from \cite{ALG} the definition of the positive excursion measure $\N^*_0$ and 
its re-rooting representation in terms of truncated snakes (Theorem \ref{re-root-rep}),
which plays an important role in some of the subsequent proofs.

The goal of Section \ref{sec:lazy} is to prove
the already mentioned convergence of the (rescaled) perimeter of the lazy hull of radius $\delta\sqrt{n}$ 
in a Boltzmann quadrangulation conditioned on the event $\dg(\rho,\xi)>\delta\sqrt{n}$. To verify that this convergence holds jointly with the 
convergence in distribution to the free Brownian map, which is crucial for our purposes, we rely on previous
results \cite{CLG-peeling} about the convergence of the perimeter and volume of
hulls in large random planar maps. 

The main objective of Section \ref{conv-glued-bdry} 
is to prove Theorem \ref{identif-glued-unpointed}, which
identifies the free Brownian disk with glued boundary (of perimeter $z$) to
the random metric space constructed under $\N^{*,z}_0$ from the pseudo-metric $D$
defined in \eqref{formulaD-i}. The general idea of the argument is the one presented in the preceding lines,
but a few technicalities are needed. 
The goal of Section \ref{sec:cons-metric} is then to prove the first assertion of Theorem \ref{main}, namely
the fact that the pseudo-metric $\Delta$ has a 
continuous extension to the boundary. The other assertions of Theorem \ref{main} are proved 
in Section \ref{sec:ident}. 

Section \ref{sec:unibdry} investigates the uniform measure on the boundary and
proves Proposition \ref{measure-boundary} in a more precise form. Here we rely on an approximation result for the
exit local time of the Brownian snake (Proposition \ref{approx-LT}). Roughly speaking,
this approximation result provides the needed asymptotics under $\N^*_0$ for the measure of the set of times where $\wh W_s<\ve$, which corresponds to the 
quantity $\int \bV(\mathrm{d}x)\,\mathbf{1}_{\{\Delta(x,\partial\Theta)<\ve\}}$ in our representation of the Brownian disk. 

Section \ref{sec:Bnet}
proves that the connected components of the complement of the
Brownian net are independent Brownian disks (Theorem \ref{Bnet}). Here we rely on results of  \cite{subor} showing that these 
connected components are in one-to-one correspondence with excursions of the Brownian snake above its minimum, as defined in
\cite{ALG}. Since it is shown in \cite{ALG} that the latter excursions are independent and distributed according to (conditional versions of)
the measure $\N^*_0$, the proof boils down to verifying that the intrinsic metric on each component coincides 
with the metric on the Brownian disk associated with the corresponding excursion.  The idea of the proof of Theorem \ref{ccBm}
in  Section \ref{sec:BdBm}
is similar, using now results of \cite{ALG} for Brownian snake excursions above $0$. There are however additional technical difficulties
because Theorem \ref{ccBm} is concerned with the measure $\N^{(1)}_0$ (that is, we consider the standard Brownian map
instead of the free Brownian map under $\N_0$), and we really want to consider Brownian snake excursions above 
the random level $W_* + r$ rather than above $0$. A key idea to circumvent this last difficulty is to use the 
invariance of $\N^{(1)}_0$ under re-rooting (formula \eqref{uniform-re-root} in Section \ref{Brsnexc}). 

Finally, the 
appendices give the proofs of two technical results.

\tableofcontents

\bigskip\medskip
{\parindent=0pt 
{\large\bf Main notation}

\medskip
\begin{tabular}{l@{\ }l}
$d^{(k)}_{GHP}$ &Gromov-Hausdorff-Prokhorov distance on $k$-pointed measure metric spaces (Sect.~\ref{GHP})\cr
$\M$ ($\M^{k\bullet}$)& space of all ($k$-pointed) measure metric spaces (Sect.~\ref{GHP})\cr
$\Q_n$ &set of all rooted planar quadrangulations with $n$ faces (Sect.~\ref{sec:quadr})\cr
$\rho$ &root vertex of $q\in \Q_n$ (Sect.~\ref{sec:quadr})\cr
$V(q)$ &vertex set of $q\in \Q_n$ (Sect.~\ref{sec:quadr})\cr
$\dg$ &graph distance on $V(q)$ (Sect.~\ref{sec:quadr})\cr
$\mu$ &counting measure on $V(q)$ (Sect.~\ref{sec:quadr})\cr
$\Q_n^\bullet$ &set of all rooted and pointed planar quadrangulations with $n$ faces (Sect.~\ref{sec:Bmap})\cr
$\xi$ &distinguished vertex of $q\in \Q_n^\bullet$ (Sect.~\ref{sec:Bmap})\cr
$\Q^{\partial,k}$ ($\Q^{\partial,k}_n$) &set of all pointed quadrangulations with a boundary of size $2k$ (Sect.~\ref{subsec:quadb})\cr
$(C^q_k)_{0\leq k\leq 2n}$&contour function of $q\in \Q_n$ (Sect.~\ref{sec:quadr})\cr
$(L^q_k)_{0\leq k\leq 2n}$&label function of $q\in \Q_n$ (Sect.~\ref{sec:quadr})\cr
$Q$&Boltzmann (rooted and pointed) planar quadrangulation (Sect.~\ref{sec:Bolquad})\cr
$B_{(k)}$&Boltzmann quadrangulation with a boundary of size $2k$ (Sect.~\ref{subsec:quadb})\cr
$(B_{(k)}^\dagger,\dg^\dagger)$&Boltzmann quadrangulation with ``glued'' boundary of size $2k$ (Sect.~\ref{sec:convBdisk})\cr
\end{tabular}

\begin{tabular}{l@{\ }l}
\noalign{\smallskip}
$\mathcal{W}$ &space of all one-dimensional stopped paths $\w=(\w(t))_{0\leq t\leq \zeta_{(\w)}}$ (Sect.~\ref{sec:snake})\cr
$\wh \w=\w(\zeta_{(\w)})$ &terminal point of $\w\in\mathcal{W}$ (Sect.~\ref{sec:snake})\cr
$\tau_b(\w)$&$=\inf\{t\geq 0: \w(t)=b\}$ (Sect.~\ref{sec:snake})\cr
$\underline{\w}$&minimal value of $\w\in\mathcal{W}$\ (Sect.~\ref{sec:snake})\cr
$\S$ &space of all snake trajectories (Sect.~\ref{sec:snake})\cr
$\sigma(\omega)$ &duration of a snake trajectory $\omega\in\S$ (Sect.~\ref{sec:snake})\cr
$(W_s)_{s\geq 0}$,$(\zeta_s)_{s\geq 0}$ &canonical process and lifetime process ($\zeta_s=\zeta_{(W_s)}$) on $\S$ (Sect.~\ref{sec:snake})\cr
$\t_\zeta$ &genealogical tree of a snake trajectory (tree coded by $(\zeta_s)_{s\geq 0}$) (Sect.~\ref{sec:snake})\\
$d_\zeta$ &metric on $\t_\zeta$ (Sect.~\ref{sec:snake})\cr
$p_\zeta$ &canonical projection from $[0,\sigma]$ onto $\t_\zeta$ (Sect.~\ref{sec:snake})\cr
$Z_a$ &labels on the tree $\t_\zeta$ (Sect.~\ref{sec:snake})\cr
$D(a,b)$ &(and $D^\circ(a,b)$) metric on $\t_\zeta$ defined using labels (Sect.~\ref{sec:snake})\cr
$\omega^{[s]}$ &snake trajectory $\omega$ re-rooted at $s\in [0,\sigma(\omega)]$ (Sect.~\ref{sec:snake})\cr
$\mathrm{tr}_y(\omega)$ &snake trajectory $\omega$ truncated at $y$ (Sect.~\ref{sec:snake})\cr
$(\mm,D,\bv)$ &measure metric space associated with a snake trajectory (Sect.~\ref{sec:consmms})\cr
$\Pi$&canonical projection from $\t_\zeta$ onto $\mm$ (Sect.~\ref{sec:consmms})\cr
$\ll,\ll^\bullet,\ll^{\bullet,\bullet}$ &mappings from snake trajectories to measure metric spaces
(Sect.~\ref{sec:consmms})\cr
\noalign{\smallskip}
$\N_x,\N_x^{(s)}$&Brownian snake excursion measures (Sect.~\ref{Brsnexc})\cr
$\N_x^{[y]}$&for $x>y$, $\N_x$ conditioned on $\{W_*< y\}$ (Sect.~\ref{sec:Bstrunc}, \ref{sec:cons-metric})\cr
$W_*,W^*$& minimal and maximal values of $\wh W_s$, $0\leq s\leq \sigma$ (Sect.~\ref{Brsnexc})\cr
$s_*,a_*$&$s_*$ time of the minimum of $s\mapsto \wh W_s$, $a_*=p_\zeta(s_*)$ (Sect.~\ref{Brsnexc}, \ref{sec:defBM})\cr
$x_*,x_0$&$x_*=\Pi(a_*)$, $x_0=\Pi(p_\zeta(0))$ (Sect. \ref{sec:defBM})\cr
$(\ell^b_s)_{s\geq 0}$&exit local time from $(b,\infty)$ under $\N_x$, $x>b$  (Sect.~\ref{sec:Bstrunc})\cr
$\z_b$ &exit measure from $(b,\infty)$ under $\N_x$, $x>b$ (Sect.~\ref{sec:Bstrunc})\cr
$\N^*_0$&positive excursion measure of the Brownian snake (Sect.~\ref{sec:posexc})\cr
$\z^*_0$&exit measure at $0$ under $\N^*_0$ (Sect.~\ref{sec:posexc})\cr
$\N^{*,z}_0$&$\N^*_0$ conditioned on $\{\z^*_0=z\}$ (Sect.~\ref{sec:posexc})\cr
\noalign{\smallskip}
$(\D^\bullet_r,D^\partial,\bv^\bullet_r)$&free pointed Brownian disk with perimeter $r$ (Sect.~\ref{sec:freeBdisk})\cr
$\F_r$ ($\F^\bullet_r$)&distribution of the free (pointed) Brownian disk with perimeter $r$ (Sect.~\ref{sec:freeBdisk})\cr
$(\D^{\bullet,\dagger}_r,D^{\partial,\dagger},\bv^{\bullet,\dagger}_r)$&free pointed Brownian disk with glued boundary and perimeter $r$ (Sect.~\ref{sec:Bdglued})\cr
$\F_r^\dagger$ ($\F_r^{\bullet,\dagger}$)&distribution of the free (pointed) Brownian disk with glued boundary (Sect.~\ref{sec:Bdglued})\cr
$\QQ_1,\ldots,\QQ_\Xi$&sequence of lazy hulls of $Q$ (Sect.~\ref{sec:peel-layers})\cr
$Q^{(\infty)}$ &UIPQ (uniform infinite planar quadrangulation) (Sect.~\ref{sec:peel})\cr
$H_i$,$V_i$, $H^{(\infty)}_i$,$V^{(\infty)}_i$&half-perimeter and volume of the lazy hull of radius $i$ of $Q$ or $Q^{(\infty)}$ (Sect.~\ref{sec:asymp-hull})\cr
$S^\downarrow_k$ or $S^\uparrow_k$&half-perimeter in the peeling algorithm of $Q$ or $Q^{(\infty)}$ (Sect.~\ref{sec:asymp-hull})\cr
$Q_n$&quadrangulation ``filling in'' the lazy hull of $Q$ of radius $\lfloor\delta\sqrt{n}\rfloor$ (Sect.~\ref{conv-glued-bdry})\cr
\noalign{\smallskip}
$\Delta(a,b)$ &(and $\Delta^\circ(a,b)$) Brownian disk metric under $\N^{[0]}_r$ or under $\N^*_0$ (Sect.~\ref{sec:cons-metric})\cr
$(\Theta,\Delta,\bV)$&Brownian disk constructed under $\N^*_0$ (Sect.~\ref{sec:cons-metric}, \ref{sec:ident})\cr
\end{tabular}

}
\section{Preliminaries on Gromov-Hausdorff-Prokhorov convergence}
\label{GHP}

\subsection{Gromov-Hausdorff-Prokhorov convergence and correspondences}

A (compact) measure metric space is a compact metric space $(X,d)$ equipped with a Borel
finite measure $\mu$. We write $\M$ for the set of all measure metric spaces, where two such spaces
$(X,d,\mu)$ and $(X',d',\mu')$ are identified if there exists an isometry $\phi$ from $X$ onto 
$X'$ such that $\phi_*\mu=\mu'$. 

The Gromov-Hausdorff-Prokhorov distance on $\M$ is then defined by
$$d_{\mathrm{GHP}}^{(0)}((X,d,\mu),(X',d',\mu'))
= \inf_{\phi:X\to E,\phi':X'\to E} \Big\{ d^E_\mathrm{H}(\phi(X),\phi'(X)) \vee d^E_\mathrm{P}(\phi_*\mu,\phi'_*\mu')\Big\},$$
where the infimum is over all isometric embeddings $\phi$ and $\phi'$ of $X$ and $X'$ into a compact
metric space $(E,d^E)$, $d^E_\mathrm{H}$ stands for the Hausdorff distance between compact subsets of $E$ and
$d^E_\mathrm{P}$ is the Prokhorov distance on the space of finite measures on $E$. Writing $C^\ve$ for
the (closed) $\ve$-enlargement of a closed subset $C$ of $E$, 
$$d^E_\mathrm{P}(\nu,\nu')=\inf\{\ve >0: \nu(C)\leq \nu'(C^\ve)+\ve\hbox{ and }\nu'(C)\leq \nu(C^\ve)+\ve\hbox{ for every closed
subset }C\hbox{ of }E\}.$$
According to \cite[Theorem 2.5]{AD}, $d_{\mathrm{GHP}}^{(0)}$ defines a metric on $\M$, and 
$(\M,d_{\mathrm{GHP}}^{(0)})$ is a Polish space. 

We will need to consider $k$-pointed measure metric spaces, for every integer $k\geq 0$. A $k$-pointed 
measure metric space is a measure metric space equipped with $k$ distinguished points $\rho_1,\ldots,\rho_k$
(the order of these points is important). We write $\M^{k\bullet}$ for the space of all (equivalence classes of)
$k$-pointed measure metric spaces. We can similarly define the 
Gromov-Hausdorff-Prokhorov distance on $\M^{k\bullet}$ by
\begin{align*}
&d_{\mathrm{GHP}}^{(k)}((X,d,\mu,(\rho_1,\ldots,\rho_k)),(X',d',\mu',(\rho'_1,\ldots,\rho'_k)))\\
&= \inf_{\phi:X\to E,\phi':X'\to E} \Big\{ d^E_\mathrm{H}(\phi(X),\phi'(X)) \vee d^E_\mathrm{P}(\phi_*\mu,\phi'_*\mu')
\vee \max_{1\leq i\leq k} d^E(\phi(\rho_i),\phi'(\rho'_i))\Big\},
\end{align*}
where, as previously, the infimum is over all isometric embeddings $\phi$ and $\phi'$ of $X$ and $X'$ into a compact
metric space $(E,d^E)$. Again, using the same arguments as in \cite{AD}, one verifies that $(\M^{k\bullet},d_{\mathrm{GHP}}^{(k)})$ is a Polish space. 

It will be convenient to have a bound on the 
Gromov-Hausdorff-Prokhorov distance in terms of 
correspondences. We consider
two $k$-pointed measure metric spaces 
$(X,d,\mu,(\rho_1,\ldots,\rho_k))$ and $(X',d',\mu',(\rho'_1,\ldots,\rho'_k))$, and write $\pi$, resp. $\pi'$, for the
canonical projection from $X\times X'$ onto $X$, resp. onto $X'$. Recall that a correspondence 
between $X$ and $X'$ is a subset $\rr$ of $X\times X'$ such that
$\pi(\rr)=X$ and $\pi'(\rr)=X'$. The distortion of the correspondence $\rr$ is
$$\mathrm{dis}(\rr):=\sup\{|d(x,y)-d'(x',y')|: (x,x')\in\rr, (y,y')\in\rr\}.$$

\begin{lemma}
\label{GHP-corrresp}
Let $\ve>0$ and suppose that there is a correspondence $\rr$ between $X$ and $X'$
with distortion bounded above by $\ve$, such that
$(\rho_i,\rho'_i)\in \rr$ for every $1\leq i\leq k$, and a finite measure $\nu$ on the product 
$X\times X'$ such that $\nu(\rr^c)<\ve$ and
$$d^X_\mathrm{P}(\pi_*\nu,\mu)<\ve\;,\quad d^{X'}_\mathrm{P}(\pi'_*\nu,\mu')<\ve.$$
Then $d_{\mathrm{GHP}}^{(k)}((X,d,\mu),(X',d',\mu'))\leq 3\ve$.
\end{lemma}

\proof Let us only sketch the proof (see \cite[Section 6]{Mie} for very similar arguments). 
We equip the disjoint union $X\sqcup X'$ with the metric $\delta(\cdot,\cdot)$ such that the
restriction of $\delta$ to $X\times X$ is $d$ and the restriction of $\delta$ to $X'\times X'$
is $d'$, and for every $(x,x')\in X\times X'$,
$$\delta(x,x')=\inf\{d(x,y)+\ve+d(y',x'):(y,y')\in\rr\}.$$
We note that $\delta(x,x')=\ve$ if $(x,x')\in\rr$. We then apply the definition of 
$d^{(k)}_{\mathrm{GHP}}$, letting $\phi$, respectively $\phi'$,
be the identity mapping from $X$, respectively $X'$, into $E:=X\sqcup X'$, which is equipped with
$d^E=\delta$. 
It is clear that $d^E_\mathrm{H}(\phi(X),\phi'(X'))\leq \ve$ and $d^E(\phi(\rho_i),\phi'(\rho'_i))=\ve$
for every $1\leq i\leq k$. So we need only check that $d^E_\mathrm{P}(\phi_*\mu,\phi'_*\mu')\leq 3\ve$. 
To this end, we note that we can view $\pi_*\nu$ and $\pi'_*\nu$ as measures on $E$, and that
our assumption gives $d^E_\mathrm{P}(\phi_*\mu,\pi_*\nu)<\ve$ and 
$d^E_\mathrm{P}(\phi'_*\mu,\pi'_*\nu)<\ve$, so that it suffices
to verify that $d^E_\mathrm{P}(\pi_*\nu,\pi'_*\nu)\leq \ve$. 

Let $\mathbf{C}$ be a closed subset of $E$
and $C=\mathbf{C}\cap X$. As previously, we write $\mathbf{C}^\ve$, resp. $C^\ve$, for the
$\ve$-enlargement of $\mathbf{C}$, resp. of $C$, in $E$. Then $\pi_*\nu(\mathbf{C})=\pi_*\nu(C)=\nu(C\times X')\leq \nu((C\times X')\cap \rr)+\ve$. 
However, if $(x,x')\in (C\times X')\cap \rr$, the fact that $\delta(x,x')=\ve$ shows that $x'\in C^\ve$, and
thus $(x,x')\in X\times (C^\ve\cap X')$. It follows that $ \nu((C\times X')\cap \rr)\leq \nu(X\times (C^\ve\cap X'))=\pi'_*\nu(C^\ve\cap X')$.
We have thus obtained $\pi_*\nu(\mathbf{C})\leq \pi'_*\nu(C^\ve\cap X')+\ve\leq \pi'_*\nu(\mathbf{C}^\ve)+\ve$. The same argument gives 
$\pi'_*\nu(\mathbf{C})\leq  \pi_*\nu(\mathbf{C}^\ve)+\ve$, and we
conclude that
$d^E_\mathrm{P}(\pi_*\nu,\pi'_*\nu)\leq \ve$ as desired. \endproof

We note that $\M^{0\bullet}=\M$ and in what follows, we write $\M^\bullet$ instead of $\M^{1\bullet}$ and $\M^{\bullet\bullet}$
instead of $\M^{2\bullet}$. 

\subsection{Measure metric spaces equipped with a uniformly distributed infinite sequence}
\label{sec:infinite}

Let $\M^{\infty\bullet}$ be the set of all (equivalence classes 
modulo measure-preserving isometries of) measure metric spaces $(X,d,\mu)$ equipped with a sequence
$(x_n)_{n\geq 1}$ of points of $X$. For every $k\geq 0$, there is an obvious 
projection $\mathbf{P}_k$ from $\M^{\infty\bullet}$ onto $\M^{k\bullet}$.
We equip $\M^{\infty\bullet}$ with the smallest $\sigma$-field 
for which all projections $\mathbf{P}_k$ are measurable. Note that a mapping
$\varphi$ with values in $\M^{\infty\bullet}$ is then measurable
if and only if all mappings $\mathbf{P}_k\circ \varphi$ are measurable.

Let $(X,d,\mu)\in\M$, with $\mu\not =0$. Consider a sequence $(U_n)_{n\geq 1}$ of
i.i.d. random variables with values in $X$ whose common distribution is $\mu/\mu(X)$. Then,
$(X,d,\mu,(U_n)_{n\geq 1})$ is a random variable with values in $\M^{\infty\bullet}$, and we write
$\mathbf{Q}((X,d,\mu),\cdot)$ for its distribution. By adapting the arguments in \cite[Lemma 13]{Mie},
one obtains that the mapping $(X,d,\mu)\mapsto \mathbf{Q}((X,d,\mu),A)$
is measurable whenever $A$ is of the form $\mathbf{P}_k^{-1}(B)$
for some measurable subset $B$ of $\M^{k\bullet}$. By standard
monotone class arguments, the same mapping is measurable whenever
$A$ is a measurable subset of $\M^{\infty\bullet}$.
Hence $\mathbf{Q}(\cdot,\cdot)$ defines a kernel from $\{(X,d,\mu)\in\M:\mu\not=0\}$ to $\M^{\infty\bullet}$.

Suppose now that $(\mathcal{X},\mathcal{D},\theta)$ is a random variable with values in $\M$, with $\theta\not = 0$ a.s.
We can then consider a random variable $(\mathcal{X},\mathcal{D},\theta,(V_n)_{n\geq 1})$ with values
in $\M^{\infty\bullet}$ such that,
conditionally on $(\mathcal{X},\mathcal{D},\theta)$, $(\mathcal{X},\mathcal{D},\theta,(V_n)_{n\geq 1})$ is distributed according to
$\mathbf{Q}((\mathcal{X},\mathcal{D},\theta),\cdot)$. In this way, we can make sense of choosing a sequence
of  independent uniformly distributed points in a random metric measure space. 

\begin{lemma}
\label{infinite-points}
Let $(\mathcal{X},\mathcal{D},\theta)$ be a random variable with values in $\M$, such that the 
topological support of $\theta$ is a.s. equal to $\mathcal{X}$, and let $(V_n)_{n\geq 1}$
be a sequence
of  independent uniformly distributed points in $\mathcal{X}$. Then the random measure metric spaces
$$\Big(\{V_1,V_2,\ldots,V_n\},\mathcal{D},\frac{1}{n}\sum_{i=1}^n \delta_{V_i}\Big)$$
converge to $(\mathcal{X},\mathcal{D},\theta/\theta(\mathcal{X}))$ a.s. in $\M$ as $n\to \infty$.
\end{lemma}

\proof It is enough to verify the convergence when $(\mathcal{X},\mathcal{D},\theta)=(X,d,\mu)$ is 
deterministic and such that $\mu$ has full topological support. Then, a.s.
$$\{V_1,V_2,\ldots,V_n\} \build{\la}_{n\to\infty}^{} X$$
in the sense of the Hausdorff distance, and
$$\frac{1}{n}\sum_{i=1}^n \delta_{V_i}  \build{\la}_{n\to\infty}^{}\frac{\mu}{\mu(X)},$$
in the sense of weak convergence of probability measures,
as an application of the law of large numbers. The desired result 
follows.
\endproof

\section{Convergence of rescaled quadrangulations}
\label{sec:Crq}

\subsection{Quadrangulations and Schaeffer's bijection}
\label{sec:quadr}

For every integer $n\geq 1$, we let $\Q_n$ stand for the set of all rooted planar quadrangulations
with $n$ faces, and, if $q\in \Q_n$, we write $|q|=n$. We recall that
\begin{equation}
\label{nbquad}
\#\Q_n=\frac{2}{n+2}\cdot 3^n\,c_n
\end{equation}
where $c_n$ is the $n$-th Catalan number.
The root vertex (tail of the root edge) of a quadrangulation 
$q\in \Q_n$ will be denoted by $\rho$. 
The graph distance 
on the vertex set $V(q)$ of $q$ is denoted by $\dg^q$, or simply $\dg$ if there is no risk of confusion. We will also use the notation $\mu_q$,
or simply $\mu$, for the
counting measure on $V(q)$. 

For every integer $n\geq 1$, Schaeffer's bijection (see \cite[Section 3]{CS}) 
is a bijection from $\Q_n$ onto the set $ \T_n$, where $\T_n$ is the set of all well-labeled plane trees with $n$ edges.
Here a well-labeled plane tree is a pair $(\tau,(\ell_u)_{u\in V(\tau)})$, where 
$\tau$ is a (rooted) plane tree, with vertex set $V(\tau)$, and $(\ell_u)_{u\in V(\tau)}$ is a collection
of integer labels assigned to the vertices of $\tau$, in such a way that labels are positive integers, the root of $\tau$ has 
label $1$ and $|\ell_u-\ell_v|\leq 1$ if the vertices $u$ and $v$ of $\tau$ are adjacent. 
It will be convenient to write $\varnothing$ for the root of the tree $\tau$. 

Occasionally, we will also use the notion of a labeled plane tree: This is
a pair $(\tau,(\ell_u)_{u\in V(\tau)})$ satisfying the same conditions as above, except that
we drop the positivity conditions on labels, and the label of the root is $0$
instead of $1$.

Although we will not need the details of Schaeffer's bijection, we record a few useful facts.
Let $q\in \Q_n$, and let $(\tau,(\ell_u)_{u\in V(\tau)})$ be the well-labeled plane tree 
corresponding to $q$ via Schaeffer's bijection. Then $V(q)\backslash \{\rho\}$ is canonically 
identified with $V(\tau)$, in such a way that 
the endpoint of the root edge of $q$ is identified to the root $\varnothing$ of $\tau$.  Using this identification we have $\dg^q(\rho,u)=\ell_u$
for every $u\in V(q)\backslash \{\rho\}$. 

Let $u$ and $v$ be two vertices of $V(\tau)$. Let $k$ be the minimal label on the geodesic between $u$ and $v$ 
in the tree $\tau$. Then  any path in $q$ between $u$ and $v$ (now viewed as vertices of $q$) must visit a vertex with label
$k$. This is the ``discrete cactus property'', see in particular Proposition 4.3 in \cite{CLM} in a more general setting.
Conversely, it is not hard to construct a path in $q$ between $u$ and $v$ that visits only vertices with labels greater than
or equal to $k-1$ (with the convention that the label of $\rho$ is $0$). We leave the proof as an exercise for the reader (observe that two adjacent vertices of 
$V(\tau)$ are connected by a path of $q$ of length at most $2$).

Notice that $\tau$
can be viewed as a genealogical tree in a population whose ancestor is the root $\varnothing$ and that the ``children'' of
every individual are ordered. We
use the notation $|u|$ for the generation (graph distance from $\varnothing$ in $\tau$) of a vertex $u\in V(\tau)$. The contour sequence of $\tau$ 
is the finite sequence $u_0,u_1,\ldots,u_{2n}$ of vertices of $\tau$ defined inductively as follows. First $u_0=\varnothing$ is the 
root of $\tau$, and then, for every $0\leq i\leq 2n-1$, $u_{i+1}$ is either the first child of $u_i$ that does not 
appear in $\{u_0,\ldots,u_i\}$, or, if there is no such child, the parent of $u_i$. We can then define
the contour function $(C^q_k)_{k\geq 0}$ and the label function $(L^q_k)_{k\geq 0}$ of the well-labeled tree $(\tau,(\ell_u)_{u\in V(\tau)})$
by setting 
$$C^q_k=|u_k|\;,\quad L^q_k=\ell_{u_k}\;,\qquad \hbox{if }0\leq k\leq 2n,$$
and $C^q_k=L^q_k=0$ if $k>2n$. 

Let $Q_{(n)}$ be distributed uniformly over $\Q_n$. We recall the following useful estimate. There exist
positive constants $K_1$ and $K_2$, which do not depend on $n$, such that, for every $y>0$,
\begin{equation}
\label{estimmin}
\P\Big(\max_{k\geq 0} L^{Q_{(n)}}_k > y\Big)\leq K_1\,\exp(-K_2 \,\frac{y}{n^{1/4}}).
\end{equation}
Note that the left-hand side in \eqref{estimmin} is also equal to
$$\P\Big(\max_{u\in V(Q_{(n)})} \dg^{Q_{(n)}}(\rho,u) > y\Big).$$
The bound in \eqref{estimmin} then follows from \cite[Proposition 4]{CS} and the version of
Schaeffer's bijection for rooted and pointed quadrangulations (see, for instance, \cite[Section 5.4]{LGM}).

\subsection{Snake trajectories and the Brownian snake}
\label{sec:snake}

It will be convenient to use the formalism of snake trajectories, which has been introduced in \cite{ALG}.
Recall that a (one-dimensional) finite path $\w$ is just a continuous mapping $\w:[0,\zeta]\la\R$, where the
number $\zeta=\zeta_{(\w)}$ is called the lifetime of $\w$. We let 
$\W$ denote the space of all finite paths. The set $\W$ is a Polish space when equipped with the
distance
$$d_\W(\w,\w')=|\zeta_{(\w)}-\zeta_{(\w')}|+\sup_{t\geq 0}|\w(t\wedge
\zeta_{(\w)})-\w'(t\wedge\zeta_{(\w')})|.$$
The endpoint or tip of the path $\w$ is denoted by $\wh \w=\w(\zeta_{(\w)})$.
For every $x\in\R$, we set $\W_x=\{\w\in\W:\w(0)=x\}$. The trivial element of $\W_x$ 
with zero lifetime is identified with the point $x$ --- in this way we view $\R$
as the subset of $\W$ consisting of all finite paths with zero lifetime. 
Occasionally we will use the notation $\underline\w=\min\{\w(t):0\leq t\leq \zeta_{(\w)}\}$.

The following definition is taken from \cite{ALG}.

\begin{definition}
\label{def:snakepaths}
Let $x\in \R$. A snake trajectory with initial point $x$ is a continuous mapping
\begin{align*}
\omega:\ \R_+&\to \W_x\\
s&\mapsto \omega_s
\end{align*}
which satisfies the following two properties:
\begin{enumerate}
\item[\rm(i)] We have $\omega_0=x$ and the number $\sigma(\omega):=\sup\{s\geq 0: \omega_s\not =x\}$,
called the duration of the snake trajectory $\omega$,
is finite (by convention $\sigma(\omega)=0$ if $\{s\geq 0: \omega_s\not =x\}$ is empty). 
\item[\rm(ii)] For every $0\leq s\leq s'$, we have
$$\omega_s(t)=\omega_{s'}(t)\;,\quad\hbox{for every } 0\leq t\leq \min_{s\leq r\leq s'} \zeta_{(\omega_r)}.$$
\end{enumerate} 
We write $\S_x$ for the set of all snake trajectories  with initial point $x$, and 
$$\S:=\bigcup_{x\in \R} \S_x$$
for the set of all snake trajectories. If $\omega\in \S$, we write $W_s(\omega)=\omega_s$ and $\zeta_s(\omega)=\zeta_{(\omega_s)}$
for every $s\geq 0$. The set $\S$ is equipped with the distance
$$d_\S(\omega,\omega')= |\sigma(\omega)-\sigma(\omega')|+ \sup_{s\geq 0} \,d_\W(W_s(\omega),W_{s'}(\omega)).$$
\end{definition}

Let $\omega\in \S$ be a snake trajectory and $\sigma=\sigma(\omega)$. The lifetime function $s\mapsto \zeta_s(\omega)$ codes a
compact $\R$-tree, which will be denoted 
by $\t_\zeta$ and called the {\it genealogical tree} of the snake trajectory. This $\R$-tree is defined rigorously as the quotient space of the interval $[0,\sigma]$ 
for the equivalence relation
$$s\sim s'\ \hbox{if and only if }\ \zeta_s=\zeta_{s'}= \min_{s\wedge s'\leq r\leq s\vee s'} \zeta_r,$$
which is equipped with the distance induced by
$$d_\zeta(s,s')= \zeta_s+\zeta_{s'}-2 \min_{s\wedge s'\leq r\leq s\vee s'} \zeta_r.$$
(see e.g.~\cite[Section 3]{LGM} for more information about the
coding of $\R$-trees by continuous functions).  Let $p_\zeta:[0,\sigma]\la \t_\zeta$ stand
for the canonical projection. By convention, $\t_\zeta$ is rooted at 
$p_\zeta(0)=p_\zeta(\sigma)$, and the volume measure on $\t_\zeta$ is defined as the push forward of
Lebesgue measure on $[0,\sigma]$ under $p_\zeta$. The ancestral line of a vertex $a$ of $\t_\zeta$ is the line segment connecting 
$a$ to the root.

We observe that  $\omega$ is completely determined 
by the knowledge of the lifetime function $s\mapsto \zeta_s(\omega)$ and of the tip function $s\mapsto \wh W_s(\omega)$ (cf. \cite[Proposition 8]{ALG}).
To explain this, we note that, for every $s\in [0,\sigma]$, $\wh W_s$ only depends on $p_\zeta(s)$, and the 
values $W_s(t)$, $0\leq t\leq \zeta_s$ are recovered from the values of $\wh W$ along the ancestral line
of $p_\zeta(s)$ in $\t_\zeta$. 

It will be convenient to use the notation $Z_a(\omega)=\wh W_s(\omega)$ whenever $a\in\t_\zeta$
and $s\in[0,\sigma]$ are such that $a=p_\zeta(s)$. We interpret $Z_a$ as a ``label'' assigned to the ``vertex'' $a$ of $\t_\zeta$. 
Notice that the mapping $a\mapsto Z_a$ is continuous on $\t_\zeta$. We also set
\begin{align*}
W_*(\omega)&=\min\{\wh W_s(\omega):0\leq s\leq \sigma\}=\min\{Z_a:a\in\t_\zeta\}\;,\\
 W^*(\omega)&=\max\{\wh W_s(\omega):0\leq s\leq \sigma\}
=\max\{Z_a:a\in\t_\zeta\}\;,
\end{align*}
and often refer to $W_*(\omega)$ as the minimal label of $\omega$.

\smallskip
We will need to define {\it lexicographical}
intervals on the tree $\t_\zeta$. If $s,t\in[0,\sigma]$ and $s>t$, we abuse notation
by writing $[s,t]=[s,\sigma]\cup[0,t]$ (and of course if $s\leq t$, $[s,t]$ is the usual interval).
Then, if $a,b\in\t_\zeta$, there is a smallest ``interval'' $[s,t]$ with $s,t\in[0,\sigma]$,
such that $p_\zeta(s)=a$ and $p_\zeta(t)=b$, and we define $[a,b]=p_\zeta([s,t])$.
We also use the notation $]a,b[=[a,b]\backslash\{a,b\}$. 

\smallskip
Let us explain the {\it re-rooting operation} on snake trajectories of $\S_0$ (see \cite[Section 2.3]{LGW}
and \cite[Section 2.2]{ALG}). Let $\omega\in \S_0$ and
$r\in[0,\sigma(\omega)]$. Then $\omega^{[r]}$ is the new snake trajectory in $\S_0$ such that
$\sigma(\omega^{[r]})=\sigma(\omega)$ and for every $s\in [0,\sigma(\omega)]$,
\begin{align*}
\zeta_s(\omega^{[r]})&= d_\zeta(r,r\oplus s),\\
\wh W_s(\omega^{[r]})&= \wh W_{r\oplus s}-\wh W_r,
\end{align*}
where we use the notation $r\oplus s=r+s$ if $r+s\leq \sigma$, and $r\oplus s=r+s-\sigma$ otherwise. 
It will be convenient to write $\zeta_s^{[r]}(\omega)=\zeta_s(\omega^{[r]})$ and $W^{[r]}_s(\omega)=W_s(\omega^{[r]})$. 
The tree $\t_{\zeta^{[r]}}$ is then interpreted as the tree $\t_\zeta$ re-rooted at the vertex $p_\zeta(r)$: More precisely,
the mapping $s\mapsto r\oplus s$ induces an isometry from $\t_{\zeta^{[r]}}$
onto $\t_\zeta$, which maps the root of $\t_{\zeta^{[r]}}$ to $p_\zeta(r)$. Furthermore, the vertices
of $\t_{\zeta^{[r]}}$ receive the ``same'' labels as in $\t_\zeta$,
shifted so that the label of the root is still $0$.

\smallskip

The notion of the {\it truncation} of snake trajectories will also play an important role in this
work. Roughly speaking, if $\omega\in \S_x$ and $y<x$, the truncation of $\omega$ at $y$ is the new snake trajectory
$\omega'$ such that the values $\omega'_s$ are exactly the values $\omega_s$ for all $s$ such that $\omega_s$ does not hit $y$, or
hits $y$ for the first time at its lifetime. Let us give a more precise definition. First, for any $\w\in\W$ and $y\in\R$, we set
$$\tau_y(\w):=\inf\{t\in[0,\zeta_{(\w)}]: \w(t)=y\}$$
with the usual convention $\inf\varnothing =\infty$. Then, if $x\in\R$ and $y\in(-\infty,x)$, and if 
$\omega\in \S_x$, we set for every $s\geq 0$,
$$\eta_s(\omega)=\inf\Big\{u\geq 0:\int_0^u \mathrm{d}r\,\mathbf{1}_{\{\zeta_{(\omega_r)}\leq\tau_y(\omega_r)\}}>s\Big\}$$
(note that the condition $\zeta_{(\omega_r)}\leq\tau_y(\omega_r)$ holds if and only if $\tau_y(\omega_r)=\infty$ or $\tau_y(\omega_r)=\zeta_{(\omega_r)}$).
Then, setting $\omega'_s=\omega_{\eta_s(\omega)}$ for every $s\geq 0$ defines an element of $\S_x$,
which will be denoted by  $\omega'={\rm tr}_y(\omega)$ and called the truncation of $\omega$ at $y$.
See \cite[Proposition 10]{ALG} for a proof. The genealogical tree of 
${\rm tr}_y(\omega)$ is canonically and isometrically identified 
with the closed subset of $\t_\zeta$ consisting of all $a$ such that
$Z_b(\omega)>y$ for every strict ancestor $b$ of $a$ (we leave the proof as an exercise 
for the reader). By abuse of notation, we often write ${\rm tr}_y(W)$ instead of ${\rm tr}_y(\omega)$.

\subsection{Constructing a measure metric space from a snake trajectory}
\label{sec:consmms}

Let us fix $\omega\in \S$. Recall the definition of the tree 
$\t_\zeta$ and of lexicographical intervals on that tree. We define, for every $a,b\in\t_\zeta$,
\begin{equation}
\label{Dzero}
D^\circ(a,b)=Z_a + Z_b -2\max\Big(\min_{c\in[a,b]} Z_c, \min_{c\in[b,a]} Z_c\Big).
\end{equation}
We record two easy but important properties of $D^\circ$. First, for every $a,b\in\t_\zeta$,
\begin{equation}
\label{low-bd-D}
D^\circ(a,b)\geq |Z_a-Z_b|.
\end{equation}
Then, if $a_*$ is such that $Z_{a_*}=W_*$, we have for every $a\in\t_\zeta$,
\begin{equation}
\label{equal-D}
D^\circ(a_*,a)=Z_a-Z_{a_*}.
\end{equation}

We let $D(a,b)$ be the largest symmetric function of the pair $(a,b)$ that
is bounded above by $D^\circ(a,b)$ and satisfies the triangle inequality: For every
$a,b\in\t_\zeta$,
\begin{equation}
\label{formulaD}
D(a,b) = \inf\Big\{ \sum_{i=1}^k D^\circ(a_{i-1},a_i)\Big\},
\end{equation}
where the infimum is over all choices of the integer $k\geq 1$ and of the
elements $a_0,a_1,\ldots,a_k$ of $\t_\zeta$ such that $a_0=a$
and $a_k=b$. Then $D$ is a pseudo-metric on $\t_\zeta$, and 
we let $\mm$ be the associated quotient space (the quotient of
$\t_\zeta$ for the equivalence relation $a\approx b$ if and only if $D(a,b)=0$)
equipped with the distance induced by $D$, for which we keep the same notation.
Then $(\mm,D)$ is a compact metric space. If $\Pi$ denotes the
canonical projection from $\t_\zeta$ onto $\mm$, we define the 
volume measure $\bv$ on $\mm$ as the push forward of the volume measure 
on $\t_\zeta$ under $\Pi$. We can therefore view 
$(\mm,D,\bv)$ as a measure metric space. There are two distinguished
points in $\mm$. One of them is $\Pi(p_\zeta(0))$ (that is, the image under
the projection $\Pi$ of the root of $\t_\zeta$). The other one 
is $\Pi(a_*)$, where $a_*$ denotes any vertex of $\t_\zeta$ such that
$Z_{a_*}=W_*$ 
(the existence of such a vertex is immediate by a compactness argument, and conversely, if
$a,a'$ are two vertices of $\t_\zeta$ such that
$Z_a=Z_{a'}=W_*$, we have $D^\circ(a,a')=0$
and therefore $D(a,a')=0$, so that $\Pi(a)=\Pi(a')$). 

Note that, as a consequence of \eqref{low-bd-D} and \eqref{equal-D}, we have 
$D(a_*,a)=Z_a-W_*$ for every $a\in\t_\zeta$.

The preceding construction obviously depends on the choice of $\omega$. We
claim that it does so in a measurable way.

\begin{lemma}
\label{measurability}
The mappings
\begin{align*}
\ll:\;&\omega\mapsto (\mm,D,\bv)\\
\ll^\bullet:\;&\omega\mapsto (\mm,D,\bv, \Pi(a_*))\\
\ll^{\bullet\bullet}:\;&\omega\mapsto (\mm,D,\bv,\Pi(a_*),\Pi(p_\zeta(0)))
\end{align*}
from $\S$ into $\M,\M^\bullet,\M^{\bullet\bullet}$ respectively,
are measurable.
\end{lemma}

\proof Let us explain the argument for the mapping $\ll:\omega\mapsto (\mm,D,\bv)$.
It is convenient to introduce the space $\dd$ 
of all continuous pseudo-metrics on a compact interval of the form
$[0,\sigma]$. In other words, a mapping $d:[0,\sigma_d]^2\la \R_+$ (where 
$\sigma_d\geq 0$) belongs to $\dd$ if it is continuous and symmetric, vanishes 
on the diagonal and 
satisfies the triangle inequality. We equip $\dd$ with the distance
$$\delta(d_1,d_2)= \Big(\sup_{s,t\geq 0} |d_1(s\wedge \sigma_{d_1},t\wedge \sigma_{d_1})
- d_2(s\wedge \sigma_{d_2},t\wedge \sigma_{d_2})|\Big) + |\sigma_{d_1}-\sigma_{d_2}|,$$
and with the associated Borel $\sigma$-field.

Then we first notice that the mapping from $\S$ into $\dd$ defined
by 
\begin{equation}
\label{measurability1}
\S\ni \omega \mapsto \Big( [0,\sigma]^2\ni (s,t) \mapsto D(p_\zeta(s),p_\zeta(t))\Big)
\end{equation}
is measurable. We leave the proof as an exercise for the reader.

We can then define a mapping from $\dd$ into $\M$ in the following way.
If $d\in\dd$, we consider the associated equivalence relation $\sim_d$ on $[0,\sigma_d]$
($s\sim_d t$ if and only if $d(s,t)=0$) and denote the associated quotient space
$[0,\sigma_d]/\sim_d$ by $M_d$. This space is equipped with the metric induced by
$d$ (still denoted by $d$) and with the volume measure $\bv_d$ which is the push forward of
Lebesgue measure under the canonical projection $\mathrm{p}_d:[0,\sigma_d]\to M_d$.
Then the mapping  $\omega\mapsto (\mm,D,\bv)$ is the composition 
of the mapping $d\mapsto (M_d,d,\bv_d)$ with the (measurable)
mapping in \eqref{measurability1}. 

So, to get the desired measurability property, we only have to verify that 
the mapping $d\mapsto (M_d,d,\bv_d)$ is measurable from $\dd$ into $\M$. In
fact, we prove that this mapping is continuous. 

To this end, let $(d_n)_{n\geq 1}$ be a sequence in $\dd$ such that $d_n$
converges to $d$ in $\dd$. For every $n\ge 1$, we define a correspondence $\cc_n$
between $M_{d_n}$ and $M_d$ by declaring that a pair $(x_n,x)$ belongs to
$\cc_n$ if and only if there exists $t\geq 0$ such that $x_n=\mathrm{p}_{d_n}(t\wedge \sigma_{d_n})$
and $x=\mathrm{p}_{d}(t\wedge \sigma_d)$. Then the distortion of $\cc_n$ is bounded above by
$$\sup_{s,t\geq 0} |d_n(s\wedge \sigma_{d_n},t\wedge \sigma_{d_n})
- d(s\wedge \sigma_{d},t\wedge \sigma_{d})| \leq \delta(d_n,d).$$
We can then apply Lemma \ref{GHP-corrresp} (in the case $k=0$)
to the correspondence $\cc_n$ and to the measure $\nu_n$ on $M_{d_n}\times M_d$
defined as the push forward of Lebesgue measure on $[0,\sigma_{d_n}\wedge \sigma_d]$
under the mapping $t\mapsto (\mathrm{p}_{d_n}(t),\mathrm{p}_{d}(t))$. It 
follows that $(M_{d_n},d_n,\bv_{d_n})$ converges to $(M_d,d,\bv_d)$ in $\M$ as desired.
\endproof

\rem The measure metric space $(\mm,D,\bv)$ and the pointed space $(\mm,D,\bv, \Pi(a_*))$
do not change if $\omega$ is replaced by the re-rooted snake trajectory $\omega^{[r]}$
for some $r\in[0,\sigma]$. To explain this, recall that the mapping $s\mapsto r\oplus s$ induces an
isometry $\mathcal{I}$ from $\t_{\zeta^{[r]}}$ onto $\t_\zeta$, and that 
$Z_a(\omega^{[r]})=Z_{\mathcal{I}(a)}(\omega)-\wh W_r(\omega)$ for every $a\in\t_{\zeta^{[r]}}$, by construction. The isometry $\mathcal{I}$ preserves 
intervals, in the sense that $\mathcal{I}([a,b])=[\mathcal{I}(a),\mathcal{I}(b)]$ for every 
$a,b\in \t_{\zeta^{[r]}}$. It follows that we have also $D_{\omega^{[r]}}(a,b)=D_\omega(\mathcal{I}(a),\mathcal{I}(b))$
for every 
$a,b\in \t_{\zeta^{[r]}}$ (with the obvious notation $D_{\omega^{[r]}},D_\omega$).
Furthermore, if $a$ is a vertex of $\t_{\zeta^{[r]}}$ with minimal label, the same is true for $\mathcal{I}(a)$
in $\t_{\zeta}$.
The first two mappings in Lemma \ref{measurability} (but not the third one) are thus invariant under the re-rooting operation.

\subsection{Brownian snake excursion measures}
\label{Brsnexc}

\smallskip
We now define excursion measures of the Brownian snake, which are ($\sigma$-finite) measures on 
$\S$ that play a fundamental role in this work. For every $x\in\R$, we define $\N_x$
as the $\sigma$-finite measure on $\S_x$ that satisfies the following two properties: Under $\N_x$,

\begin{enumerate}
\item[(i)] the distribution of the lifetime function $(\zeta_s)_{s\geq 0}$ is the It\^o 
measure of positive excursions of linear Brownian motion, normalized so that, for every $\ve>0$,
$$\N_x\Big(\sup_{s\geq 0} \zeta_s >\ve\Big)=\frac{1}{2\ve};$$
\item[(ii)] conditionally on $(\zeta_s)_{s\geq 0}$, the tip function $(\wh W_s)_{s\geq 0}$ is
a Gaussian process with mean $x$ and covariance function 
$$K(s,s')= \min_{s\wedge s'\leq r\leq s\vee s'} \zeta_r.$$
\end{enumerate}

Notice that the quantity $\sigma$ (in part (i) of Definition \ref{def:snakepaths}) corresponds under 
$\N_x$ to the duration of the excursion $(\zeta_s)_{s\geq 0}$. Under the normalization in (i),
we have for every $s>0$,
\begin{equation}
\label{lawduration}
\N_x(\sigma >s)= \frac{1}{\sqrt{2\pi s}}\,.
\end{equation}
We refer to \cite[Chapter XII]{RY} for more information about the It\^o excursion measure, and to
\cite{Zurich} for a detailed study of the Brownian snake (our presentation using snake 
trajectories is slightly different from the one in \cite{Zurich}). 

We will use the formula
\begin{equation}
\label{law-min}
\N_x(W_*<y)= \frac{3}{2(x-y)^2}\;,\qquad\hbox{for every }y\in(-\infty,x).
\end{equation}
See \cite[Chapter 6]{Zurich} for a proof. 

Let us take $x=0$. For every $r>0$, we define the probability measure $\N_0^{(r)}:=\N_0(\cdot\mid \sigma = r)$, which can be characterized
by properties exactly similar to (i) and (ii) above, with the only difference that in (i) the It\^o
measure is replaced by the law of a positive Brownian excursion conditioned to have duration $r$. From \eqref{lawduration}, we have
\begin{equation}
\label{desintN}
\N_0=\int_0^\infty \frac{\mathrm{d}s}{2\sqrt{2\pi s^3}}\,\N^{(s)}_0.
\end{equation}

One can prove (see e.g. \cite{LGW}) that $\N_0$ a.e. or $\N^{(r)}_0$ a.s. there is a unique time $s_*$
in $[0,\sigma]$ such that $\wh W_{s_*}=W_*$. We may consider the snake trajectory $W$
re-rooted at $s_*$, which is denoted by $W^{[s_*]}$ in Section \ref{sec:snake}. 
The distribution of $W^{[s_*]}$ under $\N^{(r)}_0$
 can be interpreted as $\N^{(r)}_0$ conditioned on the event
that $\wh W_s\geq 0$ for every $s$ (see \cite{LGW}). 

Let us recall the invariance property of the measures $\N^{(r)}_0$ under
re-rooting. For any $r>0$, and any nonnegative measurable function 
$G$ on $\S_0$, for every $s\in[0,r]$,
\begin{equation}
\label{uniform-re-root}
\N^{(r)}_0\Big(G(W^{[s]})\Big) = \N^{(r)}_0\Big(G(W)\Big).
\end{equation}
See formula (3) in \cite{LGW} (this result is initially due to \cite{MM1}). 

We will use the fact that, under $\N^{(r)}_0$, $s_*$ is uniformly distributed over $[0,r]$
and independent of the snake trajectory $W^{[s_*]}$ (see \cite{LGW}, noting that both
these properties follow from the invariance of $\N^{(r)}_0$ under re-rooting). 

\subsection{The definition of the Brownian map}
\label{sec:defBM}

The (standard) Brownian map is the (random) measure metric space $(\mm,D,\bv)$ constructed
as in Section \ref{sec:consmms} from a random snake trajectory distributed according
to $\N^{(1)}_0$. In what follows, we will view $(\mm,D,\bv)$ 
as a $2$-pointed space, where the distinguished points are $\Pi(a_*)$ and 
$\Pi(p_\zeta(0))$ in this order, as in the third mapping of Lemma \ref{measurability} --- here
$a_*=p_\zeta(s_*)$ in agreement with the notation of Section \ref{sec:consmms}. 
We will often write $x_*=\Pi(a_*)$ and $x_0=\Pi(p_\zeta(0))$ for the two distinguished points of $\mm$.
Notice that Lemma \ref{measurability} is used to make sense 
of the Brownian map as a random variable with values in $\M^{\bullet\bullet}$. 

We will also be interested in the ($2$-pointed) space $(\mm,D,\bv)$  under the 
infinite measure $\N_0$ and then we speak about the free Brownian map. 

A property that plays an important role in the study of the Brownian map
is the fact that, $\N_0$ or $\N^{(r)}_0$ a.e., for every $a,b\in \t_\zeta$,
$D(a,b)=0$ implies that $D^\circ(a,b)=0$ (the converse is obvious since $D\leq D^\circ$).
See \cite{Invent}. We also mention the ``continuous cactus bound'' \cite[Corollary 3.2]{Geo}:
$\N_0$ or $\N^{(r)}_0$ a.e., for every $a,b\in\t_\zeta$,
\begin{equation}
\label{ctscactus}
D(a,b)\geq Z_a+Z_b - 2\min_{c\in\llbracket a,b\rrbracket} Z_c,
\end{equation}
where $\llbracket a,b\rrbracket$ denotes the geodesic segment between $a$ and $b$
in the tree $\t_\zeta$ (not to be confused with the lexicographical interval $[a,b]$). 

\subsection{Convergence to the Brownian map}
\label{sec:Bmap}

For every $n\geq 1$, let $Q_{(n)}$ be a uniformly distributed quadrangulation in $\Q_n$, and write 
$C^{(n)}=(C^{(n)}_k)_{k\geq 0}$ and $L^{(n)}=(L^{(n)}_k)_{k\geq 0}$ for the contour and label functions of the well-labeled
tree $\tau_{Q_{(n)}}$ associated with $Q_{(n)}$ in Schaeffer's bijection. A key ingredient of the proof of the
convergence to the Brownian map is the convergence in distribution, in the functional sense on the 
Skorokhod space,
\begin{equation}
\label{conv-coding}
\Big((2n)^{-1/2}C^{(n)}_{\lfloor 2ns\rfloor},\sqrt{\frac{3}{2}}\,(2n)^{-1/4} L^{(n)}_{\lfloor 2ns\rfloor}\Big)_{s\geq 0}
\build{\la}_{n\to\infty}^{\rm(d)} (\zeta^{[s_*]}_s,\wh W^{[s_*]}_s)_{s\geq 0}\quad\hbox{under }\N^{(1)}_0,
\end{equation}
where the limit process is the pair consisting of the lifetime function
and the tip function of a snake trajectory
``re-rooted at the minimum'' under $\N^{(1)}_0$. See \cite[Theorem 2.5]{Invent}. 

In what follows, we will consider quadrangulations that are both rooted and pointed. We let
$\Q_n^\bullet$ stand for the set of all rooted and pointed planar quadrangulations
with $n$ faces. An element $q^\bullet$ of $\Q_n^\bullet$ thus consists
of a rooted quadrangulation $q\in\Q_n$ and a distinguished vertex that
we will denote by $\xi$. Note that $\#\Q_n^\bullet=(n+2)\#\Q_n$. 

Suppose now that $Q_{(n)}^\bullet=(Q_{(n)},\xi_{(n)})$ is uniformly distributed over $\Q^\bullet_n$, and define 
$C^{(n)}=(C^{(n)}_k)_{k\geq 0}$ and $L^{(n)}=(L^{(n)}_k)_{k\geq 0}$ as previously, so that the convergence
\eqref{conv-coding} holds. We can in fact reinforce this convergence as follows. Let $k_{(n)}$ be the 
first integer $k\geq 0$ such that the term of index $k$ in the contour sequence of $\tau_{Q_{(n)}}$ is equal 
to $\xi_{(n)}$ (this makes sense unless $\xi_{(n)}$ coincides with the root vertex $\rho$ of
$Q_{(n)}$, in which case we take $k_{(n)}=0$ by convention). Then, we have 
\begin{equation}
\label{conv-codingbis}
\Big(\Big((2n)^{-1/2}C^{(n)}_{\lfloor 2ns\rfloor},\sqrt{\frac{3}{2}}\,(2n)^{-1/4} L^{(n)}_{\lfloor 2ns\rfloor}\Big)_{s\geq 0}, \frac{1}{2n}k_{(n)}\Big)
\build{\la}_{n\to\infty}^{\rm(d)} \Big((\zeta^{[s_*]}_s,\wh W^{[s_*]}_s)_{s\geq 0},1-s_*\Big)\quad\hbox{under }\N^{(1)}_0. 
\end{equation}
Notice that in the above limit, the quantity $1-s_*$ is uniformly distributed over $[0,1]$ and independent of 
$(\zeta^{[s_*]}_s,\wh W^{[s_*]}_s)_{s\geq 0}$. Thus we could replace $1-s_*$ by $s_*$ or by
any random variable $U$ uniform over $[0,1]$ and independent of 
$(\zeta^{[s_*]}_s,\wh W^{[s_*]}_s)_{s\geq 0}$. The point of writing $1-s_*$ is the fact that the
vertex $p_{\zeta^{[s_*]}}(1-s_*)$ corresponds in the re-rooted tree $\t_{\zeta^{[s_*]}}$ to the root
of the tree $\t_\zeta$ (recall that the isometry from $\t_{\zeta^{[s_*]}}$ onto $\t_\zeta$
is induced by the mapping $r\mapsto s_*\oplus r$). 

This fact allows us to deduce the following convergence from \eqref{conv-codingbis}. For $0\leq j\leq 2n$, set
\begin{align*}
\wt C^{(n)}_j&=C^{(n)}_{k_n\oplus j} + C^{(n)}_{k_n} - 2\,\min_{(k_n\oplus j)\wedge k_n\leq i\leq (k_n\oplus j)\vee k_n} C^{(n)}_i\,\\
\wt L^{(n)}_j&= L^{(n)}_{k_n\oplus j} - L^{(n)}_{k_n},
\end{align*}
where, only in this formula, we use the notation $k_n\oplus j=k_n+j$ if $k_n+j\leq 2n$, and $k_n\oplus j=k_n+j-2n$ if $k_n+j> 2n$.
Also take $\wt C^{(n)}_j=\wt L^{(n)}_j=0$ for $j>2n$. 
Then (excluding the case $\xi_{(n)}=\rho$), $\wt C^{(n)}$ and $\wt L^{(n)}$ are the contour and label functions of the tree
$\wt\tau_{Q^\bullet_{(n)}}$ defined as the tree  $\tau_{Q_{(n)}}$ re-rooted at the first corner of $\xi_{(n)}$, with labels shifted so that the label of the root is $0$. 
Noting the equality $(\omega^{[r]})^{[\sigma-r]}=\omega$ for a snake trajectory 
$\omega \in \S_0$ and $r\in[0,\sigma]$, it is then straightforward to deduce from  \eqref{conv-codingbis} that we have also
the convergence in distribution
\begin{equation}
\label{conv-codingter}
\Big((2n)^{-1/2}\wt C^{(n)}_{\lfloor 2ns\rfloor},\sqrt{\frac{3}{2}}\,(2n)^{-1/4} \wt L^{(n)}_{\lfloor 2ns\rfloor}\Big)_{s\geq 0}
\build{\la}_{n\to\infty}^{\rm(d)} (\zeta_s,\wh W_s)_{s\geq 0}\quad\hbox{under }\N^{(1)}_0.
\end{equation}

\smallskip

Let us now state the convergence of rescaled quadrangulations to the 
Brownian map \cite{Uniqueness,Mie2}. 

\begin{theorem}
\label{convBmap}
For every $n\geq1$, let $Q_{(n)}^\bullet=(Q_{(n)},\xi_{(n)})$ be uniformly distributed over 
$\Q^\bullet_n$, let $\mu_{(n)}$ stand for the counting measure on $V(Q_{(n)})$
and view $V(Q_{(n)})$ as a $2$-pointed space with distinguished points
$\rho$ and $\xi_{(n)}$ in this order. Then,
\begin{equation}
\label{convBM}
\Big(V(Q_{(n)}),\sqrt{\frac{3}{2}}(2n)^{-1/4}\dg,n^{-1}\mu_{(n)}\Big)\build{\la}_{n\to\infty}^{\rm (d)} (\mm, D, \bv)\quad\hbox{under }\N^{(1)}_0,
\end{equation}
where the convergence holds in distribution in $\M^{\bullet\bullet}$. Furthermore,
this convergence holds jointly with \eqref{conv-codingter}. 
\end{theorem}

In fact the statement of Theorem \ref{convBmap} is stronger than the formulation in 
\cite{Uniqueness,Mie2} because these papers consider only convergence in
the Gromov-Hausdorff sense. In the remaining part of this section, we briefly explain how this stronger form 
can be derived.  We recall the notation $\Pi$ for the
canonical projection from $\t_\zeta$ onto $\mm$ and we write $\bp=\Pi\circ p_\zeta$ in the following lines.
The two distinguished points of $\mm$ are thus $x_*=\bp(s_*)$ and $x_0=\bp(0)$ in that order.

To begin with, we recall the arguments used to obtain the
Gromov-Hausdorff convergence
(see \cite[Section 4]{Mie-Saint-Flour} for a pedagogical presentation). 
In what follows, we implicitly exclude the case $\xi_{(n)}=\rho$, which 
occurs with vanishing probability when $n\to\infty$.
Let $(u^{(n)}_0,u^{(n)}_1,\ldots, u^{(n)}_{2n})$ be the contour sequence of the tree $\wt\tau_{Q^\bullet_{(n)}}$
(in particular, $u^{(n)}_0=\xi_{(n)}$).
Notice that
the vertex set $V(Q_{(n)})$ is identified to $V(\tau_{Q_{(n)}})\cup\{\rho\}$ in Schaeffer's bijection,
and that $V(\wt\tau_{Q^\bullet_{(n)}})=V(\tau_{Q_{(n)}})$. 
For $0\leq j\leq 2n-1$, set $\eta_n(j)= j$ if $\wt C^{(n)}_{j+1}<\wt C^{(n)}_j$ and
$\eta_n(j) = j+1$ if $\wt C^{(n)}_{j+1}>\wt C^{(n)}_j$. The point of introducing $\eta_n$
is the fact that, for every vertex $u\in V(Q_{(n)})\backslash\{\rho,\xi_{(n)}\}$ there are exactly
two values of $j$ in $\{0,\ldots,2n-1\}$ such that $u^{(n)}_{\eta_n(j)}=u$. 
Consider then the correspondence $\cc_{(n)}$ between $V(Q_{(n)})$ and $\mm$
defined by 
$$\cc_{(n)}=\{(u^{(n)}_{\eta_n(\lfloor 2nt\rfloor)},\bp(t)):t\in[0,1)\} \cup \{(\rho,x_*),(\xi_{(n)},x_0)\}.$$
By the tightness argument developed in \cite{Invent}, we may assume that, at least along
a subsequence of values of $n$,
the sequence 
$$\Big((9/8)^{1/4}n^{-1/4}\,\dg(u^{(n)}_{\lfloor 2ns\rfloor}, u^{(n)}_{\lfloor 2nt\rfloor})\Big)_{s,t\in[0,1]}$$
converges in distribution, jointly with \eqref{conv-codingter}, and the uniqueness
of the Brownian map \cite{Uniqueness,Mie2} ensures that the limiting distribution is that of $(D(\bp(s),\bp(t)))_{s,t\in[0,1]}$ independently of the
chosen subsequence. 
Using the Skorokhod representation theorem, one can
construct the whole sequence $(Q^\bullet_{(n)})_{n\geq 1}$ so that
the convergence \eqref{conv-codingter} holds a.s.,
and $((9/8)^{1/4}n^{-1/4}\dg(u^{(n)}_{\lfloor 2ns\rfloor}, u^{(n)}_{\lfloor 2nt\rfloor}))_{s,t\in[0,1]}$
converges uniformly to $(D(\bp(s),\bp(t))))_{s,t\in[0,1]}$ a.s. It easily follows that
the distortion of $\cc_n$ (when  $V(Q_{(n)})$ is equipped with $(9/8)^{1/4}n^{-1/4}\dg$ and 
$\mm$ with $D$)
tends to $0$ a.s. 
This gives the convergence in the Gromov-Hausdorff sense.

To get the convergence in the Gromov-Hausdorff-Prokhorov sense, 
we just have to apply Lemma \ref{GHP-corrresp} to the measure $\nu_{(n)}$ defined
on the product $V(Q_{(n)})\times \mm$ by the
formula
$$\langle \nu_{(n)},\varphi\rangle= \int_0^1 \varphi(u^{(n)}_{\eta_n(\lfloor 2nt\rfloor)},\bp(t))\,\mathrm{d}t.$$
By construction, $\nu$ is supported on $\cc_{(n)}$. On the other hand, using the notation in
Lemma \ref{GHP-corrresp}, we have $\pi'_*\nu_{(n)}=\bv$, and $\pi_*\nu_{(n)}$ is the uniform
probability measure on $V(Q_{(n)})\backslash\{\rho,\xi_{(n)}\}$, so that the Prokhorov distance
between $\pi_*\nu_{(n)}$ and $n^{-1}\mu_{(n)}$ clearly tends to $0$ as $n\to\infty$. 

Since $(\rho,\bp(s_*))$ and $(\xi_{(n)}, \bp(0))$ both belong to $\cc_{(n)}$, Lemma \ref{GHP-corrresp}
in fact gives the convergence in $\M^{\bullet\bullet}$ as desired. 

\subsection{Convergence of Boltzmann quadrangulations}
\label{sec:Bolquad}

We will be interested in the convergence of random quadrangulations whose size is not
fixed. In this section, we consider a Boltzmann quadrangulation $Q$: This means that $Q$
is a random rooted and pointed quadrangulation such that, for every $n\geq 1$ and every $q\in \Q_n^\bullet$,
$$\P(Q=q)= \frac{1}{2}\,12^{-n}.$$
Note that the factor $1/2$ corresponds to the formula
$$\sum_{n=1}^\infty 12^{-n}\,\#\Q_n^\bullet = 2,$$
which follows from \eqref{nbquad}
and the identity $\#\Q_n^\bullet=(n+2)\#\Q_n$. Using asymptotics for Catalan numbers, we have
\begin{equation}
\label{asympnbquad}
\P(|Q|=n)= 4^{-n}c_n \build{\sim}_{n\to\infty}^{} \frac{1}{\sqrt{\pi}}\,n^{-3/2}.
\end{equation}
In particular, there is a constant $K$ such that $\P(|Q|=n)\leq Kn^{-3/2}$ for every $n\geq 1$. 

Recall that $\xi$ denotes the distinguished vertex of $Q$, and assume that $\xi\not =\rho$.
We then consider the tree $\wt\tau_Q$, which is the tree 
associated with $Q$ in Schaeffer's bijection and re-rooted at the first corner of $\xi$, as in 
Section \ref{sec:Bmap} --- recall that labels are shifted so that the label of $\xi$ in $\wt\tau_Q$ is $0$. We
write $(\wt C_k)_{k\geq 0}$ for the contour function of $\wt\tau_Q$, and $(\wt L_k)_{k\geq 0}$
for its label function. 
We also write $\wt L_*=\min\{\wt L_k:k\geq 0\}$ for the minimal label in $\wt \tau_Q$,
and note that $\dg(\rho,\xi)=-\wt L_* +1$. As previously, $\mu$ stands for the 
counting measure on $V(Q)$ and $V(Q)$ has the two distinguished points
$\rho$ and $\xi$. We write $\D(\R_+,\R^2)$ for the classical
Skorokhod space of c\`adl\`ag functions from $\R_+$ into $\R^2$. 

\begin{corollary}
\label{convBmapB}
Let $\delta >0$. The distribution under $\P(\cdot \mid \dg(\rho,\xi)>\delta \sqrt{n})$ of 
$$\Bigg(\Big(n^{-1}\wt C_{\lfloor n^2t\rfloor},\sqrt{\frac{3}{2}}n^{-1/2} \wt L_{\lfloor n^2t\rfloor}\Big)_{t\geq 0}, \Big(V(Q),\sqrt{\frac{3}{2}}n^{-1/2} \dg, 2n^{-2}\mu\Big)\Bigg)$$
converges as $n\to \infty$ to the distribution under $\N_0(\cdot\mid W_*<-\delta\sqrt{3/2})$ of
$$\Big((\zeta_t,\wh W_t)_{t\geq 0},(\mm,D,\bv)\Big).$$
The convergence is in the sense of weak convergence of probability
measures on $\D(\R_+,\R^2)\times \M^{\bullet\bullet}$.
\end{corollary}

\proof This is basically a consequence of Theorem \ref{convBmap}, but we provide some details. Let $r>0$. 
As a consequence of Theorem \ref{convBmap} and easy scaling arguments, we get that the distribution of 
$$\Bigg(\Big(n^{-1}\wt C_{\lfloor n^2t\rfloor},\sqrt{\frac{3}{2}}n^{-1/2} \wt L_{\lfloor n^2t\rfloor}\Big)_{t\geq 0},
\Big(V(Q),\sqrt{\frac{3}{2}}n^{-1/2} \dg, 2n^{-2}\mu\Big)\Bigg)$$
under $\P(\cdot \mid |Q|=\lfloor n^2r\rfloor)$ converges as $n\to \infty$ to the distribution of
$$\Big((\zeta_t,\wh W_t)_{t\geq 0},(\mm,D,\bv)\Big)$$
under $\N_0^{(2r)}$. Then, let $\Phi$ and $\Psi$ be bounded continuous functions defined respectively on the
space $\D(\R_+,\R^2)$ and on the space $\M^{\bullet\bullet}$.
To simplify notation, set
$$\begin{array}{ll}
\Phi_n= \Phi\Big( \Big(n^{-1}\wt C_{\lfloor n^2t\rfloor},\sqrt{\frac{3}{2}}\,n^{-1/2}\wt L_{\lfloor n^2t\rfloor}\Big)_{t\geq 0}\Big)\;,\ &\Psi_n= \Psi\Big(V(Q),\sqrt{\frac{3}{2}}\,n^{-1/2} \dg,
2n^{-2}\mu\Big),\\
\Phi_\infty=\Phi((\zeta_t,\wh W_t)_{t\geq 0})\;,&\Psi_\infty=\Psi(\mm,D,\bv).
\end{array}$$
 Fix two positive constants $a$ and $A$ with $a<1<A$. Then, recalling that $\dg(\rho,\xi)=- \wt L_*+1$
 on the event $\{\xi\not =\rho\}$, we have
\begin{equation}
\label{convBm1}
\E\Big[ \Phi_n\Psi_n\,\mathbf{1}_{\{an^2\leq |Q|\leq A n^2\}\cap \{\dg(\rho,\xi)>\delta \sqrt{n}\}}\Big]
= n^2\int_a^A \mathrm{d} r\,\E\Big[ \Phi_n\Psi_n\,\mathbf{1}_{\{|Q|=\lfloor n^2 r\rfloor\}}\,\mathbf{1}_{\{ -\wt L_*+1>\delta\sqrt{n}\}}\Big]
 + O(n^{-3}),
\end{equation}
using the bound $\P(|Q|=n)\leq K\,n^{-3/2}$.

By the first observation of the proof,
we have, for every $r\in[a,A]$,
\begin{equation}
\label{convBm2}
\E\Big[ \Phi_n\Psi_n\,\mathbf{1}_{\{ -\wt L_*+1>\delta\sqrt{n}\}}\,\Big|\, |Q|=\lfloor n^2 r\rfloor\Big]
\quad\build{\la}_{n\to\infty}^{} 
\N_0^{(2r)}\Big( \Phi_\infty\Psi_\infty\,\mathbf{1}_{\{W_*< -\delta\sqrt{3/2}\}}\Big),
\end{equation}
using also the fact that the law of $W_*$ under $\N^{(2r)}_0$ has no atom at $-\delta\sqrt{3/2}$ (we omit the details, but note that
\eqref{law-min} shows that this is true for a.a. $r$, which would suffice for our argument). 

On the other hand, we observe that
$$\P(|Q| > An^2) \leq \sum_{k=\lfloor An^2\rfloor}^\infty Kk^{-3/2} \leq \frac{K'}{n\sqrt{A}},$$
for some constant $K'$, and, using the estimate \eqref{estimmin},
$$\P(\{|Q|< a n^2\}\cap \{- \wt L_*+1>\delta\sqrt{n}\})
\leq \sum_{k=1}^{\lfloor an^2\rfloor} K\,k^{-3/2} \times K_1\,\exp\Big(-K_2 \frac{\delta\sqrt{n}}{k^{1/4}}\Big)
\leq \frac{\tilde K}{n}\,\exp(-K_2\,\frac{\delta}{2a^{1/4}}),
$$
with some constant $\tilde K$. This shows that, given $\ve>0$, we can fix 
$a$ and $A$ such that, for every $n$, we have 
\begin{equation}
\label{convBm4}
n\,\P(|Q| > An^2) \leq \ve\;,\quad n\,\P(\{|Q|< a n^2\}\cap \{\dg(\rho,\xi)>\delta \sqrt{n}\})\leq \ve.
\end{equation}

By \eqref{asympnbquad} and \eqref{convBm2},
$$n^3\E\Big[ \Phi_n\Psi_n\,\mathbf{1}_{\{- \wt L_*+1>\delta\sqrt{n}\}}\,\mathbf{1}_{\{|Q|=\lfloor n^2 r\rfloor\}}\Big]
\build{\la}_{n\to\infty}^{} \frac{1}{\sqrt{\pi r^3}}\, \N_0^{(2r)}\Big( \Phi_\infty\Psi_\infty\,\mathbf{1}_{\{W_*< -\delta\sqrt{3/2}\}}\Big),$$
so that, by \eqref{convBm1} and dominated convergence (justification is easy),
$$n\,\E\Big[ \Phi_n\Psi_n\,\mathbf{1}_{\{an^2\leq |Q|\leq A n^2\}\cap \{\dg(\rho,\xi)>\delta \sqrt{n}\}}\Big]
\build{\la}_{n\to\infty}^{} \int_a^A \frac{\mathrm{d}r}{\sqrt{\pi r^3}}\,\N_0^{(2r)}\Big( \Phi_\infty\Psi_\infty\,\mathbf{1}_{\{W_*< -\delta\sqrt{3/2}\}}\Big).$$
Using the estimates \eqref{convBm4}, we conclude that
\begin{align}
\label{convBm5}
&n\,\E\Big[ \Phi_n\Psi_n\,\mathbf{1}_{\{\dg(\rho,\xi)>\delta \sqrt{n}\}}\Big]\nonumber\\
&\qquad\build{\la}_{n\to\infty}^{} \int_0^\infty \frac{\mathrm{d}r}{\sqrt{\pi r^3}}\,\N_0^{(2r)}\Big( \Phi_\infty\Psi_\infty\,\mathbf{1}_{\{W_*< -\delta\sqrt{3/2}\}}\Big)
=4\,\N_0\Big( \Phi_\infty\Psi_\infty\,\mathbf{1}_{\{W_*< -\delta\sqrt{3/2}\}}\Big).
\end{align}
We divide the asymptotics \eqref{convBm5} by the same asymptotics written with $\Phi=1$ and $\Psi=1$, and we arrive at
the desired result. \endproof

We conclude this section with a technical lemma that will be useful later. We again leave aside the case $\xi=\rho$ and, for every integer $k\geq 0$, we write 
$\mathcal{K}_k$ for the collection of those vertices $v$ of $\wt\tau_Q$
such that 
the geodesic (in the tree $\wt\tau_{Q}$) from $v$ to $\xi$ visits at least one vertex
with label at most $\wt L_*+k-1$. Since the label of
$\xi$ in $\wt\tau_Q$ is $0$, $\mathcal{K}_k=V(\wt\tau_Q)$ when $k\geq -\wt L_*+1=\dg(\rho,\xi)$. 

Recall the notation $\underline\w=\min\{\w(t):0\leq t\leq \zeta_{(\w)}\}$ for a stopped path $\w$. 

\begin{lemma}
\label{lem-tec-hull}
Let $\delta>0$ and $r\in(0,\delta]$. The distribution of $n^{-2}\#\mathcal{K}_{\lfloor r\sqrt{n}\rfloor}$ under $\P(\cdot \mid \dg(\rho,\xi)>\delta \sqrt{n})$
converges as $n\to\infty$ to the distribution of 
$$\frac{1}{2}\int_0^\sigma \mathrm{d}s\,\mathbf{1}_{\{\un{W}_s\leq W_*+r\sqrt{3/2}\}}$$
under $\N_0(\cdot\mid W_*<-\delta\sqrt{3/2})$, and this convergence holds jointly with that
of Corollary
\ref{convBmapB}.
\end{lemma} 

\proof This is a relatively straightforward consequence of Corollary
\ref{convBmapB}, and we only sketch the arguments. Write $u_0,u_1,\ldots,u_{2|Q|}$
for the contour sequence of $\wt \tau_Q$, and for every $m\in\{0,1,\ldots,2|Q|\}$, let
$\un{\wt L}_m$ be the minimal label of all ancestors of $u_m$ in $\wt\tau_Q$ 
(equivalently this is the minimal label on the geodesic from $\xi$ to $u_m$). 
Then, the distribution of 
$$\Big( \frac{1}{\sqrt{n}}\,\un{\wt L}_{\lfloor n^2 s\rfloor\wedge 2|Q|}\Big)_{s\geq 0}$$
under $\P(\cdot \mid \dg(\rho,\xi)>\delta \sqrt{n})$ converges as $n\to\infty$ to the distribution of 
$(\sqrt{2/3}\,\un{W}_s)_{s\geq 0}$
under $\N_0(\cdot\mid W_*<-\delta\sqrt{3/2})$, and this convergence
holds jointly with that of Corollary
\ref{convBmapB} (the point is that the convergence of the rescaled contour and label functions in
Corollary
\ref{convBmapB} entails the convergence of the associated discrete snakes, in the spirit of
the homeomorphism theorem of \cite{MM} --- we omit a few details here). For every $i\in\{0,1,\ldots,2|Q|-1\}$, set $\eta(i)=i+1$ if $\wt C_{i+1}>\wt C_i$
and $\eta(i)=i$ if $\wt C_{i+1}<\wt C_i$, in such a way that, for every 
vertex $v$ of $V(\wt\tau_Q)\backslash\{\xi\}$ there are exactly two values of
$i$ such that $u_{\eta(i)}=v$.  Then, for every $0\leq k\leq -\wt L_*$,
$$\#\mathcal{K}_k=\frac{1}{2}\#\{i\in\{0,1,\ldots,2|Q|-1\}: \un{\wt L}_{\eta(i)}\leq \wt L_*+k-1\}$$
and (taking $k=\lfloor r\sqrt{n}\rfloor$) it follows
from the first observation of the proof that the distribution under $\P(\cdot \mid \dg(\rho,\xi)>\delta \sqrt{n})$ of
$$n^{-2}\#\mathcal{K}_{\lfloor r\sqrt{n}\rfloor}= \frac{1}{2} \int_0^{2n^{-2}|Q|}
\mathrm{d}s\,\mathbf{1}\Big\{\frac{1}{\sqrt{n}}\un{\wt L}_{\eta(\lfloor n^2 s\rfloor)}\leq \frac{1}{\sqrt{n}}\wt L_* + \frac{\lfloor r\sqrt{n}\rfloor-1}{\sqrt{n}}\Big\},$$
converges as $n\to\infty$ to the law of
$$\frac{1}{2}\int_0^\sigma \mathrm{d}s\,\mathbf{1}_{\{\un{W}_s\leq W_*+r\sqrt{3/2}\}}$$
jointly with the convergence of Corollary
\ref{convBmapB}, as desired. \endproof

\rem In the last step of the proof, we implicitly use the fact that the set of all $s\in[0,\sigma]$
such that $\un{W}_s= W_*+r\sqrt{3/2}$ has zero Lebesgue measure, $\N_0$ a.e. This can be
derived by a scaling argument, and we leave the proof as an exercise for the reader.

\section{Quadrangulations with a boundary and Brownian disks}
\label{sec:quadbBd}

\subsection{Quadrangulations with a boundary}
\label{subsec:quadb}

\begin{figure}[!h]
 \begin{center}
 \includegraphics[width=10cm]{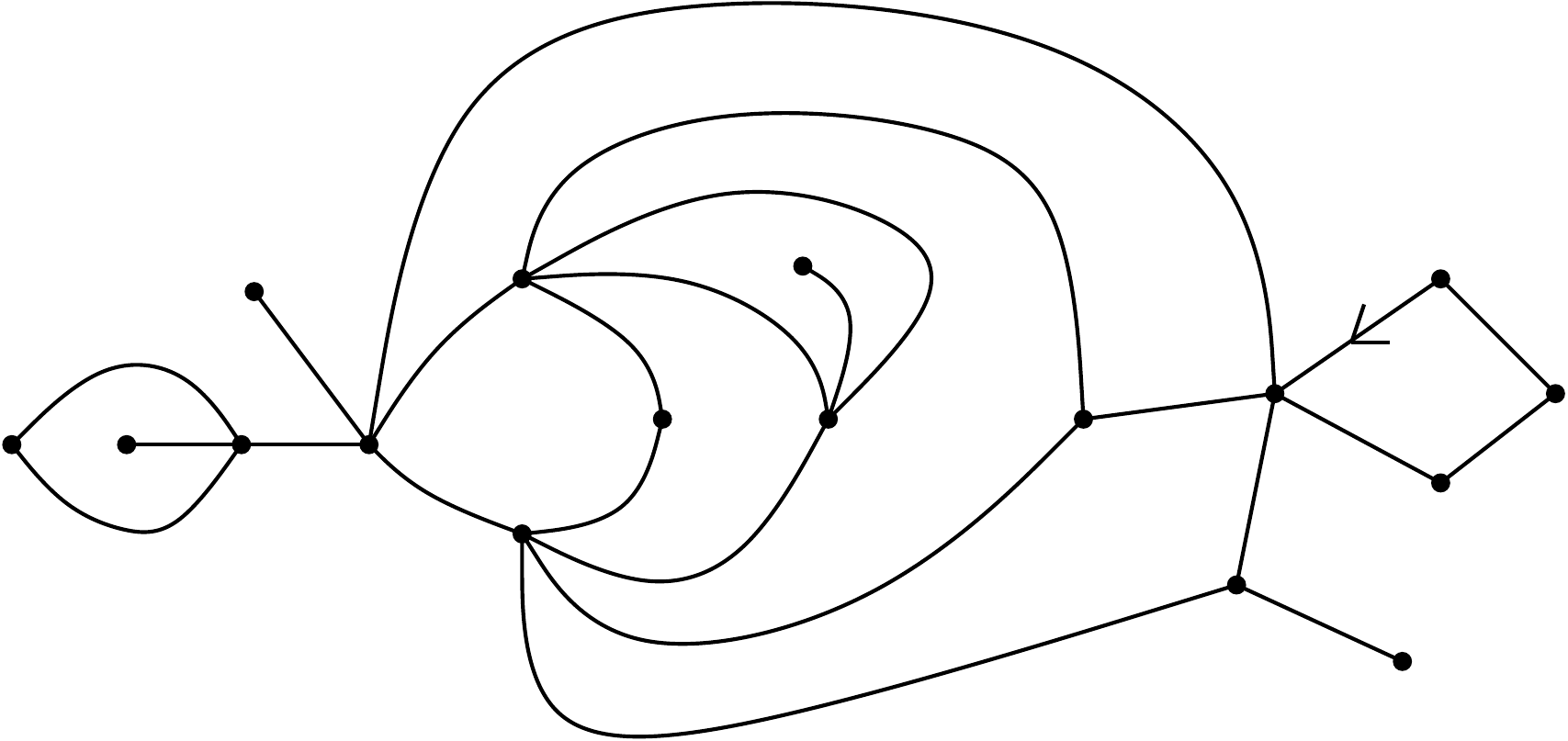}
 \caption{\label{quad-bd}
 A quadrangulation with a boundary of size 16.}
 \end{center}
 \vspace{-5mm}
 \end{figure}

Recall that a quadrangulation with a (general) boundary is a rooted planar map $q$
such that all faces but the root face (lying to the right of the root edge) have degree $4$. The root face is also called the outer face
and the other faces are inner faces. The degree of the outer face, which is an even integer, is then called the boundary size 
or the perimeter of $q$. We will use the notation $\partial q$ for the the collection of all
vertices incident to the outer face. By definition, the root corner of $q$ is the corner of the outer face that is incident to the root vertex to the right of
the root edge. See Fig.~\ref{quad-bd} for an example.

The quadrangulation with a boundary $q$ is said to be pointed if
there is a distinguished vertex, which will be denoted by $\xi$. For every integer $k\geq 1$, we denote the
set of all pointed quadrangulations with a boundary of size $2k$ by $\Q^{\partial,k}$.
For every integer $n\geq 0$, the subset of  $\Q^{\partial,k}$ consisting of those 
quadrangulations $q$ that have $n$ inner faces is denoted by  $\Q^{\partial,k}_n$.
Then, for every $k\geq 1$, there is a constant $b_k>0$ such that
$$\#\Q^{\partial,k}_n \build{\sim}_{n\to\infty}^{} b_k\,12^n\, n^{-3/2}.$$
See formula (4) in \cite{CM}, noting that this formula applies to
non-pointed quadrangulations. 

A random variable $B_{(k)}$ with values in $\Q^{\partial,k}$
is called a Boltzmann (pointed) 
quadrangulation with a boundary of size $2k$ if, for every integer $n\geq 0$
and every $q\in \Q^{\partial,k}_n$,
$$\P(B_{(k)}=q)= \wt b_k\,12^{-n},$$
where $\wt b_k>0$ is the appropriate normalizing constant.

\subsection{Free Brownian disks}
\label{sec:freeBdisk}

In this section, we recall the construction of (free pointed) Brownian disks 
from \cite{BM}. Our presentation is in fact a little simpler than the one
in \cite{BM} because we do not condition on the volume of the
Brownian disks.

We start from a standard linear Brownian motion $X=(X_t)_{t\geq 0}$ with $X_0=0$,
and, for every $r>0$, we set $T_r=\inf\{t\geq 0:X_t=-r\}$. For every $t\geq 0$, 
we also set $\un X_t=\inf_{0\leq s\leq t} X_s$. We consider a process
$Y^\circ$ such that, conditionally on $X$, $Y^\circ$ is a centered Gaussian 
process with covariance
$$\E[Y^\circ_sY^\circ_{s'}\mid X]= \inf_{s\wedge s'\leq u\leq s\vee s'} (X_u-\un X_u).$$

From now on, we fix $r>0$, and let $(\mathrm{b}(s))_{0\leq s\leq r}$ be a standard
Brownian bridge (starting and ending at $0$) with duration $r$, independent of $(X,Y^\circ)$. We set,
for every $t\in [0,T_r]$,
$$Y_t= Y^\circ_t + \sqrt{3}\,\mathrm{b}(-\un X_t).$$

We may view $(Y_t)_{0\leq t\leq T_r}$ as labels assigned to 
the vertices of the tree coded by $(X_t-\un X_t)_{0\leq t\leq T_r}$ and
then proceed in a way similar to the construction of the
Brownian map in Section \ref{sec:Bmap}. Let us explain the details.
For every $s,s'\in[0,T_r]$, we define
$$\un Y_{s,s'}=\inf\{ Y_u:u\in [s,s']\},$$
where by convention $[s,s']=[s,T_r]\cup[0,s']$ if $s>s'$. We set
$$D^{\partial,\circ}(s,s')= Y_s+Y_{s'} - 2\max\{\un Y_{s,s'},\un Y_{s',s}\}.$$
We then let $D^{\partial}$ be the maximal 
pseudo-metric on $[0,T_r]$ that is bounded above
by $D^{\partial,\circ}$ and is such that, for every $s,s'\in[0,T_r]$,
$$X_s=X_{s'}=\min_{s\wedge s'\leq u\leq s\vee s'} X_u\; \Rightarrow\; D^{\partial}(s,s')=0.$$
See formula (17) in \cite{BM} for an ``explicit'' formula for $D^{\partial}(s,s')$. To interpret the
last display, note that $(X_t)_{0\leq t\leq T_r}$ codes a forest of $\R$-trees in the way explained in \cite[Section 2.1]{BM}, and that $D^{\partial}$ can be viewed as defined on pairs of
vertices of this forest, in a way very similar
to the construction of the Brownian map. 

We let $\D^\bullet_r$ be the quotient space of $[0,T_r]$ for the equivalence 
relation associated with the pseudo-metric $D^{\partial}$, and equip
$\D^\bullet_r$ with the induced metric (still denoted by $D^{\partial}$).
We let $p^\bullet_r$ be the canonical projection from $[0,T_r]$ onto $\D^\bullet_r$.
The boundary $\partial \D^\bullet_r$ is by definition
$$\partial \D^\bullet_r=p^\bullet_r(\{s\in[0,T_r]:X_t=\un X_t\}).$$
The terminology is justified by the fact that $\D^\bullet_r$ is homeomorphic
to the unit disk in the plane, and that this homeomorphism maps
$\partial \D^\bullet_r$ to the unit circle \cite{Bet,BM}.

We view $\D^\bullet_r$  as a pointed measure metric space:
The volume measure $\bv^\bullet_r$ is the image of Lebesgue measure on $[0,T_r]$ under the
canonical projection $p^\bullet_r$, and the distinguished point is 
$p^\bullet_r(s_*)$, where $s_*$ is the (unique) time in $[0,T_r]$ 
at which $Y$ attains its minimum.

By definition, the pointed measure metric space $(\D^\bullet_r, D^{\partial}, \bv^\bullet_r)$
is the free pointed Brownian disk with perimeter $r$. We note that the total mass of
the volume measure $\bv^\bullet_r$ is $T_r$, which has density
$$v\mapsto\frac{r}{\sqrt{2\pi v^3}}\, \exp(-\frac{r^2}{2v}).$$
By scaling arguments, it is straightforward to verify that
\begin{equation}
\label{scaling-disk}
(\D^\bullet_r, D^{\partial}, \bv^\bullet_r)\build{=}_{}^{\rm(d)} (\D^\bullet_1, \sqrt{r}\,D^{\partial}, r^2\bv^\bullet_1)
\end{equation}

We denote the distribution of 
the free pointed Brownian disk with perimeter $r$ by $\F^\bullet_r$. Thus,  $\F^\bullet_r$ is a probability measure 
on $\M^\bullet$. 
If $M^\bullet\in \M^\bullet$,
we set $\kappa(M^\bullet)=M$, where $M$ is the measure metric space obtained 
by ``forgetting'' the distinguished point of $M^\bullet$. 
The distribution of the free Brownian disk can then be defined 
via its density with respect to $\kappa_*\F^\bullet_r$. 
For notational convenience,
if $M^\bullet\in \M^\bullet$, we write 
$\Sigma(M^\bullet)$ for the total mass of the volume measure of $M^\bullet$.
The distribution $\F_r$ of the free Brownian disk with perimeter $r$ is the probability measure 
on $\M$ obtained by setting, for any nonnegative measurable function $G$ on $\M$,
$$\F_r(G)= r^2\,\F_r^\bullet\Big(\frac{1}{\Sigma}\,G\circ\kappa\Big).$$
Note that this defines a probability measure since
$$\F_r^\bullet\Big(\frac{1}{\Sigma}\Big) = \int_0^\infty \mathrm{d}{v} \frac{r}{\sqrt{2\pi v^5}}\, \exp(-\frac{r^2}{2v}) = r^{-2}.$$

\subsection{Brownian disks with glued boundary}
\label{sec:Bdglued}

For technical reasons, we will also need to consider the pointed Brownian disk with ``boundary
glued into a single point''. To explain this, let $(\D^\bullet_r,D^{\partial},\bv^\bullet_r)$ be the free pointed Brownian disk
with perimeter $r>0$ constructed as in Section \ref{sec:freeBdisk}. Recalling the notation $\partial \D^\bullet_r$ for the boundary of $\D^\bullet_r$, we set
\begin{equation}
\label{metric-glued}
D^{\partial,\dagger}(x,y)= \min\{D^{\partial}(x,y),D^{\partial}(x,\partial\D^\bullet_r)+
D^{\partial}(y,\partial \D^\bullet_r)\},
\end{equation}
for every $x,y\in \D^\bullet_r$. It is easy to verify that $D^{\partial,\dagger}$ satisfies the
triangle inequality, and that $D^{\partial,\dagger}(x,y)=0$ if and only
if either $x=y$, or $x$ and $y$ both belong to $\partial \D^\bullet_r$. 
We write $\D^{\bullet,\dagger}_r$ for the set obtained from 
$\D^\bullet_r$ by identifying the whole boundary to a single point denoted by $x_\partial$, and
equip $\D^{\bullet,\dagger}_r$ with the distance induced by $D^{\partial,\dagger}$, and with the volume measure $\bv^{\bullet,\dagger}_r$
which is the image of the volume measure $\bv^\bullet_r$ on $\D^{\bullet}_r$ under the canonical projection
(we notice that $\bv^\bullet_r$ gives no mass to the boundary).
Finally we let $\F_r^{\bullet,\dagger}$ be the distribution of 
$\D^{\bullet,\dagger}_r$ viewed as a random $2$-pointed measure metric space, where the
two distinguished vertices are
the ``boundary point'' $x_\partial$ and the distinguished vertex of $\D^\bullet_r$.

As an immediate consequence of \eqref{scaling-disk}, we have also
\begin{equation}
\label{scaling-disk2}
(\D^{\bullet,\dagger}_r, D^{\partial,\dagger}, \bv^{\bullet,\dagger}_r)
\build{=}_{}^{\rm(d)} (\D^{\bullet,\dagger}_1, \sqrt{r}\,D^{\partial,\dagger}, r^2\bv^{\bullet,\dagger}_1)
\end{equation}

We can use the same method to define the (unpointed) free Brownian disk with glued boundary from the distribution 
$\F_r$. We get a probability measure $\F^\dagger_r$ on the space $\M^\bullet$ --- the distinguished point 
is the ``boundary point''.

\subsection{Convergence to the free Brownian disk}
\label{sec:convBdisk}

In this section, we recall the convergence of rescaled quadrangulations with a boundary to
the free pointed Brownian disk, which is obtained in \cite{BM}. As in the case of
the Brownian map, we need a more precise formulation than the one in \cite{BM},
which we state below with a sketch of proof (see also \cite[Section 4.1]{GM0}). 

As in Section \ref{subsec:quadb}, we let $B_{(k)}$ be a Boltzmann 
quadrangulation with a boundary of size $2k$, for every $k\geq 1$. 
We recall that $B_{(k)}$ comes with a distinguished vertex $\xi$.
As previously, the graph distance on the vertex set $V(B_{(k)})$
is denoted by $\dg$. We denote the counting measure on $V(B_{(k)})$
by $\mu_{(k)}$ and we view $(V(B_{(k)}),\dg,\mu_{(k)})$ as a random pointed
measure metric space with distinguished point $\xi$. 

\begin{theorem}
\label{conv-quadb}
We have
$$\Big(V(B_{(k)}),\sqrt{\frac{3}{2}}\,k^{-1/2}\,\dg, 2k^{-2}\mu_{(k)}\Big) \build{\la}_{k\to\infty}^{\rm(d)} 
(\D^\bullet_1, D^{\partial}, \bv^\bullet_1)$$
where the convergence holds in distribution in $\M^\bullet$.
\end{theorem}

\proof
The convergence of the theorem essentially follows from Theorem 8 in \cite{BM}, except that 
the latter result only deals with the Gromov-Hausdorff convergence,
and so some additional work is needed here to get convergence
in $\M^\bullet$. This is very similar to the argument explained 
after the statement of Theorem \ref{convBmap} for the convergence
to the Brownian map, but we will provide some details that will also be
useful later. We rely on the analog of Schaeffer's bijection for quadrangulations
with a boundary that is
presented in \cite[Section 3.3]{BM}. There exists a bijection between 
the set $\Q^{\partial,k}$ and the collection of all triples 
$(\epsilon,(\tau_1,\ldots,\tau_k), (b_1,\ldots,b_k))$ consisting of 
a number $\epsilon\in\{-1,1\}$, 
a forest of $k$ labeled plane trees $(\tau_1,\tau_2,\ldots,\tau_k)$
(in the sense of Section \ref{sec:quadr}), and a ``bridge''
$(b_1,\ldots,b_k)$, which is a finite sequence of 
integers such that $b_1=0$ and $b_{i+1}-b_i\geq -1$
for every $1\leq i\leq k$, with the convention that $b_{k+1}=0$.
Given a triple $(\epsilon,(\tau_1,\ldots,\tau_k), (b_1,\ldots,b_k))$, 
the corresponding quadrangulation with a boundary $q$
is such that its vertex set is identified canonically with
$V(\tau_1)\cup\cdot\cup V(\tau_k)\cup\{\xi\}$, where 
$\xi$ is an extra vertex which is the distinguished vertex of $q$. 

Fix such a triple $(\epsilon,(\tau_1,\ldots,\tau_k), (b_1,\ldots,b_k))$,
and, for every $i\in\{1,\dots,k\}$, let $(u^{i}_0,\ldots,u^{i}_{2n_i})$
be the contour sequence of $\tau_i$. The contour sequence of the
forest $(\tau_1,\ldots,\tau_k)$ is then the sequence
$$(u^{1}_0,\ldots,u^{1}_{2n_1},u^{2}_0,\ldots,u^{2}_{2n_2},\ldots,u^{k}_0,\ldots,u^{k}_{2n_k}).$$
Set $n=2(n_1+\cdots+n_k)+k$, so that the contour 
sequence can be written as $(u_0,\ldots,u_{n-1})$. For every $j\in\{0,1,\ldots,n-1\}$, we write $\alpha_j=i$ if the
$j$-th term in the contour sequence corresponds to a vertex of $\tau_i$
(equivalently if $2(n_1+\cdots+n_{i-1})+i-1\leq j<2(n_1+\cdots+n_{i})+i$). We define
the contour and label functions by 
$$C_j=|u_j|-\alpha_j+1\;,\quad L_j=b_{\alpha_j}+\ell_{u_j}\,,\quad \hbox{for }0\leq j\leq n-1\,,$$
where $\ell_{u_j}$ is the label of $u_j$. By convention, we also 
take $C_n=-k$ and $L_n=0$. We extend the definition of $C_j$
and $L_j$ to real values of the parameter in $[0,n]$ by linear interpolation.

Suppose now that we consider a Boltzmann (pointed)
quadrangulation $B_{(k)}$ with a boundary of size $2k$. Then the associated 
forest $(\tau_1,\tau_2,\ldots,\tau_k)$ is a forest of $k$ independent
Galton-Watson trees with geometric offspring distribution of parameter $1/2$
with labels chosen uniformly among 
admissible labels. Also, $(b_1,\ldots,b_k)$ is chosen uniformly at random
among all bridges of length $k$, independently of $(\tau_1,\ldots,\tau_k)$. 
If $(C^{(k)}_i)_{0\leq i\leq N_k}$ and $(L^{(k)}_i)_{0\leq i\leq N_k}$ are
the contour and label functions associated with $B^{(k)}$, one has
\begin{equation}
\label{conv-coding-bdry}
\Big(k^{-1}C^{(k)}_{k^2s}, \sqrt{\frac{3}{2}}\,k^{-1/2}\,L^{(k)}_{k^2s}\Big)_{0\leq s\leq k^{-2}N_k}
\build{\la}_{k\to\infty}^{\rm(d)} \Big(X_t,Y_t\Big)_{0\leq t\leq T_1}
\end{equation}
where $X,Y$ and $T_1$ are as Section \ref{sec:freeBdisk}, and the convergence 
in distribution makes sense in the space of finite paths defined in Section \ref{sec:snake}.
See \cite[Proposition 7]{Bet}, which in fact gives a stronger result than \eqref{conv-coding-bdry}.
As in Section \ref{sec:freeBdisk}, we let 
$(\D^\bullet_1, D^{\partial}, \bv^\bullet_1)$ be the free Brownian disk 
constructed from the pair $(X,Y)$, and $p^\bullet_1$ denotes the canonical projection from 
$[0,T_1]$ onto $\D^\bullet_1$. Also, the contour sequence of the forest
associated with $B_{(k)}$ is denoted by $(u^{(k)}_0,\ldots,u^{(k)}_{N_k-1})$, and 
for $0\leq i\leq N_k-1$, we set $\eta_k(i)=i$ if $C^{(k)}_{i+1}<C^{(k)}_i$, and $\eta_k(i)=i+1$
otherwise.

Following \cite[Section 5.2]{BM}, we can use the Skorokhod representation theorem
to construct the sequence $(B_{(k)})_{k\geq1}$ in such a way that
the convergence \eqref{conv-coding-bdry} holds a.s. and moreover the following holds.
If $\cc_k$ is the correspondence  between $(V(B_{(k)}),(3/2k)^{1/2}\,\dg)$
and $(\D^\bullet_1, D^{\partial})$ defined by declaring that, for every 
$s\in[0,T_1\vee k^{-2}(N_k-1)]$, $(u^{(k)}_{\eta_k(\lfloor k^2s\rfloor\wedge (N_k-1))},p^\bullet_1(s
\wedge T_1))\in\cc_k$, and furthermore the distinguished vertex $\xi$ of $B^{(k)}$ is in correspondence
with the distinguished vertex of $\D^\bullet_1$, then
the distortion of the correspondence $\cc_k$
tends to $0$ as $k\to\infty$. This property gives the convergence
 in the Gromov-Hausdorff sense.
 
To get the stronger convergence in $\M^\bullet$, we rely again on Lemma \ref{GHP-corrresp}.
With the notation of this lemma, it is easy to check that
the measure $\nu_{(k)}$ defined on the product $V(B^{(k)})\times \D^\bullet_1$ by
$$\langle \nu_{(k)},\phi\rangle = \int_0^{k^{-2}N_k\wedge T_1} \mathrm{d}t\,\phi(u^{(k)}_{\eta_k(\lfloor k^2t\rfloor)},p^\bullet_1(t))$$
is supported on $\cc_k$, and  both the Prokhorov distance
between $\pi_*\nu_{(k)}$ and $2k^{-2}\mu_{(k)}$ and the Prokhorov distance
between $\pi'_*\nu_{(k)}$ and $\bv^\bullet_1$ tend to $0$ as $k\to\infty$. 
\endproof

Let $B_{(k)}$ be as in the previous statement. We can 
view $B_{(k)}$  as a submap of a (rooted and pointed) planar
quadrangulation which is defined as follows.
We first add to $B_{(k)}$
an extra vertex $\varpi$ belonging to the 
outer face of $B_{(k)}$. Then, if $c_0,c_1,\ldots,c_{2k}=c_0$ are the corners
incident to the outer face enumerated in clockwise order starting from the root corner $c_0$, we draw an edge connecting
$\varpi$ to each of the corners $c_0,c_2,c_4,\ldots, c_{2k-2}$. In this way, we obtain a planar quadrangulation, which is rooted
at the oriented edge from $\varpi$ to the corner $c_0$.
We write $B_{(k)}^\dagger$ for this quadrangulation, and use the notation
$\dg^\dagger$ for the graph distance on $V(B_{(k)}^\dagger)=V(B_{(k)})\cup\{\varpi\}$. 
We view $(V(B_{(k)}^\dagger),\dg^\dagger,\mu_{(k)})$ as 
a random $2$-pointed measure metric space, 
whose distinguished points are $\varpi$ and $\xi$.

\begin{corollary}
\label{conv-quadb-glued}
We have
$$\Big(V(B_{(k)}^\dagger),\sqrt{\frac{3}{2}}\,k^{-1/2}\,\dg^\dagger, 2k^{-2}\mu_{(k)}\Big) \build{\la}_{k\to\infty}^{\rm(d)} 
(\D^{\bullet,\dagger}_1, D^{\partial,\dagger}, \bv^{\bullet,\dagger}_1)$$
where the convergence holds in distribution in $\M^{\bullet\bullet}$,
and the limit is the Brownian disk with glued boundary
defined in Section \ref{sec:Bdglued}. 
\end{corollary}

\proof We argue as in the proof of Theorem \ref{conv-quadb},
assuming that the convergence \eqref{conv-coding-bdry} holds a.s. and
that the distortion of the correspondences $\cc_k$ introduced in this proof converges to $0$ a.s.  
Then,
for every $k$, we define a 
 correspondence $\cc_k^\dagger$ between $V(B_{(k)}^\dagger)$ and $\D^{\bullet,\dagger}_1$
by declaring that: 
 \begin{enumerate}
 \item[$\bullet$] If $(u,x)\in \cc_k$
 and $x\notin \partial\D^{\bullet}_1$ then $(u,x)\in\cc^\dagger_k$.
 \item[$\bullet$] If $(u,x)\in \cc_k$
 and $x\in \partial\D^{\bullet}_1$ then $(u,x_\partial)\in\cc^\dagger_k$.
\item[$\bullet$] $(\varpi,x_\partial)\in\cc^\dagger_k$. 
\end{enumerate}

We claim that the distortion of $\cc_k^\dagger$ (as a correspondence between
$(V(B_{(k)}^\dagger),(3/2k)^{1/2}\,\dg^\dagger)$
and $(\D^{\bullet,\dagger}_1, D^{\partial,\dagger})$) tends to $0$. We first observe 
that, for every $u,v\in V(B_{(k)})$,
\begin{equation}
\label{easybound}
|\dg^\dagger(u,v) - \min\{\dg(u,v),\dg(u,\partial B_{(k)})+\dg(v,\partial B_{(k)})\}| \leq 4.
\end{equation}
and also recall that, by definition,
$$D^{\partial,\dagger}(x,y)= \min\{D^{\partial}(x,y),D^{\partial}(x,\partial\D^\bullet_1)+
D^{\partial}(y,\partial \D^\bullet_1)\},$$
for $x,y\in \D^\bullet_1$. 

Taking into account the last two formulas, and recalling that we already know that the distortion of $\cc_k$
tends to $0$, our claim will 
follow if we can verify that
\begin{equation}
\label{compardist0}
\sup\{|(3/2k)^{1/2} \,\dg(u,\partial B_{(k)}) - D^{\partial}(x,\partial\D^\bullet_1)|: 
(u,x)\in\cc_k\} \build{\la}_{k\to\infty}^{\rm a.s.} 0.
\end{equation}

For every $\ve>0$, write $(\partial \D^\bullet_1)_\ve$ for the set of all points of
$\D^\bullet_1$ whose $D^{\partial}$-distance from $\D^\bullet_1$
is less than $\ve$. Similarly, for every $A>0$, write  $(\partial B_{(k)})_A$
for the set of all vertices of
$B_{(k)}$ whose $\dg$-distance from $\partial B_{(k)}$
is less than $A$. We will prove that for every $\ve>0$, a.s. for $k$ large enough,
\begin{enumerate}
\item[\rm(i)] $(\partial \D^\bullet_1)_\ve$ contains all $x\in \D^\bullet_1$ such that $(u,x)\in\cc_k$
for some $u\in \partial B_{(k)}$. 
\item[\rm(ii)] $(\partial B_{(k)})_{\ve\sqrt{k}}$ contains all $u\in V(B_{(k)})$ such that $(u,x)\in\cc_k$
for some $x\in \partial \D^\bullet_1$. 
\end{enumerate}
Recalling the notation $\mathrm{dis}\,\cc_k$ for the distortion of $\cc_k$, it will follow that a.s. for $k$ large enough, for every $(u,x)\in \cc_k$,
$$(3/2k)^{1/2}\,\min_{v\in\partial B_{(k)}}\dg(u,v)
\geq \min_{y\in (\partial \D^\bullet_1)_\ve} D^{\partial}(x,y)-\mathrm{dis}\,\cc_k
\geq \min_{z\in \partial \D^\bullet_1} D^{\partial}(x,z)-\ve-\mathrm{dis}\,\cc_k$$
using (i). Similarly, using (ii), we will have for $k$ large enough, for $(u,x)\in\cc_k$,
$$\min_{z\in \partial \D^\bullet_1} D^{\partial}(x,z)
\geq \min_{v\in (\partial B_{(k)})_{\ve\sqrt{k}}} (3/2k)^{1/2}\dg(u,v) -\mathrm{dis}\,\cc_k
\geq \min_{w\in \partial B_{(k)}} (3/2k)^{1/2}\dg(u,w) -2\ve -\mathrm{dis}\,\cc_k.$$
The last two displays give \eqref{compardist0}. 

We still have to prove (i) and (ii). Let us start with (ii). Let
$u\in V(B_{(k)})$ such that $(u,x)\in\cc_k$
for some $x\in \partial \D^\bullet_1$. We note that $x=p^\bullet_1(s)$ with 
$s\in[0,T_1]$ such that $X_s=\un X_s$. By standard properties 
of linear Brownian motion, $s$ must be a time of decrease of $t\mapsto \un X_t$.
Using this remark and a compactness argument,  the a.s.
convergence \eqref{conv-coding-bdry} entails that, given $\delta>0$, we
can choose $k$ large enough (independently of $s$) such that
there exists $j\in\{0,1,\ldots,N_k-1\}$ with $|s-\frac{j}{k^2}|<\delta$
and 
$$C^{(k)}_{j+1}<\min_{0\leq i\leq j} C^{(k)}_i.$$
The last display implies that $u^{(k)}_j$ is the root of one of
the trees in the forest associated with $B_{(k)}$ and therefore belongs 
to $\partial B_{(k)}$ (see \cite[Section 3.3]{BM}). Here we recall that
the vertex set $V(B_{(k)})$ is canonically identified with the union of the vertex set
of the associated forest and the singleton $\{\xi\}$. Finally, recalling the definition of the 
correspondence $\cc_k$, we get that $(u^{(k)}_j,p^\bullet_1(\frac{j}{k^2}\wedge T_1))\in\cc_k$,
so that
$$(2k/3)^{-1/2}\,\dg(u,u^{(k)}_j)\leq D^{\partial}(s,\frac{j}{k^2}\wedge T_1) +\mathrm{dis}\,\cc_k,$$
and we just have to use the property $|s-\frac{j}{k^2}|<\delta$ together with the
continuity of $(t,t')\mapsto D^{\partial}(t,t')$ on the diagonal. 

The proof of (i) is similar, but we need an extra argument because not
every point of $\partial B_{(k)}$ is the root of one of the trees
in the forest associated with $B_{(k)}$. Still, if $(\mathrm{b}^{(k)}_1,\ldots,
\mathrm{b}^{(k)}_k)$
is the bridge corresponding to $B_{(k)}$, a close look at the 
combinatorial bijection for quadrangulations with a boundary
(see again \cite[Section 3.3]{BM}) shows that
the maximal graph distance between a point of $\partial B_{(k)}$ and
the collection of roots of all trees in the coding forest is bounded above by
$$\max_{1\leq i\leq k} (\mathrm{b}^{(k)}_{i+1}-\mathrm{b}^{(k)}_{i}),$$
where $\mathrm{b}^{(k)}_{k+1}=0$ by convention. The quantity in the last
display normalized by $k^{-1/2}$ converges to $0$ a.s., simply because
the rescaled processes $(k^{-1/2}\mathrm{b}^{(k)}_{\lfloor kt\rfloor})_{0\leq t\leq 1}$
converge to a Brownian bridge (this is indeed a consequence
of \eqref{conv-coding-bdry}). Modulo this observation, the proof of (i) is
very similar to that of (ii) --- we now use the fact that if $u^{(k)}_j$
is the root of a tree in the coding forest, then $j/k^2$ must be close to
a time $t\in[0,T_1]$ such that $\un X_t=X_t$ and therefore 
$p^\bullet_1(t)\in\partial \D^\bullet_1$ --- and we omit the details. 

This completes the proof of our claim that the distortion of 
$\cc_k^\dagger$ tends to $0$ and thus establishes the
convergence of the corollary in the Gromov-Hausdorff sense. To see that the
convergence holds in $M^{\bullet\bullet}$, we use Lemma \ref{GHP-corrresp}
with the same measure $\nu_{(k)}$ that was used in the proof
of Theorem \ref{conv-quadb}, or more precisely the image of 
$\nu_{(k)}$ under the mapping $(u,x)\mapsto (u,\tilde x)$, where
$\tilde x=x$ if $x\notin \partial \D^\bullet_1$ and 
$\tilde x=x_\partial$ if $x\in \partial \D^\bullet_1$. \endproof

\section{Some results about the Brownian snake}
\label{sec:resu-snake}

\subsection{The Brownian snake truncated at $W_*+\delta$}
\label{sec:Bstrunc}

In this section, we establish some properties of the 
Brownian snake that will play an important role
in the proof of our main results.

Recall the notation $\mathcal{S}_x$ for the space of all
snake trajectories with initial point $x$.
If $\omega\in \S_x$, we can make sense
of the truncation $\mathrm{tr}_b(\omega)$, for every $b<x$.
Let us also define the point process of excursions of $\omega$
outside $(b,\infty)$. Recalling the notation
$\tau_b(\w)=\inf\{t\in[0,\zeta_{(\w)}]:\w(t)=b\}$, we observe that the set
$$\{s\geq 0: \tau_b(\omega_s)<\zeta_s\}$$
is open and can therefore be written as a union of disjoint open intervals
$(\alpha_i,\beta_i)$, $i\in I$, where $I$ may be empty. From the fact that $\omega$
is a snake trajectory, it is not hard to verify that we must have $p_\zeta(\alpha_i)=p_\zeta(\beta_i)$ for every $i\in I$. Furthermore the path 
$\omega_{\alpha_i}=\omega_{\beta_i}$ hits $b$ exactly at its lifetime $\zeta_{\alpha_i}=\zeta_{\beta_i}$, and the paths 
$\omega_s$, $s\in[\alpha_i,\beta_i]$, coincide with $\omega_{\alpha_i}=\omega_{\beta_i}$ over the time interval $[0,\zeta_{\alpha_i}]$.
We then define the excursion $\omega_i$, for every $i\in I$, by declaring that, for every $s\geq 0$,
$W_s(\omega_i)$
is the finite path $(\omega_{(\alpha_i+s)\wedge \beta_i}(\zeta_{\alpha_i}+t)-b)_{0\leq t\leq \zeta^i(s)}$
with lifetime $\zeta^i(s)=\zeta_{(\alpha_i+s)\wedge \beta_i}-\zeta_{\alpha_i}$. 
We observe that the genealogical tree of $\omega^i$ (which is coded by $\zeta^i$)
is identified to the subtree of $\t_\zeta$ consisting of all descendants of $p_\zeta(\alpha_i)=p_\zeta(\beta_i)$. 
The point measure of excursions of $\omega$
outside $(b,\infty)$ is the point measure on $\S_0$ defined by
$$\nn_b(\omega):=\sum_{i\in I} \delta_{\omega_i}.$$
 Note that our definition is slightly different from the one in \cite{ALG} because
 we have shifted the excursions so that their starting
point is $0$ instead of $b$. 

We also introduce {\it exit measures}, following \cite[Chapter 5]{Zurich}. Let $x\in\R$ and $b<x$. There exists a continuous increasing process
$\ell^b=(\ell^b_s)_{s\geq 0}$ called the exit local time  such that, $\N_x$ a.e. for every $s\geq 0$,
$$\ell^b_s=\lim_{\ve\to 0} \frac{1}{\ve}\int_0^s \mathrm{d}r\,\mathbf{1}_{\{\tau_b(W_r)<\zeta_r
<\tau_b(W_s)+\ve\}}.$$
It is clear that the topological support $\mathrm{supp}(\mathrm{d}\ell^b_s)$ of the measure $\mathrm{d}\ell^b_s$
is contained in $\{s\geq 0:\tau_b(W_s)=\zeta_s\}$, but in fact more is true:
\begin{equation}
\label{support-exit}
\mathrm{supp}(\mathrm{d}\ell^b_s)=\{s\geq 0:\tau_b(W_s)=\zeta_s\}\;,\qquad \N_x\hbox{ a.e.}
\end{equation}
This follows from the more general results in \cite[Chapter 6]{Zurich}, see in particular the remark following the proof of
\cite[Theorem 6.9]{Zurich}.
The quantity $\z_b:=\ell^b_\sigma$ is called the exit measure from $(b,\infty)$ (our terminology
is slighly different from the one in \cite{Zurich}, where the exit measure would be the measure $\z_b\,\delta_b$). As explained in \cite[Section 2.5]{ALG}, the process 
$(\z_{b})_{b\in(-\infty,x)}$ has a c\`adl\`ag modification under $\N_x$, which
we consider from now on.

Let us fix $\delta>0$ and write $\N^{[-\delta]}_0=\N_0(\cdot\mid W_*<-\delta)$. The remaining part of this section will
provide results about the truncated snake $\mathrm{tr}_{W_*+\delta}(W)$, the point measure of excursions ${\mathcal{N}}_{W_*+\delta}$,
and the exit measure $\z_{W_*+\delta}$ under $\N^{[-\delta]}_0$. The motivation for these results is best understood 
from the construction of the free Brownian map as the $2$-pointed measure metric space $\mathcal{L}^{\bullet\bullet}(\omega)$
under $\N_0(\mathrm{d}\omega)$ (Section \ref{sec:defBM}). In this construction, the conditioning by $\{W_*<-\delta\}$ corresponds to the
property $D(x_*,x_0)>\delta$, the truncated snake $\mathrm{tr}_{W_*+\delta}(W)$ is closely related to
the connected component of the complement of the ball $B(x_*,\delta)$ that contains $x_0$, 
$\z_{W_*+\delta}$ is interpreted as a (generalized) length of the boundary of this component, and the 
point measure ${\mathcal{N}}_{W_*+\delta}$ yields information about the other components 
of the complement of the ball $B(x_*,\delta)$.

In the next proposition, we use the
following notation: If $\mathcal{N}$ is a point measure on the space of snake trajectories, $M(\mathcal{N})$
denotes the infimum of the quantities $W_*(\omega)$ over all atoms $\omega$ of $\mathcal{N}$.

\begin{proposition}
\label{condidistr}
Under the probability measure $\N^{[-\delta]}_0$, $\mathrm{tr}_{W_*+\delta}(W)$
and ${\mathcal{N}}_{W_*+\delta}$ are independent conditionally given $\z_{W_*+\delta}$,
and the conditional distribution of ${\mathcal{N}}_{W_*+\delta}$
given $\z_{W_*+\delta}$ is that of a Poisson point measure $\mathcal{N}$ with
intensity $\z_{W_*+\delta}\cdot\N_0$ conditioned to have $M(\mathcal{N})=-\delta$. If
$f$ and $F$ are nonnegative measurable functions defined respectively on
$\R_+$ and on $\mathcal{S}_0$, we have
$$\N_0\Big( \mathbf{1}_{\{W_*<-\delta\}}\,f(\z_{W_*+\delta})\,F(\mathrm{tr}_{W_*+\delta}(W))\Big)
=3\delta^{-3} \int_{-\infty}^0 \mathrm{d} b\,\N_0\Big(\z_b\,\exp(-\frac{3\z_b}{2\delta^2})\,f(\z_b)\,F(\mathrm{tr}_b(W))\Big).
$$
\end{proposition}

\rem If $\mathcal{N}$ is a Poisson point measure
with intensity $z\,\N_0$,  it is straightforward to define the conditional
distribution of $\mathcal{N}$ given that $M(\mathcal{N})=-\delta$, for instance as the limit when
$\ve \to 0$ of the law of $\mathcal{N}$ given that $-\delta-\ve<M(\mathcal{N})<-\delta$
(note that formula \eqref{law-min} readily gives the density of $M(\mathcal{N})$ --- see the proof below). 
Indeed, properties of Poisson measures show that this conditional distribution is the law of the sum of a Poisson point measure with intensity
$z\,\N_0(\cdot\cap\{W_*>-\delta\})$ and the Dirac mass at an independent random snake trajectory
distributed according to $\N_0(\cdot\mid W_*=-\delta)$ (the conditional distribution $\N_0(\cdot\mid W_*=-\delta)$ is studied in \cite{Mini}). 

\proof
Let $f$ and $F$ be as in the statement, and also assume that $f$ and $F$
are bounded and continuous, and that $F(\omega)=0$ if $W_*(\omega)>-\beta$ or $W_*(\omega)<-A$, for some $\beta,A>0$. Let $\Phi$ be a bounded nonnegative continuous function on $\mathcal{S}_0$,
such that $\Phi(\omega)=0$ if $\sup\{|\wh W_s(\omega)|:s\geq 0\}\leq \alpha$, for some $\alpha>0$. Set
$$G({\mathcal{N}}_b)=\exp\Big(-\int {\mathcal{N}}_b(\mathrm{d} \omega)\,\Phi(\omega)\Big).$$
Then,
\begin{align*}
&\N_0\Big( \mathbf{1}_{\{W_*<-\delta\}}\,f(\z_{W_*+\delta})\,F(\mathrm{tr}_{W_*+\delta}(W))\,G({\mathcal{N}}_{W_*+\delta})\Big)\\
&\quad=\lim_{\ve \to 0}\frac{1}{\ve} \N_0\Big( \mathbf{1}_{\{W_*<-\delta-\ve\}}
\int_{W_*+\delta}^{W_*+\delta+\ve} \mathrm{d} b\, f(\z_b)\,F(\mathrm{tr}_b(W))\,G({\mathcal{N}}_{b})\Big).
\end{align*}
To justify this, we observe that, $\N_0$ a.e. on $\{W_*<-\delta\}$, as $b\downarrow W_*+\delta$, we have $\z_b\la \z_{W_*+\delta}$ (by the right-continuity of $b\mapsto \z_b$),
$\mathrm{tr}_b(W)\la \mathrm{tr}_{W_*+\delta}(W)$
(see Lemma 11 in \cite{ALG}), and $G({\mathcal{N}}_{b})\la G({\mathcal{N}}_{W_*+\delta})$
(we leave the verification of the last convergence as an exercise for the reader). 

We next observe that, for $\ve>0$ small,
\begin{align*}
&\N_0\Big( \mathbf{1}_{\{W_*<-\delta-\ve\}}
\int_{W_*+\delta}^{W_*+\delta+\ve} \mathrm{d} b\, f(\z_b)\,F(\mathrm{tr}_b(W))\,G({\mathcal{N}}_{b})\Big)\\
&\quad= \int_{-\infty}^0 \mathrm{d} b\, \N_0\Big(\mathbf{1}_{\{b-\delta-\ve<W_*<b-\delta\}}\,f(\z_b)\,F(\mathrm{tr}_b(W))\,G({\mathcal{N}}_{b})\Big).
\end{align*}
Fix $b<0$. By applying the special Markov property (see e.g. \cite[Proposition 13]{ALG}) to the interval $(b,\infty)$, we get
\begin{align*}
&\N_0\Big(\mathbf{1}_{\{b-\delta-\ve<W_*<b-\delta\}}\,f(\z_b)\,F(\mathrm{tr}_b(W))\,G({\mathcal{N}}_{b})\Big)\\
&\quad= \N_0\Big(f(\z_b)\,F(\mathrm{tr}_b(W))\, \E_{(\z_b)}\Big[ G(\mathcal{N})\,\mathbf{1}_{\{-\delta-\ve<M(\mathcal{N})<-\delta\}}\Big]\Big),
\end{align*}
where, under the probability measure $\P_{(z)}$ (for any $z\geq 0$), $\mathcal{N}$ is a Poisson point measure
with intensity $z\,\N_0$.

We have thus obtained
\begin{align}
\label{conditech1}
&\N_0\Big( \mathbf{1}_{\{W_*<-\delta\}}\,f(\z_{W_*+\delta})\,F(\mathrm{tr}_{W_*+\delta}(W))\,G({\mathcal{N}}_{W_*+\delta})\Big)\nonumber\\
&\quad=\lim_{\ve \to 0}\int_{-\infty}^0 \mathrm{d} b\,\N_0\Big(f(\z_b)\,F(\mathrm{tr}_b(W))\, \frac{1}{\ve}\E_{(\z_b)}\Big[ G(\mathcal{N})\,\mathbf{1}_{\{-\delta-\ve<M(\mathcal{N})<-\delta\}}\Big]\Big).
\end{align}
At this point, we use \eqref{law-min} to get that the distribution of $M(\mathcal{N})$ under $\P_{(z)}$ is given (for $z>0$) by
$$\P_{(z)}(M(\mathcal{N})>-y)= \exp(-\frac{3z}{2y^2})$$
for every $y>0$, and thus the density of $M(\mathcal{N})$ under $\P_{(z)}$ is
$$3z\,|y|^{-3}\,\exp(-\frac{3z}{2y^2}) \,\mathbf{1}_{\{y<0\}}.$$
As explained in the remark following the statement of the proposition, we have then, for every $z>0$,
$$\frac{1}{\ve}\E_{(z)}\Big[ G(\mathcal{N})\,\mathbf{1}_{\{-\delta-\ve<M(\mathcal{N})<-\delta\}}\Big]\Big)
\build{\la}_{\ve \to 0}^{} 3z\,\delta^{-3}\,\exp(-\frac{3z}{2\delta^2})\,\E_{(z)}^\delta[G(\mathcal{N})],$$
where, under $\P_{(z)}^\delta$, $\mathcal{N}$ is a Poisson point measure with intensity $z\,\N_0$, 
conditioned to have $M(\mathcal{N})=-\delta$. 

We can now pass to the limit $\ve\to 0$ in the right-hand side of \eqref{conditech1}, recalling our assumptions on $F$
to justify dominated convergence.
It follows that
\begin{align}
\label{conditech2}
&\N_0\Big( \mathbf{1}_{\{W_*<-\delta\}}\,f(\z_{W_*+\delta})\,F(\mathrm{tr}_{W_*+\delta}(W))\,G({\mathcal{N}}_{W_*+\delta})\Big)\nonumber\\
&\quad = 3\,\delta^{-3}\int_{-\infty}^0 \mathrm{d} b\,\N_0\Big(f(\z_b)\,F(\mathrm{tr}_b(W))\, \z_b\,\exp(-\frac{3\z_b}{2\delta^2})\,\E_{(\z_b)}^\delta[G(\mathcal{N})]
\Big).
\end{align}
If we replace $G$ by $G'=1$ and $f$ by the function $z\mapsto f(z)\,\E_{(z)}^\delta[G(\mathcal{N})]$, we get that
\begin{align*}
&\N_0\Big( \mathbf{1}_{\{W_*<-\delta\}}\,f(\z_{W_*+\delta})\,F(\mathrm{tr}_{W_*+\delta}(W))\,G({\mathcal{N}}_{W_*+\delta})\Big)\\
&\quad = \N_0\Big( \mathbf{1}_{\{W_*<-\delta\}}\,f(\z_{W_*+\delta})\,F(\mathrm{tr}_{W_*+\delta}(W))\,\E_{(\z_{W_*+\delta})}^\delta[G(\mathcal{N})]\Big).
\end{align*}
This is enough to conclude that, under $\N_0(\cdot\mid W_*<-\delta)$, 
${\mathcal{N}}_{W_*+\delta}$ is independent of $\mathrm{tr}_{W_*+\delta}(W)$ conditionally given $\z_{W_*+\delta}$,
and its conditional distribution is as described in the proposition. Finally, the last assertion of the proposition is the special case 
$G=1$ in \eqref{conditech2}. \endproof

\begin{corollary}
\label{limit-law}
The quantity $\z_{W_*+\delta}$ is distributed under  $\N^{[-\delta]}_0$
according to the Gamma distribution with parameter $\frac{1}{2}$ with density
$$z\mapsto \frac{1}{\delta}\,\sqrt{\frac{3}{2\pi z}}\,\exp(-\frac{3z}{2\delta^2}).$$
\end{corollary}

\proof 
We specialize the formula given in Proposition \ref{condidistr} 
to the case $F=1$ and  $f(z)=e^{-\lambda z}$ with $\lambda >0$, to get
\begin{equation}
\label{limit-law11}
\N_0\Big( \mathbf{1}_{\{W_*<-\delta\}}\,\exp(-\lambda\z_{W_*+\delta})\Big)
=3\delta^{-3} \int_{-\infty}^0 \mathrm{d} b\,\N_0\Big(\z_b\,\exp(-\lambda\z_b-\frac{3\z_b}{2\delta^2})\Big).
\end{equation}
We compute the
right-hand side of the last display. By formula (6) in \cite{Hull}, for every 
$x>0$,
$$\N_0(1- e^{-\lambda \z_{-x}})= \Big(\lambda^{-1/2} +\sqrt{\frac{2}{3}} \,x\Big)^{-2}.$$
Differentiating with respect to $\lambda$ gives
$$\N_0(\z_{-x}e^{-\lambda \z_{-x}})= \Bigg(1+x\sqrt{\frac{2\lambda}{3}}\Bigg)^{-3}.$$
It then follows from \eqref{limit-law11} that
$$\N_0\Big( \mathbf{1}_{\{W_*<-\delta\}}\,\exp(-\lambda\z_{W_*+\delta})\Big)
=3\delta^{-3} \int_0^\infty \mathrm{d} x\,\Bigg(1+\, x\sqrt{\delta^{-2}+\frac{2\lambda}{3}}\Bigg)^{-3}=  \frac{3}{2\delta^3}\,\frac{1}{\sqrt{\delta^{-2}+\frac{2\lambda}{3}}}.$$
Since $\N_0(W_*<-\delta)= 3/(2\delta^2)$, we get
$$\N_0\Big( \exp(-\lambda\z_{W_*+\delta})\,\Big|\, W_*<-\delta\Big)=\sqrt{\frac{\delta^{-2}}{\delta^{-2}+\frac{2\lambda}{3}}},$$
giving the stated result. \endproof

\subsection{The positive excursion measure}
\label{sec:posexc}

The underlying idea of the results presented in Sections \ref{sec:Bmap}
and \ref{sec:Bolquad}  is the observation that scaling
limits of random quadrangulations of the sphere can be described
by the Brownian snake under its excursion measure. The main 
motivation of the present work is to show that the scaling limit of quadrangulations
with a boundary
can be described similarly by the Brownian snake under
its ``positive excursion measure''. We now give a brief presentation 
of this positive excursion measure, referring to \cite{ALG}
for more details.

Let $\S_0^+ $ stand for the set of all snake trajectories 
$\omega\in\S_0$ such that $\omega_s$ takes values in $\R_+$
for every $s\geq 0$, and, for every $\delta>0$, let $\S^{(\delta)}$
be the set of all $\omega\in \S$ such that 
$\sup_{s\geq 0}(\sup_{t\in[0,\zeta_s(\omega)]}|\omega_s(t)|)\geq \delta$. 
There exists a $\sigma$-finite measure $\N^*_0$ on
$\S$, which is supported on  $\S_0^+$, 
and gives finite mass to the sets $\S^{(\delta)}$ for every $\delta>0$, such that
$$\N^*_0(G)=\lim_{\ve\to 0}\frac{1}{\ve}\,\N_\ve(G(\tr_0(W))),$$
for every bounded continuous function $G$ on $\S$ that vanishes 
on $\S\backslash\S^{(\delta)}$ for some $\delta>0$ (see \cite[Theorem 23]{ALG}). 
Under $\N^*_0$, each of the paths $W_s$, $0<s<\sigma$, starts from $0$, then
stays positive during some time interval $(0,\alpha)$, and is stopped immediately
when it returns to $0$, if it does return to $0$. 

We will use the re-rooting representation of the positive excursion measure.
Recall our notation $W^{[r]}$ for the snake trajectory $W$
``re-rooted at $r$'' (see Section \ref{sec:snake}). The following
result is a consequence of \cite[Theorem 28]{ALG}.

\begin{theorem}
\label{re-root-rep}
For any nonnegative measurable function $G$ on $\S$, 
$$\N^*_0\Big(\int_0^\sigma \mathrm{d}r\,G(W^{[r]})\Big)
=2\int_{-\infty}^0 \mathrm{d}x\,\N_0\Big(\z_x\,G(\tr_x(W))\Big).$$
\end{theorem}

Since the definition of the exit measure $\z_x$ makes sense 
under $\N_0$, for every $x<0$, the identity of the previous
theorem makes it possible to also define the exit measure at $0$
under $\N^*_0$ (see
\cite[Section 6]{ALG} for details of the construction): To avoid confusion, we will write $\z^*_0$
for this exit measure at $0$ under $\N^*_0$. Roughly
speaking, $\z^*_0$ ``counts'' the number of paths $W_s$
that return to $0$. The identity of Theorem \ref{re-root-rep} then implies that we have
also, for any nonnegative measurable function 
$\varphi$ on $[0,\infty)$,
\begin{equation}
\label{re-root-rep-tec0}
\N^*_0\Big(\int_0^\sigma \mathrm{d}r\,\varphi(\z^*_0)\,G(W^{[r]})\Big)
=2\int_{-\infty}^0 \mathrm{d}x\,\N_0\Big(\z_x\,\varphi(\z_x)\,G(\tr_x(W))\Big).
\end{equation}

By \cite[Proposition 33]{ALG}, we can make sense of the conditional probability
measures $\N^{*,z}_0=\N^*_0(\cdot\mid \z^*_0=z)$, and we have
\begin{equation}
\label{desint}
\N^*_0=\sqrt{\frac{3}{2\pi}}\int_0^\infty \mathrm{d}z\,z^{-5/2}\,\N^{*,z}_0.
\end{equation}

\begin{proposition}
\label{law-trunc}
Let $\delta>0$. For every nonnegative measurable function $F$ 
on the space of snake trajectories, for every $z>0$, we have
$$\N^{[-\delta]}_0\Big( F(\mathrm{tr}_{W_*+\delta}(W)) \,\Big|\, \z_{W_*+\delta}=z\Big)
= z^{-2}\, \N^{*,z}_0\Big( \int_0^\sigma \mathrm{d}r\,F(W^{[r]})\Big).$$
\end{proposition}

In particular, and this will be important for us, the conditional
distribution of $\mathrm{tr}_{W_*+\delta}(W)$ given $\z_{W_*+\delta}=z$
does not depend on $\delta$. 

\proof Let $\varphi$ be a nonnegative measurable function on $[0,\infty)$. 
Recalling that $\N_0(W_*<-\delta)=3/(2\delta^2)$,
we use the formula of Proposition \ref{condidistr}, and then 
\eqref{re-root-rep-tec0}, to write
\begin{align}
\label{law-trunc1}
\N^{[-\delta]}_0\Big( F(\mathrm{tr}_{W_*+\delta}(W))\,\varphi(\z_{W_*+\delta})\Big)
&= \frac{2}{\delta}\int_{-\infty}^0 \mathrm{d} b\,\N_0\Big(\z_b\,\exp(-\frac{3\z_b}{2\delta^2})\,\varphi(\z_b)\,F(\mathrm{tr}_b(W))\Big)\nonumber\\
&= \frac{1}{\delta}\,\N^*_0\Big(\int_0^\sigma \mathrm{d}r\,F(W^{[r]})\,\exp(-\frac{3\z^*_0}{2\delta^2})\,\varphi(\z^*_0)\Big)\nonumber\\
&= \frac{1}{\delta}\,\sqrt{\frac{3}{2\pi}} \int_0^\infty \frac{\mathrm{d}z}{z^{5/2}}\, \exp(-\frac{3z}{2\delta^2})\varphi(z)\N^{*,z}_0\Big( \int_0^\sigma \mathrm{d}r\,F(W^{[r]})\Big)
\end{align}
using \eqref{desint}. On the other hand, if we set
$$G(\z_{W_*+\delta})=\N^{[-\delta]}_0\Big( F(\mathrm{tr}_{W_*+\delta}(W)) \,\Big|\, \z_{W_*+\delta}\Big),$$
we have
\begin{align}
\label{law-trunc2}
\N^{[-\delta]}_0\Big( F(\mathrm{tr}_{W_*+\delta}(W))\,\varphi(\z_{W_*+\delta})\Big)
&=\N^{[-\delta]}_0\Big( G(\z_{W_*+\delta})\,\varphi(\z_{W_*+\delta})\Big)\nonumber\\
&= \frac{1}{\delta} \sqrt{\frac{3}{2\pi}}\,\int_0^\infty \frac{\mathrm{d}z}{\sqrt{z}}\,\exp(-\frac{3z}{2\delta^2})\,G(z)\,\varphi(z)
\end{align}
by Corollary \ref{limit-law}. By comparing \eqref{law-trunc1}
and  \eqref{law-trunc2}, we get
$$G(z)= z^{-2}\, \N^{*,z}_0\Big( \int_0^\sigma \mathrm{d}r\,F(W^{[r]})\Big)$$
as stated. \endproof

\medskip
\rem By \cite[Proposition 31]{ALG}, the conditional density of 
$\sigma$ under $\N^*_0$ and knowing $\z^*_0=z$ is
$$f_z(s) = \frac{1}{\sqrt{2\pi}}\,z^3 \,s^{-3/2}\, \exp(-\frac{z^2}{2s})$$
so that
$$\N^{*,z}_0(\sigma)=  \frac{1}{\sqrt{2\pi}}\,z^3 \int_0^\infty \frac{\mathrm{d}s}{s^{3/2}}\,\exp(-\frac{z^2}{2s})
=\frac{z^2}{\sqrt{2\pi}} \times \sqrt{2}\int_0^\infty \frac{\mathrm{d}u}{\sqrt{u}}\,e^{-u}= z^2$$
which is consistent with the case $F=1$ in Proposition \ref{law-trunc}.

\section{The lazy hull process}
\label{sec:lazy}

\subsection{Gluing a quadrangulation with a general boundary in a face with a simple boundary}
\label{sec:glue}

\begin{figure}[!h]
 \begin{center}
 \includegraphics[width=12cm]{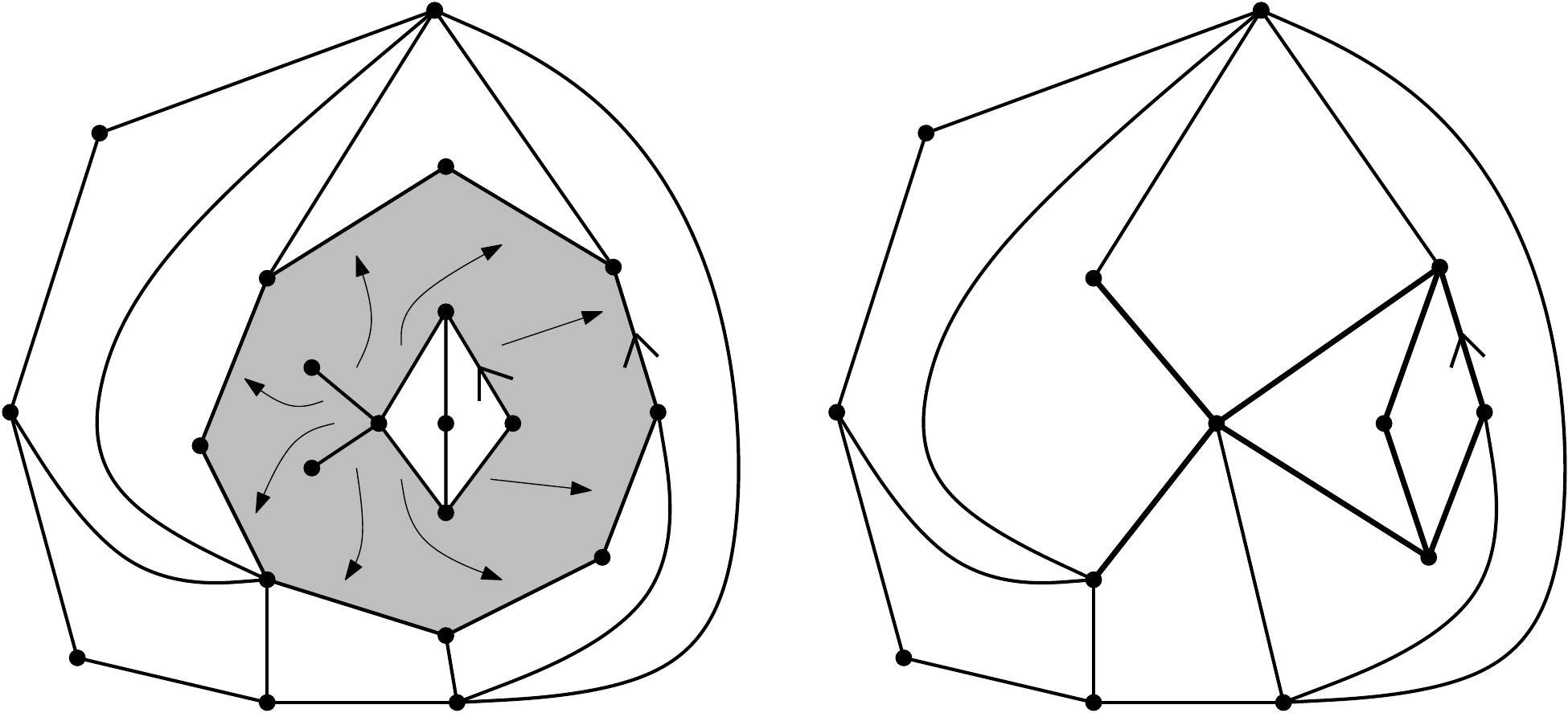}
 \caption{\label{glumap}
 Left, the planar map $\mathrm{m}$ with the distinguished face $\mathfrak{f}$ of degree $8$ in grey, and inside $\mathfrak{f}$
 the quadrangulation $\mathrm{q}$ with a boundary of perimeter $8$. The arrows indicate how the edges 
 of the boundary of $\mathrm{q}$
 are glued to the edges of the boundary of $\mathfrak{f}$. Right, the planar map obtained after the gluing operation.
}
 \end{center}
 \vspace{-2mm}
 \end{figure}

Consider a rooted planar map $\mathrm{m}$ with a distinguished face $\mathfrak{f}$ of degree  $2k$, for some $k\geq 1$. We assume that the 
boundary of $\mathfrak{f}$ is simple and that there is a distinguished oriented edge $e$ on the boundary of $\mathfrak{f}$,
such that the face $\mathfrak{f}$ lies on the left of $e$. Suppose then that
$\mathrm{q}$ is a rooted (but not pointed) quadrangulation with a (general) boundary of perimeter $2k$. Recall that the root of
$\mathrm{q}$ belongs to the boundary and that the outer face lies on the right of the root edge. There is 
then a unique way of gluing the quadrangulation $\mathrm{q}$ inside the distinguished face $\mathfrak{f}$ 
of $\mathrm{m}$, in such a way that the boundaries of $\mathfrak{f}$ and $\mathrm{q}$ are glued together, and
the root edge of $\mathrm{q}$ is glued to $e$. The construction should be 
clear from Fig.\ref{glumap}.

\subsection{The lazy peeling}
\label{sec:peel}

We will now describe an algorithm due to Budd \cite{Bud} that can be used to generate a Boltzmann quadrangulation. This algorithm
is called
the lazy peeling. We content ourselves with the properties that are needed in our applications, and
refer to \cite{Bud} and  \cite{Peccot} for more details. The main reason why we use the lazy peeling algorithm 
rather than the standard peeling (used in \cite{CLG-peeling} for instance) is the fact that we can view a Boltzmann 
quadrangulation as obtained by gluing to its lazy hull (of a certain radius) a quadrangulation with a general boundary to which we may 
apply the convergence to the Brownian disk derived in \cite{BM}. See Proposition \ref{glue-Boltzmann} below.

To describe the algorithm, we need to introduce quadrangulations with a simple 
boundary, for which we make a slightly different convention than in the case of
a general boundary: A (rooted) quadrangulation with a simple boundary of size $2k$ ($k\geq 1$) is 
a rooted planar map such that all faces have degree $4$ except for one
distinguished face of degree $2k$ (the outer face) whose 
boundary is simple --- note that, in contrast with Section \ref{subsec:quadb}, we
do not require that the root edge belongs to the boundary. We may assume that 
quadrangulations with a simple boundary are drawn in the plane so that the outer face is the infinite one, and then it makes
sense to explore the boundary in clockwise order or in counterclockwise order. 

The lazy peeling algorithm produces a finite sequence $(\qq_p)_{0\leq p\leq K}$, such that, for $0\leq p\leq K-1$, $\qq_k$
is a (rooted) quadrangulation with a {\it simple} boundary of perimeter $2L_p$, with $L_p\geq 1$, and 
$\qq_K$ is a rooted and pointed planar quadrangulation. 
To initiate the process, $\qq_0$ is the unique 
quadrangulation with a simple boundary of length $4$ and a single inner face. 

\begin{figure}[!h]
 \begin{center}
 \includegraphics[width=12cm]{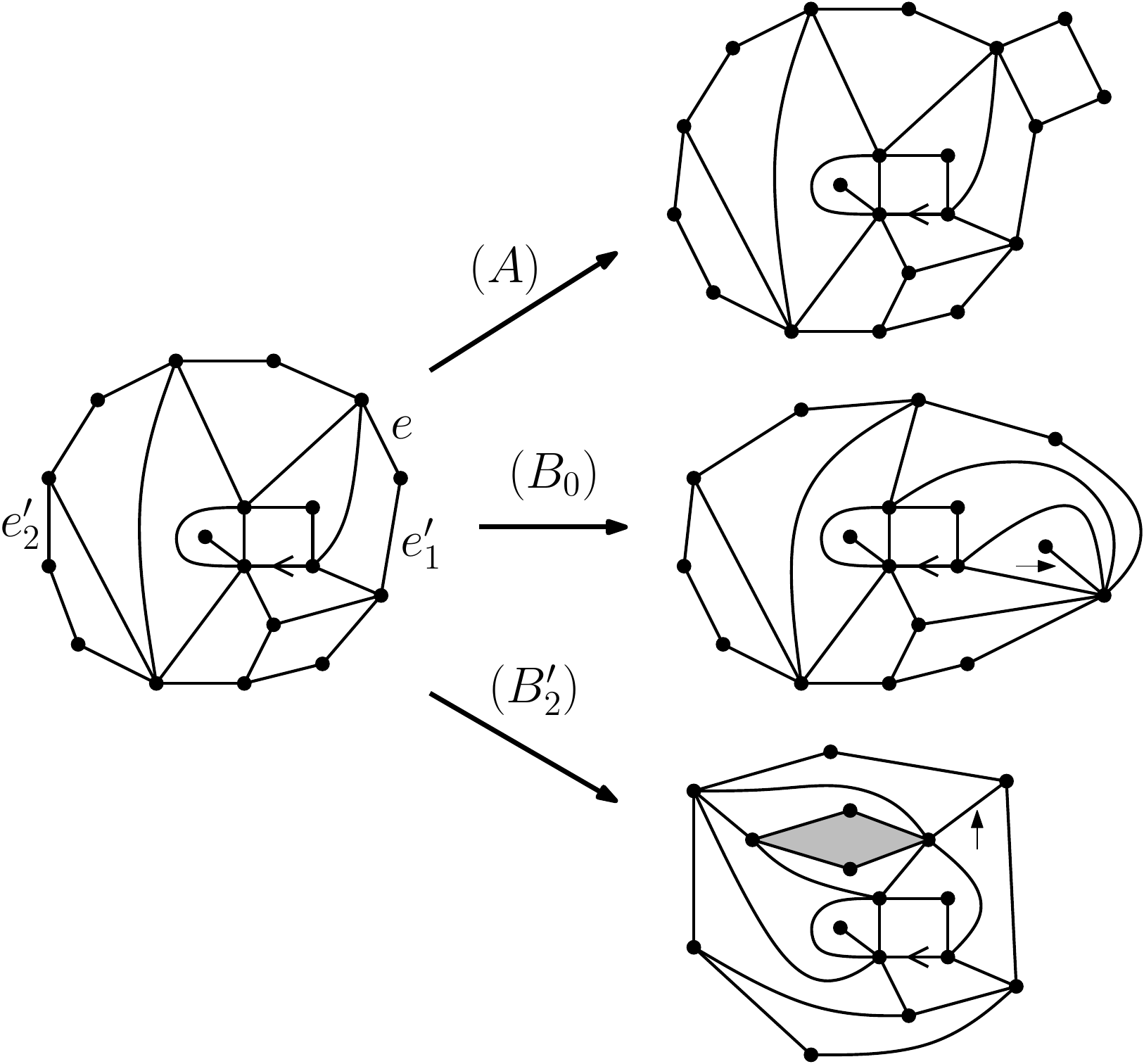}
 \caption{\label{lazy-peel}
A few possible steps in the lazy peeling algorithm. The peeled edge is the edge $e$ on the left side of the figure.
In case $(A)$ one glues a new quadrangle to $e$. In case $(B_0)$, $e$ is simply glued (in ``clockwise order'')
to the edge $e'_1$. In case $(B'_2)$, $e$ is glued in
counterclockwise order to the edge $e'_2$. In the latter case, the shaded region, whose (simple) boundary 
corresponds to the $4$ edges between $e$ and $e'_2$ in the left side, has to be filled in by
a quadrangulation with a (general) boundary of size $4$. In cases $(B_0)$ and $B'_2$, the ``glued edge''
is indicated by a small arrow in the right side. }
 \end{center}
 \vspace{-4mm}
 \end{figure}

Then suppose that at step $p$, we have $\qq_p=\qq$, where $\qq$ is 
a quadrangulation with a {\it simple} boundary of perimeter $2L$, $L\geq 1$. We construct 
$\qq_{p+1}$ in the following way. We choose an edge $e$ of the 
boundary of $\qq$, to be called the ``peeled edge'' at step $p+1$, and then:

\begin{enumerate}
\item[$(A)$] Either we glue a quadrangle to $e$, so that $\qq_{p+1}$ 
has a simple boundary of perimeter $2(L+1)$. The root
edge of $\qq_{p+1}$ is the same as the root edge of $\qq_p$.
\end{enumerate}
Or, for some $j\in\{0,1,\ldots,L-1\}$:
\begin{enumerate}
\item[$(B_j)$] We glue the edge $e$ to the edge $e'$ of the boundary of $\qq$
such that there are $2j$ edges of the boundary between $e$ and $e'$ in clockwise order. After 
this gluing the edges that were in the boundary of $\qq$ (except $e$ and $e'$) now bound two faces
of respective degrees $2j$ and $2(L-1-j)$. 
Following the device explained in Section \ref{sec:glue}, we then glue a quadrangulation $\wt q$ with a general boundary of size $2j$ in the face of degree $2j$ (bounded by the edges that were
in the  boundary of $\mathrm{q}$ between $e$ and $e'$ in clockwise order). In this way we obtain a quadrangulation with a simple boundary of size $2(L-1-j)$, which is $\qq_{p+1}$.
The root
edge of $\qq_{p+1}$ is the same as the root edge of $\qq_p$.
\item[$(B'_j)$] Or we do the same as in $(B_j)$ except that clockwise is replaced by counterclockwise. 
\end{enumerate}
See Fig.~\ref{lazy-peel} for an illustration of the different cases.

The preceding prescriptions have to be interpreted suitably in the two particular cases $j=0$ and $j=L-1$ of $(B_j)$ or $(B'_j)$.
\begin{enumerate}
\item[$\bullet$] If $j=0$, then we do not need to glue a quadrangulation $\wt q$ as explained above: in that case, we are just gluing together two adjacent edges.
\item[$\bullet$]
If $j=L-1$, then $\qq_{p+1}$ has a boundary of size $0$, which is interpreted by saying that
$\qq_{p+1}$ is a rooted planar quadrangulation. In that event, we take $K=p+1$ and the algorithm terminates. We in fact view
$q_K$ as a rooted and pointed quadrangulation, where the distinguished vertex $v$ is chosen as follows: 
in the last step of the construction, $v$ is the tail of $e$, assuming that $e$ is oriented clockwise in case $(B_j)$, counterclockwise in case $(B'_j)$
--- $v$ is also the unique vertex incident to both $e$ and $e'$, unless the length of the boundary is $2$ before the last step.
\end {enumerate}

Finally, we observe that, in case $(B_j)$ or $(B'_j)$ and when $1\leq j\leq L-1$, the gluing of the quadrangulation $\wt q$
in the face of degree $2j$ requires that we specify a distinguished edge on the boundary of this face, but this 
can be made in a prescribed manner whose choice is unimportant in the following discussion.

We now need to specify the probabilities of the different choices that are made in the preceding algorithm, and to this end
we introduce a few definitions. 
We set, for every $\ell\geq 0$,
$$h^\downarrow(\ell)= 2^{-2\ell} {2\ell\choose \ell},$$
and 
$$p_1=\frac{2}{3}\;,\quad p_{-i}= 2^{-2i+1}\,\frac{(2i-2)!}{(i-1)!(i+1)!}\;,\quad\hbox{for }i=1,2,\ldots.$$
If we also take $p_0=0$ and $p_i=0$ for every integer $i\geq 2$, elementary calculations show that the collection
$(p_i)_{i\in\Z}$ defines a probability distribution on $\Z$ with mean $0$. In particular,
\begin{equation}
\label{meanstep}
\frac{1}{2} \sum_{i=1}^\infty i\,p_{-i}=\frac{1}{3}.
\end{equation}

To specify the transition probabilities of the peeling algorithm, we first need 
to say how the peeled edge at step $p+1$ is chosen: At the present time, this choice can be made in an unimportant manner depending only on $\qq_p$. 
Then, conditionally on $\qq_0,\ldots,\qq_p$, and given that $\qq_p$ has boundary of size $2L$, the event $(A)$ occurs
with probability
$$p^{(L)}_1:=\frac{h^\downarrow(L+1)}{h^\downarrow(L)}\;p_1= \frac{1}{3}\,\frac{2L+1}{L+1},$$
and for $0\leq j\leq L-1$, either of the events $(B_j)$ or $(B'_j)$ occurs with probability
$$\frac{1}{2}\,p^{(L)}_{-j-1}:=\frac{1}{2}\,\frac{h^\downarrow(L-j-1)}{h^\downarrow(L)}\,p_{-j-1}$$
(one can check \cite{Bud,Peccot} that $p^{(L)}_1+\sum_{j=0}^{L-1} p^{(L)}_{-j-1}=1$). 
Furthermore, the quadrangulation $\wt q$ is chosen uniformly at random in the 
set of all (unpointed) rooted quadrangulations with a boundary of size $2j$. 
It follows from these prescriptions that the half-perimeter of
$\qq_p$ evolves like a Markov chain on $\Z_+$, which is stopped at the time $K$ when it hits $0$. The property $K<\infty$ a.s. (the algorithm terminates in finite time) can be derived from the fact that the 
half-perimeter process is a random walk with jump distribution $(p_i)_{i\in\Z}$ 
conditioned to hit zero before taking negative values (and stopped at that time): See the comments following Lemma 1 in \cite{Bud}.

If we perform the lazy peeling algorithm according to these probabilities,  the 
quadrangulation $\qq_K$ that we obtain  is a (rooted and pointed) Boltzmann quadrangulation \cite{Bud,Peccot}.
The key property that we will use is the following proposition, which is essentially
a special case of \cite[Proposition 1]{Bud}. 

\begin{proposition}
\label{glue-Boltzmann}
Let $T$ be a random variable with nonnegative integer values. Assume that $T\leq K$ and that
$T$ is a stopping time of the filtration generated
by $(\mathrm{q}_{n\wedge K})_{n\geq 0}$. 

Then, conditionally on the event $\{T<K\}$ and on $\qq_T$,  $\mathrm{q}_K$ has the same
distribution as the quadrangulation obtained by gluing to the outer face of $\mathrm{q}_T$ an 
independent Boltzmann 
(pointed) quadrangulation $\wt{\mathrm{q}}$ with a (general) boundary of size equal to the perimeter of
$\mathrm{q}_T$, with the convention that the root edge of the resulting map is the same
as the root edge of $\mathrm{q}_T$ and the distinguished vertex is the same as in
the quadrangulation $\wt{\mathrm{q}}$.
\end{proposition}

Again, for the gluing mentioned in the proposition to make sense, we should specify an
edge of the boundary of $\mathrm{q}_T$. The choice of this edge can however be made
in an arbitrary deterministic manner given $\mathrm{q}_T$, and then $\wt{\mathrm{q}}$ is
uniquely determined, which will be important for our purposes.

\medskip
\noindent{\bf Peeling algorithm for the UIPQ.} A variant of the preceding algorithm gives rise to the UIPQ or Uniform Infinite Planar 
Quadrangulation. We just replace the function $h^\downarrow(\ell)$ by
$$h^\uparrow(\ell)=\ell\, h^\downarrow(\ell)$$
in the definition of the probabilities of the different steps. Then the algorithm never stops (this 
is clear since $h^\uparrow(0)=0$ and thus, with the previous notation, 
the cases $(B_j)$ or $(B'_j)$ with $j=L-1$ never occur). Furthermore, the 
quadrangulation $\mathrm{q}_p$ converges locally a.s. to 
an infinite random planar map $Q^{(\infty)}$ which is the UIPQ. See \cite[Section 5]{Bud} or
\cite[Section 4.2]{Peccot} for more details in a more general setting. 
Note that $Q^{(\infty)}$ has a root edge but no distinguished vertex.
In that case, the Markov chain corresponding to the half-perimeter of $\qq_p$ is transient. 

\subsection{Peeling by layers}
\label{sec:peel-layers}

In the peeling algorithm described above, the peeled edge at step $p+1$ can be chosen
on the boundary of $\qq_p$ in an arbitrary way (depending on $\qq_p$). We will now describe 
a specific choice of the peeled edges, which produces the so-called peeling by
layers. This will lead us to define a sequence $\QQ_1,\QQ_2,\ldots,\QQ_\Xi$
of quadrangulations with a simple boundary, where $\Xi\geq 0$
is a (random) integer. 

To simplify the presentation in this section, we call {\it label} of a vertex of the boundary of $\qq_p$ its 
graph distance (in $\qq_p$) from
the root vertex.
We say that an edge $e$ of the boundary  of $\qq_p$ is of type $(i,i+1)$, resp.~$(i+1,i)$, if the labels of its ends listed in clockwise order around the boundary are $i$
and $i+1$, resp.~$i+1$ and $i$.

\begin{figure}[!h]
 \begin{center}
 \includegraphics[width=12cm]{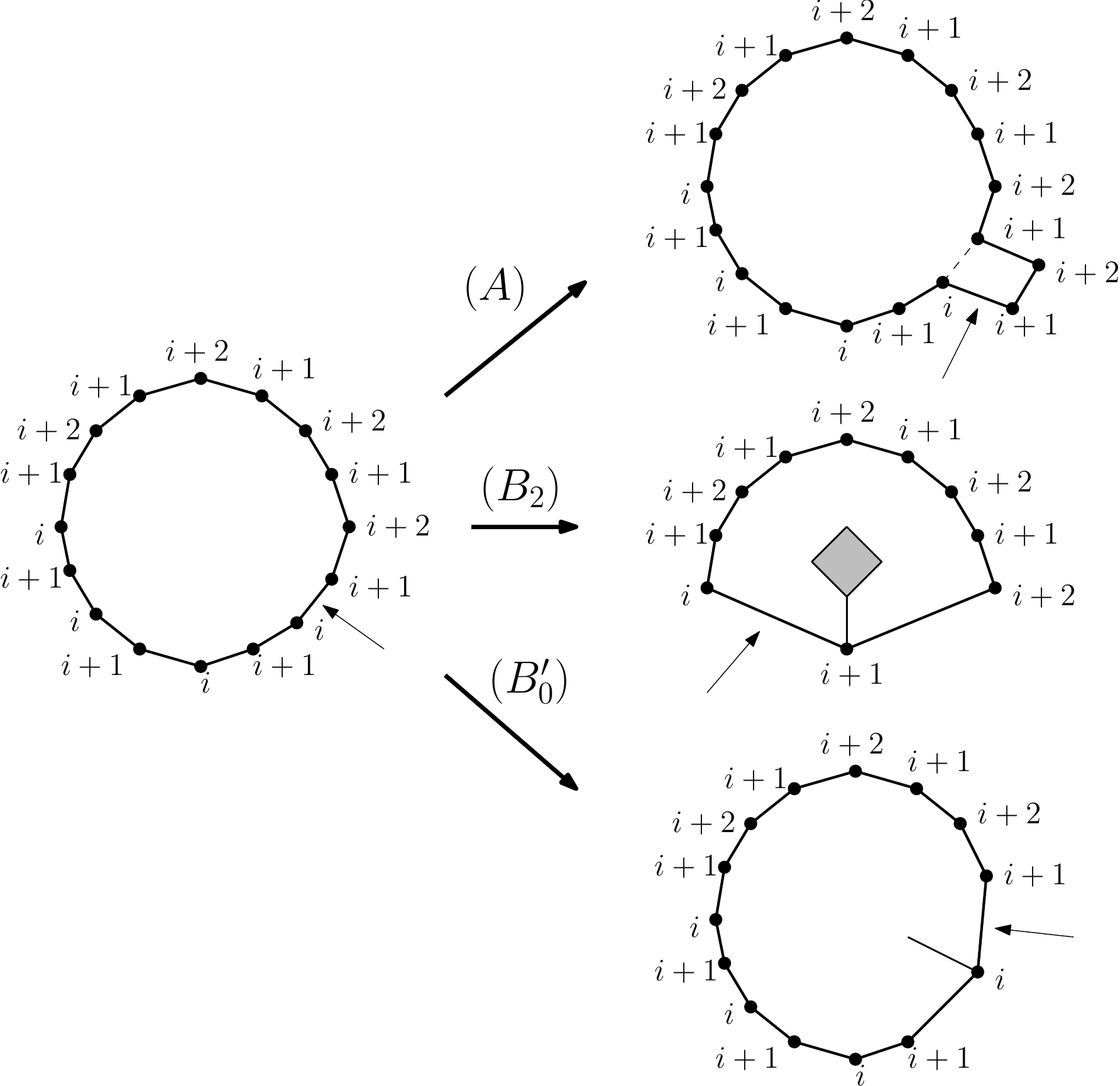}
 \caption{\label{peelayers}
Evolution of labels on the boundary in the peeling by layers in a few cases.
The peeled edge is indicated by an arrow in the left side (and in each case
in the right side the next edge to be peeled is also marked by an arrow).}
 \end{center}
 \vspace{-5mm}
 \end{figure}

Let us turn to the description of the peeling by layers algorithm. The choice of the peeled edge
at each step will be designed in such a way that, for every $0\leq p<K$, there exists an integer
$i\geq 0$ such that
\begin{enumerate}
\item[(a)] all vertices of the boundary of $\qq_p$ have label  $i,i+1$ or $i+2$.
\item[(b)] there are vertices of the boundary with label $i$;
\item[(c)] the edges of the boundary of type $(i+1,i+2)$ or $(i+2,i+1)$, if any, form a connected 
subset of the boundary.
\end{enumerate}
Supposing that these properties hold 
at step $p$, we choose the peeled edge at step $p+1$ in the following way:
\begin{enumerate}
\item[$\bullet$] If there is at least one vertex of the boundary of $\qq_p$ with label $i+2$, we peel the first edge
of type $(i+1,i)$ coming after the last vertex with label $i+2$ in clockwise order;
\item[$\bullet$] If there is no vertex of the boundary of $\qq_p$ with label $i+2$, we peel an 
edge of type $(i+1,i)$ chosen according to some rule given $\qq_p$. 
\end{enumerate}
See Fig.~\ref{peelayers} for an example.

It is clear that properties (a), (b), (c) stated above will be preserved at
each step of the algorithm (though the value of $i$ may increase to $i+1$
between steps $p$ and $p+1$) since they  hold at the initial step.
See Fig.~\ref{peelayers} for an example of the evolution of labels on the
boundary in the peeling by layers.

Since we start by peeling all edges of type $(1,0)$, there will
exist a first time $R_1$ at which the boundary contains only
edges of type $(1,2)$ or $(2,1)$ (unless the algorithm stops before this occurs,
in which case we take $R_1=\infty$). Similarly, for every $i\geq1$, we let
$R_i$ be the first time at which the boundary only contains edges of
type $(i,i+1)$ or $(i+1,i)$, with the convention that $R_i=\infty$
if this never occurs. Then there exists a (random) integer $\Xi\geq 0$ such that
$R_i<\infty$ if and only if $i\leq \Xi$. For $1\leq i\leq \Xi$, we set $\QQ_i=\qq_{R_i}$,
and we also let $Q=\qq_K$ be the rooted and pointed Boltzmann quadrangulation
obtained at the end of the algorithm. 
As previously $\rho$ and $\xi$ are respectively the
root vertex and the distinguished vertex of $Q$.
The following properties then hold.
\begin{enumerate}
\item[(a)] The sequence $\QQ_1,\QQ_2,\ldots,\QQ_\Xi$  is a deterministic function of $Q$. 
\item[(b)] We have $\Xi\leq \dg^Q(\rho,\xi)\leq \Xi+1$.
\item[(c)] For every $1\leq i\leq \Xi$, the set of all vertices of $\QQ_i$ that do not lie
on the boundary is identified canonically to a subset $\mathcal{J}_i$
of $V(Q)$. For any $v\in \mathcal{J}_i$, any path from $v$ to $\xi$ in $Q$
must visit a vertex at graph distance at most $i$ from $\rho$. Conversely,
for any vertex $v$ of $Q$ that does not belong to $\mathcal{J}_i$, there is a path from
$v$ to $\xi$ that visits only vertices whose graph distance from $\rho$
is at least $i$.
\end{enumerate}
Let us briefly explain why properties (a) and (b) are satisfied (we omit the argument for (c)). Let $\mathcal{A}$ be the set of all 
dual edges of $Q$ (each edge $e$ of $Q$ corresponds to a dual edge between the faces of
$Q$ that are incident to $e$). Fix $i\geq 1$, and let $\mathcal{A}_i$
consist of all dual edges associated with a primal edge connecting two vertices at
graph distance less than or equal to $i$ from the root vertex. 
Let $\mathcal{A}'_i$ be obtained by adding to $\mathcal{A}_i$ all
edges dual to a primal edge that lies in a connected component of the complement of
$\mathcal{A}_i$ not containing the distinguished vertex of $Q$. Then
one checks that $i\leq\Xi$ if and only $\mathcal{A}'_i\not =\mathcal{A}$. Furthermore, if this
property holds, $\QQ_i$ can be obtained informally by starting from the root face (to the right of the root edge) and 
then gluing quadrangles along the dual edges in $\mathcal{A}'_i$, and
keeping the ``same'' root as in $Q$  (as explained in \cite[Chapter 3]{Peccot}, with each connected subset 
of $\mathcal{A}$ containing the dual root edge one can associate a
rooted planar map by this gluing procedure). To give a less informal construction, view $\mathcal{A}'_i$ as the edge set of 
a planar map, and let $\mathfrak{f}$ be the face of this planar map containing $\xi$.  Note that the 
boundary of $\mathfrak{f}$ is a simple cycle $\mathcal{C}$,
and that vertices of this cycle have degree $2$ or $4$. Choose a vertex $\mathrm{v}_0$ in the face $\mathfrak{f}$, and let
$\mathcal{A}''_i$ be obtained from $\mathcal{A}'_i$ by adding,
for each vertex $\mathrm{v}$  of degree $2$ in $\mathcal{C}$, two edges 
connecting $\mathrm{v}$ to $\mathrm{v}_0$ in $\mathfrak{f}$ (of course in such a way that these edges do not cross). Then 
$\QQ_i$ is the dual planar map to the planar map with edge set $\mathcal{A}''_i$. 

\begin{figure}[!h]
 \begin{center}
 \includegraphics[width=14cm]{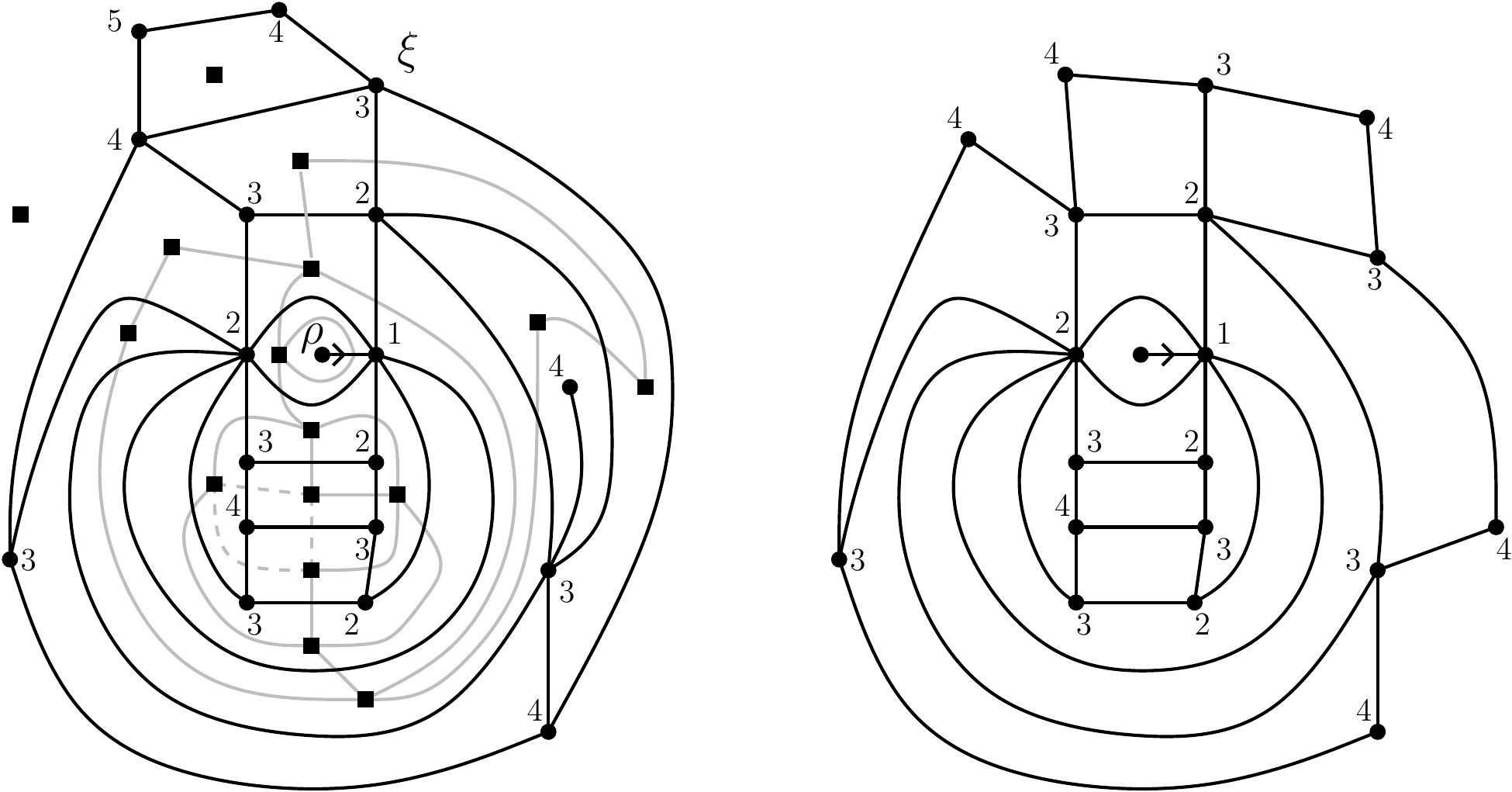}
 \caption{\label{lazy/hull}
 Left, the rooted and pointed planar quadrangulation $Q$. Vertices of $Q$ are represented by small black disks, while vertices
 of the dual map (one in each face of $Q$) are represented by small black squares. The figures correspond to
 distances from the root vertex of $Q$. The dual edges in grey form the set $\mathcal{A}_3$ and $\mathcal{A}'_3$
 is obtained by adding to $\mathcal{A}_3$ the three dashed dual edges. Right, the lazy hull $\QQ_3$ of radius $3$. 
The reader may verify that $\QQ_3$ is determined by the map $\mathcal{A}'_3$ (and the knowledge of the dual root)
as explained above.
 }
 \end{center}
 \vspace{-5mm}
 \end{figure}

We call $\QQ_i$ the lazy hull of radius $i$ in $Q$. See Fig.\ref{lazy/hull} for an illustration with $i=3$. 

\bigskip

\medskip
\noindent{\bf The UIPQ case.} The above considerations also make sense for the peeling algorithm associated
with the UIPQ $Q^{(\infty)}$. In that case however, we can make sense of the whole sequence 
$\QQ_1,\QQ_2,\ldots$ of lazy hulls. The analogs of properties (a) and (c) above
remain valid (in (c), ``path from $v$ to $\xi$'' has to be replaced by 
``path from $v$ to infinity'').

\subsection{Asymptotics for the perimeter and volume of hulls}
\label{sec:asymp-hull}

Let $Q$ and $Q^{(\infty)}$ be as previously. For every $1\leq i\leq \Xi$, we write $H_i$
and $V_i$ respectively for the half-perimeter and volume (number of inner faces) of the
lazy hull of radius $i$ in $Q$. Similarly, for every $i\geq 1$, we write $H^{(\infty)}_i$
and $V^{(\infty)}_i$ respectively for the half-perimeter and volume of the
lazy hull of radius $i$ in $Q^{(\infty)}$. By convention, 
$H_0=H^{(\infty)}_0=0$ and $V_0=V^{(\infty)}_0=0$. 

We let $(\Upsilon_t)_{t\geq 0}$ denote a centered stable L\'evy process
with index $3/2$ and no positive jumps, whose law
is specified by the equality
$$\E[\exp(\lambda \Upsilon_t)] = \exp(2\,t\,\lambda^{3/2}),$$
for every $\lambda,t\geq 0$. Write $(\Upsilon^\uparrow_t)_{t\geq 0}$
for the process $\Upsilon$
 conditioned to stay positive (see \cite[Chapter VII]{Be}).
 Let $\tau_1,\tau_2,\ldots$ be a measurable enumeration of the jump times
 of $\Upsilon^\uparrow$, and let $\theta_1,\theta_2,\ldots$ be a sequence of 
 nonnegative i.i.d. random variables with density 
\begin{equation}
\label{densi-theta}
\frac{1}{\sqrt{2\pi x^5}}\;\exp(-\frac{1}{2x})
\end{equation}
 on $(0,\infty)$. Assume that this sequence is independent of 
 $\Upsilon^\uparrow$.
 We then define a process 
 $(\mathcal{V}^\uparrow_t)_{t\geq 0}$ by setting, for every $t\geq 0$,
 $$\mathcal{V}^\uparrow_t=\frac{1}{2}\sum_{\{i:\tau_i\leq t\}} (\Delta \Upsilon^\uparrow_{\tau_i})^2\;\theta_i.$$

\begin{proposition}
\label{conv-peri-vol}
We have the convergence in distribution in the Skorokhod topology
$$\Big(n^{-2}H^{(\infty)}_{\lfloor n t\rfloor}, n^{-4} V^{(\infty)}_{\lfloor n t\rfloor}\Big)_{t\geq 0}
\build{\la}_{n\to\infty}^{\rm (d)} \Big(\mathcal{Y}^\uparrow_t,\mathcal{U}^\uparrow_t\Big)_{t\geq 0},$$
where the limiting processes are defined by
$$\mathcal{Y}^\uparrow_t:=\Upsilon^\uparrow_{\Psi_t}\;,\quad \mathcal{U}^\uparrow_t:= \mathcal{V}^\uparrow_{\Psi_t}\,,$$
with
$$\Psi_t=\inf\Big\{s\geq 0: \int_0^s \frac{\mathrm{d}r}{\Upsilon^\uparrow_r} >t\Big\}.$$
\end{proposition}

The result of the proposition is very close to similar results proved in \cite{CLG-peeling}
(see in particular Theorem 2 in \cite{CLG-peeling}). Unfortunately, it does not seem
easy to deduce Proposition \ref{conv-peri-vol} from the results of \cite{CLG-peeling},
and we postpone the details
of the proof to Appendix B below. 

Our next goal is to obtain an analog of the previous proposition for the 
perimeter and volume process of the lazy hull of a Boltzmann quadrangulation
under a suitable conditioning. We will derive this result from Proposition \ref{conv-peri-vol}
by an absolute continuity argument. 

Write $(S^\downarrow_k)_{0\leq k\leq K}$ for the Markov chain giving for every $0\leq k\leq K$ the half-perimeter 
of $\qq_k$ in the peeling algorithm of  Section \ref{sec:peel},
with the convention $S^\downarrow_K=0$, and $(S^\uparrow_k)_{k\geq 0}$
for the analogous Markov chain for the peeling algorithm of the UIPQ described at the
end of Section \ref{sec:peel}. We also let $(S_k)_{k\geq 0}$ be
a random walk with jump distribution $(p_i)_{i\in\Z}$ started from $S_0=2$,
and set $T=\inf\{k\geq 0:S_k\leq 0\}$. Then both functions $h^\downarrow$
and $h^\uparrow$ are harmonic on $\{1,2,\ldots\}$ for the random walk $S$, as it
was first observed by Budd \cite[Proposition 3, Corollary 1]{Bud} (see also formula (3.5) and Lemma 11 in \cite{Peccot}).
Futhermore, the form
of the transition probabilities in the peeling algorithm shows that $S^\downarrow$
can be viewed as the $h^\downarrow$-transform of the random walk $S$
killed when it hits $\Z_{-}$,
meaning more precisely that, for any $n\geq 0$ and any function $F:\Z^{n+1}\la \R_+$,
$$\E[F(S_0^\downarrow,S^\downarrow_1,\ldots,S^\downarrow_n)\;\mathbf{1}_{\{n\leq K\}}]
= \E\Big[ \frac{h^\downarrow(S_n)}{h^\downarrow(S_0)}\,F(S_0,S_1,\ldots,S_n)\;
\mathbf{1}_{\{n\leq T\}}\Big],$$
with the convention $h^\downarrow(i)=0$ for $i<0$. Similarly, $S^\uparrow$
can be viewed as the $h^\uparrow$-transform of the random walk $S$
killed when it hits $\Z_{-}$, so that
$$\E[F(S_0^\uparrow,S^\uparrow_1,\ldots,S^\uparrow_n)]
= \E\Big[ \frac{h^\uparrow(S_n)}{h^\uparrow(S_0)}\,F(S_0,S_1,\ldots,S_n)\;
\mathbf{1}_{\{n\leq T\}}\Big],$$
with the same convention $h^\uparrow(i)=0$ if $i<0$ (and we recall that
$h^\uparrow (0)=0$). By comparing the last two displays, 
and recalling that $h^\uparrow(\ell)=\ell h^\downarrow(\ell)$, we obtain that
\begin{equation}
\label{up-down}
\E[F(S_0^\downarrow,S^\downarrow_1,\ldots,S^\downarrow_n)\;\mathbf{1}_{\{n< K\}}]
=\E\Big[ \frac{2}{S_n^\uparrow} \,F(S_0^\uparrow,S^\uparrow_1,\ldots,S^\uparrow_n)\Big].
\end{equation}
In fact, we can reinforce this identity, noting that apart from the choices of 
cases $(A)$, $(B_j)$ or $(B'_j)$, the other random choices (in particular for the
quadrangulations with a simple boundary ``filling in the holes'') are made 
in exactly the same way in the algorithm
generating a Boltzmann quadrangulation and in the algorithm for the UIPQ.
We get more generally that, for any $n\geq 1$, for any function $F:\Z^{n+1}\la \R_+$,
and for any nonnegative function $G_n$ of the quadrangulations $\qq_1,\ldots,\qq_n$
obtained in the first $n$ steps,
$$\E[G_n\,F(S_0^\downarrow,S^\downarrow_1,\ldots,S^\downarrow_n)\;\mathbf{1}_{\{n< K\}}]
=\E\Big[ \frac{2}{S_n^\uparrow} \,G_n\,F(S_0^\uparrow,S^\uparrow_1,\ldots,S^\uparrow_n)\Big].$$
By a standard argument, we can generalize this identity to the case where the integer $n$
is replaced 
by a stopping time of the discrete filtration generated by $(\qq_{n\wedge K})_{n\geq 0}$. 

The preceding discussion is valid for any peeling algorithm of the type described in Section
\ref{sec:peel}, but we now specialize to the peeling by layers of Section \ref{sec:peel-layers},
and we consider the processes $(H_n)_{0\leq n\leq \Xi}$ and $(V_n)_{0\leq n\leq \Xi}$.
By construction, we have  $H_n=S^\downarrow_{R_n}$ for every $0\leq n\leq \Xi$, and a similar relation links $H^{(\infty)}$
to $S^\uparrow$. The previous considerations, and the
fact that $\{n\leq \Xi\}=\{R_n<K\}$, entail that, for any nonnegative
function $F$ on $(\Z^2)^n$, 
\begin{equation}
\label{abs-cont}
\E[F((H_1,V_1),\ldots,(H_n,V_n))\,\mathbf{1}_{\{n\leq\Xi\}}]
= \E\Big[\frac{2}{H_n^{(\infty)}}\;
F((H^{(\infty)}_1,V^{(\infty)}_1),\ldots,(H^{(\infty)}_n,V^{(\infty)}_n))\Big].
\end{equation}
Combining this equality with Proposition \ref{conv-peri-vol} leads to the following result.

\begin{proposition}
\label{conv-peri-vol2}
Let $t>0$. The distribution of the pair of processes
$$\Big(n^{-2}H_{\lfloor n s\rfloor}, n^{-4} V_{\lfloor n s\rfloor}\Big)_{0\leq s\leq t}$$
under the conditional probability $\P(\cdot\mid \Xi\geq \lfloor n t\rfloor)$
converges as $n\to\infty$ to the law of a pair of processes 
$$\Big(\mathcal{Y}^\downarrow_s,\mathcal{U}^\downarrow_s\Big)_{0\leq s\leq t}$$
such that, for any nonnegative measurable function $F$ on the
Skorokhod space $\D([0,t],\R^2)$,
$$\E\Big[ F\Big(\Big(\mathcal{Y}^\downarrow_s,\mathcal{U}^\downarrow_s\Big)_{0\leq s\leq t}\Big)\Big]=
\frac{t^2}{2}\,\E\Big[ \frac{1}{\mathcal{Y}^\uparrow_t}\,F\Big(\Big(\mathcal{Y}^\uparrow_s,\mathcal{U}^\uparrow_s\Big)_{0\leq s\leq t}\Big)\Big].$$
\end{proposition}

\proof Let $F$ be a bounded continuous function on $\D([0,t],\R^2)$ such that $0\leq F\leq 1$. 
By \eqref{abs-cont},
$$\E\Big[F\Big(\Big(n^{-2}H_{\lfloor n s\rfloor}, n^{-4} V_{\lfloor n s\rfloor}\Big)_{0\leq s\leq t}\Big)\,
\mathbf{1}_{\{\Xi\geq \lfloor n t\rfloor\}}\Big]
= \E\Big[ \frac{2}{H_{\lfloor n t\rfloor}^{(\infty)}}\,F\Big(\Big(n^{-2}H^{(\infty)}_{\lfloor n s\rfloor}, n^{-4} V^{(\infty)}_{\lfloor n s\rfloor}\Big)_{0\leq s\leq t}\Big)\Big].$$
From Proposition \ref{conv-peri-vol}, we get
\begin{equation}
\label{tec-conv-peri-vol}
\liminf_{n\to\infty}
n^2\,\E\Big[F\Big(\Big(n^{-2}H_{\lfloor n s\rfloor}, n^{-4} V_{\lfloor n s\rfloor}\Big)_{0\leq s\leq t}\Big)\,
\mathbf{1}_{\{\Xi\geq \lfloor n t\rfloor\}}\Big]
\geq \E\Big[ \frac{2}{\mathcal{Y}^\uparrow_t}\,F\Big(\Big(\mathcal{Y}^\uparrow_s,\mathcal{U}^\uparrow_s\Big)_{0\leq s\leq t}\Big)\Big].
\end{equation}
On the other hand, we claim that
\begin{equation}
\label{claim29}
\lim_{n\to\infty} n^2\,\P(\Xi\geq \lfloor n t\rfloor)=  \frac{4}{t^2}= \E\Big[\frac{2}{\mathcal{Y}^\uparrow_t}\Big]
\end{equation}
To verify this claim, first note that the distribution of $\mathcal{Y}^\uparrow_t$ has density
\begin{equation}
\label{distri-Y}
h(z)=\frac{2}{t^3}\,\sqrt{\frac{z}{\pi}}\, \exp(-\frac{z}{t^2})
\end{equation}
by \cite[Proposition 1.2]{Hull} (beware that the normalization in \cite{Hull} is different). 
The second equality in \eqref{claim29} immediately follows.
We then use the following lemma, whose proof is postponed after the end of the proof of
Proposition \ref{conv-peri-vol2}.

\begin{lemma}
\label{asymp-dist-root}
For a rooted and pointed quadrangulation with Boltzmann distribution,
\begin{align*}
\P(\dg(\rho,\xi)\geq n)&=\frac{4}{n(n+2)},\qquad \hbox{if }n\geq 2,\\
\P(\dg(\rho,\xi)\geq 1)&=\frac{5}{6}.
\end{align*}
\end{lemma}

Recalling that $\Xi\leq \dg(\rho,\xi)\leq \Xi+1$, we immediately get from the lemma that we have also 
$$\lim_{n\to\infty} n^2\,\P(\Xi\geq n)= 4,$$
which gives the first equality in  \eqref{claim29}. 

Finally, if we combine \eqref{tec-conv-peri-vol} with the same result with $F$ replaced by $1-F$, using \eqref{claim29}, we get that 
$$\lim_{n\to\infty}
n^2\,\E\Big[F\Big(\Big(n^{-2}H_{\lfloor n s\rfloor}, n^{-4} V_{\lfloor n s\rfloor}\Big)_{0\leq s\leq t}\Big)\,
\mathbf{1}_{\{\Xi\geq \lfloor n t\rfloor\}}\Big]
= \E\Big[ \frac{2}{\mathcal{Y}^\uparrow_t}\,F\Big(\Big(\mathcal{Y}^\uparrow_s,\mathcal{U}^\uparrow_s\Big)_{0\leq s\leq t}\Big)\Big].$$
Then we just have to divide by the first equality in \eqref{claim29} to get the desired result. \endproof

\smallskip
\noindent{\it Proof of Lemma \ref{asymp-dist-root}.}
We use a version of Schaeffer's bijection which is different from the one presented
briefly in Section \ref{sec:quadr} (see e.g. \cite{LGM} for this other version). Recall from Section \ref{sec:quadr}
the definition of a labeled plane tree, and write $\T^\circ$
for the collection of all labeled plane trees having at least one edge (we exclude the case 
where the tree consists only of its root). Then there is a bijection between the set
$$\bigcup_{n=1}^\infty \Q_n^\bullet$$
of all rooted and pointed quadrangulations of the sphere and $\T^\circ \times\{0,1\}$. 
Furthermore, if $((\tau,(\ell_u)_{u\in V(\tau)}),\epsilon)\in \T^\circ \times\{0,1\}$
and $q$ is the associated rooted and pointed quadrangulation, the
distance between the root vertex and the distinguished vertex
of $q$ is equal to
$$-\min_{u\in V(\tau)}\ell_u+\epsilon.$$
Suppose now that $q$ is a Boltzmann quadrangulation. Then one easily
verifies that the associated tree $\tau$ is a Galton-Watson tree with geometric
offspring distribution with parameter $1/2$ conditioned to have at least one edge, 
that, conditionally on $\tau$, the labels are chosen uniformly among possible choices,
and that $\epsilon=0$ or $1$ with probability $1/2$ independently of $(\tau,(\ell_u)_{u\in V(\tau)})$.
By \cite[Proposition 2.4]{CD}, we have, for every integer $n\geq 1$,
$$\P\Big(\min_{u\in V(\tau)}\ell_u \leq -n\Big)=\frac{4}{(n+1)(n+2)}.$$

It follows that, for $n\geq2$,
$$\P(\dg(\rho,\xi)\geq n)=\frac{1}{2}\,\P\Big(\min_{u\in V(\tau)}\ell_u \leq -n\Big)+\frac{1}{2}\, \P\Big(\min_{u\in V(\tau)}\ell_u \leq -n+1\Big)
=\frac{4}{n(n+2)}$$
and similarly $\P(\dg(\rho,\xi)\geq 1)=(\frac{1}{2} \times \frac{2}{3}) + \frac{1}{2}=\frac{5}{6}$.  \hfill$\square$

\smallskip
\rem Since $\Xi\leq \dg(\rho,\xi)\leq \Xi+1$ and we know that $\P(\dg(\rho,\xi)\geq n)\sim 4/n^2$
as $n\to\infty$, the result of Proposition \ref{conv-peri-vol2}
remains valid if we replace $\P(\cdot\mid \Xi\geq \lfloor n t\rfloor)$ by 
$\P(\cdot\mid \dg(\rho,\xi)\geq \lfloor n t\rfloor)$. 
We will use such remarks implicitly in what follows. 

\subsection{Asymptotics for Boltzmann quadrangulations}
\label{subsec:asBol}

Our goal in this section is to restate Proposition \ref{conv-peri-vol2} in a different manner
showing that the limiting process can be written as a functional 
of the Brownian snake, and that the convergence holds jointly with that of Corollary \ref{convBmapB}.
We start with the convergence of the process of volumes of hulls.

\begin{lemma}
\label{conv-volumes}
Let $\delta>0$. The distribution of 
$(n^{-2}V_{\lfloor r\sqrt{n}\rfloor})_{0\leq r\leq \delta}$
under $\P(\cdot \mid \dg(\rho,\xi)>\delta \sqrt{n})$ converges 
in the Skorokhod sense to the distribution of the process
$$\Big(\frac{1}{2}\int_0^\sigma \mathrm{d} s\,\mathbf{1}_{\{\un W_s \leq W_* +r\sqrt{3/2}\}}\Big)_{0\leq r\leq \delta}$$
under $\N_0(\cdot\mid W_*<-\delta\sqrt{3/2})$. Furthermore, this convergence holds jointly
with that of Corollary \ref{convBmapB}.
\end{lemma}

\proof For every integer $k\leq \dg(\rho,\xi)$
write $\mathcal{H}_k$ for the collection of all
vertices $v$ of $Q$ such that any path from $v$ to $\xi$
visits a vertex whose graph distance from $\rho$ is
(less than or)
equal to $k$. Recall also the notation $\mathcal{J}_k$ introduced in property (c)
of the lazy hulls in Section \ref{sec:peel-layers}. From the latter property, we have, for
every $1\leq k\leq \Xi$,
$$\mathcal{H}_{k-1}\subset \mathcal{J}_k\subset \mathcal{H}_{k}.$$

On the other hand, recall the notation $\wt\tau_Q$ for the (re-rooted at $\xi$) tree
associated with $Q$ via Schaeffer's bijection, and notice
that vertices of $\wt\tau_Q$ are identified with vertices of $Q$. 
As in Section \ref{sec:Bolquad}, write $\mathcal{K}_k$ for the collection of those vertices $v$ of $\wt\tau_Q$
such that 
the geodesic (in the tree $\wt\tau_{Q}$) from $v$ to $\xi$ visits at least one vertex
with label at most $\wt L_*+k-1$ (recall the notation $\wt L_*$ for the minimal label on
$\wt\tau_Q$). The properties of Schaeffer's bijection recalled in Section \ref{sec:quadr} imply 
that $\mathcal{K}_k\subset \mathcal{H}_k\subset \mathcal{K}_{k+1}
\cup\{\rho\}$, for every $1\leq k\leq \dg(\rho,\xi)$.

Let $r\in(0,\delta]$. By Lemma \ref{lem-tec-hull}, the distribution of $n^{-2}\#\mathcal{K}_{\lfloor r\sqrt{n}\rfloor}$
under $\P(\cdot\mid \dg(\rho,\xi)
>\lfloor \delta\sqrt{n}\rfloor)$ converges as $n\to\infty$
to the distribution of
$$\frac{1}{2}\int_0^\sigma \mathrm{d}s\,\mathbf{1}_{\{\un{W}_s\leq W_*+r\sqrt{3/2}\}}$$
under $\N_0(\cdot\mid W_*<-\delta\sqrt{3/2})$, and this convergence holds jointly with that
of Corollary
\ref{convBmapB}. By the preceding considerations, this holds also if we replace
$\#\mathcal{K}_{\lfloor r\sqrt{n}\rfloor}$ by $\#\mathcal{J}_{\lfloor r\sqrt{n}\rfloor}$. 
On the other hand, Euler's formula, and the fact that the size of the boundary of 
the lazy hull of radius $\lfloor r\sqrt{n}\rfloor$ is of order $n$ (by Proposition \ref{conv-peri-vol2})
show that $n^{-2}(\#\mathcal{J}_{\lfloor r\sqrt{n}\rfloor} - V_{\lfloor r\sqrt{n}\rfloor})$
under $\P(\cdot\mid \dg(\rho,\xi)
>\lfloor \delta\sqrt{n}\rfloor)$
tends to $0$ in probability. It follows that, for every fixed $r\geq 0$, the distribution of $n^{-2}V_{\lfloor r\sqrt{n}\rfloor}$
under $\P(\cdot\mid \dg(\rho,\xi)
>\lfloor \delta\sqrt{n}\rfloor)$ converges as $n\to\infty$
to the distribution of
$$\frac{1}{2}\int_0^\sigma \mathrm{d}s\,\mathbf{1}_{\{\un{W}_s\leq W_*+r\sqrt{3/2}\}}$$
under $\N_0(\cdot\mid W_*<-\delta\sqrt{3/2})$, and that this convergence holds jointly with that
of Corollary
\ref{convBmapB}. Since we already know from Proposition \ref{conv-peri-vol2} that the distribution of $(n^{-2}V_{\lfloor r\sqrt{n}\rfloor})_{0\leq r\leq \delta}$
under $\P(\cdot \mid \dg(\rho,\xi)>\delta \sqrt{n})$ converges 
in the Skorokhod sense, the desired result follows easily. \endproof

\begin{proposition}
\label{conv-exit}
The distribution of 
$(n^{-1} H_{\lfloor r\sqrt{n}\rfloor},n^{-2}V_{\lfloor r\sqrt{n}\rfloor})_{0\leq r\leq \delta}$ 
under $\P(\cdot \mid \dg(\rho,\xi)>\delta \sqrt{n})$
converges to the distribution of 
$$\Big(\z_{W_*+r\sqrt{3/2}},\frac{1}{2}\int_0^\sigma \mathrm{d} s\,\mathbf{1}_{\{\un W_s \leq W_* +r\sqrt{3/2}\}}\Big)_{0\leq r\leq \delta}
$$ under $\N_0(\cdot\mid W_*<-\delta\sqrt{3/2})$.
Furthermore, this convergence in distribution holds jointly with that of Corollary \ref{convBmapB}.
\end{proposition}

\proof For every $r\in[0,\delta]$, set 
$$\mathcal{X}_r:=\frac{1}{2}\int_0^\sigma \mathrm{d} s\,\mathbf{1}_{\{\un W_s \leq W_* +r\sqrt{3/2}\}}.$$
We first verify that the distribution of 
$(\z_{W_*+r\sqrt{3/2}},\mathcal{X}_r)_{0\leq r\leq \delta}$ 
under $\N_0(\cdot\mid W_*<-\delta\sqrt{3/2})$ is the same as the distribution
of the process $(\mathcal{Y}^\downarrow_r,\mathcal{U}^\downarrow_r)_{0\leq r\leq \delta}$. 
By comparing Proposition \ref{conv-peri-vol2} and Lemma \ref{conv-volumes}, we already know from that 
the distribution of $(\mathcal{X}_r)_{0\leq r\leq \delta}$ 
under $\N_0(\cdot\mid W_*<-\delta\sqrt{3/2})$  is the same as the distribution of
$(\mathcal{U}^\downarrow_r)_{0\leq r\leq \delta}$. 

Let us start by verifying that $(\z_{W_*+r\sqrt{3/2}})_{0\leq r\leq \delta}$ 
(under $\N_0(\cdot\!\mid\! W_*<-\delta\sqrt{3/2})$) has the same distribution as 
$(\mathcal{Y}^\downarrow_r)_{0\leq r\leq \delta}$. We observe that the 
process $(\z_{-r})_{r>0}$ is Markov under $\N_0$ with the transition probabilities 
of the continuous-state branching process (CSBP) with branching
mechanism $\psi(u)=\sqrt{8/3}\,u^{3/2}$ (see e.g. \cite[Section 2.5]{ALG}). By scaling, the same property holds 
for  $(\z_{-r\sqrt{3/2}})_{r>0}$ with $\psi(u)=2\,u^{3/2}$. Notice that 
$-W_*=\inf\{r>0: \z_{-r}=0\}$. It follows (we omit
a few details here) that, under $\N_0(\cdot\mid W_*<\delta\sqrt{3/2})$
and conditionally on $\z_{W_*+\delta\sqrt{3/2}}=z$, the
process $(\z_{W_*+(\delta-r)\sqrt{3/2}})_{0\leq r\leq \delta}$ is distributed as
the CSBP with branching mechanism $\psi(u)=2u^{3/2}$ started at $z$ and conditioned 
to become extinct at time $\delta$. Note that, by Corollary \ref{limit-law},
the law of $\z_{W_*+\delta\sqrt{3/2}}$ under $\N_0(\cdot\mid W_*<\delta\sqrt{3/2})$ 
has density 
$$g(z)=\frac{1}{\delta}\sqrt{\frac{1}{\pi z}}\,\exp(-\frac{z^2}{\delta^2}).$$
On the other hand, from \cite[Proposition 4.4]{Hull}, we also know 
that conditionally on $\mathcal{Y}^\uparrow_\delta=z$, the process 
$(\mathcal{Y}^\uparrow_{\delta-r})_{0\leq r\leq \delta}$ is distributed as a 
CSBP with branching mechanism $\psi(u)=2u^{3/2}$ started at $z$ and conditioned 
to become extinct at time $\delta$. From the absolute 
continuity relation in Proposition \ref{conv-peri-vol2}, the same
property holds for the process
$(\mathcal{Y}^\downarrow_{\delta-r})_{0\leq r\leq \delta}$. 
So to prove that $(\z_{W_*+(\delta-r)\sqrt{3/2}})_{0\leq r\leq \delta}$ 
(under $\N_0(\cdot\mid W_*<-\delta\sqrt{3/2})$) has the same distribution as 
$(\mathcal{Y}^\downarrow_{\delta-r})_{0\leq r\leq \delta}$, it only remains to
verify that the density of $\mathcal{Y}^\downarrow_{\delta}$ is also 
given by the function $g$ of the last display. But this follows from 
\eqref{distri-Y} and the absolute 
continuity relation in Proposition \ref{conv-peri-vol2}. 

From the construction of $\mathcal{U}^\uparrow$, and Proposition \ref{conv-peri-vol2},
we know that the conditional distribution of $(\mathcal{U}^\downarrow_r)_{0\leq r\leq \delta}$
given $(\mathcal{Y}^\downarrow_r)_{0\leq r\leq \delta}$ is that of the process
$$[0,\delta]\ni r\mapsto \frac{1}{2}\sum_{\{i:r_i\in[0,r]\}} (\Delta \mathcal{Y}^\downarrow_{r_i})^2\;\theta_i,$$
where $r_1,r_2,\ldots$ is a measurable enumeration of the jumps of 
$\mathcal{Y}^\downarrow$, and $\theta_1,\theta_2,\ldots$ is a
sequence of i.i.d. random variables with density given in \eqref{densi-theta},
which is independent of $\mathcal{Y}^\downarrow$. We then note that the same property holds
for the conditional distribution of $(\mathcal{X}_r)_{0\leq r\leq \delta}$
given $(\z_{W_*+r\sqrt{3/2}})_{0\leq r\leq \delta}$
under $\N_0(\cdot\mid W_*<-\delta\sqrt{3/2})$, as a straightforward
consequence of Corollary 4.9 in \cite{Hull} (alternatively, one could also use
\cite[Theorem 40]{ALG}).

We conclude that  the distribution of 
$(\z_{W_*+r\sqrt{3/2}},\mathcal{X}_r)_{0\leq r\leq \delta}$ 
under $\N_0(\cdot\mid W_*<-\delta\sqrt{3/2})$ is the same as the distribution
of the process $(\mathcal{Y}^\downarrow_r,\mathcal{U}^\downarrow_r)_{0\leq r\leq \delta}$,
and thus (by Proposition \ref{conv-peri-vol2}) the convergence stated in Proposition
\ref{conv-exit} holds. 

It remains to 
prove that this convergence holds jointly with that of Corollary \ref{convBmapB}.
To this end, it is enough to show that, for every $r\in(0,\delta]$, there is a measurable
function $\Phi$ from $\D([0,\delta],\R)$ into $\R$ such that, 
a.s.,
$$\mathcal{Y}^\downarrow_r=\Phi\Big((\mathcal{U}^\downarrow_s)_{0\leq s\leq \delta}\Big).$$
Indeed the same representation then holds for $\z_{W_*+r\sqrt{3/2}}$
as a function of $(\mathcal{X}_s)_{0\leq s\leq \delta}$, and
since we know that the convergence of Lemma \ref{conv-volumes}
holds jointly with that of Corollary \ref{convBmapB}, a simple tightness argument
will show that the same holds for the convergence of the proposition.

To begin with, consider a process $(\mathcal{V}_t)_{t\geq 0}$ 
defined exactly as the process $(\mathcal{V}^\uparrow_t)_{t\geq 0}$
of Section \ref{sec:asymp-hull}, except that $\Upsilon^\uparrow$
is replaced by the L\'evy process $\Upsilon$ (with L\'evy measure
$(3\sqrt{\pi}/2)\,|x|^{-5/2}\,\mathbf{1}_{\{x<0\}}\,\mathrm{d}x$). Clearly,
$(\mathcal{V}_t)_{t\geq 0}$ is a subordinator, and 
a few lines of calculations show that its L\'evy measure
is $c_0\,y^{-7/4}\,\mathrm{d}y$, where $c_0:=3.2^{-7/4}\,\Gamma(3/4)$.
Write 
$$\phi(\ve)=\int_\ve^\infty c_0\,y^{-7/4}\,\mathrm{d}y=\frac{4}{3}\,c_0\,\ve^{-3/4},$$
for every $\ve >0$. Then, for every fixed $r> 0$ and $\alpha\in(0,r)$,
$$\lim_{\ve\to 0} \phi(\ve)^{-1}\#\{s\in[r-\alpha,r]:\Delta\mathcal{V}_s>\ve\}=\alpha,\quad\hbox{a.s.}$$
It is not hard to verify that the above limit holds simultaneously for every $r\geq 0$
and every $\alpha\in(0,r)$ outside a single set of probability zero. Using the absolute continuity relations
between the laws of $\Upsilon$ and of $\Upsilon^\uparrow$, it follows that the previous
limit also holds if $\mathcal{V}$ is replaced by $\mathcal{V}^\uparrow$, simultaneously for every $r> 0$
and every $\alpha\in(0,r)$, a.s. Next recall from Proposition \ref{conv-peri-vol} the definition of $\mathcal{U}^\uparrow$ as a time change
of $\mathcal{V}^\uparrow$, with a time change given by
$$\Psi_t=\inf\{s\geq 0: \int_0^s \frac{\mathrm{d}r}{\Upsilon^\uparrow_r} >t\}=\int_0^t \mathrm{d}r\,\mathcal{Y}^\uparrow_r.$$
It follows that we have also
$$\lim_{\ve\to 0} \phi(\ve)^{-1}\#\{s\in[r-\alpha,r]:\Delta\mathcal{U}^\uparrow_s>\ve\}=\int_{r-\alpha}^r\mathrm{d}s\,\mathcal{Y}^\uparrow_s,$$
for every $r>0$ and $\alpha\in(0,r)$, a.s.
Thanks to the absolute continuity relation in Proposition \ref{conv-peri-vol2}, the preceding limit still
holds if the pair $(\mathcal{Y}^\uparrow,\mathcal{U}^\uparrow)$ is replaced by $(\mathcal{Y}^\downarrow,\mathcal{U}^\downarrow)$.
Hence, for every $r\in(0,\delta]$, a.s.,
$$\mathcal{Y}^\downarrow_r=\lim_{\alpha\to 0}\frac{1}{\alpha}\Bigg(\lim_{\ve\to 0}\phi(\ve)^{-1}\#\{s\in[r-\alpha,r]:\Delta\mathcal{U}^\downarrow_s>\ve\}\Bigg),$$
which gives the desired representation of $\mathcal{Y}^\downarrow_r$
as a function of $(\mathcal{U}^\downarrow_s)_{0\leq s\leq \delta}$. This completes the proof.
\endproof

\smallskip
\rem An argument similar to the end of the preceding proof 
is used in \cite[Section 2.3]{Budz}.

\section{Convergence of quadrangulations with glued boundary}
\label{conv-glued-bdry}

The goal of this section is to give a representation of the free Brownian disk
with glued boundary in terms of the measure  $\N^{*,z}_0$ introduced in Section \ref{sec:posexc}.
We recall from Section \ref{sec:Bdglued} the notation $\F_z^\dagger$  for the
 distribution of the free Brownian disk with perimeter $z$ and boundary glued into a
 single point. Also recall from Section \ref{sec:consmms} the notation $\mathcal{L}^\bullet(\omega)$ for the pointed
 measure metric space associated with a snake trajectory.
 
 \begin{theorem}
 \label{identif-glued-unpointed}
 For any nonnegative measurable function $G$ on
 $\M^\bullet$,
 $$\F_z^\dagger(G)= \N^{*,z}_0\Big(G(\ll^\bullet(W))\Big).$$
 \end{theorem}

The proof of this theorem, which is given at the very end of this section, requires a few preliminary results. 
Let us outline our general strategy. We start from a a Boltzmann quadrangulation $Q$. For $\delta>0$, conditionally on the event where the
perimeter of the lazy hull of radius $\lfloor \delta\sqrt{n}\rfloor$ is of order $n$, Proposition \ref{glue-Boltzmann} allows us to ``embed'' a quadrangulation $Q_n$
with a boundary  of size of order $n$ in $Q$. When $\delta$ is small,
the quadrangulation $Q^\dagger_n$ with ``glued boundary'', which is derived from $Q_n$ as in Section \ref{sec:convBdisk},
is close to $Q$ in the Gromov-Hausdorff-Prokhorov sense after rescaling distances by $n^{-1/2}$ (Lemma  \ref{GHdistance-discrete}). Then, on one hand, the scaling limit of $Q^\dagger_n$ is given 
in terms of the measures $\F^{\bullet,\dagger}_z$ (Proposition \ref{conv-glued}, which is basically a consequence
of Corollary \ref{conv-quadb-glued}). On the other hand, the convergence of 
rescaled Boltzmann quadrangulations to the free Brownian map shows that the Boltzmann quadrangulation $Q$,
under the preceding conditioning depending on $n$ and the usual rescaling, is close to the free Brownian map $(\mm,D,\bv)$ under the
conditioning that $\z_{W_*+\delta\sqrt{3/2}}$ is of order $1$, which in turn (when $\delta$ is small) is close to the random
measure metric space $\ll^{\bullet\bullet}(\mathrm{tr}_{W_*+\delta}(\omega))$ under the same conditioning
(Lemma \ref{GHdistance-truncated}). The distribution of $\mathrm{tr}_{W_*+\delta}(\omega)$ can be written 
in terms of the measures $\N^{*,z}_0$ via Proposition \ref{law-trunc}. By combining all these observations, we
arrive at an expression of $\F^{\bullet,\dagger}_z$ in terms of the measures $\N^{*,z}_0$ (Proposition \ref{identif-glued})
from which it is then easy to derive Theorem \ref{identif-glued-unpointed}.

We keep the notation introduced in the previous section. In particular, $\mathcal{Q}_1,\mathcal{Q}_2,\ldots,\mathcal{Q}_\Xi$ 
is the lazy hull process associated with the Boltzmann (pointed) quadrangulation $Q$.
We  fix $\delta>0$, and we will argue conditionally on the event $\{\Xi\geq \lfloor\delta \sqrt{n}\rfloor\}$.
By Proposition \ref{glue-Boltzmann}, conditionally on the lazy hull of radius $\lfloor \delta\sqrt{n}\rfloor$,
$Q$ is obtained by gluing to this hull an independent (pointed) 
Boltzmann quadrangulation
with boundary of size equal to $2H_{\lfloor \delta\sqrt{n}\rfloor}$, and we denote
this quadrangulation with a boundary by $Q_n$ (as already mentioned, the gluing requires that 
we specify an edge on the boundary of the hull, which can be done
in a deterministic way given the hull, and then $Q_n$ is
uniquely determined). Notice that in this gluing procedure, every vertex of the boundary of the hull corresponds
to a vertex incident to the outer face of $Q_n$, but this correspondence is not one-to-one since several vertices 
of the boundary of the hull may correspond to the same vertex in the boundary of the outer face of $Q_n$.
However, after gluing $Q_n$ to the hull, these vertices will correspond to a single vertex of $Q$.
Also, each vertex of $Q_n$ not on the boundary corresponds to a (unique) vertex of $Q$ that does not belong to the 
hull.  We write $d_n$ for the graph distance on the vertex set
$V(Q_n)$, and $|Q_n|$ for the number of inner faces of $Q_n$.  
We also write $\mu_n$ for the counting measure on $V(Q_n)$.

For notational convenience, we agree that $H_k=0$ if $k> \Xi$, so that an event 
of the type $\{H_{\lfloor \delta\sqrt{n}\rfloor}\in [n,(1+\ve) n]\}$ is contained in $\{\Xi\geq \lfloor \delta\sqrt{n}\rfloor\}$. 

We can view $Q_n$ as a 
submap of a rooted and pointed planar quadrangulation $Q_n^\dagger$, 
which is constructed in exactly in the same way as $B^\dagger_{(k)}$
was constructed from $B_{(k)}$ in Section \ref{sec:convBdisk}.
We use the notation $d_n^\dagger$ for the graph distance on
$V(Q_n^\dagger)$, and also notice that $\mu_n$
can be viewed as a measure on $V(Q_n^\dagger)$. As was the case for $B^\dagger_{(k)}$
 in Section \ref{sec:convBdisk}, $Q_n^\dagger$ has two distinguished points, namely
 the extra vertex $\varpi$ and the distinguished vertex of $Q_n$. Recall the notation $\F_r^{\bullet,\dagger}$  for the
 distribution of the free pointed Brownian disk with perimeter $r$ and boundary glued into a
 single point.

\begin{proposition}
\label{conv-glued}
The sequence of random $2$-pointed measure metric spaces
$$\Big(V(Q^\dagger_n), \sqrt{\frac{3}{2}}\,n^{-1/2} \,d^\dagger_n, 2n^{-2}\mu_n\Big)\;,\ \hbox{under } \P(\cdot\mid H_{\lfloor \delta\sqrt{n}\rfloor}\in [n,(1+\ve) n]),$$
converges in distribution  as $n\to\infty$ in $\M^{\bullet\bullet}$.
The limiting distribution is
$$c_{\delta,\ve}\,\int_1^{1+\ve} \frac{\mathrm{d}r}{\sqrt{r}}\,e^{-r/\delta^2}\, \F^{\bullet,\dagger}_r.$$
where $c_{\delta,\ve}$ is the appropriate normalizing constant. 
\end{proposition}

\proof 
If $r>0$, and $(k_n)$ is a sequence of integers such that
$n^{-1}k_n$ tends to $r$ as $n\to\infty$, Corollary \ref{conv-quadb-glued}
and scaling arguments relying on \eqref{scaling-disk2} show that
$$\Big(V(B^\dagger_{(k_n)}),\sqrt{\frac{3}{2}}\,n^{-1/2}\,\dg^\dagger, 2n^{-2}\mu_{(k_n)}\Big) \build{\la}_{k\to\infty}^{\rm(d)} 
(\D^{\bullet,\dagger}_r, D^{\partial,\dagger}, \bv^{\bullet,\dagger}_r)
$$
and the limit depends continuously on $r$. It follows that, if $\Phi$ is bounded and continuous
on $\M^{\bullet\bullet}$, 
\begin{equation}
\label{conv-glued1}
\lim_{n\to\infty} \E[\Phi(V(B^\dagger_{(\lfloor nr\rfloor)}),\sqrt{\frac{3}{2}}\,n^{-1/2}\,\dg^\dagger, 2n^{-2}\mu_{(\lfloor nr\rfloor)})]
= \E[\Phi(\D^{\bullet,\dagger}_r, D^{\partial,\dagger}, \bv^{\bullet,\dagger}_r)]= \F^{\bullet,\dagger}_r(\Phi)
\end{equation}
and the limit is uniform when $r$ varies over a compact subset of $(0,\infty)$. 
Then observe that, conditionally on $H_{\lfloor \delta\sqrt{n}\rfloor}=k$, the
random $2$-pointed measure metric space
$$\Big(V(Q_n^\dagger), \sqrt{\frac{3}{2}}\,n^{-1/2} \,d_n^\dagger,2n^{-2} \mu_n\Big)$$
has the same distribution as 
$$\Big(V(B^\dagger_{(k)}),\sqrt{\frac{3}{2}}\,n^{-1/2}\,\dg^\dagger, 2n^{-2}\mu_{(k)}\Big).$$
The desired result follows by writing
\begin{align*}
&\E\Big[\Phi\Big(V(Q^\dagger_n), \sqrt{\frac{3}{2}}\,n^{-1/2} \,d^\dagger_n, 2n^{-2}\mu_n\Big)\,\mathbf{1}_{\{
H_{\lfloor \delta\sqrt{n}\rfloor}\in [n,(1+\ve) n] \}}\Big]\\
&\qquad=n \int_1^{n^{-1}(\lfloor(1+\ve) n\rfloor+1)} \mathrm{d}r\, \P(H_{\lfloor \delta\sqrt{n}\rfloor}=\lfloor nr\rfloor)
\times \E[\Phi(V(B^\dagger_{(\lfloor nr\rfloor)}),\sqrt{\frac{3}{2}}\,n^{-1/2}\,\dg^\dagger, 2n^{-2}\mu_{(\lfloor nr\rfloor)})]
\end{align*}
and using \eqref{conv-glued1} together with Proposition \ref{conv-exit} and Corollary \ref{limit-law}. \endproof

\medskip
The next step is to introduce a random metric measure space that is a
candidate for describing the Brownian disk with boundary glued into 
a single point. With the notation of
Section \ref{sec:consmms}, this is the random ($2$-pointed) metric 
measure space  $\ll^{\bullet\bullet}(\mathrm{tr}_{W_*+\delta}(\omega))$ under 
the probability measure $\N^{[-\delta]}_0=\N_0(\cdot\mid W_*<-\delta)$. 

To give a somewhat different presentation of this space, we introduce (under 
$\N^{[-\delta]}_0$)  the subtree $\t^{(\delta)}_\zeta$ of
$\t_\zeta$ that consists of all vertices $a\in \t_\zeta$ such that
$Z_b>W_*+\delta$ for every strict ancestor $b$ of $a$.
By definition, the
``boundary'' $\partial \t^{(\delta)}_\zeta$ consists of all 
$a\in \t^{(\delta)}_\zeta$ such that $Z_a=W_*+\delta$.
We define a pseudo-metric
$D^{(\delta)}$ on $\t^{(\delta)}_\zeta$ by setting, for every $a,b\in \t^{(\delta)}_\zeta$,
$$D^{(\delta)}(a,b)=\min\{D(a,b),Z_a+Z_b-2(W_*+\delta)\}.$$
It is easy to verify that $D^{(\delta)}$ satisfies the triangle inequality. 
We also observe that $D^{(\delta)}(a,b)=0$ if and only if at least one of the following two
conditions holds:
\begin{enumerate}
\item[(i)] $D(a,b)=0$;
\item[(ii)] $a$ and $b$ both belong to the boundary $\partial \t^{(\delta)}_\zeta$.
\end{enumerate}
Similarly as in Section \ref{sec:consmms}, we let
$\mm^{(\delta)}$ be the quotient space of $\t^{(\delta)}_\zeta$ 
for the equivalence relation $a\sim b$ if and only if $D^{(\delta)}(a,b)=0$
and we equip $\mm^{(\delta)}$ with the metric induced by $D^{(\delta)}$.
Writing $\Pi^{(\delta)}$
for the canonical projection from $\t^{(\delta)}_\zeta$ onto $\mm^{(\delta)}$
we define the volume measure $\bv^{(\delta)}$ on 
$\mm^{(\delta)}$ as the image of the restriction of the volume measure on $\t_\delta$
to $\t^{(\delta)}_\zeta$ under $\Pi^{(\delta)}$. We note that 
$\Pi^{(\delta)}(\partial \t^{(\delta)}_\zeta)$ consists of a single 
point denoted by $x_{(\delta)}$ and called the boundary point
of $\mm^{(\delta)}$. We also notice that if $x= \Pi^{(\delta)}(a)$
we have $D^{(\delta)}(x,x_{(\delta)})=Z_a-(W_*+\delta)$. 

We view $(\mm^{(\delta)},D^{(\delta)},\bv^{(\delta)})$ as a $2$-pointed measure metric space,
whose distinguished points are  the boundary point $x_{(\delta)}$ and $\Pi^{(\delta)}(p_\zeta(0))$ (in this order). 

\begin{lemma}
\label{identi-mms-trunc} $\N^{[-\delta]}_0$ a.s., the $2$-pointed measure metric space 
$(\mm^{(\delta)},D^{(\delta)},\bv^{(\delta)})$ coincides with $\ll^{\bullet\bullet}(\mathrm{tr}_{W_*+\delta}(\omega))$.
\end{lemma}

\proof This is a straightforward consequence of our definitions. We note that the
genealogical tree of the snake trajectory $\mathrm{tr}_{W_*+\delta}(\omega)$ is
identified canonically with $\t^{(\delta)}_\zeta$. Modulo this identification, we write 
$D'(a,b)$, $a,b\in \t^{(\delta)}_\zeta$ for the pseudo-metric defined
via formula \eqref{formulaD} applied to $\mathrm{tr}_{W_*+\delta}(\omega)$ (instead of $\omega$).
The proof then boils down to verifying that $D'(a,b)=D^{(\delta)}(a,b)$ 
for every $a,b\in \t^{(\delta)}_\zeta$. The verification of this identity is left to the reader. \endproof

\rem Define the hull of radius $\delta$ in $(\mm,D)$ as the complement of
the connected component of the complement of the ball $B(x_*,\delta)$ that contains the second distinguished point 
$x_0$ (this makes sense under $\N^{[-\delta]}_0$, since
$D(x_*,x_0)=-W_*$). The metric space $(\mm^{(\delta)},D^{(\delta)})$ can be viewed
as the closure of the complement of the hull of radius $\delta$ in $(\mm,D)$, with 
boundary identified to a single point.

\smallskip
The next lemma bounds  the Gromov-Hausdorff-Prokhorov
distance between $(\mm,D,\bv)$ and the ``truncated space'' $(\mm^{(\delta)},D^{(\delta)},\bv^{(\delta)})$.

\begin{lemma}
\label{GHdistance-truncated}
Under $\N^{[-\delta]}_0$, we have
$${d}^{(2)}_{\rm GHP}\Big((\mm,D,\bv),(\mm^{(\delta)},D^{(\delta)},\bv^{(\delta)})\Big)
\leq 3(2(\delta + K_\delta)+\kappa_\delta),$$
where
$$K_\delta:=\sup\{ Z_a-W_*: a\in \t_\zeta\backslash \t_\zeta^{(\delta)}\}\,,
\quad\kappa_\delta:=\int_0^\sigma \mathrm{d}s\,\mathbf{1}_{\{p_\zeta(s)\notin \t^{(\delta)}_\zeta\}}.$$
\end{lemma}

The quantities $K_\delta$ and $\kappa_\delta$ have the following geometric interpretation in terms of the
$2$-pointed measure metric space $(\mm,D,\bv)$: $K_\delta$ is the maximal distance from $x_*$ in the hull of radius $\delta$ (defined in the preceding
remark), and
$\kappa_\delta$ is the volume of this hull. 

\proof We first construct a correspondence $\mathcal{C}_{(\delta)}$ between 
$(\mm,D)$ and $(\mm^{(\delta)},D^{(\delta)})$
whose distortion is bounded above by $2(\delta + K_\delta)$. We note that
two vertices of $\t^{(\delta)}_\zeta$ not belonging to $\partial \t^{(\delta)}_\zeta$ are identified in
$\mm^{(\delta)}$ if and only if they are identified in $\mm$. We
then define the correspondence  $\mathcal{C}_{(\delta)}$  as follows. First if $x\in\mm$ is of the form $x=\Pi(a)$ with $a\in \t^{(\delta)}_\zeta
\backslash \partial \t^{(\delta)}_\zeta$, then $x$ is in correspondence with the
``same'' point $\Pi^{(\delta)}(a)$ of $\mm^{(\delta)}$. Then, if 
$x=\Pi(a)$ with $a\in \t_\zeta\backslash \t^{(\delta)}_\zeta$ or $a\in\partial \t^{(\delta)}_\zeta$,
$x$ is in correspondence with the boundary point $x_{(\delta)}$. 

Let us bound the distortion of $\mathcal{C}_{(\delta)}$. We need to bound
$|D(x,y) - D^{(\delta)}(x',y')|$ when $(x,x')$ and $(y,y')$ belong
to $\mathcal{C}_{(\delta)}$. 
Consider first the case when 
$x=\Pi(a)$ and $y=\Pi(b)$ with $a,b\in \t^{(\delta)}_\zeta
\backslash \partial \t^{(\delta)}_\zeta$. Then, by construction,
$D^{(\delta)}(x',y')=D^{(\delta)}(a,b)\leq D(a,b)=D(x,y)$. On the
other hand, we have 
$$Z_a+Z_b-2(W_*+\delta)= (Z_a-W_*)+(Z_b-W_*)-2\delta \geq D(a,b) - 2\delta$$
because $Z_a-W_*=D(a,a_*)$ and $Z_b-W_*=D(b,a_*)$. It follows that
we have also $D^{(\delta)}(x',y')\geq D(x,y)-2\delta$, so 
that $|D(x,y) - D^{(\delta)}(x',y')|\leq 2\delta$ in that case. 

Suppose then that $x=\Pi(a)$ and $y=\Pi(b)$ with $a\in \t^{(\delta)}_\zeta
\backslash \partial \t^{(\delta)}_\zeta$ and 
$b\in (\t_\zeta \backslash \t^{(\delta)}_\zeta)\cup \partial \t^{(\delta)}_\zeta$. 
In that case $y'=x_{(\delta)}$.
Then we have again $D^{(\delta)}(x',x_{(\delta)})=Z_a-(W_*+\delta)\leq D(x,y)$
by the continuous cactus bound \eqref{ctscactus}.
Furthermore,
$$D(x,y)\leq D(a_*,a)+D(a_*,b)=Z_a-W_*+Z_b-W_*
\leq (D^{(\delta)}(x',x_{(\delta)}) +\delta)+K_\delta.$$
 We therefore get the
bound $|D(x,y) - D^{(\delta)}(x',y')|\leq K_\delta +\delta$. 

Finally, in the case where $a$ and $b$ both belong
to $(\t_\zeta \backslash \t^{(\delta)}_\zeta)\cup \partial \t^{(\delta)}_\zeta$,
we have $D^{(\delta)}(x,y)=0$ and 
$$D(x,y)\leq D(a_*,a)+D(a_*,b)=Z_a-W_*+Z_b-W_*\leq 2K_\delta.$$

Combining all three cases, we have always $|D(x,y) - D^{(\delta)}(x',y')|\leq 2(\delta + K_\delta)$
which gives the bound $\mathrm{dis}(\cc_{(\delta)})\leq 2(\delta + K_\delta)$. 

In order to apply Lemma \ref{GHP-corrresp}, we consider the measure 
$\nu$ on $\mm \times \mm^{(\delta)}$ defined by
$$\int \varphi(x,y)\,\nu(\mathrm{d}x,\mathrm{d}y)= \int \varphi(x,\theta(x))\,\bv(\mathrm{d}x),$$
where $\theta(x)=\Pi^{(\delta)}(a)$ if $x=\Pi(a)$ with $a\in \t^{(\delta)}_\zeta\backslash \partial\t^{(\delta)}_\zeta$,
and $\theta(x)=x_{(\delta)}$ otherwise
(we use the fact that two vertices of $\t^{(\delta)}_\zeta\backslash\partial \t^{(\delta)}_\zeta$ are identified in
$\mm^{(\delta)}$ if and only if they are identified in $\mm$). Clearly, 
$\nu$ is supported on $\cc_{(\delta)}$. Moreover,
with the notation of Lemma \ref{GHP-corrresp}, we have $\pi_*\nu=\bv$ and 
$\pi'_*\nu=\bv^{(\delta)} + \kappa_\delta\,\delta_{(x_{(\delta)})}$. An application
of Lemma \ref{GHP-corrresp} then completes the proof. 
\endproof

\medskip
We now give a discrete version of the preceding lemma. Recall
from the beginning
of this section
that the quadrangulations $Q_n$ and $Q_n^\dagger$ are well-defined on the event $\{ \Xi
\geq \lfloor \delta\sqrt{n}\rfloor\}$. The vertex set $V(Q_n)$ is identified to
a subset of $V(Q)$, and $V(Q_n^\dagger)=V(Q_n)\cup\{\varpi)$. We note that, by the construction of 
hulls in the preceding section, any vertex of $\partial Q_n$ is identified in the gluing procedure
to a vertex of $Q$ which is at graph distance either $\lfloor \delta\sqrt{n}\rfloor$ or
$\lfloor \delta\sqrt{n}\rfloor+1$ from the root vertex $\rho$. 

\begin{lemma}
\label{GHdistance-discrete}
On the event $\{ \Xi
\geq \lfloor \delta\sqrt{n}\rfloor\}$,
\begin{align*}
&{d}^{(2)}_{\mathrm{GHP}}\Big((V(Q_n^\dagger),\sqrt{\frac{3}{2}}\,n^{-1/2} \,d_n^\dagger,2n^{-2}\mu_n),(V(Q),\sqrt{\frac{3}{2}}\,n^{-1/2} \,\mathrm{d}_{\mathrm gr},2n^{-2}\mu)\Big)\\
&\qquad\leq 3\Big(2\sqrt{\frac{3}{2}}\,n^{-1/2} (\lfloor \delta\sqrt{n}\rfloor+ K_{n,\delta}+1)+ 2n^{-2}\,\kappa_{n,\delta}\Big),
\end{align*}
where
$$K_{n,\delta}=\max\{\mathrm{d}_{\mathrm gr}(\rho,v):v\in V(Q)\backslash V(Q_n)\}\,,\quad
\kappa_{n,\delta}=\#V(Q)-\#V(Q_n) + \#\partial Q_n.$$
\end{lemma}

\proof This is similar to the preceding proof. We construct a correspondence $\cc_n$ between
$V(Q)$ and $V(Q_n^\dagger)$ by declaring that every vertex of $V(Q_n)$
(viewed as a subset of $V(Q)$) is in correspondence with itself in $V(Q_n^\dagger)$,
and all vertices of $V(Q)\backslash  V(Q_n)$ are in correspondence with $\varpi$.
To bound the distortion of this correspondence, we note that, for
every $v,v'\in V(Q)$,
$$\begin{array}{ll}
\mathrm{d}_{\mathrm gr}(v,v')-2(\lfloor \delta\sqrt{n}\rfloor+1)   \leq  
d_n^\dagger(v,v')  \leq \mathrm{d}_{\mathrm gr}(v,v') +2\qquad&\hbox{if }v,v'\in V(Q_n)\\
\noalign{\smallskip}
\mathrm{d}_{\mathrm gr}(v,v')-2\,K_{n,\delta} -2 \leq  
d_n^\dagger(v,\varpi)  \leq \mathrm{d}_{\mathrm gr}(v,v') +2\qquad
&\hbox{if }v\in V(Q_n), v'\in V(Q)\backslash V(Q_n)\\
\noalign{\smallskip}
\mathrm{d}_{\mathrm gr}(v,v') \leq 2\,K_{n,\delta}\qquad&\hbox{if }v,v'\in 
V(Q)\backslash V(Q_n).
\end{array}$$
This shows that the distortion of the correspondence $\cc_n$ is bounded above by 
$2(2n/3)^{-1/2} (\lfloor \delta\sqrt{n}\rfloor+ K_{n,\delta}+1)$. On the 
other hand, we can construct a measure $\nu_n$ on the product 
$V(Q)\times V(Q_n^\dagger)$ as the image of the counting mesure $\mu$ on $V(Q)$
under the mapping $u\mapsto (u,\theta_n(u))$, where $\theta_n(u)=u$ if $u$ is a vertex 
of $Q_n$ not lying on the boundary $\partial Q_n$, and $\theta_n(u)=\varpi$
otherwise. An application of Lemma \ref{GHP-corrresp} then gives the desired result.
\endproof

In the next lemma, and in what follows, we consider implicitly that an event
of the form $\{\z_{W_*+\delta}\in[1,1+\alpha]\}$ is contained in
$\{W_*\leq -\delta\}$. 

\begin{lemma}
\label{tech-GHdistance}
Let $\beta>0$. We have for every $\alpha >0$,
$$\limsup_{n\to\infty} \P(2n^{-2}\,\kappa_{n,\delta} >\beta\mid H_{\lfloor \delta\sqrt{n}\rfloor}\in [n,(1+\alpha) n]) 
\leq \N_0(\kappa_{\delta\sqrt{3/2}}\geq \beta
\mid \z_{W_*+\delta\sqrt{3/2}}\in[1,1+\alpha]),$$
and similarly
$$\limsup_{n\to\infty} \P(n^{-1/2}\,K_{n,\delta}>\beta\mid H_{\lfloor \delta\sqrt{n}\rfloor}\in [n,(1+\alpha) n]) 
\leq \N_0(K_{\delta\sqrt{3/2}}\geq \beta
\mid \z_{W_*+\delta\sqrt{3/2}}\in[1,1+\alpha]).$$
\end{lemma}

\proof As we already noticed in the proof of Lemma \ref{conv-volumes}, the properties of the lazy hull ensure
that $V(Q)\backslash V(Q_n)\subset \mathcal{H}_{\lfloor\delta\sqrt{n}\rfloor}$,
where we recall that the notation  $\mathcal{H}_k$ refers to the collection of all
vertices $v$ of $Q$ such that any path from $v$ to $\xi$
visits a vertex whose graph distance from $\rho$ is
(less than or)
equal to $k$. We also saw that
 $\mathcal{H}_k\subset \mathcal{K}_{k+1}
\cup\{\rho\}$, and Lemma \ref{lem-tec-hull} implies that
 the distribution of $(2n^2)^{-1}\#\mathcal{K}_{\lfloor\delta\sqrt{n}\rfloor}$
under $\P(\cdot\mid \Xi
\geq \lfloor \delta\sqrt{n}\rfloor)$ converges as $n\to\infty$
to the distribution of
$$\int_0^\sigma \mathrm{d}s\,\mathbf{1}_{\{\un{W}_s\leq W_*+\delta\sqrt{3/2}\}}
=\kappa_{\delta\sqrt{3/2}}$$
under $\N_0(\cdot\mid W_*<-\delta\sqrt{3/2})$.  On the other hand, arguments very similar to the proof
of Lemma \ref{lem-tec-hull} show that
the distribution of 
$$n^{-1/2} \max\{\mathrm{d}_{\mathrm gr}(\rho,v):
v\in  \mathcal{K}_{\lfloor\delta\sqrt{n}\rfloor}(Q)\}$$
under $\P(\cdot\mid \Xi
\geq \lfloor \delta\sqrt{n}\rfloor)$ converges as $n\to\infty$
to the distribution of
$$\sup\{\wh W_s-W_*: s\in[0,\sigma], \un{W}_s\leq W_*+\delta\sqrt{3/2}\}=K_{\delta\sqrt{3/2}}$$
under $\N_0(\cdot\mid W_*<-\delta\sqrt{3/2})$.
The latter two convergences in distribution hold jointly with the 
convergence of $n^{-1}H_{\lfloor \delta\sqrt{n}\rfloor}$ to $\z_{W_*+\delta\sqrt{3/2}}$
in Proposition \ref{conv-exit}. 
Both assertions of the lemma now follow.  \endproof

\begin{proposition}
\label{conv-truncated}
Let $\alpha>0$ and $\ve>0$. Let $F$ be a bounded Lipschitz
continuous function on $\M^{\bullet\bullet}$. 
Then, for all $\delta>0$ small enough,
\begin{align*}
\limsup_{n\to\infty}
&\Big|\E\Big[F\Big(V(Q^\dagger_n), \sqrt{\frac{3}{2}}\,n^{-1/2} \,d^\dagger_n,2n^{-2}\mu_n\Big)\,\Big|\, H_{\lfloor \delta\sqrt{n}\rfloor}\in [n,(1+\alpha) n]\Big]\\
&\qquad\qquad- \N_0\Big(F(\mm^{(\delta)},D^{(\delta)},\bv^{(\delta)})\,\Big|\,
 \z_{W_*+\delta\sqrt{3/2}}\in[1,1+\alpha]\Big)\Big|\leq \ve.
 \end{align*}
\end{proposition}

\proof We may assume that $F$ is bounded by $1$ and $1$-Lipschitz. 
Recall the notation ${\mathcal{N}}_{W_*+\delta}$ for the point measure of 
excursions outside $(W_*+\delta,\infty)$ (see the beginning of Section \ref{sec:resu-snake}).
If we write
$${\mathcal{N}}_{W_*+\delta}=\sum_{i\in I_\delta} \delta_{\omega_{i,\delta}}$$
then it follows from our definitions that
$$K_\delta=\delta + \sup_{i\in I_\delta}\,W^*(\omega_{i,\delta})\,,\quad\kappa_\delta= \sum_{i\in I_\delta} \sigma(\omega_{i,\delta}).$$
From the description of the conditional law of ${\mathcal{N}}_{W_*+\delta}$ given $\z_{W_*+\delta}$ in Proposition \ref{condidistr}
and the subsequent remark,
it is not hard 
to verify that, for every $\beta>0$,
\begin{equation}
\label{conv-trunc1}
\N_0(K_\delta>\beta\mid \z_{W_*+\delta}\in[1,1+\alpha])\build{\la}_{\delta\to 0}^{}0\,,\quad\N_0(\kappa_\delta>\beta\mid \z_{W_*+\delta}\in[1,1+\alpha])\build{\la}_{\delta\to 0}^{}0.
\end{equation}
Let us explain how the first convergence in \eqref{conv-trunc1} is derived (the proof of the second one is analogous).
If $z\in[1,1+\alpha]$, the conditional law of ${\mathcal{N}}_{W_*+\delta}$ given $\z_{W_*+\delta}=z$  is the sum of a Poisson point measure 
$\mathcal{N}_{(\delta)}$ with intensity $z\,\N_0(\cdot\cap\{W_*>-\delta\})$ and the Dirac mass at an independent 
random snake trajectory $\omega_{(\delta)}$ distributed according to $\N_0(\cdot\mid W_*=-\delta)$. For $\delta<\beta/2$, we have thus
$$\N_0(K_\delta>\beta\mid \z_{W_*+\delta}=z)\leq \P(\sup\{W^*(\omega):\omega\hbox{ atom of }\mathcal{N}_{(\delta)}\} >\frac{\beta}{4}) + \P(W^*(\omega_{(\delta)})>\frac{\beta}{4}).$$
A scaling argument shows that $W^*(\omega_{(\delta)})$ has the same distribution as $\delta W^*(\omega_{(1)})$ and thus tends to $0$ in probability
when $\delta\to 0$. On the other hand, we may assume that $\mathcal{N}_{(\delta)}$ is the restriction of a Poisson point measure $\mathcal{N}$ with intensity $z\,\N_0$
to the set $\{\omega:W_*(\omega)>-\delta\}$. Now observe that $\mathcal{N}$ has only finitely many atoms $\omega$ such that $W^*(\omega)>\beta/4$ and 
that, for each of these atoms, we have $W_*(\omega)<0$. It follows that, if $\delta>0$ is small enough, no atom $\omega$ of $\mathcal{N}_{(\delta)}$
satisfies $W^*(\omega)>\beta/4$, and hence $ \P(\sup\{W^*(\omega):\omega\hbox{ atom of }\mathcal{N}_{(\delta)}\} >\frac{\beta}{4})$ tends to $0$ as $\delta\to 0$.

It follows from \eqref{conv-trunc1} that, if $\delta<\ve/100$ is small enough, 
we have both
\begin{align*}
\N_0(K_{\delta\sqrt{3/2}}>\ve/100\mid \z_{W_*+\delta\sqrt{3/2}}\in[1,1+\alpha])&<\ve/100,\\
\N_0(\kappa_{\delta\sqrt{3/2}}>\ve/100\mid \z_{W_*+\delta\sqrt{3/2}}\in[1,1+\alpha])&<\ve/100,
\end{align*}
and (using also Lemma \ref{tech-GHdistance}) for all $n$ sufficiently large,
\begin{align*}
\P(2n^{-2}\,\kappa_{n,\delta}>\ve/100\mid H_{\lfloor \delta\sqrt{n}\rfloor}\in [n,(1+\alpha) n]) 
&<\ve/100,\\
\P(n^{-1/2}K_{n,\delta}>\ve/100\mid H_{\lfloor \delta\sqrt{n}\rfloor}\in [n,(1+\alpha) n])
&<\ve/100.
\end{align*}

By Lemma \ref{GHdistance-truncated}, we have then
\begin{align*}
&\N_0\Big(|F(\mm^{(\delta)},D^{(\delta)},\bv^{(\delta)})-F(\mm,D,\bv)|\,\Big|\,
 \z_{W_*+\delta\sqrt{3/2}}\in[1,1+\alpha]\Big)\\
 &\qquad \leq 2\, \N_0(3(2(\delta + K_\delta)+\kappa_\delta)\,\mathbf{1}_{\{K_\delta\leq\ve/100,
 \kappa_\delta\leq\ve/100\}}\mid \z_{W_*+\delta\sqrt{3/2}}\in[1,1+\alpha]) 
+  4\times (\ve/100)\\
&\qquad \leq 40\times (\ve/100).
 \end{align*}
Similarly, from Lemma \ref{GHdistance-discrete}, we get for $n$ sufficiently large,
\begin{align*}
&\E_0\Big[|F(V(Q^\dagger_n), \sqrt{\frac{3}{2}}\,n^{-1/2} \,d^\dagger_n,2n^{-2}\mu_n)
-F(V(Q),\sqrt{\frac{3}{2}}\,n^{-1/2}\mathrm{d}_{\mathrm{gr}},2n^{-2}\mu)|\,\Big|\,H_{\lfloor \delta\sqrt{n}\rfloor}\in [n,(1+\alpha) n])\Big]\\
&\qquad
\leq 50\times (\ve/100).
\end{align*}
Finally the convergence to the Brownian map (Proposition \ref{convBmap}) and
Proposition \ref{conv-exit} also show that for $n$ large enough we have
\begin{align*}
&|\N_0\Big(F(\mm,D,\bv)  \,\Big|\,
 \z_{W_*+\delta\sqrt{3/2}}\in[1,1+\alpha]\Big)\\
 &\qquad-\E\Big[F(V(Q),\sqrt{\frac{3}{2}}\,n^{-1/2}\mathrm{d}_{\mathrm{gr}},2n^{-2}\mu)
 \,\Big|\,H_{\lfloor \delta\sqrt{n}\rfloor}\in [n,(1+\alpha) n])\Big]\Big| < \ve/20.
 \end{align*}
 We get the desired result by combining the last three displays. \endproof
 
 \medskip
 We will now combine Propositions \ref{conv-glued}   and \ref{conv-truncated}
 to get a representation of $\F^{\bullet,\dagger}_z$ (defined in Section \ref{sec:Bdglued}) in terms of
 the measure $\N^{*,z}_0$ introduced in Section \ref{sec:posexc}. Theorem \ref{identif-glued-unpointed} will be an easy consequence 
 of this representation. Recall the notation $\mathcal{L}^{\bullet\bullet}(\omega)$
 from Section \ref{sec:consmms}.

\begin{proposition}
\label{identif-glued}
For every $z>0$, for every nonnegative measurable function
$G$ on $\M^{\bullet\bullet}$, we have
$$\F_z^{\bullet,\dagger}(G)= z^{-2}\,\N^{*,z}_0\Big(\int_0^\sigma \mathrm{d}r\,G(\ll^{\bullet\bullet}(W^{[r]}))\Big).$$
\end{proposition}

\proof Let $\Gamma_z$ denote the conditional
distribution of $\mathrm{tr}_{W_*+\delta}(W)$ given $\z_{W_*+\delta}=z$ under $\N_0^{[-\delta]}$. By 
Proposition \ref{law-trunc}, $\Gamma_z$ does not depend on $\delta$, and moreover, for every 
nonnegative measurable function $F$ on the space of snake trajectories, we have
$$\int \Gamma_z(\mathrm{d}\omega)\,F(\omega)= z^{-2}\,\N^{*,z}_0\Big(\int_0^\sigma \mathrm{d}r\,F(W^{[r]})\Big).$$
Hence the identity of Proposition \ref{identif-glued} is equivalent to $\F_z^{\bullet,\dagger}(G)= \int \Gamma_z(\mathrm{d}\omega)\,G(\ll^{\bullet\bullet}(\omega))$.

We prove this identity for $z=1$.
We may assume that $G$ is $1$-Lipschitz and bounded by $1$. 
Let us fix $\ve>0$. We first choose $\alpha_0>0$ so that  $|\F^{\bullet,\dagger}_z(G)-\F^{\bullet,\dagger}_1(G)|\leq \ve$ if $1\leq z\leq 1+\alpha_0$.
Using Lemma \ref{identi-mms-trunc} and Corollary \ref{limit-law}, we have for every $\alpha>0$,
$$
\N_0\Big(G(\mm^{(\delta)},D^{(\delta)},\bv^{(\delta)}) \,\Big|\,
 \z_{W_*+\delta\sqrt{3/2}}\in[1,1+\alpha]\Big)
= c_{\delta,\alpha}\int_1^{1+\alpha} 
 \frac{\mathrm{d}z}{\sqrt{z}}\,\exp(-\frac{z}{\delta^2})\int \Gamma_z(\mathrm{d}\omega)\,G(\ll^{\bullet\bullet}(\omega)),
$$
 where the normalizing constant $c_{\delta,\alpha}$ is such that
 $$c_{\delta,\alpha}\int_1^{1+\alpha} 
 \frac{\mathrm{d}z}{\sqrt{z}}\,\exp(-\frac{z}{\delta^2})=1.$$
It follows that
\begin{align*}
&\Big| \N_0\Big(G(\mm^{(\delta)},D^{(\delta)},\bv^{(\delta)}) \,\Big|\,
 \z_{W_*+\delta\sqrt{3/2}}\in[1,1+\alpha]\Big) - \int \Gamma_1(\mathrm{d}\omega)\,G(\ll^{\bullet\bullet}(\omega))\Big|\\
 &\qquad \leq \sup_{1\leq z\leq 1+\alpha} \Big|\int \Gamma_z(\mathrm{d}\omega)\,G(\ll^{\bullet\bullet}(\omega))
 -\int \Gamma_1(\mathrm{d}\omega)\,G(\ll^{\bullet\bullet}(\omega))\Big|.
 \end{align*}
Thanks to the continuity of $G$ and the scaling properties linking the measures $\N^{*,z}_0$
(and hence the measures $\Gamma_z$ according to Proposition \ref{law-trunc}), it is not hard to verify that the
right-hand side tends to $0$ as $\alpha\to 0$.

We fix $\alpha\in(0,\alpha_0)$ such that the right-hand side of the last display is 
smaller than $\ve$, and note that this bound holds for
any $\delta>0$. 
We then fix $\delta>0$ small enough so that we can combine the latter bound with 
Proposition \ref{conv-truncated} to get that, for all 
large enough $n$,
$$\Big|\E\Big[G(V(Q^\dagger_n), \sqrt{\frac{3}{2}}\,n^{-1/2} \,d^\dagger_n,2n^{-2}\mu_n)\,\Big|\, H_{\lfloor \delta\sqrt{n}\rfloor}\in [n,(1+\alpha) n]\Big]
- \int \Gamma_1(\mathrm{d}\omega)\,G(\ll^{\bullet\bullet}(\omega))\Big|\leq 2\ve.$$

On the other hand, thanks to Proposition \ref{conv-glued}, we have also for 
$n$ large,
$$\Big|\E\Big[G(V(Q^\dagger_n), \sqrt{\frac{3}{2}}\,n^{-1/2} \,d^\dagger_n,2n^{-2}\mu_n)\,\Big|\, H_{\lfloor \delta\sqrt{n}\rfloor}\in [n,(1+\alpha) n])\Big]
-c_{\delta,\alpha}\int_1^{1+\alpha} 
 \frac{\mathrm{d}z}{\sqrt{z}}\,\exp(-\frac{z}{\delta^2})\,\F^{\bullet,\dagger}_z(G)\Big|\leq \ve.$$
and from the fact that $\alpha<\alpha_0$,
$$\Big| c_{\delta,\alpha}\int_1^{1+\alpha} 
 \frac{\mathrm{d}z}{\sqrt{z}}\,\exp(-\frac{z}{\delta^2})\,\F^{\bullet,\dagger}_z(G) - \F^{\bullet,\dagger}_1(G)\Big|\leq \ve.$$
 By combining all these bounds, we have
 $$\Big| \int \Gamma_1(\mathrm{d}\omega)\,G(\ll^{\bullet\bullet}(\omega)) - \F^{\bullet,\dagger}_1(G)\Big| \leq 4\ve.$$
 Since $\ve$ was arbitrary, we get
 $$\F^{\bullet,\dagger}_1(G)=\int \Gamma_1(\mathrm{d}\omega)\,G(\ll^{\bullet\bullet}(\omega)),$$
 which was the desired result. \endproof

\noindent{\it Proof of Theorem \ref{identif-glued-unpointed}.}
Recall from the remark at the end of Section \ref{sec:consmms} that
 the pointed measure metric space $\ll^\bullet(W^{[s]})$ does not depend on the
 choice of  $s\in[0,\sigma]$.

 From the relation between $\F_z$ and $\F^\bullet_z$ (Section \ref{sec:freeBdisk}), we
 have
 $$\F_z^\dagger(G)=z^2\,\F_z^{\bullet,\dagger}\Big(\frac{1}{\Sigma}\,G\circ\kappa^\bullet\Big),$$
 where $\kappa^\bullet$ is the obvious projection from $\M^{\bullet\bullet}$ onto $\M^\bullet$
 that consists in ``forgetting'' the second distinguished point, and $\Sigma$
 denotes the total mass of the volume measure. 
 By Proposition \ref{identif-glued}, the right-hand side is equal to
 $$\N^{*,z}_0\Big(\frac{1}{\sigma}\int_0^\sigma \mathrm{d}r\,G(\ll^\bullet(W^{[r]}))\Big)$$
 and this is equal to $\N^{*,z}_0(G(\ll^\bullet(W)))$ by the first observation of the proof. \hfill$\square$

 \section{Construction of the metric under $\N^*_0$}
 \label{sec:cons-metric}
 
 Theorem \ref{identif-glued-unpointed} shows that the law of the Brownian disk with perimeter $z$
 and glued boundary can be identified as the law of $\ll^\bullet(W)$ under $\N^{*,z}_0$.
 Our objective is now to recover the distribution of the Brownian disk (without gluing), by
 reconstructing the original metric from the metric with glued boundary.
 We aim at constructing this original metric under $\N^{*,z}_0$, but it will be more
 convenient to argue first under a different measure, and then to use the 
 re-rooting representation of $\N^*_0$ in Theorem \ref{re-root-rep}.
 
 In the first part of this section, we fix $r>0$, and 
 we will argue under the measure $\N_r^{[0]}:=\N_r(\cdot\mid W_*\leq 0)$. 
 To simplify
notation, we write $W^{\tr}=\mathrm{tr}_0(W)$, so that $W^\tr_s= W_{\eta_s}$, with
\begin{equation}
\label{def-eta}
\eta_s=\inf\{t\geq 0: \int_0^t \mathrm{d}u\,\mathbf{1}_{\{\tau_0(W_u)\geq \zeta_u\}}>s\}.
\end{equation}
It will be convenient to write $\t^\tr$ for the genealogical tree of 
$W^\tr$, and $p_\tr$ for the canonical projection from
$[0,\sigma(W^\tr)]$ onto $\t^\tr$. Recall that 
$\t^\tr$ is canonically (and isometrically) identified 
with the closed subset of $\t_\zeta$ consisting of all $a$ such that
$Z_b>0$ for every strict ancestor $b$ of $a$. Without risk of confusion, we keep the
notation $(Z_a)_{a\in\t^\tr}$ for the labels on $\t^\tr$ ($Z_a=\wh W^\tr_s$
if $a=p_\tr(s)$). Note that, by the definition of the truncation $\mathrm{tr}_0$, the labels $Z_a$, $a\in\t^\tr$,
are nonnegative. By definition, the boundary $\partial\t^\tr$ is
$$\partial\t^\tr=\{a\in \t^\tr: Z_a=0\},$$
and we set $\t^{\tr,\circ}=\t^\tr\backslash \partial\t^\tr$.

For $a,b\in \t^\tr$,  
 $[a,b]$ stands for the ``lexicographical'' interval in  $\t^\tr$, as in Section \ref{sec:snake}, and as previously
 we define for every $a,b\in\t^\tr$,
$$D^{\tr,\circ}(a,b)= Z_a+ Z_b -2\max\Big( \min_{c\in [a,b]} Z_c, \min_{c\in [ b,a ]} Z_c\Big).
$$
 We then set, for every $a,b\in \t^{\tr,\circ}$,
\begin{equation}
\label{Dtrunc}
\Delta^{\circ}(a,b)= \left\{ \begin{array}{ll}
D^{\tr,\circ}(a,b)\qquad&\hbox{if }\max\Big( \min_{c\in [a,b]} Z_c, \min_{c\in [b,a]} Z_c\Big)>0\\
 +\infty&\hbox{otherwise.}
 \end{array}
 \right.
 \end{equation}
 and,
 \begin{equation}
 \label{defDelta}
 \Delta(a,b)=\inf_{a=a_0,a_1,\ldots,a_{k-1},a_k=b}\  \sum_{i=1}^k \Delta^{\circ}(a_{i-1},a_i),
 \end{equation}
 where the infimum is over all choices of the integer $k\geq 1$ and of the
elements $a_0,a_1,\ldots,a_k$ of $\t^\tr$ such that $a_0=a$
and $a_k=b$. Notice that $\Delta(a,b)\geq |Z_a-Z_b|$ by construction.

\smallskip
\rem If $a,b\in \t^{\tr,\circ}$, we can also define $D^\circ(a,b)$ by viewing $a$ and $b$ as elements of
$\t_\zeta$, using \eqref{Dzero}. If $a$ and $b$  are such that $\max( \min_{c\in [a,b]} Z_c, \min_{c\in [b,a]} Z_c)>0$ in \eqref{Dtrunc},
then $\Delta^{\circ}(a,b)=D^{\tr,\circ}(a,b)=D^\circ(a,b)$. The point is that, if $\min_{c\in [a,b]} Z_c>0$, the points of the interval $[a,b]$ of $\t^\tr$
correspond  to the ``same'' interval of $\t_\zeta$ modulo the preceding identification 
of $\t^\tr$ as a subset of $\t_\zeta$ (and similarly if $[a,b]$
is replaced by $[b,a]$). 
 
 \begin{proposition} 
 \label{propDelta}
 The following properties hold $\N^{[0]}_r$ a.s.
 \begin{enumerate}
 \item[\rm(i)] $\Delta(a,b)<\infty$ for every $a,b\in \t^{\tr,\circ}$.
 \item[\rm(ii)] $\Delta$ defines a pseudo-metric on $\t^{\tr,\circ}$.
 \item[\rm(iii)] The mapping $(a,b)\mapsto \Delta(a,b)$ is continuous on $\t^{\tr,\circ}\times\t^{\tr,\circ}$.
 \item[\rm(iv)] For every $a,b\in \t^{\tr,\circ}$, $\Delta(a,b)=0$ if and only if $\Delta^\circ (a,b)=0$. 
 \end{enumerate}
 \end{proposition}
 
 We omit the easy proof of these properties (for (i), observe that, if $\llbracket a,b\rrbracket$ denotes the 
 geodesic segment between $a$ and $b$ in  $\t^{\tr}$, $\t^{\tr}\backslash \llbracket a,b\rrbracket$
 has only finitely many connected components $C_1,\ldots,C_k$ that intersect $\partial\t^{\tr}$, and,
 if $a_i$ denotes the unique point of $\ov{C}_i\cap \llbracket a,b\rrbracket$, then, assuming that the components
 have been ranked so that $a_1,\ldots,a_k$ come in this order from $a$ to $b$ along $\llbracket a,b\rrbracket$, we have 
 $\Delta^\circ(a_{i-1},a_i)<\infty$ for every $i=1,\ldots,k+1$, with the convention $a_0=a$ and $a_{k+1}=b$). 
 In order to verify (iv), we use the analogous property for the Brownian map recalled in Section \ref{sec:defBM}.
 
 \begin{proposition}
 \label{extenD}
 $\N^{[0]}_r$ a.s., the mapping $(a,b)\mapsto \Delta(a,b)$ has a continuous extension to $\t^\tr\times\t^\tr$.
 \end{proposition}
 
 \proof
 An important ingredient of the proof is the the reduced tree of $\t^\tr$, which consists
 of all vertices of $\t^\tr$ that have at least one descendant with label $0$. Let us briefly explain why 
 this reduced tree, which is denoted by $\t^\bigtriangledown$, is relevant for our purposes. In the infimum in the 
 definition \eqref{defDelta} of $\Delta(a,b)$, we interpolate between $a$ and $b$
 using points $b_1,\ldots,b_{k-1}$ such that $ \Delta^{\circ}(b_{i-1},b_i)$ is finite for every $i=2,\ldots,k-1$, and thus 
 the minimum of $Z$ on the interval $[b_{i-1},b_i]$ (or on $[b_i,b_{i-1}]$) is positive. A way of finding such 
 points is to let $b_1,\ldots,b_{k-1}$ be successive branching points of the reduced tree $\t^\bigtriangledown$,
 meaning that there is no branching point of $\t^\bigtriangledown$ in the interior of the line segment 
 between $b_{i-1}$ and $b_{i}$. To get bounds on the quantities $\Delta^\circ(b_{i-1},b_i)$, we then need to control the
 labels $Z_{a}$ when $a$ varies over branching points of $\t^\bigtriangledown$. 
 In order to obtain the result of the proposition, it will be enough to prove that 
 the maximum of these labels over branching points of $\t^\bigtriangledown$ that are ``close'' to the boundary $\partial\t^\tr$ tends to $0$ sufficiently fast 
 (see the bound \eqref{estimdep15} below).

 The probabilistic structure of $\t^\bigtriangledown$
 is described in \cite[Theorem 4.7.2]{DLG} (in a much more general setting). Let us mention only
 the properties that are relevant for our purposes. The tree $\t^\bigtriangledown$ is a binary $\R$-tree, which can be constructed by induction 
 as follows. One starts 
from a line segment connecting the root of $\t^\tr$ to a first branching point 
 $a_\varnothing$. To this branching point are attached two other line segments connecting
 $a_\varnothing$ respectively to branching points $a_1$ and $a_2$, listed in the
 lexicographical order of $\t^\tr$. To $a_1$ (respectively to $a_2$) are then attached two line segments
 connecting $a_1$ (resp. $a_2$) to branching points $a_{(1,1)}$ and $a_{(1,2)}$ (resp.
 $a_{(2,1)}$ and $a_{(2,2)}$) and so on. The tree $\t^\bigtriangledown$ has the following recursive structure:
 conditionally on the line segment joining the root to $a_\varnothing$, and on the labels 
 along this line segment, the two (labeled) subtrees
 rooted at $a_\varnothing$ are distributed independently according to 
 the law of $\t^\bigtriangledown$ under $\N^{[0]}_{Z_{a_\varnothing}})$. 
 
 Thanks to this recursive property, the distribution of $\t^\bigtriangledown$ (and of labels along $\t^\bigtriangledown$) is characterized by the law of the (length of the)
 segment connecting the root to $a_\varnothing$, and the distribution of the labels along this line segment.
 These labels can be represented by a stopped path $W^\bigtriangledown=(W^\bigtriangledown(t):0\leq t\leq \zeta^\bigtriangledown)$, with 
 $W^\bigtriangledown(0)=r$ and $W^\bigtriangledown(\zeta^\bigtriangledown)=Z_{a_\varnothing}$ --- in particular $\zeta^\bigtriangledown$ is the
 length of the line segment from the root to $a_\varnothing$. The distribution of $W^\bigtriangledown$
 is then given by the formula
 $$\N^{[0]}_r(\Phi(W^\bigtriangledown))
 = 2\,u(r)^{-1}\,\int_0^\infty \mathrm{d}s\,\mathbf{E}_r\Big[ \mathbf{1}_{\{s<\tau_0\}}\,u(B_s)^2\exp\Big(-4\int_0^s \mathrm{d}t\,u(B_t)\Big)
 \Phi(B_t,0\leq t\leq s)\Big],$$
 where $\Phi$ is a nonnegative measurable function on the space $\mathcal{W}$, $(B_t)_{t\geq 0}$ is a linear Brownian motion that starts from $r$ under the probability measure $\mathbf{P}_r$,
 $\tau_0=\inf\{t\geq 0: B_t=0\}$, and for every $y>0$, $u(y)=\N_y(W_*\leq 0)= 3/(2y^2)$. 
 
 We are in fact interested in the distribution of $Z_{a_\varnothing}$ (the terminal label of $W^\bigtriangledown$). From now on, we take $r=1$
 to simplify the presentation.
 From the preceding formula, we get, for every nonnegative measurable function $\phi$ on $[0,\infty)$,
 $$\N^{[0]}_1(\phi(Z_{a_\varnothing}))
 = 3\,\int_0^\infty \mathrm{d}s\,\mathbf{E}_1\Big[ \mathbf{1}_{\{s<\tau_0\}}\,(B_s)^{-4}\exp\Big(-6\int_0^s \frac{\mathrm{d}t}{(B_t)^{2}}\Big)
 \phi(B_s)\Big].$$
 At this point we use classical absolute continuity relations between Bessel processes
 (see e.g. \cite[Proposition 2.6]{LGW}), which show that, for 
 every $s>0$,
 $$\mathbf{E}_1\Big[ \mathbf{1}_{\{s<\tau_0\}}\,(B_s)^{-4}\exp\Big(-6\int_0^s \frac{\mathrm{d}t}{(B_t)^{2}}\Big)
 \phi(B_s)\Big]= \mathbf{E}_1^{(9)}\Big[ (R_s)^{-8}\,\phi(R_s)\Big],$$
 where $R=(R_t)_{t\geq 0}$ stands for a nine-dimensional Bessel process that starts at $1$ under the probability measure
 $\mathbf{P}^{(9)}_1$. We are thus led to the calculation of
 $$3 \,\int_0^\infty \mathrm{d}s\,\mathbf{E}_1^{(9)}\Big[ (R_s)^{-8}\,\phi(R_s)\Big],$$
 and we specialize to the case $\phi(x)=x^p$, for $p\geq 0$. For $y,z\in \R^9$, let $G(y,z)=c_9|z-y|^{-7}$ be the Green function 
 of nine-dimensional Brownian motion, with $c_9=\Gamma(7/2)/(2\pi^{9/2})$. We also write
 $\mathbf{1}$ for a fixed point of the unit sphere of $\R^9$. Since a nine-dimensional Bessel process has
 the distribution of the norm of a nine-dimensional Brownian motion, we have
 $$3 \,\int_0^\infty \mathrm{d}s\,\mathbf{E}_1^{(9)}\Big[ (R_s)^{-8+p}\Big]
 = 3\int_{\R^9} \mathrm{d}z\, G(\mathbf{1},z)\,|z|^{-8+p}
 =\frac{6}{7} \int_0^\infty \mathrm{d}\rho \,\rho^{p}\,\int \pi_\rho(\mathrm{d}z)\,|z-\mathbf{1}|^{-7},$$
 where $\pi_\rho$ stands for the uniform probability measure on the sphere $\{|z|=\rho\}$.
 In the last equality, we also note that, if $c'_9= 2\pi^{9/2}/\Gamma(9/2)$ is the volume of
 the unit sphere of $\R^9$, we have $c_9c'_9=2/7$. Finally, we observe that
 $$\int \pi_\rho(\mathrm{d}z)\,|z-\mathbf{1}|^{-7}= 1\wedge \rho^{-7}$$
 and thus
 $$\frac{6}{7} \int_0^\infty \mathrm{d}\rho \,\rho^{p}\,\int \pi_\rho(\mathrm{d}z)\,|z-\mathbf{1}|^{-7}
 =\frac{6}{7} \int_0^\infty \mathrm{d}\rho \,\rho^{p}\,(1\wedge \rho^{-7})= \frac{6}{7}\Big(\frac{1}{p+1} + \frac{1}{6-p}\Big),$$
 provided that $0\leq p<6$. Summarizing, we have obtained that, for $0\leq p<6$, we have
 $$\N^{[0]}_1((Z_{a_\varnothing})^p)= \frac{6}{7}\Big(\frac{1}{p+1} + \frac{1}{6-p}\Big).$$
 As a function of $p$, the right-hand side attains its minimal value for $p=5/2$,
 and
 $$\N^{[0]}_1((Z_{a_\varnothing})^{5/2})=\frac{24}{49}<\frac{1}{2}.$$
 
 At this point we use the recursive structure of the labeled tree $\t^\bigtriangledown$, together with scaling 
 properties. If
 $$I=\bigcup_{n=0}^\infty \{1,2\}^n \quad(\hbox{where }\{1,2\}^0=\{\varnothing\})$$
 then for every $(i_1,\ldots,i_n)\in I$, $Z_{a_{(i_1,\ldots,i_n)}}$ has the distribution of
 the product $\xi_0\xi_1\cdots\xi_n$, where $\xi_0,\ldots,\xi_n$ are independent and have
 the distribution of $Z_{a_\varnothing}$ under $\N^{[0]}_1$. 
 Consequently,
$$\N^{[0]}_1\Big((Z_{a_{(i_1,\ldots,i_n)}})^{5/2}\Big)=(\frac{24}{49})^{n+1}.$$
Fix a real $\alpha\in (0,1)$ such that $2\,\alpha^{-5/2}\,\frac{24}{49} < 1$.
Then 
$$\N^{[0]}_1(Z_{a_{(i_1,\ldots,i_n)}}>\alpha^n) \leq \alpha^{-5n/2}\,\N^{[0]}_1\Big((Z_{a_{(i_1,\ldots,i_n)}})^{5/2}\Big)=\alpha^{-5n/2}\,(\frac{24}{49})^{n+1}.$$
If we sum this bound over the $2^n$ choices of $(i_1,\ldots,i_n)\in\{1,2\}^n$, we obtain that
$$\sum_{n=0}^\infty\N^{[0]}_1\Bigg(\bigcup_{(i_1,\ldots,i_n)\in\{1,2\}^n} \{Z_{a_{(i_1,\ldots,i_n)}}>\alpha^n\}
\Bigg) < \infty,$$
and, by the Borel-Cantelli lemma, we get that $\N^{[0]}_1$ a.s. for $n$ large enough, for every $(i_1,\ldots,i_n)\in\{1,2\}^n$,
\begin{equation}
\label{estimdep15}
Z_{a_{(i_1,\ldots,i_n)}}\leq \alpha^n.
\end{equation}
Let $J$ be the set of all infinite sequences $u=(i_1,i_2,\ldots)$ where $i_j=1$ or $2$. For $u=(i_1,i_2,\ldots)\in J$ and $k\geq 1$ write $[u]_k=(i_1,\ldots,i_k)$
for the truncated sequence at order $k$. It follows from the preceding considerations that
\begin{equation}
\label{estimdep}
\lim_{\ell\to\infty} \Bigg(\sup_{u\in J} \sum_{k=\ell}^\infty Z_{a_{[u]_k}} \Bigg) =0,\qquad \hbox{a.s.}
\end{equation}

Let us use \eqref{estimdep} to prove the statement of the proposition. Let $\delta\in(0,1)$.
We note that, under $\N^{[0]}_1$, each excursion of $W^{\tr}$ outside $(\delta,\infty)$ 
that hits $0$ (in the sense of the beginning of Section \ref{sec:Bstrunc}) corresponds
to a subtree of $\t^\tr$ rooted at a vertex $a$ of $\t^\tr$ such that $Z_a=\delta$.
We let $\t_{1,\delta},\ldots,\t_{N_\delta,\delta}$ be the subtrees of $\t^\tr$
obtained in this way. 
We observe that it is enough to prove that
\begin{equation}
\label{unifdep}
\lim_{\delta\to 0} \Bigg( \sup_{1\leq i\leq N_\delta} \;\sup_{x,y\in \t_{i,\delta}, Z_x\wedge Z_y>0} \Delta(x,y) \Bigg) =0.
\end{equation}
Indeed, assume that \eqref{unifdep} holds. Then, if we are given a sequence $(x_n,y_n)_{n\geq 1}$ in $\t^\tr\times\t^\tr$
such that $Z_{x_n}\wedge Z_{y_n}>0$ and if $(x_n,y_n)$ converges to $(x,y)$, with $Z_x=Z_y=0$, then, for every $\delta >0$, 
 there exist $i,j\in \{1,\ldots,N_\delta\}$ such that, for all 
$n$ large enough, say for $n\geq n_\delta$, we have $x_n,x\in \t_{i,\delta}$ and $y_n,y\in\t_{j,\delta}$. Therefore we also
get that, if $m,n\geq n_\delta$, we have $\Delta(x_n,x_m)\leq r_\delta$ and $\Delta(y_n,y_m)\leq r_\delta$, where 
$r_\delta$ denotes the supremum over $1\leq i\leq N_\delta$ appearing in the last display. Consequently, 
$|\Delta(x_n,y_n)-\Delta(x_m,y_m)|\leq 2\,r_\delta$ for all $n,m\geq n_\delta$, and we obtain that $\Delta(x_n,y_n)$ is
a Cauchy sequence. We can therefore set $\Delta(x,y)=\lim \Delta(x_n,y_n)$ as $n\to\infty$, and this gives the
desired continuous extension. 

It remains to prove \eqref{unifdep}. For every $k\geq 1$, let $\t^\bigtriangledown_{(k)}$ be the subtree of $\t^\bigtriangledown$ that consists of
all ancestors of the vertices $a_{(i_1,\ldots,i_k)}$ for all $(i_1,\ldots,i_k)\in \{1,2\}^k$. By compactness,
$$\ve_k:=\inf_{a\in \t^\bigtriangledown_{(k)}} Z_a >0.$$
Suppose that $\delta<\ve_k$, and let $x,y\in \t_{i,\delta}$ for some $1\leq i\leq \t_{i,\delta}$, with $Z_x\wedge Z_y>0$. 
Let $\rho_{i,\delta}$ stand for the root of $\t_{i,\delta}$. Then,
$$\Delta(x,y)\leq \Delta(x,\rho_{i,\delta})+ \Delta(y,\rho_{i,\delta})$$
and in order to bound $\Delta(x,\rho_{i,\delta})$ (or $\Delta(y,\rho_{i,\delta})$), we may proceed as
follows. Write $x_\bullet$ for the most recent ancestor of $x$ that belongs
to $\t^\bigtriangledown$. Note that $x_\bullet$ belongs to $\t_{i,\delta}$ (because $\rho_{i,\delta}\in\t^\bigtriangledown$ and so $x_\bullet$ is
a descendant of $\rho_{i,\delta}$). We can find integers
$m$ and $\ell$, with $k\leq m\leq \ell$, such that $x_\bullet\in \llbracket a_{(i_1,\ldots,i_\ell)}, a_{(i_1,\ldots,i_{\ell+1})}\rrbracket $
and $\rho_{i_,\delta}\in \llbracket a_{(i_1,\ldots,i_{m})}, a_{(i_1,\ldots,i_{m+1})}\rrbracket$
for some $(i_1,\ldots,i_{\ell+1})\in \{1,2\}^{m+1}$ --- here we recall the notation 
$\llbracket a,b\rrbracket$ for the geodesic segment
between two points $a,b\in\t^\bigtriangledown$. Then,
$$\Delta^\circ(x,x_\bullet)=D^{\tr,\circ}(a,b)\leq Z_x + Z_{x_\bullet}$$
where the first equality holds because no vertex $y$ of the lexicographical interval $[x,x_\bullet]$ of $\t^\tr$ can be such that $Z_y=0$
(otherwise, this would contradict the fact that $x_\bullet$ is the most recent ancestor of $x$ that belongs
to $\t^\bigtriangledown$). For similar reasons, we have
\begin{align*}
\Delta^\circ(x_\bullet,a_{(i_1,\ldots,i_{\ell})})&\leq Z_{x_\bullet}+ Z_{a_{(i_1,\ldots,i_{\ell})}}\\
\Delta^\circ(a_{(i_1,\ldots,i_{m+1})},\rho_{i_,\delta})&\leq Z_{a_{(i_1,\ldots,i_{m+1})}}+Z_{\rho_{i_,\delta}}
\end{align*}
and, for every $j$ such that $m\leq j\leq \ell$,
$$\Delta^\circ(a_{(i_1,\ldots,i_{j})}, a_{(i_1,\ldots,i_{j+1})})\leq Z_{a_{(i_1,\ldots,i_{j})}}+Z_{a_{(i_1,\ldots,i_{j+1})}}.$$
By combining these bounds, we get
$$\Delta(x,\rho_{i,\delta})\leq 2\sum_{j=k}^{\ell+1} Z_{a_{(i_1,\ldots,i_{j})}} + 3 \,\sup_{1\leq n\leq N_\delta} \sup_{a\in \t_{n,\delta}} Z_a.$$
The supremum in the last display tends to $0$ when $\delta\to 0$, by a simple
uniform continuity argument, and we just have to use 
\eqref{estimdep} to complete the proof of \eqref{unifdep}. \endproof

\noindent {\bf Remark.} Instead of the Borel-Cantelli argument used to derive \eqref{estimdep}, we could have 
referred to standard large deviation estimates for tree-indexed random walk, see in particular \cite[Section 18]{Peres}.

\smallskip
We keep the notation $\Delta$ for the continuous extension of $\Delta$ to $\t^\tr\times \t^\tr$. 

\begin{proposition}
\label{properDelta}
The following properties hold $\N^{[0]}_r$ a.s.
\begin{enumerate}
\item[\rm(i)] For every $a\in\t^\tr$,
$$\inf\{ \Delta(a,c):c\in \partial\t^{\tr}\}= Z_a.$$
\item[\rm(ii)] For every $b,b'\in \partial\t^\tr$,
$\Delta(b,b')=0$ if and only if $Z_c>0$ for every $c\in]b,b'[$ or the same holds 
with $]b,b'[$ replaced by $]b',b[$.
\end{enumerate}
\end{proposition}

\proof (i) Let $a\in\t^\tr$. We can assume that $Z_a>0$. There exists $s\in[0,\sigma(W^\tr)]$
such that $p_\tr(s)=a$. Then set $s':=\inf\{t\in(s,\sigma(W^\tr)]: \wh W^\tr_t =0\}$
and $b=p_\tr(s')$ (this definition 
only makes sense if $\{t\in(s,\sigma(W^\tr)]: \wh W^\tr_t =0\}\not =\varnothing$, but if not we take 
$s'=\sup\{t\in[0,s):\wh W^\tr_t =0\}$ and make the necessary adaptations in the following lines).
Then $Z_b=0$ and $b\in\partial\t^\tr$. For every $n\geq 1$ such that $n^{-1}<Z_a$, set
$s_n=\inf\{t\in(s,\sigma(W^\tr)]: \wh W^\tr_t =n^{-1}\}$, and $b_n=p_\tr(s_n)$. Then 
$s_n\uparrow s'$ and thus $b_n\to b$ as $n\to\infty$. On the other hand,
$\Delta(a,b_n)=Z_a-Z_{b_n}$ , so that we get $\Delta(a,b)=Z_a$. This shows that
$\inf\{ \Delta(a,c):c\in \partial\t^{\tr}\}\leq Z_a$. The reverse bound is obvious
from the inequality $\Delta(a,c)\geq |Z_a-Z_c|$.

\noindent(ii) The fact that $Z_c>0$ for every $c\in]b,b'[$ implies $\Delta(b,b')=0$
is easy and left to the reader. Let us prove only the reverse implication. If $a\in\t^\tr$, we use the notation $\wt a$ for the  ``same'' point in $\t_\zeta$
(recall that
$\t^\tr$ is canonically identified to a closed subset of $\t_\zeta$).

Let $b,b'\in \partial\t^\tr$, with $b\not =b'$, and assume that $\Delta(b,b')=0$.
Write $b=p_\tr(s)$ and $b'=p_\tr(t)$ with $s,t\in [0,\sigma(W^\tr)]$. For definiteness,
we may assume that $s<t$. For every $n\geq 1$ such that $n^{-1}<r$, set
$s_n=\inf\{u>s: \wh W^\tr_u\geq n^{-1}\}$, $b_n=p_\tr(s_n)$, and similarly
$t_n=\sup\{u<t:\wh W^\tr_u\geq n^{-1}\}$, $b'_n=p_\tr(t_n)$. Then $b_n\to b$
and $b'_n\to b'$ as $n\to\infty$, so that $\Delta(b_n,b'_n)$ tends to $0$
as $n\to\infty$. Now recall formula \eqref{defDelta} defining $\Delta(b_n,b'_n)$.
In this formula we can restrict our attention to choices of $b_n=a_0,a_1,\ldots,a_k=b'_n$
such that, for every $i=1,\ldots,k$, $\Delta^\circ(a_{i-1},a_i)<\infty$, and therefore
$$\max\Big( \min_{c\in [a_{i-1},a_i]} Z_c, \min_{c\in [a_i,a_{i-1}]} Z_c\Big)>0.$$
By the remark preceding Proposition \ref{propDelta}, this implies that
$$\Delta^\circ(a_{i-1},a_i)=D^\circ(\wt a_{i-1},\wt a_i).$$
It follows from these observations 
that $\Delta(b_n,b'_n)\geq D(\wt b_n,\wt b'_n)$. By passing to the limit $n\to\infty$
we get that $D(\wt b,\wt b')=0$. By the property recalled at the end of
Section \ref{sec:defBM}, this implies that $D^\circ(\wt b,\wt b')=0$
and thus $Z_c\geq 0$ for every $c$ in the interval $[\wt b,\wt b']$ of $\t_\zeta$ (or the same with $[\wt b,\wt b']$
replaced by $[\wt b',\wt b]$).  The latter property can hold only
if $Z_c>0$ for every $c\in]\wt b,\wt b'[$, since otherwise this would mean that 
the mapping $s\mapsto \wh W_s$ has a local minimum 
at $0$, which does not occur $\N^{[0]}_1$ a.s. (this would contradict \eqref{support-exit}).
From the construction of truncations, this also implies that
$Z_c>0$ for every $c$ belonging to the interval $]b,b'[$ of $\t^\tr$.
\endproof

Then, $\N^{[0]}_r$ a.s.,
$\Delta$ is a pseudo-metric on $\t^\tr$, and we can consider the associated 
quotient metric space, which we denote by $(\Theta',\Delta)$. We equip $\Theta'$ with the volume measure
$\bV'$, which is the image of the volume measure on $\t^\tr$ under the canonical projection. The image 
of $\partial \t^\tr$ under the canonical projection is the boundary $\partial \Theta'$ (by definition).

We now observe that, thanks to the re-rooting representation of $\N^*_0$
(Theorem \ref{re-root-rep}), the same 
construction works as well under $\N^*_0$. Arguing now under $\N^*_0$, we set
$$ \t_\zeta^\circ:=\{c\in \t_\zeta:Z_c>0\}.$$
For $a,b\in \t_\zeta$, we define $D^\circ(a,b)$  by
\eqref{Dzero}, then, for 
$a,b\in \t_\zeta^\circ$, 
we define $\Delta^\circ(a,b)$ as in \eqref{Dtrunc},
replacing $D^{\tr,\circ}$ by $D^\circ$, and finally, for 
$a,b\in \t_\zeta^\circ$,  we define $\Delta(a,b)$ by the formula analogous to
\eqref{defDelta}, where $a_1,\ldots,a_{k-1}$ now vary in $\t_\zeta$
(under $\N^*_0$, we keep the same notation $\Delta^\circ$, $\Delta$ as in \eqref{Dtrunc} and \eqref{defDelta}
under $\N_r^{[0]}$, but this should
create no confusion). 
These definitions are consistent with those given in the introduction.

\begin{corollary}
\label{consDisk}
The results of the three previous propositions remain valid if 
the measure $\N^{[0]}_r$ is replaced by $\N^*_0$,
provided
$\t^\tr$ is replaced by $\t_\zeta$, $\t^{\tr,\circ}$
is replaced by $\t_\zeta^\circ$, and $\partial \t^\tr$ is replaced by
$\t_\zeta\backslash \t_\zeta^\circ$.
\end{corollary}
As already mentioned, Corollary \ref{consDisk} follows from
the re-rooting representation of $\N^*_0$
(Theorem \ref{re-root-rep}).

Then, under the measure $\N^*_0$, we can define a measure metric space
$(\Theta,\Delta,\bV)$ by exactly the same procedure as before
to define $(\Theta',\Delta,\bV')$ under $\N_r^{[0]}$: $\Theta$
is the quotient space of $\t_\zeta$ for the equivalence relation 
associated with the pseudo-metric $\Delta$, and 
$\bV$ is the image of the volume measure on $\t_\zeta$ under the canonical projection. By definition, the boundary 
$\partial\Theta$ is the image  of the 
set $\t_\zeta\backslash \t_\zeta^\circ$ under the canonical projection. By simple scaling arguments, the
construction also works 
under $\N_0^{*,z}$ for every $z>0$, and we recover the first part of Theorem \ref{main}.

\section{Identification of the Brownian disk}
\label{sec:ident}
In this section, we prove Theorem \ref{main}, which was stated in the introduction.
Let $(\D_z,D^\partial,\bv_z)$ be a free Brownian disk with perimeter $z$. As in Section \ref{sec:Bdglued} above, we can
construct  $\D_z^\dagger$ by gluing the boundary of $\D_z$  into
a single point, and equipping the resulting quotient space with the metric 
$D^{\partial,\dagger}$ defined as in \eqref{metric-glued}. Let $\bv_z^\dagger$ be the volume measure on $\D_z^\dagger$, and
denote the total mass of $\bv_z^\dagger$ by $|\bv_z^\dagger|$.
By Theorem \ref{identif-glued-unpointed}, we know that the distribution $\F_z^\dagger$ of $\D_z^\dagger$,
as a random pointed metric measure space,
is the same as the distribution of $\ll^\bullet(W)$ under $\N^{*,z}_0$. According to Section \ref{sec:infinite}, we may
consider a sequence $U_1,\ldots,U_n,\ldots$ of
independent uniformly distributed points in $(\D_z^\dagger,D^{\partial,\dagger},\bv_z^\dagger)$, and similarly a sequence
$V_1,\ldots,V_n,\ldots$ of
independent uniformly distributed points in $\ll^\bullet(W)$. Then, the distribution of
$$(|\bv_z^\dagger|, (\D_z^\dagger,\partial, U_1,\ldots,U_n,\ldots))$$
coincides with the distribution of 
$$(\sigma, (\ll^\bullet(W),\partial, V_1,\ldots,V_n,\ldots))$$
under $\N^{*,z}_0$. Here, we use the notation $\partial$ for the 
distinguished point (boundary point) of $\D_z^\dagger$ or of $\ll^\bullet(W)$, and  both
$(\D_z^\dagger,\partial,U_1,\ldots,U_n,\ldots)$ and $(\ll^\bullet(W),\partial,V_1,\ldots,V_n,\ldots)$ are viewed as a random variables taking values in the
space $\M^{\infty \bullet}$ of  metric measure spaces equipped with an infinite sequence 
of distinguished points. 

We then claim that, if $x,y\in \D_z\backslash \partial \D_z$ are such that there exists a 
$D^\partial$-geodesic from $x$ to $y$ that does not intersect $\partial \D_z$, 
\begin{equation}
\label{ega-dist1}
D^\partial(x,y)= \build{\inf_{y_0=x,y_1,\ldots,y_p=y}}_{D^{\partial,\dagger}(y_{i-1},y_i)<D^{\partial,\dagger}(\partial, y_{i-1})+D^{\partial,\dagger}(\partial,y_i)}^{}
\sum_{i=1}^p D^{\partial,\dagger}(y_{i-1},y_i),
\end{equation}
where the infimum is over all choices of the integer $p\geq 1$ and of the points $y_0,\ldots,y_p\in \D_z\backslash \partial \D_z$
such that $y_0=x$ and $y_p=y$, and the condition $D^{\partial,\dagger}(y_{i-1},y_i)<D^{\partial,\dagger}(\partial, y_{i-1})+D^{\partial,\dagger}(\partial,y_i)$
holds for $1\leq i\leq p$. The inequality $\leq$ in \eqref{ega-dist1} is clear from the fact that 
$D^{\partial,\dagger}(y_{i-1},y_i)=D^\partial(y_{i-1},y_i)$ under the condition 
$D^{\partial,\dagger}(y_{i-1},y_i)<D^{\partial,\dagger}(\partial, y_{i-1})+D^{\partial,\dagger}(\partial,y_i)$,
so that we can just use the triangle inequality for $D^\partial$. For the other inequality, we pick a $D^\partial$-geodesic $\gamma$
from $x$ to $y$ that does not intersect $\partial \D_z$
(so that it stays at a positive distance from $\partial \D_z$), and points $y_0=x,y_1,\ldots,y_p=y$ coming in this order along the geodesic $\gamma$ and such that 
$D^\partial(y_{i-1},y_i)<\min\{ D^\partial(u,\partial\D_z):u\in \gamma\}$
for $1\leq i\leq p$. It follows that $D^\partial(y_{i-1},y_i)=D^{\partial,\dagger}(y_{i-1},y_i)$, and 
$$\sum_{i=1}^p D^{\partial,\dagger}(y_{i-1},y_i)=\sum_{i=1}^p D^\partial(y_{i-1},y_i)=D^\partial(x,y).$$

We apply \eqref{ega-dist1} to each pair $(U_i,U_j)$, noting that a.s. the (unique) geodesic between $U_i$ and $U_j$
in $\D_z$ does not intersect $\partial \D_z$ (see Lemmas 17 and 18 in \cite{BM}). It follows that, a.s. 
for every $i,j\geq 1$,
\begin{equation}
\label{ega-dist2}
D^\partial(U_i,U_j)= \build{\inf_{i_0=i,i_1,\ldots,i_p=j}}_{D^{\partial,\dagger}(U_{i_{k-1}},U_{i_k})<D^{\partial,\dagger}(\partial, U_{i_{k-1}})
+D^{\partial,\dagger}(\partial,U_{i_k})}^{}
\sum_{k=1}^p D^{\partial,\dagger}(U_{i_{k-1}},U_{i_k}),
\end{equation}
using also the fact that any finite sequence $y_1,\ldots,y_{p-1}$ can be approximated by $U_{i_1},\ldots,U_{i_{p-1}}$
for suitable choices of $i_1,\ldots,i_{p-1}$, since the sequence $U_1,U_2,\ldots$ is (a.s.) dense in $\D_z$. 

Next recall the definition of the pseudo-metric $\Delta$ on $\t_\zeta^\circ$ under $\N^{*,z}_0$
(formulas \eqref{Dtrunc} and \eqref{defDelta} with $\t^\tr,\t^{\tr,\circ}$ and $D^{\tr,\circ}$
replaced respectively by $\t_\zeta,\t_\zeta^\circ$ and $D^\circ$)
and the definition \eqref{formulaD} of $D(a,b)$. Also note that,
$\N^{*,z}_0$ a.s., $D(a_*,a)=Z_a$ for every $a\in\t_\zeta$, where $a_*$ is any
point of $\t_\zeta$ such that $Z_{a_*}=0$.
We claim that we have, 
$\N^{*,z}_0$ a.s. for every
$a,b\in \t_\zeta^\circ$,
\begin{equation}
\label{ega-dist3}
\Delta(a,b)= \build{\inf_{a_0=a,a_1,\ldots,a_p=b}}_{D(a_{i-1},a_i)<D(a_*, a_{i-1})+D(a_*,a_i)}^{}
\sum_{i=1}^p D(a_{i-1},a_i),
\end{equation}
where the points $a_1,\ldots,a_{p-1}$ vary in $\t_\zeta^\circ$. 
To obtain the inequality $\leq$ in \eqref{ega-dist3}, we first note that, for every
$a,b\in\t_\zeta^\circ$, the condition $D(a,b)<D(a_*,a)+D(a_*,b)=Z_a+Z_b$
implies $D(a,b)=\Delta(a,b)$. Let us justify this. For $a,b\in \t^\circ_\zeta$, 
we have $D(a,b)\leq \Delta(a,b)$ by construction. On the other hand,
if, in the definition of  $D(a,b)$ as an infimum of sums, one of the terms is such that $\min_{[a_{i-1},a_i]}Z_c=\min_{[a_i,a_{i-1}]}Z_c=0$, 
we get $D^\circ(a_{i-1},a_i)=Z_{a_i}+Z_{a_i-1}$ and (using the property $D^\circ(u,v)\geq |Z_u-Z_v|$)
we obtain that the corresponding sum is greater than or equal to $Z_a+Z_b$. Thus, if $D(a,b)<Z_a+Z_b$, this means that such cases can be discarded, so that $D(a,b)=\Delta(a,b)$.

The preceding discussion shows that we can replace $\sum_{i=1}^p D(a_{i-1},a_i)$ by $\sum_{i=1}^p \Delta(a_{i-1},a_i)$ in the right-hand side 
of \eqref{ega-dist3}, and we 
use the triangle inequality for $\Delta$ to obtain the inequality $\leq$ in \eqref{ega-dist3}. 

The reverse inequality in \eqref{ega-dist3} follows from the definition \eqref{defDelta} of $\Delta(a,b)$
(with $D^{\tr,\circ}$ replaced by $D^\circ$ as explained before Corollary \ref{consDisk}),  noting 
that in this definition we need only consider the case where $\Delta^{\circ}(a_{i-1},a_i)<\infty$
for every $1\leq i\leq k$, which implies 
$D^\circ(a_{i-1},a_i)< Z_{a_{i-1}}+Z_{a_i}=D(a_*, a_{i-1})+D(a_*,a_i)$
and also
 $$D(a_{i-1},a_i)\leq D^\circ(a_{i-1},a_i)=\Delta^\circ(a_{i-1},a_i).
$$

Then, we can apply \eqref{ega-dist3} to each pair $(V_i,V_j)$ and we get, similarly as in \eqref{ega-dist2},
\begin{equation}
\label{ega-dist4}
\Delta(V_i,V_j)= \build{\inf_{i_0=i,i_1,\ldots,i_p=b}}_{D(V_{i_{k-1}},V_{i_k})<D(\partial, V_{i_{k-1}})+D(\partial,V_{i_k})}^{}
\sum_{k=1}^p D(V_{i_{k-1}},V_{i_k}).
\end{equation}

Write $U_0=\partial$ and $V_0=\partial$ for convenience.
Thanks to the identity in distribution mentioned at the beginning of this section
(following from Theorem \ref{identif-glued-unpointed}), we know that the collection
$$\Big(D(V_i,V_j)\Big)_{i,j\geq 0}\quad\hbox{under}\quad \N^{*,z}_0,$$
has the same distribution as $(D^{\partial,\dagger}(U_i,U_j))_{i,j\geq 0}$. Looking at \eqref{ega-dist2}
and \eqref{ega-dist4}, we infer that the collection
$$\Big(\Delta(V_i,V_j)\Big)_{i,j\geq 1}\quad\hbox{under}\quad \N^{*,z}_0,$$
has the same distribution as $(D^{\partial}(U_i,U_j))_{i,j\geq 0}$,
and that this identity in distribution holds jointly with that 
of $|\bv_z|$ and $\sigma$. 

By Lemma \ref{infinite-points}, we know that 
$$\Big(\{U_1,\ldots,U_n\}, D^{\partial}, \frac{1}{n}\sum_{i=1}^n \delta_{U_i}\Big)$$
converges a.s. in $\M$ to $(\D_z,D^{\partial},\wh\bv_z)$, where $\wh\bv_z=|\bv_z|^{-1}\,\bv_z$, and similarly
$$\Big(\{V_1,\ldots,V_n\}, \Delta, \frac{1}{n}\sum_{i=1}^n \delta_{V_i}\Big)$$
converges a.s. in $\M$ to $(\Theta,\Delta,\wh\bV)$, where $\wh\bV$
is the normalized volume measure on $\Theta$. Notice that
the assumption about the topological support in Lemma \ref{infinite-points} is easily verified. We then conclude that
$(|\bv_z|,(\D_z,D^{\partial},\wh\bv_z))$ and $(\sigma,(\Theta,\Delta,\wh\bV))$
have the same distribution. Theorem \ref{main} stated in the introduction now follows.

\medskip
\rem By \cite{Bet}, we know that the Brownian disk is
a.s. homeomorphic to the closed unit disk of the plane.
The boundary of $\Theta$ may then be defined as the
set of all points $x\in\Theta$ having no neighborhood 
homeomorphic to an open disk. It will follow from the
results of the subsequent sections that this boundary coincides 
with $\partial\Theta$ as defined above. Proposition \ref{loop-trunc} below shows that
$\partial \Theta$ is the range of a simple loop in $\Theta$, and Theorem
\ref{Bnet} provides an embedding of a whole sequence of Brownian disks $(\Theta^{(j)}, \Delta^{(j)}, \bV^{(j)})$ (each 
of which is a copy of $(\Theta, \Delta, \bV)$ under $\N^{*,z^{(j)}}_0$,
for a certain random value of $z^{(j)}$) in the Brownian map, in such a way that $\partial\Theta^{(j)}$
appears as 
the topological boundary of $\Theta^{(j)}$ in the Brownian map, which itself is known to be homeomorphic to
the two-dimensional sphere.

\section{The uniform measure on the boundary}
\label{sec:unibdry}

Let us fix $r>0$ and argue under the measure $\N_r$. Recall from the
beginning of Section \ref{sec:resu-snake} the definition of the
exit local time $(\ell^0_s)_{s\geq 0}$, and the fact
that $\ell^0_\sigma=\z_0$. We will need a different
approximation of this process. Recall the notation $\tau_0(\w)=\inf\{t\in[0,\zeta_{(\w)}]:\w(t)=0\}$ for $\w\in\mathcal{W}$.

\begin{proposition}
\label{approx-LT}
We have $\N_r$ a.e., for every $t\geq 0$,
$$\ell^0_t=\lim_{\ve\to 0} \frac{1}{\ve^2}\int_0^t \mathrm{d}s\,\mathbf{1}_{\{\tau_0(W_s)\geq \zeta_s,\wh W_s<\ve\}}.$$
\end{proposition}

We note that the condition $\tau_0(W_s)\geq \zeta_s$ holds only if $\tau_0(W_s)= \zeta_s$ or $\tau_0(W_s)=\infty$. We refer to Appendix A below for a proof of Proposition
\ref{approx-LT}, which is a refinement of \cite[Lemma 14]{ALG}. The motivation for this result is the fact that
it will allow us to reinterpret the exit local time in terms of $W$ truncated at $0$. As in Section
\ref{sec:cons-metric}, we write $W^{\tr}=\mathrm{tr}_0(W)$, so that $W^\tr_s= W_{\eta_s}$, with
$(\eta_s)_{s\geq 0}$ defined in \eqref{def-eta}.
We set, for every $s\geq 0$,
$$L^0_s=\ell^0_{\eta_s}.$$
Then as a consequence of the preceding proposition, we have,
$\N_r$ a.e., for every $t\geq 0$,
\begin{equation}
\label{approx-exit}
L^0_t=\lim_{\ve\to 0} \frac{1}{\ve^2}\int_0^t \mathrm{d}s\,\mathbf{1}_{\{\wh W^\tr_s<\ve\}}.
\end{equation}
We note that $(L^0_s)_{s\geq 0}$ has continuous sample paths, $\N_0$ a.e.: The point is that,
if $\eta_{s-}<\eta_s$, then necessarily $\tau_0(W_u)\leq\zeta_u$ for every 
$u\in(\eta_{s-},\eta_s)$ and this implies that
$\ell^0_{\eta_{s-}}=\ell^0_{\eta_s}$ (for instance by using the preceding proposition, together
with the fact that $\{s:\tau_0(W_s)=\zeta_s\}$ has zero Lebesgue measure).

Recall the notation $\t^\tr$, $p_\tr$ introduced in Section \ref{sec:cons-metric}.
Under the probability measure $\N_r^{[0]}=\N_r(\cdot\mid W_*\leq0)$,
the construction of Section \ref{sec:cons-metric}
yields a pointed measure metric space $(\Theta',\Delta,\bV')$ associated with 
$W^\tr$. Recall that $\Theta'$ is obtained as a quotient space of $\t^\tr$ and write $\Pi_\tr$
for the canonical projection from $\t_\tr$ onto $\Theta'$. Also recall that, by definition,
$$\partial \Theta'= \Pi_\tr(\{a\in \t^\tr: Z_a=0\}).$$

\begin{proposition}
\label{loop-trunc}
For every $t\in[0,\z_0)$, set
$$\gamma^0_t:=\inf\{s\geq 0: L^0_s>t\}$$
and $\gamma^0_{\z_0}:= \gamma^0_0$. Then, $\N^{[0]}_r$ a.s., the mapping
$$[0,\z_0]\ni t \mapsto \Pi_\tr(p_\tr(\gamma^0_t))$$
defines a simple continuous loop in $\Theta'$ whose range is $\partial\Theta'$. 
\end{proposition}

\proof From \eqref{support-exit}, we have
\begin{equation}
\label{supp-exit}
\mathrm{supp}(\mathrm{d}L^0_s)=\{s\geq 0: \tau_0(W^\tr_s)=\zeta^\tr_s\},
\end{equation}
with the obvious notation $\zeta^\tr_s$. 

Let us start by showing that the mapping of the proposition is
continuous. Since $t\mapsto \gamma^0_t$ is right-continuous with left
limits on $[0,\z_0)$, the desired continuity on $[0,\z_0)$  will follow
if we can prove that $ \Pi_\tr(p_\tr(\gamma^0_t))= \Pi_\tr(p_\tr(\gamma^0_{t-}))$
whenever $t\in[0,\z_0)$ is such that $\gamma^0_{t-}<\gamma^0_t$. 
If this occurs, we have by definition
$$L^0_{\gamma^0_t}=L^0_{\gamma^0_{t-}}$$
showing that, for $\gamma^0_{t-}<u<\gamma^0_{t}$ we must have
$\tau_0(W^\tr_u)=\infty$ (use \eqref{supp-exit}). But we know
(Proposition \ref{properDelta}) that this implies $\Delta(p_\tr(\gamma^0_t),p_\tr(\gamma^0_{t-}))=0$ and thus
$ \Pi_\tr(p_\tr(\gamma^0_t))= \Pi_\tr(p_\tr(\gamma^0_{t-}))$. 
A similar argument (noting that $\gamma^0_{\z_0-}=\sup\{s\in[0,\sigma(W^\tr)]:\wh W^\tr_s=0\}$)
implies that $ \Pi_\tr(p_\tr(\gamma^0_{\z_0-}))= \Pi_\tr(p_\tr(\gamma^0_0))$, giving the desired 
continuity at $t=\z_0$. 

We next prove that the range of the mapping of the proposition
is $\partial\Theta'$. It follows from \eqref{supp-exit} that any
$u\geq 0$ such that $u=\gamma^0_t$ for some $t\in[0,\z_0]$ 
is such that $\tau_0(W^\tr_u)=\zeta^\tr_u$, which implies
$p_\tr(u)\in\{a\in\t^\tr:Z_a=0\}$, and $ \Pi_\tr(p_\tr(\gamma^0_u))\in\partial\Theta'$.

Conversely, if $a\in\t^\tr$ and $Z_a=0$, then $a=p_\tr(s)$ for some
$s\in[0,\sigma(W^\tr)]$ such that $\wh W^\tr_s=0$, hence necessarily
$\tau_0(W^\tr_s)=\zeta^\tr_s$. By \eqref{supp-exit}, $s$ must be an increase time of $L^0$
and so $s=\gamma^0_t$ or $s=\gamma^0_{t-}$ for some $t\in [0,\z_\sigma]$. But
we saw that $\Pi_\tr(p_\tr(\gamma^0_t))=\Pi_\tr(p_\tr(\gamma^0_{t-}))$, so that,
in either case, we obtain that $\Pi_\tr(a)$ is in the range of the mapping of the proposition.

It only remains to verify that the mapping $[0,\z_0)\ni t \mapsto \Pi_\tr(p_\tr(\gamma^0_t))$
is injective. However, if $0\leq s<t<\z_0$, there exist times $u$ with $\gamma^0_s<u<\gamma^0_t$
such that $\wh W^\tr_u=0$, and, by Proposition \ref{properDelta} (ii), we know that $p_\tr(\gamma^0_s)$ and $p_\tr(\gamma^0_t)$
are then not identified in $\Theta'$. \endproof

Write $(\Lambda_t)_{0\leq t\leq \z_0}$ for the loop with range 
$\partial\Theta'$ defined in Proposition \ref{loop-trunc}. The image of
Lebesgue measure on $[0,\z_0]$ under $t\mapsto \Lambda_t$ then defines a finite measure 
$\nu$ on $\partial\Theta'$ with total mass $\z_0$. This is the uniform
measure on $\partial \Theta'$ in the sense of the following proposition.

\begin{proposition}
\label{mesure-bdry}
Almost surely under $\N^{[0]}_r$, for every bounded
continuous function $\varphi$ on $\Theta'$,
$$\langle \nu, \varphi\rangle = \lim_{\ve\to 0} \frac{1}{\ve^2}\int_{\Theta'} \bV'(\mathrm{d}x)\,\varphi(x)\,\mathbf{1}_{\{\Delta(x,\partial\Theta')<\ve\}}.$$
\end{proposition}

\proof By definition,
$$
\langle \nu, \varphi\rangle=\int_0^{L^0_\sigma} \mathrm{d}t\,\varphi(\Lambda_t)
=\int_0^{L^0_\sigma} \mathrm{d}t\,\varphi(\Pi_\tr(p_\tr(\gamma^0_t)))
=\int_0^{\sigma_\tr} \mathrm{d}L^0_s\,\varphi(\Pi_\tr(p_\tr(s))),
$$
where $\sigma_\tr=\sigma(W^\tr)$. By \eqref{approx-exit}, we know
that the measures 
$$\frac{1}{\ve^2}\,\mathbf{1}_{\{\wh W^\tr_t<\ve\}}\,\mathbf{1}_{[0,\sigma_\tr]}(t)\,\mathrm{d}t$$
converge weakly to $\mathrm{d}L^0_t$ as $\ve\to 0$. Consequently,
$$\langle \nu, \varphi\rangle
= \lim_{\ve\to 0} \frac{1}{\ve^2}\int_0^{\sigma_\tr} \mathrm{d}t\,\mathbf{1}_{\{\wh W^\tr_t<\ve\}}\,\varphi(\Pi_\tr(p_\tr(t))).$$
The desired result follows since, by Proposition \ref{properDelta} (i) and the definition of the volume measure on $\Theta'$,
$$\int_0^{\sigma_\tr} \mathrm{d}t\,\mathbf{1}_{\{\wh W^\tr_t<\ve\}}\,\varphi(\Pi_\tr(p_\tr(t)))
=\int_{\Theta'} \bV'(\mathrm{d}x)\,\varphi(x)\,\mathbf{1}_{\{\Delta(x,\partial\Theta')<\ve\}}.$$
\endproof

We can restate the last two propositions in terms of the Brownian disk.

\smallskip
\begin{corollary}
\label{bdry-disk}
Let $z>0$. The following properties hold $\N^{*,z}_0$ a.s.
\begin{enumerate}
\item[\rm(i)] The limit
$$L_s:=\lim_{\ve \to 0} \frac{1}{\ve^2}\int_0^s \mathrm{d}t\,\mathbf{1}_{\{\wh W_t<\ve\}}$$
exists for every $s\in[0,\sigma]$, and defines a continuous increasing function
with $L_\sigma=z$.
\item[\rm(ii)] Set $\gamma_t=\inf\{s\in[0,\sigma]:L_s>t\}$ for every $t\in[0,z]$, and $\gamma_z=\sigma$. 
The image of $(\gamma_t)_{0\leq t\leq z}$ under the canonical projection 
from $[0,\sigma]$ onto $\Theta$ defines a simple continuous loop whose range is $\partial\Theta$.
\item[\rm(iii)] Let $(\Lambda_t)_{0\leq t\leq z}$ be the continuous loop obtained in {\rm(ii)}, and let
$\nu_z$ be the image of Lebesgue measure on $[0,z]$ under the mapping $t\mapsto \Lambda_t$.
Then,  for every bounded
continuous function $\varphi$ on $\Theta$,
$$\langle \nu_z, \varphi\rangle = \lim_{\ve\to 0} \frac{1}{\ve^2}\int_{\Theta} \bV(\mathrm{d}x)\,\varphi(x)\,\mathbf{1}_{\{\Delta(x,\partial\Theta)<\ve\}}.$$
\end{enumerate}
\end{corollary}

All these properties follow from the previous propositions and the
re-rooted representation theorem (Theorem \ref{re-root-rep}),
which links the measure
$\N^*_0$ with the laws of  truncated snakes. 
To be more precise, one first obtains that the properties of the corollary hold under $\N^*_0$.
It follows that they hold under $\N^{*,z}_0$ for a.a. values of $z>0$, but then
a scaling argument shows that they must hold for every value of $z>0$. 

Proposition \ref{measure-boundary} stated in the introduction follows from part (iii) of Corollary \ref{bdry-disk}. 

\section{Brownian disks filling in the holes of the Brownian net}
\label{sec:Bnet}

In this section, we deal with the free Brownian map $\mm$
under the measure $\N_0$, as defined in Section \ref{sec:defBM}. Recall that
$\mm$ has the two distinguished points $x_*$
and $x_0$, and that $D(x_*,x_0)=-W_*$.

The metric net $\mathbf{N}$ is the closed subset of $\mm$
defined as the closure of
$$\bigcup_{0<r<-W_*} \partial \mathfrak{H}_r,$$
where $\mathfrak{H}_r$ denotes the hull of radius $r$, which is defined for $0<r<-W_*$
(recall from Section \ref{conv-glued-bdry} that this hull is the complement of 
the connected component of the complement of the ball $B(x_*,r)$ containing $x_0$).
The connected components of
the complement
of the metric net are called Brownian disks in \cite{MS0}. The goal of the 
next theorem is to provide a precise justification to
this terminology.

Before stating this result, recall from \cite[Section 8]{subor}
that connected components of $\mm\backslash\mathbf{N}$
are in one-to-one correspondence with 
jump times of the exit measure process $(\z_r)_{r<0}$
(we noticed in Section \ref{sec:resu-snake} that this process
has a c\`adl\`ag modification). We 
denote the sequence of these jump times (ordered by decreasing size
of the jumps) by $r_1,r_2,\ldots$ and the associated connected components 
of $\mm\backslash\mathbf{N}$ by $\mathbf{D}^{(1)},\mathbf{D}^{(2)},\ldots$. As usual,
$\ov{\mathbf{D}}^{(j)}$ denotes the closure of $\mathbf{D}^{(j)}$ in $\mm$. For every
$j\geq 1$, we can equip $\mathbf{D}^{(j)}$ with the intrinsic metric $d^{(j)}_{\mathrm{intr}}$,
such that $d^{(j)}_{\mathrm{intr}}(x,y)$ is the infimum of the lengths (computed with
respect to $D$) of continuous paths in $\mathbf{D}^{(j)}$ connecting $x$ to $y$.
We also write $\bv^{(j)}$ for the restriction of the volume measure on 
$\mm$ to $\mathbf{D}^{(j)}$.

\begin{theorem}
\label{Bnet}
$\N_0$ a.e., for every $j=1,2,\ldots$, the metric $d^{(j)}_{\mathrm{intr}}$ has
a continuous extension to $\ov{\mathbf{D}}^{(j)}$, which is a metric on $\ov{\mathbf{D}}^{(j)}$. Furthermore,
conditionally on the exit measure process $(\z_r)_{r<0}$, 
the measure metric spaces $(\ov{\mathbf{D}}^{(j)},d^{(j)}_{\mathrm{intr}},\bv^{(j)})$, $j=1,2,\ldots$
are independent free Brownian disks with respective perimeters $\Delta\z_{r_j}$, $j=1,2,\ldots$
\end{theorem} 

\proof We rely strongly on the results of \cite[Section 8]{subor} and \cite{ALG}. According to
\cite[Section 8]{subor}, the connected components $\mathbf{D}^{(1)},\mathbf{D}^{(2)},\ldots$
are also in one-to-one correspondence with excursions of the Brownian snake above
its minimum, as defined in \cite[Section 3]{ALG}. For the sake of completeness, we briefly
recall the definition of these excursions (notice that this definition is quite different from the
definition of excursions outside $(b,\infty)$ in Section \ref{sec:Bstrunc}). We first introduce the
closed subset $\mathbf{F}$ of $\t_\zeta$ defined by
$$\mathbf{F}:=p_\zeta(\{s\in [0,\sigma]: \un W_s=\wh W_s\}).$$
Roughly speaking, each excursion above the minimum describes the process of labels restricted
to a connected component of the open set $\t_\zeta\backslash\mathbf{F}$. To explain this, say that
$a\in \mathbf{F}$ is
an excursion debut if $a$ has a strict descendant $b$ such that $Z_b>Z_a$ and $Z_c > Z_a$ for every interior
point $c$ of the geodesic segment between $a$ and $b$ in $\t_\zeta$. For a given excursion debut $a$, the collection of such points $b$ then forms a
connected component of $\t_\zeta\backslash\mathbf{F}$ denoted by $\mathbf{C}_{(a)}$. Moreover, $\Pi(\mathbf{C}_{(a)})$ 
is a connected component of $\mm\backslash\mathbf{N}$. In this way, one obtains one-to-one correspondences
between excursion debuts, connected components of $\t_\zeta\backslash\mathbf{F}$ and 
connected components of $\mm\backslash\mathbf{N}$ (see \cite[Section 8]{subor}).
Next, if $a$ is a fixed excursion debut, there are exactly two times $0<s_1<s_2<\sigma$ such that $p_\zeta(s_1)=p_\zeta(s_2)=a$, and the labels of descendants
of $a$ are described by the snake trajectory $\omega^{(a)}$ defined by saying that $\wh \omega^{(a)}_s=\wh\omega_{s_1+s} - Z_a$
and $\zeta_{(\omega^{(a)}_s)}=\zeta_{s_1+s}-\zeta_{s_1}$ for $0\leq s\leq \sigma(\omega^{(a)}):=s_2-s_1$. 
The excursion above the minimum associated with $a$ (or with the connected component $\Pi(\mathbf{C}^{(a)})$) is then the truncation at $0$ of the snake trajectory $\omega^{(a)}$ (this requires a 
minor extension of the definition given in Section \ref{sec:snake}, where the truncation level had to be different
from the initial point of the snake trajectory --- see the next section for an analogous 
definition of excursions above $0$).

For every $j=1,2,\ldots$, we write $W^{(j)}$ for the 
excursion above the minimum associated with $\mathbf{D}^{(j)}$,
$\t^{(j)}$ for the genealogical tree of $W^{(j)}$ and $(Z^{(j)}_b)_{b\in\t^{(j)}}$ for the labels on $\t^{(j)}$,
and we also set $\t^{(j),\circ}:=\{b\in \t^{(j)}:Z^{(j)}_b>0\}$ and $\partial \t^{(j)}=\t^{(j)}\backslash \t^{(j),\circ}$.  

By \cite[Theorem 40]{ALG}, 
conditionally on $(\z_r)_{r<0}$, $W^{(1)}, W^{(2)},\ldots$ are independent 
and distributed according to $\N^{*,\Delta\z_{r_1}}_0,\N^{*,\Delta\z_{r_2}}_0,\ldots$
respectively. We can thus define, for every $j=1,2,\ldots$, a measure metric space
$(\Theta^{(j)},\Delta^{(j)},\bV^{(j)})$ constructed from $W^{(j)}$ via the 
procedure of Section \ref{sec:cons-metric}. Thanks to Theorem \ref{main}, the proof of Theorem \ref{Bnet} then
boils down to verifying that, $\N_0$ a.e. for every $j=1,2,\ldots$,
$d^{(j)}_{\mathrm{intr}}$ can be extended continuously to a metric on $\ov{\mathbf{D}}^{(j)}$,
in such a way that
\begin{equation}
\label{Bnet1}
(\ov{\mathbf{D}}^{(j)},d^{(j)}_{\mathrm{intr}},\bv^{(j)})=(\Theta^{(j)},\Delta^{(j)},\bV^{(j)})
\end{equation}
in the sense of equality in $\M$. 

Let us proceed to this verification. We note that $\t^{(j)}$ is identified canonically to the closure $\ov{C}^{(j)}$  of the connected component 
$C^{(j)}$ of $\t_\zeta\backslash\mathbf{F}$
corresponding to $\mathbf{D}^{(j)}$ (this basically follows from the construction of excursions
above the minimum, see \cite[Section 3]{ALG}). 
More precisely, there is an isometric bijection $\varphi_j:\t^{(j)}\la \ov{C}^{(j)}$
which preserves labels up to a shift (to be specific, $Z_{\varphi_j(b)}=Z^{(j)}_b + z_j$, where
$z_j=Z_{a_j}$ if $a_j$ is the excursion debut corresponding to $C^{(j)}$) and is such that
the following holds. We have 
$\varphi_j(\t^{(j),\circ})=C^{(j)}$, and $\varphi_j$
maps the volume measure on $\t^{(j)}$ to the volume measure on $\ov{C}^{(j)}$, and preserves 
intervals, in the sense that, for every $a,b\in\t^{(j)}$,
$$\varphi_j([a,b])=[\varphi_j(a),\varphi_j(b)]\cap \ov C^{(j)}.$$
In particular, if $a,b\in \t^{(j),\circ}$, we have $[a,b]\cap \partial \t^{(j)}=\varnothing$
if and only if  $[\varphi_j(a),\varphi_j(b)]\cap \partial C^{(j)}=\varnothing$ and then
$\varphi_j([a,b])=[\varphi_j(a),\varphi_j(b)]$. 

As already mentioned, we have $\Pi(C^{(j)})=\mathbf{D}^{(j)}$ and moreover $\Pi(\partial C^{(j)})=
\partial \mathbf{D}^{(j)}$ (see \cite[Lemma 14]{subor}).
Writing $\Pi^{(j)}$ for the canonical projection from $\t^{(j)}$ onto $\Theta^{(j)}$, we next observe 
that, for $a,b\in \t^{(j),\circ}$, we have
\begin{equation}
\label{Bnet2}
\Pi^{(j)}(a)=\Pi^{(j)}(b)\quad\hbox{if and only if}\quad \Pi(\varphi_j(a))=\Pi(\varphi_j(b)).
\end{equation}
This 
is a straightforward consequence of the preceding remarks: If $\Pi^{(j)}(a)=\Pi^{(j)}(b)$,
then $Z^{(j)}_a=Z^{(j)}_b$ and at least one of the intervals $[a,b]$ or $[b,a]$ does not intersect $\partial\t^{(j),\circ}$
and is such that labels on this interval remain greater than or equal to $Z^{(j)}_a$.
It follows that a similar property holds for $\varphi_j(a)$ and $\varphi_j(b)$ so that $ \Pi(\varphi_j(a))=\Pi(\varphi_j(b))$.
The converse is obtained in the same way.

Thanks to \eqref{Bnet2}, the restriction of $\varphi_j$ to $\t^{(j),\circ}$ induces a bijection
$\Phi^{(j)}$ from $\Theta^{(j)}\backslash \partial \Theta^{(j)}=\Pi^{(j)}(\t^{(j),\circ})$ onto $\mathbf{D}^{(j)}=\Pi(C^{(j)})$, and
we have $\Phi^{(j)}\circ\Pi^{(j)}(a)=\Pi \circ \varphi_j(a)$ for every $a\in\t^{(j),\circ}$. 
This bijection is an isometry:  
$$\Delta^{(j)}(x,y)=d^{(j)}_{\mathrm{intr}}(\Phi^{(j)}(x),\Phi^{(j)}(y))$$
for every $x,y\in \Theta^{(j),\circ}$. Again this is a simple consequence
of the preceding remarks. To prove that  $\Delta^{(j)}(x,y)\geq d^{(j)}_{\mathrm{intr}}(\Phi^{(j)}(x),\Phi^{(j)}(y))$,
we note that, in the definition (\ref{defDelta}) of $\Delta(a,b)$, each term
$\Delta^\circ(a_{i-1},a_i)$ can be interpreted (provided it is finite) as the length of a curve connecting
$a_{i-1}$ to $a_i$, namely the union of two simple geodesics starting from $a_{i-1}$
and from $a_i$ up to their merging time (see \cite{Geo} for the definition and properties
of simple geodesics). To get the reverse inequality, we use the definition 
of the intrinsic distance, formula \eqref{formulaD} expressing $D(a,b)$ as an infimum
and the fact that, with an obvious notation, $D^\circ(a,b)=\Delta^{(j),\circ}(\varphi_j^{-1}(a),\varphi_j^{-1}(b))$
if $a,b\in C^{(j)}$ are such that $D^\circ(a,b)<Z_a+Z_b-2z_j$
(see the remark before Proposition \ref{propDelta}). Details are left to the reader. 

Since $\Phi^{(j)}$ is isometric, simple arguments show that it has a unique extension to a mapping $\bar \Phi^{(j)}$ from $\Theta^{(j)}$ onto 
$\ov{\mathbf{D}}^{(j)}$, such that $\bar\Phi^{(j)}\circ\Pi^{(j)}=\Pi \circ \varphi_j$. The fact that $\bar \Phi^{(j)}$ is onto is
immediate from the last equality.
We claim that
$\bar\Phi^{(j)}$ is also one-to-one. This will complete the proof of the proposition, since 
then the bijection $\bar\Phi^{(j)}$ will map $\Delta^{(j)}$ to the (unique) continuous
extension of $d^{(j)}_{\mathrm{intr}}$ to $\ov{\mathbf{D}}^{(j)}$, thus showing that
\eqref{Bnet1} holds (the fact that $\bar\Phi^{(j)}$ also maps the volume 
measure on $\Theta^{(j)}$ to the volume measure on $\ov{\mathbf{D}}^{(j)}$
is immediate).

It remains to verify that $\bar\Phi^{(j)}$ is one-to-one. So let
$x,y\in \partial \Theta^{(j)}$ with $x\not=y$. Write $x=\Pi^{(j)}(a)$ and 
$y=\Pi^{(j)}(b)$, with $a,b\in\partial \t^{(j)}$. By Proposition \ref{properDelta} (ii) (and Corollary \ref{consDisk}), the fact that $x\not =y$
implies that there exists $c\in]a,b[$ such that $Z^{(j)}_c=0$ (and the same
if $]a,b[$ is replaced by $]b,a[$). In fact, there must
exist infinity many such values of $c$, because otherwise an application of the
triangle inequality for $\Delta^{(j)}$ would imply that $\Delta^{(j)}(a,b)=0$ and so
$\Pi^{(j)}(a)=\Pi^{(j)}(b)$. Next, from the fact that $\varphi_j$ preserves intervals (and labels up to a shift) we also
get that both $]\varphi_j(a),\varphi_j(b)[$ and $]\varphi_j(b),\varphi_j(a)[$
contain infinitely many values of $c$ such that $Z_c=Z_{\varphi_j(a)}=Z_{\varphi_j(b)}$.
This implies that $\Pi\circ \varphi_j(a)\not =\Pi\circ \varphi_j(b)$ so that $\bar\Phi^{(j)}(x)\not = \bar\Phi^{(j)}(y)$.
This completes the proof. \endproof

\section{Brownian disks in the Brownian map}
\label{sec:BdBm}

This section is devoted to the proof of Theorem \ref{ccBm}, which is somewhat more delicate than the proof
of Theorem \ref{Bnet} in the preceding section. In several parts of the proof, we argue under the 
measure $\N_0$ rather than under $\N^{(1)}_0$.  

For the reader's convenience, we have broken the
proof in several steps. The first three steps are devoted to checking that connected components of
$\mm\backslash B(x_*,r)$ are in one-to-one correspondence with connected components of $\{a\in\t_\zeta: Z_a>W_*+r\}$,
then using this to define the boundary lengths of components of $\mm\backslash B(x_*,r)$ (via Corollary \ref{bdry-disk}), and
to verify that the intrinsic distance on each component $\mathbf{C}$ has a continous extension to the closure $\ov{\mathbf{C}}$.
It turns out that it is easier to derive the latter properties for the random value $r=-W_*$ under $\N_0$, and then a re-rooting 
invariance argument allows us to obtain that they hold for a fixed value of $r$, still under $\N_0$ (Step 2). We in fact want these properties
to also hold  under $\N_0^{(s)}=\N_0(\cdot\mid \sigma = s)$, for every $s>0$, and then we need an extra absolute continuity argument (Step 3). 
Step 4 is the heart of the proof. Again, we start by considering the special case $r=-W_*$, and we use the excursion theory of \cite{ALG}
to derive the conditional independence and the law of the connected components of $\mm\backslash B(x_*,-W_*)$, and a re-rooting 
invariance argument then allows us to prove that the same properties hold under $\N^{(s)}_0$ for almost all values of $r$. Since we want the
result to hold for {\it every} value of $r$, we need another argument (Step 5) relying on the absolute continuity properties already 
used in Step 3.

We will use the local times of the
process $(\wh W_s)_{s\geq 0}$. There exists a process $(L^{(u)})_{u\in\R}$ with continuous sample paths
and nonnegative values, such
that, $\N_0$ a.e. for every nonnegative measurable function $\varphi$ on $\R$,
$$\int_0^\sigma \mathrm{d}t\,\varphi(\wh W_t) =\int_\R \mathrm{d}u\,\varphi(u)\,L^{(u)}.$$
The same result holds if we replace $\N_0$ by $\N^{(s)}_0$ for any $s>0$.
See \cite{BMJ} for the case of $\N^{(1)}_0$, from which the other cases follow immediately.
Notice that $L^{(r)}=0$ if $r\notin(W_*,W^*)$. On the other hand,
the explicit distribution of $L^{(0)}$ under $\N^{(1)}_0$ found in \cite[Corollary 3.4]{BMJ} shows that $L^{(0)}>0$, $\N_0$ a.e. or $\N_0^{(s)}$ a.s.
We also set
$$\sigma_-:=\int_0^\sigma \mathrm{d}t\,\mathbf{1}_{\{\wh W_t<0\}}.$$

\noindent{\it Step 1.} $\N_0$ a.e., or $\N_0^{(s)}$ a.s., for every $r\in(W_*,W^*)$, the mapping $\cc\mapsto \Pi(\cc)$ induces a one-to-one 
correspondence between connected components of $\{a\in\t_\zeta: Z_a>W_*+r\}$ and connected components of $\mm\backslash B(x_*,r)$.

Indeed, let $\cc$ be a connected component of $\{a\in\t_\zeta: Z_a>W_*+r\}$. Then $\Pi(\cc)$ is connected and it
is easy to verify that $\Pi^{-1}(\Pi(\cc))=\cc$
(the continuous cactus bound \eqref{ctscactus} shows that if $a\in \cc$ the property $\Pi(a)=\Pi(b)$
may only hold if $b\in \cc$). Since the topology on $\mm$ is the quotient topology, it follows that
$\Pi(\cc)$ is open. Finally, if $b\notin \cc$ the geodesic path in $\t_\zeta$
from $b$ to an arbitrary point $a$ of $\cc$ visits a point $c$ such that $Z_c\leq W_*+r$, so that by 
\cite[Proposition 3.1]{Geo}, any path in $\mm$ from $\Pi(b)$ to $\Pi(a)$ must visit a point of $B(x_*,r)$, thus proving that
$\Pi(b)$ does not belong to the connected component of $\mm\backslash B(x_*,r)$ containing $\Pi(\cc)$, and we conclude that this connected component is indeed $\Pi(\cc)$. 
For future use, we also notice that the boundary of $\Pi(\cc)$ is equal to $\Pi(\partial\cc)$ (we leave the proof as an easy exercise).

\medskip
\noindent{\it Step 2.} We claim that, for every fixed $r>0$, the following properties hold $\N_0$ a.e. on the event $\{r<W^*-W_*\}$:
\begin{enumerate}
\item[(a)] For every connected component $\mathbf{C}$ of $\mm\backslash B(x_*,r)$, the limit
$$|\partial\mathbf{C}|:=\lim_{\ve\to 0} \frac{1}{\ve^2}\,\bv(\{x\in\mathbf{C}: D(x,\partial \mathbf{C})<\ve\})$$
exists in $(0,\infty)$. 
\item[(b)] The quantities $|\partial\mathbf{C}|$ when $\mathbf{C}$ varies among the connected components of $\mm\backslash B(x_*,r)$
are distinct.
\item[(c)] For every connected component $\mathbf{C}$ of $\mm\backslash B(x_*,r)$, the intrinsic metric on $\mathbf{C}$ has a continuous extension to the closure $\ov{\mathbf{C}}$, which is a metric on $\ov{\mathbf{C}}$.
\end{enumerate}

To obtain the preceding properties when $r$ is a fixed (deterministic) quantity, we will in fact start with the case where 
$r=-W_*$ is random. Note that, by Step 1, the connected components of $\mm\backslash B(x_*,-W_*)$ are 
in one-to-one 
correspondence with the connected components of $\{a\in\t_\zeta: Z_a>0\}$. 
Let 
$\cc^{(1)},\cc^{(2)},\ldots$ be the connected components of the
set $\{a\in\t_\zeta:Z_a>0\}$ (we will explain later how they are ordered), and let
$\mathbf{C}^1=\Pi(\cc^{(1)}), \mathbf{C}^2=\Pi(\cc^{(2)}),\ldots$ be the
corresponding connected components of $\mm\backslash B(x_*,-W_*)$. These objects are also in one-to-one correspondence with
the Brownian snake excursions above $0$, which will be denoted by $W^{(1)},W^{(2)},\ldots$. 
Let us recall the definition of these excursions, following closely
\cite{ALG}. We first notice that, for every $j=1,2,\ldots$, there
is a unique $a_j\in\t_\zeta$ that belongs to the closure of $\cc^{(j)}$ 
and is such that $\cc^{(j)}$ is contained in the set of descendants of $a_j$ in $\t_\zeta$. 
Furthermore, there exist $0<u_j<v_j<\sigma$ such the latter set of
descendants is equal to $p_\zeta([u_j,v_j])$, and we define a snake trajectory $\wt W^{(j)}\in\S_0$ with duration $v_j-u_j$
by setting 
$$\wt W^{(j)}_s(t):= W_{(u_j+s)\wedge v_j}(t)\,,\quad 0\leq t\leq \wt\zeta^{(j)}_s:=\zeta_{(u_j+s)\wedge v_j}.$$
Since we are only interested in values $s\in[0,v_j-u_j]$ such that $p_\zeta(u_j+s)$ belongs to
(the closure of) $\cc^{(j)}$, we introduce the time change
$$\pi^{(j)}_s:=\inf\{u\geq 0:\int_0^u \mathrm{d}t\,\mathbf{1}_{\{\tau^*_0(\wt W^{(j)}_t)>\wt\zeta^{(j)}_t\}}>s\},$$
with the notation $\tau_0^*(\w)=\inf\{t>0: \w(t)=0\}$ for $\w\in\mathcal{W}$. 
The effect of this time change is to eliminate the paths $W^{(j)}_s$ that return to $0$ and then survive for a positive amount of time
(in fact, the closure of $\cc^{(j)}$ is equal to $\{p_\zeta((u_j+\pi^{(j)}_s)\wedge v_j):s\geq 0\}$).
We define another snake trajectory by setting
$$W^{(j)}_s=\wt W^{(j)}_{\pi^{(j)}_s},$$
for every $s\geq 0$. We call $W^{(1)},W^{(2)},\ldots$ the Brownian snake excursions above $0$
--- note that we can similarly define the excursions below $0$ corresponding to the connected components of 
$\{a\in\t_\zeta:Z_a<0\}$. 

It follows from \cite[Theorem 4]{ALG} (together with Corollary \ref{bdry-disk} (i)) that, $\N_0$ a.e. for every $j=1,2,\ldots$,
the boundary size $|\partial\cc^{(j)}|$ can be defined by
\begin{equation}
\label{approbdry}
|\partial\cc^{(j)}|:=\lim_{\ve \to 0} \frac{1}{\ve^2}\mathrm{Vol}(\{a\in \cc^{(j)}:Z_a<\ve\})
= \lim_{\ve \to 0} \frac{1}{\ve^2}\int_0^{\sigma(W^{(j)})} \mathrm{d}s\,\mathbf{1}_{\{\wh W^{(j)}_s <\ve\}},
\end{equation}
where we write $\mathrm{Vol}$ for the
volume measure on $\t_\zeta$.

The quantities $|\partial\cc^{(1)}|,|\partial\cc^{(2)}|,\ldots$ are distinct $\N_0$ a.e. (they indeed correspond to jumps of
a continuous-state branching process with stable branching mechanism) and therefore we may and will assume that the components 
$\cc^{(1)},\cc^{(2)},\ldots$ (and consequently the excursions $W^{(1)},W^{(2)},\ldots$) have been ranked so that $|\partial\cc^{(1)}|>|\partial\cc^{(2)}|>\cdots$.
Furthermore, \cite[Theorem 4]{ALG} implies that, under $\N_0$ and conditionally on the sequence
of boundary sizes
$(|\partial\cc^{(1)}|,|\partial\cc^{(2)}|,\ldots)$, the excursions $W^{(1)},W^{(2)},\ldots$
are independent and distributed respectively according to 
$$\N_0^{*,|\partial\cc^{(1)}|},
\N_0^{*,|\partial\cc^{(2)}|},\ldots$$ and are also independent of the excursions below $0$.
 

It now follows from \eqref{approbdry} that property (a) stated above holds when $r=-W_*$,
and $|\partial \mathbf{C}^j|=|\partial \cc^{(j)}|$
(notice that  $D(x,\partial\mathbf{C}^{j})=Z_a$ if $x=\Pi(a)$ and $a\in \cc^{(j)}$).
By a preceding remark, property (b) also holds. As for (c), 
a straightforward adaptation of the arguments of the proof of Theorem \ref{Bnet} allows us
to verify that, for every $j=1,2,\ldots$, the intrinsic metric on $\mathbf{C}^j$
can be extended continuously to a metric on the closure $\ov{\mathbf{C}}^j$, and 
the resulting metric space equipped with the restriction of the volume measure $\bv$
coincides (as an element of  $\M$) with the Brownian disk associated with the excursion $W^{(j)}$ ---
this makes sense since we know
that, conditionally on $|\partial\cc^{(j)}|$,
$W^{(j)}$ is distributed according to $\N_0^{*,|\partial\cc^{(j)}|}$. 
In what follows, we simply write $\ov{\mathbf{C}}^j$ for this random measure metric space. 

So properties (a),(b),(c) hold if $r$ is replaced by the random quantity $-W_*$.
By a re-rooting invariance argument (see \cite[Theorem 2.3]{LGW}), we then 
obtain that these properties hold for $r=\wh W_t-W_*$, for a.a. $t\in(0,\sigma)$, $\N_0$ a.e.
It follows that they also hold for a.a. $r\in(0,W^*-W_*)$, $\N_0$ a.e. From a Fubini argument, one
can pick one value of $r_0>0$ such that (a),(b),(c) hold for $r=r_0$, $\N_0$ a.e. 
on the event $\{r_0<W^*-W_*\}$. But then the scaling property of $\N_0$ shows that
(a),(b),(c) hold for every $r>0$ (on the event $\{r<W^*-W_*\}$). This completes Step 2.

\medskip
\noindent{\it Step 3.} Let $s>0$ and $r>0$. The properties (a),(b),(c) stated in Step 2 also hold 
$\N_0^{(s)}$ a.e. on the event $\{r<W^*-W_*\}$.

The fact that this is  true for a.a. $s>0$ follows from Step 2, but it does not seem
easy to get the result for {\it every} $s>0$. However, we may use the following absolute continuity 
argument that will also be useful later. We fix $\eta>0$, and, on the 
event $\{\zeta_{\sigma/2}>\eta\}$, we define
$$R_{(\eta)}:=\sup\{t\in(0,\sigma/2): \zeta_t=\eta\}\;,\quad S_{(\eta)}:=\inf\{t\in(\sigma/2,\sigma):\zeta_t=\eta\}.$$
Still on the 
event $\{\zeta_{\sigma/2}>\eta\}$, we consider the snake trajectory $(W^\eta_u)_{u\geq 0}$ defined by
$$W^\eta_u(t):=W_{(R_{(\eta)}+u)\wedge S_{(\eta)}}(\eta+ t) -\wh W_{R_{(\eta)}}\,,\quad 0\leq t\leq \zeta^\eta_u:=\zeta_{(R_{(\eta)}+u)\wedge S_{(\eta)}} -\eta,$$
and we also set 
$$m^\eta:=\min\{\wh W_u: u\in[0,R_{(\eta)}]\cup [S_{(\eta)},1]\}\,\quad 
M^\eta:=\max\{\wh W_u: u\in[0,R_{(\eta)}]\cup [S_{(\eta)},1]\}.$$
Roughly speaking, the definition of $W^\eta$ means that we remove a part of the genealogical tree near the root and shift the labels so that 
the label of the new
root is again zero. It is not hard to prove (see the proof of Proposition 10 in \cite{Cactus2}) that the distribution 
of $W^\eta$ under $\N^{(s)}_0(\cdot \cap \{\zeta_{\sigma/2}>\eta\})$ is absolutely continuous with respect to $\N_0$. 

Recall the notation $\ll^\bullet$ in Lemma \ref{measurability}.
On the event $\{\zeta_{\sigma/2}>\eta\}\cap \{W_*<m^\eta<M^\eta<W_*+r\}$,
each connected component of the complement of the ball $B(x_*,r)$ in $\ll^\bullet(W)$
is identified  with a corresponding connected component of the complement of
the same ball in $\ll^\bullet(W^\eta)$. Furthermore this identification is isometric if either
connected component is equipped with its intrinsic metric. We can thus use the preceding absolute continuity property and Step 2 to see that properties (a),(b),(c)
hold $\N^{(s)}_0$ a.e. on the event $\{W_*<m^\eta<M^\eta<W_*+r\}\cap\{\zeta_{s/2}>\eta\}\cap \{r<W^*-W_*\}$. 

Clearly, if $\eta$ is small, $\N^{(s)}_0(\zeta_{s/2}>\eta)$ is close to $1$. This is not the case for $\N^{(s)}_0(W_*<m^\eta<M^\eta<W_*+r)$,
but we may use the re-rooting invariance property \eqref{uniform-re-root}, recalling that the space $\ll^\bullet(W)$
does not change if we replace $W$ by $W^{[t]}$ for some $t\in[0,\sigma]$.
 If $\eta>0$ is small enough then, with $\N^{(s)}_0$-probability close to one, we will be able to find a rational
 $t\in[0,s]$ such that the property $W_*<m^\eta<M^\eta<W_*+r$ holds when $W$
 is replaced by $W^{[t]}$. This completes Step 3.
 
 \medskip
 \noindent{\it Step 4.} Thanks to Step 3, we can define, $\N_0^{(s)}$ a.e. on the event $\{r<W^*-W_*\}$, the sequence 
 $\mathbf{C}^{r,j}$, $j=1,2,\ldots$ of connected components of $\mm\backslash B(x_*,r)$ and their boundary sizes. 
This sequence is ordered
 in such a way that $|\partial \mathbf{C}^{r,1}|>|\partial \mathbf{C}^{r,2}|>\cdots$. Also for every $j=1,2,\ldots$, we know 
 that the intrinsic metric $d^{r,j}_{\rm intr}$ on $\mathbf{C}^{r,j}$ has a continuous extension to a metric on $\ov{\mathbf{C}}^{r,j}$.
 For the sake of simplicity, in the following lines, the notation $\ov{\mathbf{C}}^{r,j}$ will be used
 to represent the measure metric space $(\ov{\mathbf{C}}^{r,j},d^{r,j}_{\rm intr},\bv^{r,j})$.
 
 We now aim at establishing the last part of Theorem \ref{ccBm}.
Let us come back to the excursions $W^{(1)},W^{(2)},\ldots$ above $0$ under $\N_0$.
By preceding observations, under $\N_0$ and conditionally on the sequence
$(|\partial\cc^{(1)}|,|\partial\cc^{(2)}|,\ldots)$, the excursions $W^{(1)},W^{(2)},\ldots$ are 
independent of the triple
$$\Big(W_*, \sigma_-,L^{(0)}\Big),$$
because $W_*$ is just the minimal value attained in all excursions below $0$,
$\sigma_-$ is the sum of the sizes of all the latter excursions
and 
the local time $L^{(0)}$ is a measurable function of the collection of 
excursions below $0$, since
$$L^{(0)}=\lim_{\ve\to 0} \frac{1}{\ve} \int_0^\sigma \mathrm{d}t\,\mathbf{1}_{\{-\ve\leq\wh W_t<0\}}.$$

Let $n\geq1$ and let $G$ and $H$ be nonnegative measurable
functions defined respectively on $\R^2$ and on $\R\times \R^\N$. Also let
$\Phi_1,\ldots,\Phi_n$ be nonnegative measurable
functions on $\M$. Recall that, for every $j=1,2,\ldots$, 
$\ov{\mathbf{C}}^j$ is the Brownian disk associated with $W^{(j)}$. Using the fact that, conditionally on the sequence
of boundary sizes
$(|\partial\cc^{(1)}|,|\partial\cc^{(2)}|,\ldots)$, the excursions $W^{(1)},W^{(2)},\ldots$
are independent and distributed respectively according to 
$\N_0^{*,|\partial\cc^{(1)}|},
\N_0^{*,|\partial\cc^{(2)}|},\ldots$, together with Theorem \ref{main} and the observation of the preceding paragraph, we have
\begin{align*}
&\N_0\Big(G(W_*,L^{(0)})\,H\Big(\sigma_-,(|\partial\mathbf{C}^{i}|)_{i\geq 1}\Big)\;
\prod_{i=1}^n \Phi_i(\ov{\mathbf{C}}^i)\Big)\\
&\qquad=\N_0\Big(G(W_*,L^{(0)})\,H\Big(\sigma_-,(|\partial\mathbf{C}^{i}|)_{i\geq 1}\Big)\;
\prod_{i=1}^n \F_{|\partial\mathbf{C}^{i}|}(\Phi_i)\Big).
\end{align*}

We slightly modify the last formula by conditioning also on the volumes
$\bv(\mathbf{C}^{1}),\bv(\mathbf{C}^{(2)}),\ldots$
Recalling the notation $\F_{z,v}$ for the distribution of the
Brownian disk with perimeter $z$ and volume $v$, we get from the previous display and standard arguments that
\begin{align*}
&\N_0\Big(G(W_*,L^{(0)})\,H\Big(\sigma_-,(|\partial\mathbf{C}^{i}|,\bv(\mathbf{C}^{i}))_{i\geq 1}\Big)\;
\prod_{i=1}^n\Phi_i(\ov{\mathbf{C}}^i)\Big)\\
&\quad=\N_0\Big(G(W_*,L^{(0)})\,H\Big(\sigma_-,(|\partial\mathbf{C}^{i}|,\bv(\mathbf{C}^{i}))_{i\geq 1}\Big)\;
\prod_{i=1}^n\F_{|\partial\mathbf{C}^{i}|,\bv(\mathbf{C}^{i})}(\Phi_i)\Big),
\end{align*}
where the function $H$ is now nonnegative and measurable on $\R\times (\R^2)^\N$. Noting that
$$\sigma =\sigma_- + \sum_{j=1}^\infty \bv(\mathbf{C}^{j}),$$
we get from the preceding display that, for any nonnegative measurable functions $G$, $g$ and $H$
defined respectively on $\R^2$, on $\R$ and on $(\R^2)^\N$, we have
\begin{align*}
&\N_0\Big(G(W_*,L^{(0)})\,g(\sigma) \,H\Big((|\partial\mathbf{C}^{i}|,\bv(\mathbf{C}^{i}))_{i\geq 1}\Big)\;
\prod_{i=1}^n\Phi_i(\ov{\mathbf{C}}^i)\Big)\\
&\quad=\N_0\Big(G(W_*,L^{(0)})\,g(\sigma)\, H\Big((|\partial\mathbf{C}^{i}|,\bv(\mathbf{C}^{i}))_{i\geq 1}\Big)\;
\prod_{i=1}^n\F_{|\partial\mathbf{C}^{i}|,\bv(\mathbf{C}^{i})}(\Phi_i)\Big),
\end{align*}

The last formula still holds if we replace $\N_0$ by $\N_0^{(s)}$ for any $s>0$. Indeed, using \eqref{desintN} and the fact that
$g$ can be any nonnegative measurable function, we get that the formula must hold under $\N_0^{(s)}$ simultaneously for 
any choice of $G,g,H,\Phi_1,\ldots\Phi_n$,
for a.a.~$s$, but then a scaling argument shows that it holds for all $s>0$. 

We next apply the formula of the last display (with $\N_0$ replaced by $\N^{(s)}_0$) with $g=1$ and 
$G(W_*,L^{(0)})=h(-W_*)/L^{(0)}$, where $h$ is a nonnegative measurable function 
on $\R_+$. It follows that
\begin{align}
\label{BdBm-tech}
&\N^{(s)}_0\Big(\frac{h(-W_*)}{L^{(0)}}\, H\Big((|\partial\mathbf{C}^{i}|,\bv(\mathbf{C}^{i}))_{i\geq 1}\Big)\;
\prod_{i=1}^n\Phi_i(\ov{\mathbf{C}}^i)\Big)\nonumber\\
&\quad=\N^{(s)}_0\Big(\frac{h(-W_*)}{L^{(0)}}\, H\Big((|\partial\mathbf{C}^{i}|,\bv(\mathbf{C}^{i}))_{i\geq 1}\Big)\;
\prod_{i=1}^n\F_{|\partial\mathbf{C}^{i}|,\bv(\mathbf{C}^{i})}(\Phi_i)\Big).
\end{align}

We rewrite the left-hand side of \eqref{BdBm-tech} by using 
the invariance of $\N^{(s)}_0$ under re-rooting. Note that, if $W$ is replaced by $W^{[t]}$, for some fixed
$t\in[0,s]$, then $W_*$ is replaced by $W_* - \wh W_t$, $L^{(0)}$ is replaced 
by $L^{(\wh W_t)}$, and $\mathbf{C}^{1}, \mathbf{C}^{2},\ldots$ are replaced 
by the connected components $\mathbf{C}^{\wh W_t-W_*,1}, \mathbf{C}^{\wh W_t-W_*,2},\ldots$
of $\mm\backslash B(x_*, \wh W_t-W_*)$. Hence, we get that the left-hand side of \eqref{BdBm-tech} is equal to
\begin{align*}
&\frac{1}{s}\,\N^{(s)}_0\Big(\int_0^s \mathrm{d}t\,\frac{h(\wh W_t-W_*)}{L^{(\wh W_t)}}\,
H\Big((|\partial\mathbf{C}^{\wh W_t-W_*,i}|,\bv(\mathbf{C}^{\wh W_t-W_*,i}))_{i\geq 1}\Big)\;
\prod_{i=1}^n \Phi_i(\ov{\mathbf{C}}^{\wh W_t-W_*,i})\Big)\\
&\quad=\frac{1}{s}\,\N^{(s)}_0\Big(\int_{W_*}^{W^*} \mathrm{d}u\,h(u-W_*)\,
H\Big((|\partial\mathbf{C}^{u-W_*,i}|,\bv(\mathbf{C}^{u-W_*,i}))_{i\geq 1}\Big)\;
\prod_{i=1}^n \Phi_i(\ov{\mathbf{C}}^{u-W_*,i})\Big)\\
&\quad =\frac{1}{s}\,\N^{(s)}_0\Big(\int_{0}^{W^*-W_*} \mathrm{d}r\,h(r)\,
H\Big((|\partial\mathbf{C}^{r,i}|,\bv(\mathbf{C}^{r,i}))_{i\geq 1}\Big)\;
\prod_{i=1}^n \Phi_i(\ov{\mathbf{C}}^{r,i})\Big)\\
&\quad =\frac{1}{s}\,\int_0^\infty \mathrm{d}r\,h(r)\,\N^{(s)}_0\Big(\mathbf{1}_{\{r<W^*-W_*\}}\,
H\Big((|\partial\mathbf{C}^{r,i}|,\bv(\mathbf{C}^{r,i}))_{i\geq 1}\Big)\;
\prod_{i=1}^n \Phi_i(\ov{\mathbf{C}}^{r,i})\Big).
\end{align*}
In the first equality, we used the definition of the local times. 

Note that, in the previous lines we left aside certain measurability issues, concerning in particular
the measurability of the mapping $r\mapsto (\ov{\mathbf{C}}^{r,1},\ldots,\ov{\mathbf{C}}^{r,n})$.
To deal with these issues, one may observe that the connected components of
$\{a\in \t_\zeta:Z_a>W_*+r\}$ can be represented by Brownian snake excursions above $0$ (in the case
$r=-W_*$ these are the excursions $W^{(1)},W^{(2)},\ldots$ introduced above) and the collection of these excursions depends on $r$ in a measurable way. 
Furthermore, the measure metric spaces $\ov{\mathbf{C}}^{r,1},\ov{\mathbf{C}}^{r,2},\ldots$ are obtained as measurable
functions of these excursions (this measurability property can be verified by an adaptation of the arguments of Section
\ref{sec:consmms}, but we omit the details). 

By the same manipulations, we obtain that the right-hand side of \eqref{BdBm-tech} is equal to
$$\frac{1}{s}\,\int_0^\infty \mathrm{d}r\,h(r)\,\N^{(s)}_0\Big(\mathbf{1}_{\{r<W^*-W_*\}}\,
H\Big((|\partial\mathbf{C}^{r,i}|,\bv(\mathbf{C}^{r,i}))_{i\geq 1}\Big)\;
\prod_{i=1}^n \F_{|\partial\mathbf{C}^{r,i}|,\bv(\mathbf{C}^{r,i})}(\Phi_i)\Big).$$

Since the function $h$ was arbitrary, the equality of the quantities in the
last two displays means that for a.a. $r>0$ we have
\begin{align}
\label{ccBmtec1}
&\N^{(s)}_0\Big(\mathbf{1}_{\{r<W^*-W_*\}}\,
H\Big((|\partial\mathbf{C}^{r,i}|,\bv(\mathbf{C}^{r,i}))_{i\geq 1}\Big)\;\prod_{i=1}^n \Phi_i(\ov{\mathbf{C}}^{r,i})\Big)\nonumber\\
&\quad=\N^{(s)}_0\Big(\mathbf{1}_{\{r<W^*-W_*\}}\,
H\Big((|\partial\mathbf{C}^{r,i}|,\bv(\mathbf{C}^{r,i}))_{i\geq 1}\Big)\;
\prod_{i=1}^n \F_{|\partial\mathbf{C}^{r,i}|,\bv(\mathbf{C}^{r,i})}(\Phi_i)\Big).
\end{align}

\noindent{\it Step 5.} 
We now prove that \eqref{ccBmtec1} holds for {\it every} $r>0$, which 
will complete the proof of Theorem \ref{ccBm}.
We may assume
that $H,\Phi_1,\ldots,\Phi_n$ are bounded and continuous. It is then enough 
to prove that both sides of \eqref{ccBmtec1} are continuous functions of $r$ for a fixed value of $s>0$. 
By a simple scaling argument, this is equivalent to showing continuity in $s$ for a fixed value of $r>0$.

It turns out to be easier to deal first with the quantities
\begin{equation}
\label{ccBmtec2}
\N^{(s)}_0\Big(\mathbf{1}_{\{-W_*<r\}}\,\mathbf{1}_{\{r<W^*-W_*\}}\,
H\Big((|\partial\mathbf{C}^{r,i}|,\bv(\mathbf{C}^{r,i}))_{i\geq 1}\Big)\;\prod_{i=1}^n \Phi_i(\ov{\mathbf{C}}^{r,i})\Big)
\end{equation}
and 
\begin{equation}
\label{ccBmtec3}
\N^{(s)}_0\Big(\mathbf{1}_{\{-W_*<r\}}\,\mathbf{1}_{\{r<W^*-W_*\}}\,
H\Big((|\partial\mathbf{C}^{r,i}|,\bv(\mathbf{C}^{r,i}))_{i\geq 1}\Big)\;
\prod_{i=1}^n \F_{|\partial\mathbf{C}^{r,i}|,\bv(\mathbf{C}^{r,i})}(\Phi_i)\Big).
\end{equation}

So let us fix $r>0$ and discuss continuity in the variable $s$ of the quantities \eqref{ccBmtec2} and \eqref{ccBmtec3}. We rely on the same absolute continuity argument as in Step 3 above. Recall 
the notation $W^\eta, m^\eta, M^\eta$ introduced in Step 3 (these random quantities are defined 
under $\N^{(s)}_0$ on the event $\{\zeta_{s/2}>\eta\}$). We also note that,
as $s'\to s$, the distribution of $W^\eta$ under $\N^{(s')}_0(\cdot \cap \{\zeta_{s'/2}>\eta\})$ converges to the distribution
of $W^\eta$ under $\N^{(s)}_0(\cdot \cap \{\zeta_{s/2}>\eta\})$
in variation distance (recalling our notation $\zeta^\eta_u=\zeta_{(W^\eta_u)}$, it suffices to prove the analogous
convergence for the distribution of $\zeta^\eta$, and one can use the explicit formulas available under the Brownian excursion measure --- we
leave the details as an exercise for the reader). 

Let us argue on $\{-W_*<r\}$ but discard the event  $\{\zeta_{s/2}\leq \eta\}\cup \{-W_*< r \leq M^\eta-W_*\}\cup \{W_*=m^\eta\}$, which has small $\N^{(s)}_0$-probability 
when $\eta$ is small, uniformly when $s$ varies over a compact subset of $(0,\infty)$. 
We already noticed in Step 3 that, on the complement of the latter event in $\{-W_*<r\}$,
the connected components of the complement of the ball $B(x_*,r)$ in $\ll^\bullet(W)$
are identified isometrically with the connected components of the complement of
the same ball in $\ll^\bullet(W^\eta)$. It follows that the quantity 
$$\mathbf{1}_{\{-W_*<r\}}\,\mathbf{1}_{\{r<W^*-W_*\}}\,H\Big((|\partial\mathbf{C}^{r,i}|,\bv(\mathbf{C}^{r,i}))_{1\leq i\leq n}\Big)\;\prod_{i=1}^n \Phi_i(\ov{\mathbf{C}}^{r,i})$$
coincides with a measurable function of $W^\eta$, except on a set of $\N^{(s)}_0$-probability small when $\eta\to 0$ (uniformly when $s$ varies in 
a compact subset of $\R_+$). 
Since we know that (for a fixed value of $\eta$), the distribution of $W^\eta$ under $\N^{(s')}_0(\cdot \cap \{\zeta_{\sigma/2}>\eta\})$
converges in variation distance to the distribution of $W^\eta$ under $\N^{(s)}_0(\cdot \cap \{\zeta_{\sigma/2}>\eta\})$ when $s\to s'$, this suffices to
give the desired continuity in $s$ of the quantity \eqref{ccBmtec2}, and the same argument applies to the
quantity \eqref{ccBmtec3}. 

We then use a re-rooting invariance argument, recalling the remark at the end 
of Section \ref{sec:consmms}. If $W$ is replaced by $W^{[t]}$ for some $t\in[0,s]$ the quantity inside
the expectation in \eqref{ccBmtec2} does not change, except that $\mathbf{1}_{\{-W_*<r\}}$ is replaced by $\mathbf{1}_{\{\wh W_t-W_*<r\}}$. Noting that
$$\int_0^s \mathbf{1}_{\{\wh W_t-W_*<r\}}\,\mathrm{d}t= \bv(B(x_*,r))\,,\quad \N^{(s)}_0\hbox{ a.s.},$$
we see that the quantity \eqref{ccBmtec2} is also equal to 
$$\N^{(s)}_0\Big(\mathbf{1}_{\{r<W^*-W_*\}}\,\frac{\bv(B(x_*,r))}{s}\,
H\Big((|\partial\mathbf{C}^{r,i}|,\bv(\mathbf{C}^{r,i}))_{i\geq 1}\Big)\;\prod_{i=1}^n \Phi_i(\ov{\mathbf{C}}^{r,i})\Big)$$
Next since 
$$\bv(B(x_*,r))= s - \sum_{i=1}^\infty \bv(\mathbf{C}^{r,i})\,,\quad \N^{(s)}_0\hbox{ a.s.},$$
by changing the function $H$ in an appropriate way, we derive that, for every $\delta>0$, the
quantity
$$\N^{(s)}_0\Big(\mathbf{1}_{\{r<W^*-W_*\}}\,\frac{\bv(B(x_*,r))}{\bv(B(x_*,r))+\delta}\,
H\Big((|\partial\mathbf{C}^{r,i}|,\bv(\mathbf{C}^{r,i}))_{i\geq 1}\Big)\;\prod_{i=1}^n \Phi_i(\ov{\mathbf{C}}^{r,i})\Big)$$
depends continuously on $s$. When $\delta\to 0$, the quantity in the last display converges to the left-hand side of \eqref{ccBmtec1}
uniformly when $s$ varies in a compact subset of $(0,\infty)$, so that we also obtain that the left-hand side of \eqref{ccBmtec1}
depends continuously on $s$. The same argument applies to the right-hand side, and this completes the proof
of Theorem \ref{ccBm}.

\section*{Appendix A. Proof of Proposition \ref{approx-LT}.}

We will use properties of the 
Brownian snake that can be found in \cite{Zurich}. We fix $r>0$ and
view the Brownian snake as a strong Markov process with values in
$\mathcal{W}_r$, which is 
 killed at the first time $\sigma$ when its lifetime hits $0$ (the measure $\N_r$ is then viewed as the excursion measure
of this Markov process away from the trivial path consisting
only of the point $r$).
We keep the notation $(W_s)_{s\geq 0}$ for this Markov process, and 
$(\zeta_s)_{s\geq 0}$ for the associated lifetime process
 --- we slightly abuse notation since in the previous
sections $(W_s)_{s\geq 0}$ stood for the canonical process on 
the space of snake trajectories.   For $\w\in\mathcal{W}_r$, we write 
$\P_\w$ for the probability measure under which $W_0=\w$. 
As in Section \ref{sec:cons-metric}, we write $(B_t)_{t\geq 0}$ for  a standard linear Brownian motion
 that starts at $x$ under the probability measure $\mathbf{P}_x$, and 
 $\tau_y:=\inf\{t\geq 0:B_t=y\}$
for every $y\in \R$. 

The definition of the exit local time $(\ell^0_s)_{s\geq 0}$ makes sense under
each probability measure $\P_\w$, $\w\in \mathcal{W}_r$. We recall the
first moment formulas
\begin{equation}
\label{1-mom}
\N_r\Big(\int_0^\sigma \mathrm{d}s\,F(W_s)\Big) = \int_0^\infty\mathrm{d}u\,\mathbf{E}_r[F((B_t)_{0\leq t\leq u})],
\end{equation}
and 
\begin{equation}
\label{1-mom-exit}
\N_r\Big(\int_0^\sigma \mathrm{d}\ell^0_s\,F(W_s)\Big) = \mathbf{E}_r[F((B_t)_{0\leq t\leq \tau_0})],
\end{equation}
which hold for any nonnegative measurable function $F$ on $\mathcal{W}_r$ (see Propositions IV.2 and V.3
in \cite{Zurich}). 

We fix $K>1$ and, for every $\ve>0$, we set
$$A^{\ve,K}_t:=\frac{1}{\ve^2}\int_0^t \mathrm{d}s\,\mathbf{1}_{\{\tau_0(W_s)\geq \zeta_s,
\tau_K(W_s)=\infty,\wh W_s<\ve,
\}},$$
for every $t\geq 0$. We also define
$$\ell^{0,K}_t:=\int_0^t \mathrm{d}\ell^0_s \,\mathbf{1}_{\{\tau_K(W_s)=\infty\}}.$$
Then both $A^{\ve,K}_t$ and $\ell^{0,K}_t$ are additive functionals of the Brownian snake. The potential
of these additive functionals is easy to compute using \cite[Lemma V.5]{Zurich}
and formulas \eqref{1-mom} and \eqref{1-mom-exit}. For every
$\w\in\mathcal{W}_r$,
\begin{align}
h(\w)&:= \E_\w[\ell^{0,K}_\sigma]=2\int_0^{\zeta_{(\w)}\wedge \tau_0(\w)\wedge \tau_K(\w)} \mathrm{d}t\,\psi(\w(t)),\label{pot-1}\\
h_\ve(\w)&:= \E_\w[A^{\ve,K}_\sigma]= 2 \int_0^{\zeta_{(\w)}\wedge \tau_0(\w)\wedge \tau_K(\w)}
\mathrm{d}t\,\varphi_\ve(\w(t)),\label{pot-2}
\end{align}
where, for every $x\in(0,K)$,
$$\psi(x):=\mathbf{P}_{x}(\tau_0<\tau_K)= \frac{K-x}{K},$$
and
$$\varphi_\ve(x):= 
\mathbf{E}_x\Big[\frac{1}{\ve^2}\int_0^{\tau_0\wedge \tau_K} \mathrm{d}s\,\mathbf{1}_{\{B_s<\ve\}}\Big]=
\frac{1}{\ve^2}\int_{(0,\ve)} \mathrm{d}y\,G(x,y),$$
where we write $G$ for the Green function of Brownian motion killed at $\tau_0\wedge \tau_K$. From the
explicit formula $G(x,y)=2\,K^{-1}(x\wedge y)(K-(x\vee y))$, we readily get that we also
$\varphi_\ve(x)= \psi(x)$ if $x\geq \ve$, and $\varphi_\ve(x)\leq \psi(x)\leq 1$ for every $x\in(0,K)$.

Then,
\begin{align*}
\N_r((\ell^{0,K}_\sigma)^2)=2\N_r\Big(\int_0^\sigma \mathrm{d}\ell^{0,K}_s\,\E_{W_s}[\ell^{0,K}_\sigma]\Big)
&=4\,\mathbf{E}_r\Big[ \mathbf{1}_{\{\tau_0<\tau_K\}}\int_0^{\tau_0\wedge \tau_K} \mathrm{d}t\,\psi(B_t)\Big]\\
&=4\int_0^\infty \mathrm{d}t\,\mathbf{E}_r\Big[ \mathbf{1}_{\{t<\tau_0\wedge \tau_K\}}\,\psi(B_t)\,\mathbf{P}_{B_t}(\tau_0<\tau_k)\Big]\\
&= 4\,\mathbf{E}_r\Big[\int_0^{\tau_0\wedge \tau_K} \mathrm{d}t\,\psi(B_t)^2\Big] 
\end{align*}
using \eqref{pot-1} and \eqref{1-mom-exit} in the second equality. Similarly,
\begin{align*}
\N_r(\ell^{0,K}_\sigma\,A^{\ve,K}_\sigma)=2\N_r\Big(\int_0^\sigma \mathrm{d}\ell^{0,K}_s\,
\E_{W_s}[A^{\ve,K}_s]\Big)&= 4\,\mathbf{E}_r\Big[ \mathbf{1}_{\{\tau_0<\tau_K\}}\int_0^{\tau_0\wedge \tau_K} \mathrm{d}t\,\varphi_\ve(B_t)\Big]\\
&= 4\,\mathbf{E}_r\Big[ \int_0^{\tau_0\wedge \tau_K} \mathrm{d}t\,\psi(B_t)\varphi_\ve(B_t)\Big],
\end{align*}
using \eqref{pot-2} and \eqref{1-mom-exit}.
Furthermore, using \eqref{pot-2} and \eqref{1-mom}, it is also
easy to verify that
$$\N_r((A^{\ve,K}_\sigma)^2)
=4\,\mathbf{E}_r\Big[ \int_0^{\tau_0\wedge \tau_K} \mathrm{d}t\,\varphi_\ve(B_t)^2\Big],$$
in such a way that we obtain
$$\N_r((\ell^{0,K}_\sigma-A^{\ve,K}_\sigma)^2)
=4\,\mathbf{E}_r\Big[ \int_0^{\tau_0\wedge \tau_K} \mathrm{d}t\,(\psi(B_t)-\varphi_\ve(B_t))^2\Big].$$
Take $\ve\in(0,r)$. Since $\psi(x)=\varphi_\ve(x)$ when $x\in(\ve,K)$, we get
$$\N_r((\ell^{0,K}_\sigma-A^{\ve,K}_\sigma)^2)
\leq 4\,\mathbf{E}_r\Big[ \int_0^{\tau_0\wedge \tau_K} \mathrm{d}t\,\mathbf{1}_{\{B_t<\ve\}}\Big]\leq 4\ve^2.$$

Let us fix $\delta>0$.
We note that both $(A^{\ve,K}_{t\wedge \sigma}+h_\ve(W_{t\wedge\sigma}))_{t\geq \delta}$
and $(\ell^{0,K}_{t\wedge \sigma}+h(W_{t\wedge\sigma}))_{t\geq \delta}$ are uniformly integrable martingales 
under $\N_0(\cdot\mid \sigma\geq \delta)$. Since the terminal 
values of these martingales are $A^{\ve,K}_\sigma$ and $\ell^{0,K}_\sigma$ respectively, we deduce from
the last display and Doob's maximal inequality for martingales that
$$\N_r\Big(\sup_{t\geq \delta} \Big(\ell^{0,K}_{t\wedge \sigma}-A^{\ve,K}_{t\wedge \sigma}
+ h(W_{t\wedge\sigma})-h_\ve(W_{t\wedge\sigma})\Big)^2\,\Big|\, \sigma\geq\delta\Big)
\leq 4\,\N_r((\ell^{0,K}_\sigma-A^{\ve,K}_\sigma)^2\mid \sigma \geq\delta) \leq c_\delta\,\ve^2,$$
where $c_\delta>0$ is a constant depending on $\delta$. Set $\ve_n=n^{-1}$
for every $n\geq 1$. It follows from the preceding display that
$$\lim_{n\to\infty} \sup_{t\geq \delta}\Big|\ell^{0,K}_{t\wedge \sigma}-A^{\ve_n,K}_{t\wedge \sigma}
+ h(W_{t\wedge\sigma})-h_{\ve_n}(W_{t\wedge\sigma})\Big| = 0,\quad\N_r\hbox{ a.e.}$$
Since $\delta$ is arbitrary, and since  $h_{\ve_n}(\w)\la h(\w)$ as $n\to\infty$, for any $\w\in \mathcal{W}_r$, it follows that
$$\lim_{n\to\infty} A^{\ve_n,K}_{t\wedge \sigma}= \ell^{0,K}_{t\wedge \sigma}
\quad\hbox{ for every }t\geq 0,\ \N_r\hbox{ a.e.}$$
We apply this result to a sequence of values of $K$ tending
to $+\infty$, noting that we have $ A^{\ve_n,K}_t=\ve_n^{-2}\int_0^t \mathrm{d}s\,\mathbf{1}_{\{\tau_0(W_s)\geq \zeta_s, \wh W_s<\ve_n\}}$ and $\ell^{0,K}_t=\ell^0_t$ for every $t\geq 0$ and every $n$,
except on a set of $\N_r$-measure tending to $0$ as $K\to\infty$. This
gives the result of the proposition, except that we restricted ourselves to the
sequence $(\ve_n)_{n\geq 1}$. The proof is however completed by
a straightforward monotonicity argument. \hfill$\square$

\section*{Appendix B. Proof of Proposition \ref{conv-peri-vol}.}
We write $\mathrm{q}^\uparrow_k$ for the quadrangulation with a boundary obtained at the 
$k$-th step of the (lazy) peeling algorithm of the UIPQ. Recall that $S^\uparrow_k$
is the half-perimeter of $\mathrm{q}^\uparrow_k$ and write $M^\uparrow_k$ for
the volume (number of faces) of $\mathrm{q}^\uparrow_k$. We have then the
following convergence in distribution in the Skorokhod sense
\begin{equation}
\label{App-tec1}
\Big(n^{-2/3} S^\uparrow_{\lfloor nt\rfloor},n^{-4/3} M^\uparrow_{\lfloor nt\rfloor}\Big)_{t\geq 0}
\build{\la}_{n\to\infty}^{\rm(d)} \Big( \Upsilon^\uparrow_{t/3}, \mathcal{V}^\uparrow_{t/3}\Big)_{t\geq 0},
\end{equation}
where the processes $\Upsilon^\uparrow$ and $\mathcal{V}^\uparrow$ are
as in Section \ref{sec:asymp-hull}. See \cite[Theorem 3]{Bud} or \cite[Section 5.3]{Peccot},
noting that these references deal with a more general setting, and that one needs
to compute the value of the relevant constants in the case of quadrangulations. At least for the first component,
the convergence \eqref{App-tec1} is basically a consequence of invariance principles for
random walks conditioned to stay positive \cite{CC}.

On the other hand, we have, for every $k\geq 1$,
\begin{equation}
\label{App-tec10}
H^{(\infty)}_k=S^\uparrow_{R_k}\;,\quad V^{(\infty)}_k= M^\uparrow_{R_k},
\end{equation}
where we abuse notation by still writing $R_k$ for the number of steps of the peeling
algorithm needed for the $k$ first layers of the UIPQ (in Section \ref{sec:peel-layers},
we used the same notation but for the peeling of a Boltzmann quadrangulation instead of the UIPQ, this
should however create no confusion here). In order to deduce Proposition \ref{conv-peri-vol}
from \eqref{App-tec1}, we thus need to control the time change $(R_k)$. To this end, the rough 
idea is that $R_{k+1}-R_k$ (the number of steps needed to discover the $k+1$-st layer)
is close to $3$ times the half-perimeter $S^\uparrow_{R_k}$ of the hull of radius $k$. The appearance of the
constant $3$ comes from the property \eqref{meanstep}. 

In order to make the preceding idea more precise, we 
introduce a probability measure $\P_\ell$ under which we run the peeling by layers algorithm of the UIPQ
starting from a boundary of length $2\ell$. Under $\P_\ell$, vertices of the initial boundary receive labels $0$ or $1$ 
alternatively. When running the peeling by layers algorithm, newly created vertices receive labels
equal to their graph distances from the set of vertices with label $0$. Note that the label of a vertex may change 
in the case of an event $(B_j)$ or $(B'_j)$. We still write $R_1$ for the number of steps needed to
complete the first layer, meaning that the boundary only contains vertices labeled $1$ or $2$, and 
$(S^\uparrow_k)_{k\geq 0}$ for the half-perimeter proces under $\P_\ell$
(in particular $\P_\ell(S^\uparrow_0=\ell)=1$).

Before we state our first lemma, we notice that, for any $1\leq \ell'\leq \ell$, we have
\begin{equation}
\label{mini-peri}
\P_\ell\Big(\min_{k\geq 0} S^\uparrow_k\leq \ell'\Big) \leq \frac{\ell'}{\ell}.
\end{equation}
This is a straightforward consequence of the identity \eqref{up-down}
relating the laws of $S^\uparrow$ and of $S^\downarrow$, or
rather of the extension of this identity to the case where $S^\uparrow$
and $S^\downarrow$ both start from $\ell$ (and then the constant $2$
has to be replaced by $\ell$). We leave the details to the reader.

\begin{lemma}
\label{App-est}
Let $\ve>0$. Then,
\begin{equation}
\label{App1}
\P_\ell(R_1 \geq (3-\ve)\ell) \build{\la}_{\ell\to \infty}^{} 1.
\end{equation}
Furthermore, we can find constants $A_\ve$, $C_\ve$ and $\rho_\ve>0$ such that,
for every $\ell\geq 1$,
\begin{equation}
\label{App2}
\P_\ell\Big(R_1 > (3+\ve)\ell ,\;\inf_{k\geq 0} S^\uparrow_k > A_\ve\Big) \leq C_\ve\,e^{-\rho_\ve\ell}.
\end{equation}
\end{lemma}

\proof We start by proving \eqref{App2}. It is convenient to introduce a sequence $(\xi_k)_{k\geq 1}$
of integer-valued random variables defined as follows. If, at the $k$-th step of the peeling algorithm,
the event $(B_j)$ occurs (for some $j\in\{0,1,\ldots\}$), we take $\xi_k=j+1$, and otherwise we take $\xi_k=0$. 
Writing $\mathcal{F}_k$ for the $\sigma$-field generated by the first $k$ steps of the algorithm, the conditional
distribution of $\xi_{k+1}$ knowing $\mathcal{F}_k$ is derived from the prescriptions in Section \ref{sec:peel}:
$$\P_\ell(\xi_{k+1}=j\mid \mathcal{F}_k)=\left\{
\begin{array}{ll}
\frac{1}{2} \,p^{(S^\uparrow_k)}_{-j}\quad&\hbox{if }j=1,2,\ldots,S^\uparrow_k\,,\\
1-\frac{1}{2} \,\sum_{i=1}^{S^\uparrow_k} p^{(S^\uparrow_k)}_{-i}\quad&\hbox{if } j=0,
\end{array}
\right.
$$
where the weights $p^{(L)}_{-j}$, for $L\geq 1$ and $1\leq j\leq L$, are defined in Section \ref{sec:peel}.

Let $N_k$ be the number of vertices labeled $0$ at the $k$-th step of the algorithm (and $N_0=\ell$), so that
$R_1=\inf\{k\geq 1:N_k=0\}$. From the definition of the algorithm, one easily gets that
$N_{k+1}\leq (N_k-\xi_{k+1})^+$, for every $k=0,1,\ldots,R_1-1$. It follows that
$N_k\leq (\ell -(\xi_1+\cdots+\xi_k))^+$ for $k=0,1,\ldots,R_1$, and therefore,
$$R_1\leq \min\{k\geq 1: \xi_1+\cdots+\xi_k\geq \ell\}.$$
For every integer $m\in\{1,\ldots,\ell\}$, set $\Gamma_m:=\min\{k\geq 0:S^\uparrow_k\leq m\}$. 
Using the explicit form of $h^\uparrow$ (and in particular the fact that $h^\uparrow(\ell+1)/h^\uparrow(\ell)$
is a decreasing function of $\ell$), we can, for every integer $k\geq1$, couple $\xi_1,\ldots,\xi_k$
with i.i.d. random variables $\xi^{(m)}_1,\ldots,\xi^{(m)}_k$ taking values in $\{0,1,\ldots\}$
and with common distribution determined by $\P(\xi^{(m)}_i=j)=p^{(m)}_{-j}$ for $j\geq 1$, 
in such a way that $\xi^{(m)}_i\leq \xi_i$ for every $i\in\{1,\ldots,k\}$ on the event $\{\Gamma_m>k\}$. 
Using also the last display, we get that
$$\Big(\{R_1 > k\}\cap \{\Gamma_m>k\}\Big) \subset \{\xi^{(m)}_1+\cdots+\xi^{(m)}_k < \ell\}.$$
Recalling \eqref{meanstep}, we can fix $m=A_\ve$ large enough so that the mean value
of the variables $\xi^{(A_\ve)}_i$ is (strictly) greater than $(3+\ve)^{-1}$. We have then 
$$\P_\ell\Big(R_1 > (3+\ve)\ell,\;\inf_{k\geq 0} S^\uparrow_k > A_\ve\Big) \leq
\P(\xi^{(A_\ve)}_1+\cdots+\xi^{(A_\ve)}_{\lfloor (3+\ve)\ell\rfloor} < \ell),$$
and standard large deviation results for sums of i.i.d. random variables give \eqref{App2}. 

Let us turn to the proof of \eqref{App2}. The argument is slightly more involved since 
the upper bound $N_{k+1}\leq (N_k-\xi_{k+1})^+$ does not help us to find 
a lower bound for $R_1$ (the point is the fact that events of the type $(B'_j)$
may also lead to a decrease of $N_k$). To begin with, we observe that we may
run the peeling by layers algorithm on a model with infinite boundary represented by $\Z$,
in such a way that even integers are labeled $1$ and odd integers are labeled $0$. We start 
the algorithm from the edge $(0,1)$ on the boundary, using the probabilities 
$p_1$ and $\frac{1}{2}p_{-j-1}$ for steps of type $(A)$ and $(B_j)$ (or $(B'_j)$)
respectively. We notice that a step of type $(A)$ occurs with 
probability $2/3$ and creates a vertex labeled $2$. An event of type $(B_j)$
does not change the number of vertices labeled $2$, whereas an event of
type $(B'_j)$ can decrease the number of vertices labeled $2$ by at most $j+1$. 
Recalling \eqref{meanstep}, we easily deduce from the law of large numbers 
that the number of vertices labeled $2$ converges to infinity a.s., 
which also means that after a certain (random) number of steps, the vertices
that lie to the left of vertices labeled $2$ are no longer affected by the algorithm.
The same conclusion holds a fortiori if instead of using the probability weights
$p_1,p_{-j-1}\; (j=0,1,\ldots)$ we use $p^{(L)}_1,p^{(L)}_{-j-1}\; (j=0,1,\ldots,L-1)$.
Indeed we just have to observe that $p_1\leq p^{(L)}_1$ and $p_{-j-1}\geq p^{(L)}_{-j-1}$.

Let us return to the case with finite boundary. By the first part of the proof
and \eqref{mini-peri}, we know that $\P_\ell(R_1\leq 4\ell)$ tends to $1$
as $\ell\to\infty$, and it easily follows that, for every $\delta>0$,
$$\P_\ell\Big(\sup_{k\leq R_1} |S^\uparrow_k-\ell| >\delta\ell\Big) \build{\la}_{\ell\to\infty}^{} 0.$$
A coupling argument relying on the law of large numbers also shows that
$\P_\ell(R_1\geq \sqrt{\ell})$ tends to $1$
as $\ell\to\infty$. Therefore, if we set
$$E_\ell:=\{R_1\geq \sqrt{\ell}\}\cap \Big\{\sup_{k\leq R_1} |S^\uparrow_k-\ell| \leq\ell/2\Big\},$$
we have $\P_\ell(E_\ell)\la 1$ as $\ell\to\infty$. A comparison argument with the case 
of the infinite boundary discussed above now shows that the following holds on the event $E_\ell$ 
except maybe on a set of probability vanishing when $\ell\to\infty$: for every $k$ such that 
$\lfloor\sqrt{\ell}\rfloor< k\leq R_1$, the $k$-th step
of the peeling algorithm affects
the number of vertices labeled $0$ only if it is of type $(B_j)$. In other words,
we can find an event $E'_\ell\subset E_\ell$ with $\P_\ell(E_\ell\backslash E'_\ell)\la 0$, such that
on $E'_\ell$ we have
$$N_{k+1}= (N_k -\xi_{k+1})^+$$
for every $\lfloor\sqrt{\ell}\rfloor\leq k< R_1$. In a way similar to the first part of the proof,
it follows that, still on the event $E'_\ell$,
$$R_1-\lfloor\sqrt{\ell}\rfloor =\min\{j\geq 1: \xi_{\lfloor\sqrt{\ell}\rfloor+1}+\cdots+ \xi_{\lfloor\sqrt{\ell}\rfloor+j} \geq N_{\lfloor\sqrt{\ell}\rfloor}\}.$$
Since it is also clear that $\P_\ell(\ell- N_{\lfloor\sqrt{\ell}\rfloor}>\ve \ell)\la 0$ as $\ell\to\infty$, for every $\ve>0$, and since we may couple 
the variables $\xi_k$ with i.i.d.~variables $\xi'_k$ taking values in $\{0,1,\ldots\}$ with distribution given by $\P(\xi'_k=j)=p_{-j}$
for $j\geq 1$, in such a way that $\xi_k\leq \xi'_k$ for every $k$, a renewal-type argument 
using \eqref{meanstep} gives \eqref{App1}. \endproof

In the next lemma, we come back to the peeling by layers algorithm of the UIPQ. We consider 
the  ``inverse process'' of $(R_k)_{k\geq 1}$ defined by
$$T_n=k\quad\hbox{ if and only if }\quad R_k\leq n<R_{k+1},$$
with the convention $R_0=0$. 

\begin{lemma}
\label{layers0}
Let $\ve \in(0,1)$ and $0<s<t$.
We have
\begin{equation}
\label{layers1}
\P\Big((3-\ve)(T_{\lfloor nt\rfloor}-T_{\lfloor ns\rfloor})\times \min_{{\lfloor ns\rfloor}\leq i\leq {\lfloor nt\rfloor}} S^\uparrow_i \leq {\lfloor nt\rfloor}-{\lfloor ns\rfloor}+\ve n\Big)
\build{\la}_{n\to\infty}^{} 1,
\end{equation}
and 
\begin{equation}
\label{layers2}
\P\Big((3+\ve)(T_{\lfloor nt\rfloor}-T_{\lfloor ns\rfloor})\times \max_{{\lfloor ns\rfloor}\leq i\leq {\lfloor nt\rfloor}} S^\uparrow_i \geq {\lfloor nt\rfloor}-{\lfloor ns\rfloor}\Big)
\build{\la}_{n\to\infty}^{} 1.
\end{equation}
\end{lemma}

\proof We start by proving \eqref{layers1}. Since $R_{T_{\lfloor ns\rfloor}+1}\geq \lfloor ns\rfloor$ and $R_{T_{\lfloor nt\rfloor}}\leq \lfloor nt\rfloor$, we have
\begin{equation}
\label{App-tec4}
{\lfloor nt\rfloor}-{\lfloor ns\rfloor}\geq \sum_{j=1}^{ T_{\lfloor nt\rfloor}-T_{\lfloor ns\rfloor}-1} \Big(R_{T_{\lfloor ns\rfloor}+j+1}-R_{T_{\lfloor ns\rfloor}+j}\Big).
\end{equation}
Now note that, for every $j\geq 1$, conditionally on the event $S^\uparrow_{R_{T_{\lfloor ns\rfloor}+j}}=\ell$, the distribution of
$$R_{T_{\lfloor ns\rfloor}+j+1}-R_{T_{\lfloor ns\rfloor}+j}$$
coincides with the distribution of $R_1$ under $\P_\ell$. Recalling that $S^\uparrow_k\to \infty$ 
a.s. as $k\to\infty$, we then deduce from 
\eqref{App1} that, for every $j\geq 1$,
$$\P\Big( R_{T_{\lfloor ns\rfloor}+j+1}-R_{T_{\lfloor ns\rfloor}+j} > (3-\ve)\,S^\uparrow_{R_{T_{\lfloor ns\rfloor}+j}}\Big) \build{\la}_{n\to\infty}^{} 1.$$

Let $\eta>0$. Consider also a constant $K>1$. Then for any fixed $\delta\in(0,1)$, the bound 
$$R_{T_{\lfloor ns\rfloor}+j+1}-R_{T_{\lfloor ns\rfloor}+j} > (3-\ve)\,S^\uparrow_{R_{T_{\lfloor ns\rfloor}+j}}$$
holds for at least $\lfloor Kn^{1/3}\rfloor-\lfloor \delta n^{1/3}\rfloor$ values of $j\in\{1,2,\ldots,\lfloor K n^{1/3}\rfloor\}$, except 
possibly on an event of probability bounded by $\eta$ when $n$ is large. Using \eqref{App-tec4}, we obtain that the bound 
$$\Big(((T_{\lfloor nt\rfloor}-T_{\lfloor ns\rfloor}-1)\wedge \lfloor K n^{1/3}\rfloor) - \delta n^{1/3}\Big) \times (3-\ve) \min_{{\lfloor ns\rfloor}\leq i\leq {\lfloor nt\rfloor}} S^\uparrow_i 
\leq {\lfloor nt\rfloor}-{\lfloor ns\rfloor}$$
holds except on an event of probability bounded by $\eta$ when $n$ is large. 

Now use \eqref{App-tec1} to observe that, if $\delta$ has been chosen small enough, we have
$$\delta n^{1/3} \times 3\min_{{\lfloor ns\rfloor}\leq i\leq {\lfloor nt\rfloor}} S^\uparrow_i \leq \ve n,$$
except on a set of probability bounded by $\eta$ when $n$ is large, and
therefore the bound
$$\Big((T_{\lfloor nt\rfloor}-T_{\lfloor ns\rfloor}-1)\wedge \lfloor K n^{1/3}\rfloor\Big) \times (3-\ve) \min_{{\lfloor ns\rfloor}\leq i\leq {\lfloor nt\rfloor}} S^\uparrow_i 
\leq {\lfloor nt\rfloor}-{\lfloor ns\rfloor}+\ve n$$
holds except on a set of probability bounded by $2\eta$ when $n$ is large. Furthermore, we may also assume
that $K$ has been chosen large enough so that
$$\lfloor K n^{1/3}\rfloor \times (3-\ve) \min_{{\lfloor ns\rfloor}\leq i\leq {\lfloor nt\rfloor}} S^\uparrow_i > \lfloor nt\rfloor - \lfloor ns\rfloor + \ve n $$
except on a set of probability bounded by $\eta$ when $n$ is large. By putting together the previous observations, we conclude that
$$(T_{\lfloor nt\rfloor}-T_{\lfloor ns\rfloor}-1) \times (3-\ve) \min_{{\lfloor ns\rfloor}\leq i\leq {\lfloor nt\rfloor}} S^\uparrow_i 
\leq {\lfloor nt\rfloor}-{\lfloor ns\rfloor}+\ve n$$
except on a set of probability bounded by $3\eta$ when $n$ is large. This proves \eqref{layers1}.

Let us turn to the proof of \eqref{layers2}. In a way similar to  \eqref{App-tec4} we now 
write
\begin{equation}
\label{App-tec5}
{\lfloor nt\rfloor}-{\lfloor ns\rfloor}\leq 
R_{T_{\lfloor ns\rfloor} +1} -\lfloor ns\rfloor +
\sum_{j=1}^{ T_{\lfloor nt\rfloor}-T_{\lfloor ns\rfloor}} \Big(R_{T_{\lfloor ns\rfloor}+j+1}-R_{T_{\lfloor ns\rfloor}+j}\Big).
\end{equation}
Let $\delta>0$. Using \eqref{App-tec1} and \eqref{mini-peri}, we can fix $\alpha>0$ small enough
so that, for every $n$,
\begin{equation}
\label{App-tec6}
\P\Big(\min_{k\geq \lfloor ns\rfloor} S_k^\uparrow \leq \alpha n\Big) \leq \delta.
\end{equation}
Then, using the Markovian properties of the peeling algorithm and 
\eqref{App2}, we get that for $n$ large enough (such that $\alpha n>A_\ve$), for every $j\geq 1$,
$$\P\Big( R_{T_{\lfloor ns\rfloor}+j+1}-R_{T_{\lfloor ns\rfloor}+j} \geq (3+\ve)\,S^\uparrow_{R_{T_{\lfloor ns\rfloor}+j}},\;\min_{k\geq \lfloor ns\rfloor} S_k^\uparrow \geq \alpha n\Big)
\leq C_\ve\,e^{-\rho_\ve\alpha n}.$$
The same bound holds for the probability of the event
$$\Big\{ R_{T_{\lfloor ns\rfloor} +1} -\lfloor ns\rfloor \geq (3+\ve)S^\uparrow_{\lfloor ns\rfloor},\;
\min_{k\geq \lfloor ns\rfloor} S_k^\uparrow \geq \alpha n\Big\}.$$
By \eqref{layers1}, we have $T_{\lfloor nt\rfloor}- T_{\lfloor ns\rfloor}\leq n$ with probability at least
$1-\delta$ when $n$ large. Using also \eqref{App-tec6} and the preceding estimates, we
deduce from \eqref{App-tec5} that we
have for all large $n$, 
$${\lfloor nt\rfloor}-{\lfloor ns\rfloor} \leq (T_{\lfloor nt\rfloor}-T_{\lfloor ns\rfloor}+1)
\times (3+\ve) \max_{{\lfloor ns\rfloor}\leq i\leq {\lfloor nt\rfloor}} S^\uparrow_i,$$
except on a set of probability bounded above by $2\delta+ C_\ve(n+1)e^{-\rho_\ve\alpha n}$.
This completes the proof of \eqref{layers2}. \endproof

We need a last lemma before we proceed to the proof of Proposition \ref{conv-peri-vol}.

\begin{lemma}
\label{tightness}
The sequence of the distributions of $n^{-1/3}T_n$, $n\geq 1$, is tight.
\end{lemma}

\proof We need to verify that $\P(T_n>An^{1/3})$ can be made small
uniformly in $n$ by choosing $A>0$ large. We observe that
$$\P(T_n>A n^{1/3}) \leq \P(R_{\lfloor A n^{1/3}\rfloor}\leq n)
\leq \P(V^{(\infty)}_{\lfloor A n^{1/3}\rfloor}\leq M^\uparrow_n)$$
since $V^{(\infty)}_k=M^\uparrow_{R_k}$ for every $k\geq 1$. By
\eqref{App-tec1}, we have
$$n^{-4/3}\,M_n^\uparrow \build{\la}_{n\to\infty}^{\rm(d)} \mathcal{V}^{\uparrow}_{1/3}.$$
On the other hand, the known results about the convergence of hulls in the
UIPQ (see \cite[Section 5]{Hull} or \cite[Section 6.2]{CLG-peeling}) 
imply that
$$n^{-4/3}\,V^{(\infty)}_{\lfloor A n^{1/3}\rfloor} \build{\la}_{n\to\infty}^{\rm(d)} \chi_A,$$
with a limiting process $(\chi_t)_{t\geq 0}$ such that $\chi_t\uparrow \infty$
as $t\uparrow \infty$. Notice that \cite{Hull,CLG-peeling} deal with ``simple'' hulls instead
of the lazy hulls we consider here, but as far as volumes are concerned this makes
no difference (we leave it as an exercise for the reader to check that the volume 
of the lazy hull of radius $r\geq 1$ is bounded above by the volume of the simple
hull of the same radius, and bounded below by the volume of the simple
hull of radius $r-1$). It follows that
$$\limsup_{n\to\infty} \P(V^{(\infty)}_{\lfloor A n^{1/3}\rfloor}\leq M^\uparrow_n) \leq \P(\chi_A\leq 
\mathcal{V}^{\uparrow}_{1/3})$$
and the right-hand side can be made arbitrarily small by choosing $A$ large. This completes the proof.
\endproof

Let us turn to the proof of Proposition \ref{conv-peri-vol}. We follow the method of \cite[Proof of Proposition 10]{CLG-peeling}. We fix $0<s<t$. From  Lemma 
\ref{tightness} (or from the bound \eqref{layers1}), we may assume that along a suitable sequence of values of $n$, we have, for every integer $k\geq 0$ and every $1\leq i\leq 2^k$,
$$n^{-1/3}\Big( T_{\lfloor n(s+i2^{-k}(t-s))\rfloor} - T_{\lfloor n(s+(i-1)2^{-k}(t-s))\rfloor}\Big)
\build{\la}_{n\to\infty}^{\rm (d)} \Lambda^{(s,t)}_{k,i}$$
and these convergences hold jointly, and jointly with \eqref{App-tec1}. From
Lemma \ref{layers0}, we then get
$$\frac{2^{-k}(t-s)} {3\;{\displaystyle \sup_{s+(i-1)2^{-k}(t-s)\leq r \leq s+i2^{-k}(t-s)} \Upsilon^\uparrow_{r/3}}}
\leq \Lambda^{(s,t)}_{k,i} \leq 
\frac{2^{-k}(t-s)} {3\;{\displaystyle \inf_{s+(i-1)2^{-k}(t-s)\leq r \leq s+i2^{-k}(t-s)} \Upsilon^\uparrow_{r/3}}}.
$$
By summing these bounds over $i$ and then letting $k\to \infty$, we get that necessarily
$$\Lambda^{(s,t)}_{0,1} = \frac{1}{3} \int_s^t \frac{\mathrm{d}r}{\Upsilon^\uparrow_{r/3}}.$$
Since the limit does not depend on the sequence of values of $n$
we were considering, we have thus proved that, for every $0<s<t$,
$$n^{-1/3}(T_{\lfloor nt\rfloor}-T_{\lfloor ns\rfloor}) \build{\la}_{n\to\infty}^{\rm (d)}
\frac{1}{3} \int_s^t \frac{\mathrm{d}r}{\Upsilon^\uparrow_{r/3}}$$
and this convergence holds jointly with \eqref{App-tec1}. 

By Lemma \ref{tightness}, $n^{-1/3}T_{\lfloor ns\rfloor}$ is small in probability when $s$
is small, uniformly in $n$. The convergence in the last display then implies that
we have also, for every $t>0$,
$$n^{-1/3}\,T_{\lfloor nt\rfloor}\build{\la}_{n\to\infty}^{\rm (d)}
\frac{1}{3} \int_0^t \frac{\mathrm{d}r}{\Upsilon^\uparrow_{r/3}}=\int_0^{t/3}\frac{\mathrm{d}r}{\Upsilon^\uparrow_{r}}$$
and this convergence holds jointly with \eqref{App-tec1}. By a monotonicity argument, the 
convergence in the last display holds in the functional sense (in the Skorokhod space), if both sides are
viewed as processes indexed by $t\geq 0$. It follows that we have also
$$(n^{-3}R_{\lfloor nt\rfloor})_{t\geq 0} \build{\la}_{n\to\infty}^{(d)} (3\Psi_t)_{t\geq 0},$$
where 
$$\Psi_t=\inf\Big\{s\geq 0: \int_0^s \frac{\mathrm{d}r}{\Upsilon^\uparrow_{r}} >t\Big\}$$
in agreement with the notation introduced in Proposition \ref{conv-peri-vol}. 
The preceding converges again holds in the functional sense and jointly with
\eqref{App-tec1}. Finally, we use \eqref{App-tec10} and \eqref{App-tec1}
to conclude that
$$\Big(n^{-2}H^{(\infty)}_{\lfloor n t\rfloor}, n^{-4} V^{(\infty)}_{\lfloor n t\rfloor}\Big)_{t\geq 0}
= \Big(n^{-2}S^\uparrow_{n^3\times(n^{-3}R_{\lfloor n t\rfloor})}, n^{-4} M^\uparrow
_{n^3\times(n^{-3}R_{\lfloor n t\rfloor})}\Big)_{t\geq 0}$$
converges in distribution to
$$\Big(\Upsilon^\uparrow_{\Psi_t},\mathcal{V}^\uparrow_{\Psi_t}\Big)_{t\geq 0}.$$
This completes the proof of Proposition \ref{conv-peri-vol}. \endproof

\medskip
\noindent{\bf Acknowledgments.} I wish to thank the referee for a very careful reading of the manuscript and several useful
remarks. I am indebted to Nicolas Curien for many helpful discussions about
the lazy peeling. I also thank Gr\'egory Miermont for keeping me informed of his ongoing work
on Brownian disks.

\end{document}